\documentclass[10pt,reqno]{amsart}
\usepackage[a4paper,margin=.8in]{geometry}
\usepackage[utf8]{inputenc}
\usepackage{enumerate}
\usepackage{float}
\usepackage{amsmath}
\usepackage{amsfonts}
\usepackage{amssymb}
\usepackage{physics}
\usepackage{cancel}
\usepackage{graphicx}
\usepackage{caption}
\usepackage{subcaption}
\usepackage[showonlyrefs]{mathtools}
\usepackage{dsfont}
\usepackage[english]{babel}
\usepackage{amsthm}
\usepackage{mathrsfs}
\usepackage{tikz}
\usetikzlibrary{decorations.pathmorphing,calc}
\usetikzlibrary{arrows.meta}
\usepackage{blindtext}
\usepackage[colorinlistoftodos]{todonotes}
\usepackage{upgreek}
\usepackage[usenames,dvipsnames]{xcolor}
\usepackage[colorlinks=true,linkcolor=Red,citecolor=Green]{hyperref}
\usepackage[T1]{fontenc}
\usepackage{xcolor}
\usepackage{listings}
\usepackage{textcomp}

\definecolor{mmaBlue}{RGB}{38,72,145}
\definecolor{mmaComment}{RGB}{105,105,105}
\definecolor{mmaString}{RGB}{115,62,100}
\definecolor{mmaBackground}{RGB}{248,249,252}

\usetikzlibrary{arrows.meta,decorations.markings}

\lstdefinestyle{MPMathematica}{
  language=Mathematica,
  basicstyle=\ttfamily\footnotesize,
  keywordstyle=\color{mmaBlue},
  commentstyle=\color{mmaComment},
  stringstyle=\color{mmaString},
  backgroundcolor=\color{mmaBackground},
  frame=l,
  rulecolor=\color{mmaBlue!45},
  framerule=0.7pt,
  xleftmargin=8pt,
  framexleftmargin=5pt,
  aboveskip=8pt,
  belowskip=8pt,
  columns=fullflexible,
  keepspaces=true,
  showstringspaces=false,
  upquote=true,
  breaklines=true,
  breakatwhitespace=true,
  tabsize=2,
  title={\color{mmaBlue}\small\ttfamily In[1]:=},
  morekeywords={
    FullSimplify,Manipulate,Plot,Text,Scaled,
    Sqrt,Abs,Epilog,PlotRange
  },
  literate=
    {Δ}{{\ensuremath{\Delta}}}1
    {δ}{{\ensuremath{\delta}}}1
}

\mathtoolsset{showonlyrefs,showmanualtags}

\makeatletter
\let\oldlabel@in@display\label@in@display
\def\label@in@display#1{%
    \oldlabel@in@display{#1}%
    \MT@newlabel{#1}%
}
\makeatother

\theoremstyle{plain}
\newtheorem{theorem}{Theorem}[section]
\newtheorem{lemma}[theorem]{Lemma}
\newtheorem{proposition}[theorem]{Proposition}
\newtheorem{corollary}[theorem]{Corollary}

\theoremstyle{definition}
\newtheorem{definition}[theorem]{Definition}
\newtheorem{remark}[theorem]{Remark}

\DeclareMathOperator\supp{supp}
\DeclareMathOperator\ric{Ric}
\DeclareMathOperator\diag{diag}
\DeclareMathOperator\Char{Char}

\DeclareMathOperator\ind{ind}
\DeclareMathOperator\Op{Op}
\DeclareMathOperator\esssup{ess\;sup}

\DeclareMathOperator\Ell{Ell}
\DeclareMathOperator\Diff{Diff}

\newcommand{\scxi}{{}^{sc}\xi}
\newcommand{\sceta}{{}^{sc}\eta}

\newcommand{\scT}{{}^{sc}T}

\newcommand{\bxi}{{}^{b}\xi}
\newcommand{\beeta}{{}^{b}\eta}

\newcommand{\bT}{{}^{b}T}

\newcommand{\WF}{\text{WF}}

\newcommand{\mfrm}{\mathfrak{m}}
\newcommand{\bfa}{\mathbf{a}}
\newcommand{\bfb}{\mathbf{b}}
\newcommand{\bfDelta}{\Delta_{r}}
\newcommand{\bfr}{r_{\mathbf{a},\mathbf{b}}}
\newcommand{\bfp}{p_{\mathbf{a},\mathbf{b}}}
\newcommand{\sgn}{\text{sgn}}
\newcommand{\chr}{\text{Char}}

\numberwithin{equation}{section}

\title[Price's law on 1+4 dimensional Myers--Perry]
{Price's Law for scalar waves on 1+4 dimensional subextremal Myers--Perry spacetimes}
\author{Pieter-Bart Peters}
\address{Department of Mathematics, ETH Z\"urich,
R\"amistrasse 101, 8092 Z\"urich, Switzerland}
\email{pieterbart.peters@math.ethz.ch}
\date{\today}

\begin{document}

\begin{abstract}
We consider Myers--Perry spacetimes
$(M,g_{\bfa,\bfb,\mfrm})$ in $1+4$ dimensions with two spin
parameters $(\bfa,\bfb)\in\mathbb{R}^2$, which are a direct
generalization of the well known Kerr spacetime in $1+3$
dimensions. We prove that the Hamiltonian flow associated to
$g_{\bfa,\bfb,\mfrm}^{-1}$ is normally hyperbolically trapped
for the full subextremal range
\[
    (|\bfa|+|\bfb|)^2<2\mfrm.
\]
As an application, we prove a Price law with explicit leading
order terms for solutions $\phi(t,x)$ of the scalar wave equation
\[
    \Box_{g_{\bfa,\bfb,\mfrm}}\phi=0
\]
under the assumption of real mode stability. Our precise
asymptotics imply that
$\phi(t,x)=ct^{-3}+\mathcal{O}(t^{-4+})$, where
$c\in\mathbb{C}$ is an explicitly computable constant.
\end{abstract}

\maketitle

\tableofcontents

\newpage
\section{Introduction}

\subsection{$1+d$-dimensional Myers-Perry spacetimes}

In $1986$ and $1987$ R.C. Myers and M.J. Perry proposed the generalization of the Kerr black hole to $1+d$-dimensions with $d\geq 3$ \cite{myers1986black},\cite{myers1987higher}, that are now known as Myers-Perry spacetimes. The metric $g_{\bfa_i,\mu}$, describing a spinning black hole with $N:=\left\lfloor \frac{d}{2}\right\rfloor$ with angular momenta (rotational axes) $\{\bfa_i\}$ and mass $\mu$ is a stationary solution to the the vacuum Einstein field equations 
\begin{align}\label{eq:GeneralRicci}
    \ric{(g_{\bfa_i,\mu})}=0.
\end{align}
To describe the metrics $g_{\bfa_i,\mu}$ in Boyer-Lindquist coordinates for a given spatial dimension $d$, let us first define on $\mathbb{S}^{d-1}\subset \mathbb{R}^d$ the directional cosines by taking a point $\omega\in \mathbb{S}^{d-1}\subset \mathbb{R}^d$ and writing 
\begin{align}
    (\omega_{2i-1},\omega_{2i})=\mu_i(\cos \phi_i,\sin\phi_i), \; i=1,...,N,
\end{align}
and imposing 
\begin{align}
    \begin{cases}
        &\sum_{i=1}^{N} \mu^2_i=1\;\;\;\;\;\;\;\;\;\;\;\;\;d=2N\\
        & \sum^N_{i=1}\mu_i^2+\alpha^2=1\;\;\;\;\;d=2N+1,
    \end{cases}
\end{align}
where $\alpha = \omega_d\in [-1,1]$ is the remaining component in the second case. With this convention we can define the domain of outer communication $(\tilde M,g_{\bfa_i,\mu})$ both for even and odd dimensions $d$ in Boyer-Lindquist coordinates, i.e. on a chart 
\begin{align}
    \tilde M =(-\infty,\infty)_t\times (r_+,\infty)_r \times \mathbb{S}^{d-1}\subset \mathbb{R}^{1+d}
\end{align}
where $r_+$ is the largest positive root of 
\begin{align}\label{eq:GeneralDeltar}
    \Delta_r := \frac{\Pi-\mu r^\gamma}{r^{2N-2}}, \text{ with } \gamma =
\begin{cases}
2, & d \text{ even},\\
1, & d \text{ odd}.
\end{cases}
\end{align}

For even dimensional $d=2N$ the Myers-Perry metric is then given by $(-,\underbrace{+,+,...,+}_d)$ convention \cite{horowitz2012black} by
\begin{align}
    g_{\bfa_i,\mu}=-dt^2+ \frac{2\mu}{\Pi F}\left(dt+\sum^{N}_{i=1} \bfa_i \mu_i d\phi_i\right)^2 +\frac{\Pi F}{\Pi-\mu r^2}dr^2 + \sum^{N}_{i=1} (r^2+\bfa_i^2)(d\mu_i^2+\mu_i^2 d\phi_i^2).
\end{align}

Here the quantities $\Pi$ and $F$ are given by 
\begin{align}
    F&= 1-\sum^{N}_{i=1}\frac{\bfa_i^2\mu_i^2}{r^2+\bfa_i^2} \text{ and } \Pi=\prod_{i=1}^{N}(r^2+\bfa_i^2).
\end{align}

For odd dimensional Myers-Perry metrics $d=2N+1$ (of which Kerr is a member for $d=3$) also generalizes as 
\begin{align}
    g_{\bfa_i,\mu}=-dt^2+\frac{\mu r}{\Pi F}\left(dt  +\sum_{i=1}^{\lfloor\frac{d}{2}\rfloor} \bfa_i \mu_i^2 d\phi_i \right)^2&+ \frac{\Pi F}{\Pi -\mu r}dr^2 \nonumber\\
    &+\sum_{i=1}^{\lfloor\frac{d}{2}\rfloor} (r^2+\bfa_i^2)(d\mu_i^2+\mu_i^2d\phi_i^2) +r^2d\alpha^2
\end{align}
The ADM mass of the black hole for Myers-Perry is given in both cases through the relation
\begin{align}
    \mfrm = \frac{(d-1)\Omega_{d-1}}{16\pi G}\mu,\; \text{ and } \; \Omega_{d-1}:= \frac{2\pi^{\frac{d}{2}}}{\Gamma\left(\frac{d}{2}\right)}, \label{eq:BlackHoleMassRelation}
\end{align}
where we set $G=1$ (note that in the Kerr-case $d=3$ we have $\mu = 2\mfrm$) and $\Gamma$ denotes the Euler gamma function.
For generic values of the spin parameters $\{\bfa_i\}$ the connected isometry group of the Myers-Perry spacetimes solution to (\ref{eq:GeneralRicci}) is given by
\begin{align}
    \mathbb{R}_t \times U(1)^{N},
\end{align}
where $U(1)^N$ is the maximal torus of the spatial rotation group $\text{SO}(d)$. Equivalently, the rotational symmetry of the symmetry of the static Schwarzschild-Tangherlini solution ($\bfa_i=0$ for all $1\leq i \leq N$) is reduced to its maximal commuting subgroup once the independent rotation planes are specified.

\subsection{Extremal Myers-Perry black holes and dimensionality.}

Our main theorem involves Price's law in the special case $d=4$ for the wave equation 
\begin{align}\label{eq:WaveEquationGeneral}
    \Box_{g_{\bfa,\bfb,\mfrm}}\phi = -|g|^{-1/2}\partial_\mu (|g|^{1/2}g^{\mu\nu}\partial_\nu\phi) =0,\; \phi(0,x)=\phi_0(x), \; \partial_t \phi (0,x) =\phi_1(x),
\end{align}

where to keep consistency with literature on the case $d=4$ we rename the momenta $\{\bfa_i\}_i$ for $i=1,2$ in this case to $\bfa$ and $\bfb$ and also renormalize $G$ such that $\mu=2\mfrm$. The theorem is as of writing still condition on the assumption of real mode stability for the wave equation.

\begin{theorem}\label{thm:MainPriceLaw}
    (Price's Law on $1+4$-dimenional Myers-Perry spacetime in the full subextremal range.). Let $(M,g_{\bfa,\bfb,\mfrm})$ be a $1+4$-dimensional Myers-Perry spacetime and $\bfa,\bfb\in \mathbb{R}^2$ satisfy the subextremality constraint
    \begin{align}\label{eq:FirstExtremalityConstraint}
        (|\bfa|+|\bfb|)^2< 2\mfrm
    \end{align}
    Let $\phi$ denote the associated solution to (\ref{eq:WaveEquationGeneral}) with compactly supported initial date $\phi_0,\phi_1\in C_c^\infty(X^\circ)$ on the hypersurface $X^\circ=\{t_*=0\}$. Then, for a constant $c$, computable in terms of $\phi_0,\phi_1 \in C_c^\infty(X^\circ)$ we have that the solution $\phi$ decays according to 
    \begin{align}\label{eq:MyersPerryDecayRemainder}
        \left|\phi - ct_*^{-3}\left(\frac{t_*}{t_*+2r}\right)^{3/2}\right|\leq C_\epsilon t_*^{-4+\epsilon}\left(\frac{t_*}{t_*+2r}\right)^{3/2}.
    \end{align}
    Moreover the radiation field $F(t_*, \omega) :=\lim_{r\to \infty} r^{3/2}\phi(t_*,r,\omega)$ has a leading order term with remainder 
    \begin{align}\label{eq:RadiationFieldDecayRemainder}
        \left|F-\frac{1}{2\sqrt{2}}ct_*^{-3/2}\right|\leq C_\epsilon t_*^{-5/2+\epsilon}.
    \end{align}
    We note that the decay rates in (\ref{eq:MyersPerryDecayRemainder}) and (\ref{eq:RadiationFieldDecayRemainder}) hold for $\phi - ct_*^{-3}\left(\frac{t_*}{t_*+2r}\right)^{3/2}$ for any finite derivatives along $t_*\partial_{t_*}$, $r\partial_r$ and rotational vector fields on $\mathbb{S}^3$.
\end{theorem}

In order to prove theorem \ref{thm:MainPriceLaw} we will establish high energy bounds on the spectral family associated by the wave equation, given by

\begin{align}
    \hat{\Box}_{g_{\bfa,\bfb,\mfrm}}(\sigma):=e^{i\sigma t_*}\circ \Box_{g_{\bfa,\bfb,\mfrm}}\circ e^{-i\sigma t_*}\qquad \sigma \in \overline{\mathbb{H}}.
\end{align} 
The key result leading this is the following normally hyperbolic trapping result, valid in the full subextremal range $(|\bfa|+|\bfb|)^2<2\mfrm$.

\begin{theorem}[Normally hyperbolic trapping on $1+4$ dimensional Myers-Perry, full extremal range]\label{thm:Theorem1}
Let $(M,g_{\bfa,\bfb,\mfrm})$ be a subextremal $1+4$ dimensional Myers-Perry spacetime, with $\mfrm>0$ and $\bfa,\bfb\in\mathbb{R}^2$ satisfy
\begin{align}
    (|\bfa|+|\bfb|)^2<2\mfrm.
\end{align}
Let $\tilde M$ denote the domain of outer communication and define the dual metric function $G$ and dual light-cone $G\subset T^*\tilde{M}$ associated to $g_{\bfa,\bfb,\mfrm}$ via
\begin{align}
    G=g_{a,b,m}^{-1}(\upalpha,\upalpha),
    \qquad\text{ and }\qquad
    \Sigma_G=\{G=0\}\backslash O,
\end{align}
where $O$ denotes the zero section of $T^*\widetilde M$. Then there exists a trapped set $\Gamma\subset\Sigma$ of the null geodesic flow as a
smooth invariant submanifold of codimension two in $\Sigma$.
Near $\Gamma$, there exist smooth invariant hypersurfaces
$\Gamma^\pm\subset\Sigma$ which intersect transversely,
\begin{align}
    \Gamma=\Gamma^+\cap\Gamma^-.
\end{align}
These hypersurfaces admit smooth defining functions $\varphi_\pm$
satisfying
\begin{align}\label{eq:NormallyHyperbolicTrappingRelations}
    H_G\varphi_+=\nu_+\varphi_+,
    \qquad
    H_G\varphi_-=-\nu_-\varphi_-,
    \qquad
    \nu_\pm|_\Gamma>0,
\end{align}
and
\begin{align}
    \{\varphi_-,\varphi_+\}>0
    \quad\text{on }\Gamma.
\end{align}
In particular, the trapped set is normally hyperbolic in the sense
of Definition~\ref{def:TrappedSet}.
\end{theorem}

Theorem \ref{thm:Theorem1} concerns the exact Myers-Perry family throughout the full subextremal range. It asserts that the Hamiltonian flow $H_G$ associated to the dual metric function splits on $\Sigma_G$ in a stable manifold $\Gamma^+$ and an unstable manifold $\Gamma^-$ that intersect in the trapped set $\Gamma$ which contains null-geodesics that are bounded. We do not establish the stronger $r$-normal hyperbolicity estimates needed for the perturbative
extension of our argument. Consequently, the results below do not assert stability of the trapping structure or persistence of Price's law under general small perturbations of the metric.

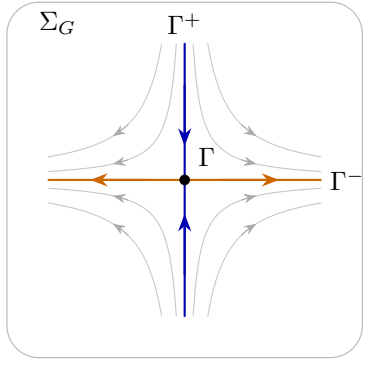
\begin{figure}[H]
    \centering

\begin{tikzpicture}[>=Stealth, line cap=round]
    \draw[gray!60, rounded corners=12pt]
        (-2.35,-2.35) rectangle (2.35,2.35);
    \node[anchor=north west] at (-2.05,2.35) {$\Sigma_G$};

    \foreach \sx in {-1,1} {
        \foreach \sy in {-1,1} {
            \begin{scope}[xscale=\sx, yscale=\sy]
                \foreach \c in {0.2,0.55} {
                    \draw[gray!45, thin]
                        plot[domain=\c/1.8:1.8, samples=60]
                        (\x,{\c/\x});
                    \draw[gray!65, ->]
                        plot[domain=0.75:0.95, samples=12]
                        (\x,{\c/\x});
                }
            \end{scope}
        }
    }

    \draw[thick, blue!70!black] (0,-1.8) -- (0,1.8);
    \draw[thick, orange!80!black] (-1.8,0) -- (1.8,0);

    \foreach \s in {-1,1} {
        \draw[thick, blue!70!black, ->]
            (0,{\s*1.25}) -- (0,{\s*0.45});

        \draw[thick, orange!80!black, ->]
            ({\s*0.45},0) -- ({\s*1.25},0);
    }

    \fill (0,0) circle (2pt);
    \node[above right] at (0.06,0.06) {$\Gamma$};
    \node[above] at (0,1.8) {$\Gamma^+$};
    \node[right] at (1.8,0) {$\Gamma^-$};
\end{tikzpicture}
    \caption{Sketch of the flow of the $H_G$ inside of $\Sigma_G\subset T^*M$ splititing into a stable manifold $\Gamma^+$ and unstable manifold $\Gamma^-$. The flow persists for all $(\bfa,\bfb)\in \mathbb{R}^2$ satisfying $(|\bfa|+|\bfb|)^2<2\mfrm$.}
\end{figure}

To give motivation to the case why $1+4$-dimensional Myers-Perry spacetimes are, together with Kerr, a special case in the class of Myers-Perry spacetimes of arbitary dimension we need to highlight the notion of extremality and its relation to dimension. For $d=4$ only, Myers-Perry spacetimes exhibit extremality in a very similar fashion to the standard Kerr black hole $d=3$ in the case when $\bfa\neq 0$, $\bfb\neq 0$. The $\Delta_r$ as defined by equation (\ref{eq:GeneralDeltar}) is in this case a quartic polynomial, but has for the entire subextremal range (\ref{eq:FirstExtremalityConstraint}) two negative and two positive real roots and two distinct real roots. The largest root $r_+$ can be used to identify the event horizon as the hypersurface $\{r=r_+\}$ and its smaller root $r_-$ can be used to find a regular Cauchy horizon at $\{r=r_-\}$. Whenever extremality is reached for any configuration $(\bfa,\bfb)\in \mathbb{R}^2$ satisfying $(|\bfa|+|\bfb|)^2=2\mfrm$ there is equality $r_+=r_-$ and for any $(\bfa,\bfb)$ with $(|\bfa|+|\bfb|)^2>2\mfrm$ this structure collapses again into a spacetime with a naked ring singularity no positive real roots $r_+$ and $r_-$ remain, see Appendix B in \cite{horowitz2012black}. An important difference with the Kerr case is that for whenever one spin turns zero, say $\bfb=0$, the single spinning Myers-Perry black hole retains the causal structure of Schwarzschild-Tangherlini and it will not exhibit a regular Cauchy horizon for $|\bfa|^2< 2\mfrm$ and not even a regular event horizon "at the vertices" $|\bfa|^2=2\mfrm$.
For the $1+4$-dimensional case we will revisit (\ref{eq:FirstExtremalityConstraint}) and give explicit notions to $r_+,r_-$ and associated concepts like black-hole temperature $T_g$ and surface gravity $\kappa$ in Section \ref{sec:Myers-Perry spacetimes} following the literature. \\
For any $d>4$ the situation described above changes dramatically. To illustrate this, we follow  \cite{emparan2008black} and consider black holes with a single spinning parameter (\cite{horowitz2012black} also considers this as an example). 
In this case the metric (in both even and odd dimensional $d$) reduces and can be written as 
\begin{align}
    g_{\bfa,\mfrm}= -dt^2 + \frac{2\mfrm}{r^{d-4}r^2_{\bfa}}&(dt-\bfa \sin^2\theta)^2+\frac{r_{\bfa}^2}{\Delta_r}dr^2 +r^2_{\bfa}d\theta^2 +(r^2+\bfa^2)\sin^2\theta d\phi^2\nonumber\\
    &+r^2\cos^2\theta d\Omega^2_{d-3},
\end{align}

where we introduced 
\begin{align}
    r_{\bfa}^2=r^2+\bfa^2\cos^2\theta.
\end{align}
Whilst this metric looks deceptively similar to Kerr, one of its key differences lies in expression of $\Delta_r$, which depends on the dimension $d$:
\begin{align}
    \Delta_r = r^2 +\bfa^2 -\frac{2\mfrm}{r^{d-4}}, \label{eq:GeneralDeltaSingleSpinning}
\end{align}
where we note that $\mu = 2\mfrm$ for both $d=3$. Now for $d=3$ and for $d=4$\footnote{For general $(\bfa,\bfb)$ the equation for $d=4$ is quartic.} we obtain in this case quadratic equation in $r$ that are easily solved for solutions $r_+$ and $r_-$ that when real, can be linked to the event horizon and the Cauchy horizon as before. For $d>4$ however, there is \textit{always} an event horizon as when $r\to 0$ the term $-\frac{\mu}{r^{d-4}}$ starts dominating (\ref{eq:GeneralDeltaSingleSpinning}) and because of its signature will always guarantee a real root for $\Delta_r$. Thus, the existence of an event horizon no longer requires that the momentum $\bfa$ always be bounded by the mass term $\mu$ and $d=4$ is the highest dimension where a momentum bound by the mass $\mfrm$ is a necessary condition for the existence of a horizon. \\
To completely visualize this, consider the cases $d=4$ and $d=5$. In this case we have via (\ref{eq:GeneralDeltar}) that
\begin{align}
    \Delta^{(4)}_r &= r^2+\bfa^2+\bfb^2 +\frac{\bfa^2\bfb^2}{r^2}-2\mfrm \\
    \Delta_r^{(5)}&= r^2+\bfa^2+\bfb^2+\frac{\bfa^2\bfb^2}{r^2}-\frac{2\mfrm}{r}.\label{eq:SixdimensionalDeltar}
\end{align}
The existence for a double root in general requires 
\begin{align}
    \Delta_r =0 \;\text{ and }\;\frac{d}{dr}\Delta_r := \Delta'_r =0. \label{eq:HorizonEquation1}
\end{align}
Solving these equations indeed gives for $\Delta^{(4)}_r$ that $(|\bfa|+|\bfb|)^2= 2\mfrm$ which defines a bounded diamond n $\mathbb{R}^2$. For $\Delta_r^{(5)}(r)$ we have that the $4$ curves that guarantee the existence of an event horizon are given by 
\begin{align}\label{eq:6DimensionalExtremalityCurves}
    (\bfa,\bfb)= \left(\left(\frac{\pi}{2\sqrt{3}}\right)^{1/3}\sqrt{\frac{1-4\nu^3\pm \sqrt{1-16\nu^3}}{4\nu}},\left(\frac{\pi}{2\sqrt{3}}\right)^{1/3}\sqrt{\frac{1-4\nu^3\mp \sqrt{1-16\nu^3}}{4\nu}}\right)
\end{align}
where $\nu=\frac{r_+}{(2\mfrm)^{1/3}}$ is a dimensionless constant depending on $\mfrm$ and the largest root $r_+$ of (\ref{eq:SixdimensionalDeltar}) in (\ref{eq:HorizonEquation1}). This represents a diamond in $\mathbb{R}^2$ that is ''opened up'' at its vertices, giving access to the unbounded regimes in $\mathbb{R}^2$ where either $\bfa=0$ or $\bfb=0$, which are referred to as the ultra-spinning regimes. These ultra-spinning regimes for $d\geq 5$ exhibit instabilities very similar to Gregory-Laflamme instabilities owing to the event horizon "pinching" off into a black membrane \cite{emparan2003instability}, which is also supported by numerical evidence \cite{dias2009instability},\cite{dias2011ultraspinning},\cite{shibata2010bar}. This leaves $d=3$ and $d=4$ as the special cases to investigate when it comes to stability of the wave equation and Einstein field equations.

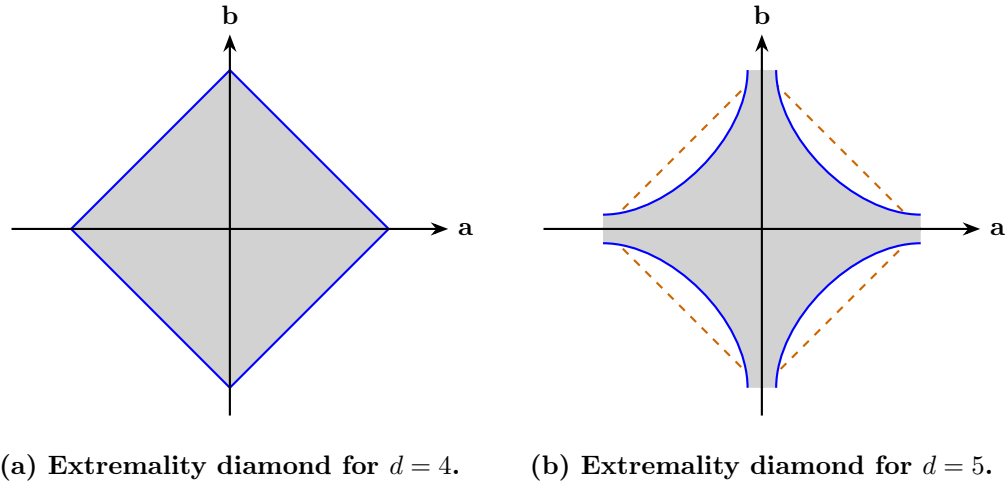
\begin{figure}[H]
\centering
\begin{tikzpicture}[
  >=Stealth,
  scale=1.05,
  axis/.style={->, black, thick},
  ext/.style={blue, thick},
  sing/.style={orange!80!black, dashed, line width=0.8pt},,
  ultra/.style={green!50!black, dotted, very thick}
]

\def\J{2.0}
\def\g{0.18}
\def\c{0.85}

\begin{scope}[shift={(0,0)}]
\node at (0,-3) { \textbf{(a) Extremality diamond for $d=4$.}};

\fill[gray!35] (0,\J)--(\J,0)--(0,-\J)--(-\J,0)--cycle;
\draw[ext] (0,\J)--(\J,0)--(0,-\J)--(-\J,0)--cycle;

\draw[axis] (-2.75,0)--(2.75,0) node[right] {$\mathbf a$};
\draw[axis] (0,-2.35)--(0,2.45) node[above] {$\mathbf b$};
\end{scope}

\begin{scope}[shift={(6.7,0)}]
\node at (0,-3) { \textbf{(b) Extremality diamond for $d=5$.}};

\draw[sing] (0,\J)--(\J,0)--(0,-\J)--(-\J,0)--cycle;

\coordinate (Tleft)  at (-\g,\J);
\coordinate (Tright) at ( \g,\J);
\coordinate (Rtop)   at ( \J,\g);
\coordinate (Rbot)   at ( \J,-\g);
\coordinate (Bright) at ( \g,-\J);
\coordinate (Bleft)  at (-\g,-\J);
\coordinate (Lbot)   at (-\J,-\g);
\coordinate (Ltop)   at (-\J,\g);

\fill[gray!35]
  (Tright)
  .. controls +(0,-\c) and +(-\c,0) .. (Rtop)
  -- (Rbot)
  .. controls +(-\c,0) and +(0,\c) .. (Bright)
  -- (Bleft)
  .. controls +(0,\c) and +(\c,0) .. (Lbot)
  -- (Ltop)
  .. controls +(\c,0) and +(0,-\c) .. (Tleft)
  -- cycle;

\draw[ext]
  (Tright) .. controls +(0,-\c) and +(-\c,0) .. (Rtop);

\draw[ext]
  (Rbot) .. controls +(-\c,0) and +(0,\c) .. (Bright);

\draw[ext]
  (Bleft) .. controls +(0,\c) and +(\c,0) .. (Lbot);

\draw[ext]
  (Ltop) .. controls +(\c,0) and +(0,-\c) .. (Tleft);

\draw[axis] (-2.75,0)--(2.75,0) node[right] {$\mathbf a$};
\draw[axis] (0,-2.35)--(0,2.45) node[above] {$\mathbf b$};
\end{scope}

\end{tikzpicture}
\caption{(a) The Myers-Perry extremality diamond for $d=4$, satisfying $(|\bfa|+|\bfb|)^2<2\mfrm$. (b) The extremality diamond for $d=5$ the bounds of which satisfy (\ref{eq:6DimensionalExtremalityCurves}) with the $d=4$ case as an orange inset. It should be noted that for $d>5$ one finds higher dimensional generalizations of this figure in $\mathbb{R}^{\lfloor \frac{d}{2}\rfloor}$, see e.g. \cite{emparan2008black} Figure 3 for a visualization.}
\end{figure}

\subsection{Comparison with the Kerr literature}

To sketch further parallels between Kerr and $1+4$-dimensional Myers-Perry we focus on the existing literature. In this comparison we will mostly restrict ourselves to results on Kerr pertaining to the wave equation only. For $1+3$-dimensional Kerr, one could argue that the history of stability results dates back to $1989$, where Whiting \cite{whiting1989mode} proved a mode-stability result for solutions $e^{-i\sigma t} e^{i m\phi} R(r) S_{\kappa l m}(\theta)$ to the separated wave equation in the upper half plane $\text{Im }(\sigma)>0$. For $1+4$ and higher dimensional Myers-Perry, separation of variables and integrability of the geodesic flow was an upcoming research topic in the early $2000$'s. Notable mentions especially in the $1+4$ dimensional case regarding this topic on Myers-Perry black holes \cite{frolov2003quantum}, \cite{frolov2003particle} and also \cite{wu2008separability} where both authors show separability of the wave equation and construct a corresponding Killing-Yano tensor for both the scalar wave-equation as well as the massive Dirac equation which can be used to derive analogons of Carter's constant \cite{carter1968hamilton} for this particular case. Later these arguments were extended by \cite{vasudevan2005particle} for black holes with equal sets of spinning parameters as well as \cite{frolov2007hidden} for Myers-Perry of arbitrary dimension and arbitrary spin. Eventually full Liouville integrability of the geodesic flow for Kerr-NUT-(A)DS black holes (with as special case the Myers-Perry black holes) was proven in \cite{page2007complete}.\\
Only recently \cite{cvetivc2022mode}, Cvetic, Pope, Whiting and Zhang proved an analogous result for the separated wave equation on Myers-Perry using an integration method very reminiscent of \cite{whiting1989mode}, partially due to the fact that the radial part of the wave equation in both cases reduces to an equation of confluent Heun type (standard vs reduced). Real mode stability for the wave equation on Kerr has been shown in \cite{shlapentokh2015quantitative} and \cite{andersson2017mode} but no equivalent result exists yet for the wave equation on $1+4$ dimensional Myers-Perry. It will serve as an assumption for the purposes of this paper. \\
Beyond the mode stability problem, the history of stability results on Kerr for the wave equation remains very rich. For the scalar wave equation in particular \cite{dafermos2010decay} established uniform boundedness for restricted rotation angles, whilst \cite{tataru2011local} established integrated energy decay estimates on Kerr. \cite{dafermos2012black}, \cite{dafermos2016decay} moreover established polynomial decay in the full subextremal range, using the fact that superradiant modes remain non-trapped on the Kerr interior. As far as results on $1+4$ dimensional Myers-Perry go, \cite{laul2012localized} \cite{laul2015localized} established local integrated energy decay on $1+4$ dimensional Myers-Perry black holes for small $(\bfa,\bfb)\in \mathbb{R}^2$ only. \\
When it comes to a sharp Price law for the wave equation, this goes back to \cite{price1972nonspherical} predicting $\mathcal{O}(t^{-2\ell -3})$ decay on Schwarzschild spacetimes and $\mathcal{O}(t^{-3})$ on Kerr spacetimes. \cite{tataru2013local}  established an $\mathcal{O}(t_*^{-3})$ upper decay bound for solutions on a large class of asymptotically flat stationary spacetimes under the assumption of local energy decay allowing even for long range pertubations to be conormal following a sharp $t^{-3}$ scalar decay bound without symmetry assumptions for initial data on Schwarzschild by \cite{donninger2012pointwise}. This was later generalized by \cite{metcalfe2012price} to non-stationary spacetimes with suitable decay towards stationarity. The main inspiration for the second part of this article comes from \cite{hintz2022sharp} which established sharp $t^{-3}$ with a leading order term for Kerr in the full extremal range $\bfa<\mfrm$, specifically for conormal solutions to the wave equation with compact initial data using spectral methods. Pushing this method further, and incorporating Kerr as a special case of $s=0$ for the Teukolsky equation, Millet \cite{millet2023optimal} pushed this method even further by proving a spin dependent Price law for the Teukolsky equation on Kerr spacetimes.  The direct counterpart to these results were provided in \cite{angelopoulos2023late} using physical space approaches using conserved quantities at infinity and time inversion to establish precise scalar-wave asymptotics through the subextremal range, including angular projections and the effects of mode coupling, building on a preceding results for Schwarschild and Reissner/Nordstr\"om in \cite{angelopoulos2018late}. Likewise physical space methods were used by \cite{ma2023sharp} to derive Price's law for the Teukolsky equation as well.
Contrasting this with the analogous results on Myers-Perry, the main significant result to mention here is by Luk-Oh \cite{luk2024late} which establishes a Price law for solutions to the wave equation on $1+d$-dimensional Schwarzschild-Tangherlini spacetimes with $d$ odd, employing a method that relies on the strong Huygens principle, which is absent for $1+4$-dimensinal Myers-Perry. Another important reference is \cite{hintz2023linear} which (among other result) establishes a Price law on $1+d$-dimensional asymptotically flat spacetimes with $d$ even under the assumption that the wave equation satisfies high energy estimates. Its companion paper \cite{hintz2026non} further generalized these results to non-linear wave equations on spacetimes that settle into Kerr over time.

\subsection{Outline of proof.}\label{subsec:OutlineofProof}

To strengthen the position of our paper in the context of the existing literature even further we will summarily discuss the method of proof. When it comes to the method of proof, this article will be subdivided in $3$ different parts. 

\begin{itemize}
    \item Part I: Hyperbolic trapping and horizon dynamics of $1+4$-dimensional Myers-Perry spacetimes in the full subextremal range. its most closely related references are \cite{dyatlov2019mathematical}, \cite{dyatlov2015asymptotics}, \cite{dyatlov2016spectral}, and \cite{hintz2026non}.
    \item Part II (bridge): Full dynamical picture of the Hamiltonian flow associated to $\hat{\Box}_{g_{\bfa,\bfb,\mfrm}}(\sigma)$ and corrsponding microlocal estimates leading to a high energy estimate for the wave equation, as well as Fredholm estimates down to $\sigma=0$. As a consequence it shown here that $\hat{\Box}_g(\sigma):\mathcal{X}^{s,r}\to \mathcal{Y}^{s',r'}$ is Fredholm of index $0$ between appropriate Sobolev spaces. Its most closely related references are   \cite{hintz2026non}, \cite{vasy2021limiting}, and \cite{hintz2025introduction}.
    \item  Part III: Late-time asymptotics in the full subextremal range $(\bfa,\bfb)\in D_{\mfrm}$ using a resolvent expansion for $\hat{\Box}_{g_{\bfa,\bfb,\mfrm}}(\sigma)$ near $\sigma =0$. The most closely related reference for the resolvent expansion is \cite{hintz2022sharp}, and this chapter is written with the same notation to make straightforward comparison possible. When it comes to comparing the actual asymptotics, in view of the weak Huygens principle, \cite{hintz2023linear} has been an important point of orientation for comparing the results. 
  
\end{itemize}

The main goal of this article is to eventually solve forward problems wave problems of
\begin{align}
\Box_g \phi = f \in C_c^\infty(\mathbb{R}_{t_*}\times X^\circ),
\end{align}
here $X$ is the spatial manifold associated to $M$, that is if $(M,g)$ represents the total Myers-Perry spacetime we will then decompose $M$ as $\mathbb{R}_{t_*}\times X$ where $X$ corresponds to $\{t_*=0\}$ (using again stationarity of $g$).
We will proceed by means of the Fourier transform in $t_*$ (using stationarity of the metric $g$ in $t_*$)
\begin{align}\label{eq:ForwardWaveSolution}
\phi(t_*,x) =\frac{1}{2\pi}\int_{\mathbb{R}}e^{-i\sigma t_*} \hat{\Box}_g(\sigma)^{-1}\hat{f}(\sigma,x) \,d\sigma
\end{align}
where $\hat{\phi}(\sigma,x) = \int_{\mathbb{R}} e^{i\sigma t_*}\phi(t_*,x) \, dt_*$ and $\hat{\Box}_g(\sigma)$ is obtained from $\Box_g$ by conjugation with $e^{i\sigma t_*}$, i.e. replacing all $\partial_{t_*}$ derivatives by $-i\sigma$. Once the existence of $\hat{\Box}_g(\sigma)^{-1}:\mathcal{Y}^{s',r'}\to\mathcal{X}^{s,r}$ is properly justified down to the real axis the solution in (\ref{eq:ForwardWaveSolution}) is a proper solution to 
\begin{align}
    \Box_g \phi =f, \; \text{ on the set } \{t_*\geq 0\}\cap \supp(h)\cap \supp(f),
\end{align}
where $f$ has compact support in $t_*$ and suitable decay $r^{-2-\alpha}$ where $\alpha\in (0,2)$.
In order to get to this forward solution we need to establish a lot of dynamical properties intrinsic to the Myers-Perry spacetime $(M,g)$ that are well known to Kerr but to our knowledge have not been generalized as of yet. This then leads to appropriate corresponding operator properties for the meromorphic family $\hat{\Box}_g(\sigma)$.\\
In general we follow the procedure outlined in \cite{hafner2021linear}. That is we start with a representation of the forward solution as 
\begin{align}\label{eq:ForwardWaveSolutionUpperPlane}
\phi(t_*,x) =\frac{1}{2\pi}\int_{\mathbb{R}+iC}e^{-i\sigma t_*} \hat{\Box}_g(\sigma)^{-1}\hat{\phi}(\sigma,x) \,d\sigma
\end{align}
with $C>>0$ a large constant. Equation (\ref{eq:ForwardWaveSolutionUpperPlane}) is well-defined as soon as one establishes polynomial high energy bounds for the resolvent,
\begin{align}
    \|\hat{\Box}_g(\sigma)^{-1}\|\lesssim \langle \sigma \rangle^N
\end{align}
in an upper half-plane strip for large enough $C>0$.
We then want to shift this contour integral downwards using the Payley-Wiener theorem. To be able to do this, we need to establish several properties for $\hat{\Box}_g(\sigma)$. 
\begin{enumerate}
    \item The operator $\hat{\Box}_g(\sigma):\mathcal{X}^{s,r}\to \mathcal{Y}^{s',r'}$ is Fredholm of index $0$ between suitable Sobelev spaces. For Myers-Perry, this will be dealt with in part II, especially Fredholmness in the case of the limiting absorption principle (non-zero $\sigma$).
     \item Mode stability for $\hat{\Box}_g(\sigma)$ i.e. that 
    \begin{align}
        \ker \hat{\Box}_g(\sigma)=\{0\},
    \end{align}
    on suitable function spaces. Combined with the statement that $\text{ind }\hat{\Box}(\sigma) =0$ for all $\sigma \in \overline{\mathbb{H}}$ this immediately leads to invertibility and existence of the resolvent $\hat{\Box}(\sigma)^{-1}$, as well as the
    For $\text{Im }\sigma>0$ this is exactly the content of 
    \cite{cvetivc2022mode}, however it has not been shown to hold yet down to the real axis $\text{Im }\sigma=0$. For the purposes of this paper we will take the latter as a mild conditional assumption, although a real mode stability result obtained by means analogous to \cite{shlapentokh2015quantitative},\cite{andersson2017mode} are expected to hold in the case of $1+4$ dimensional Myers-Perry as well.    
    \item Following existence of the resolvent $\hat{\Box}^{-1}_g(\sigma)$ we also need to establish sufficient control in the parameter $\sigma$ to justify the contourshift. We will thus need to show that $\hat{\Box}_g(\sigma)$ satisfies high energy estimates for energy strips above the real line, that is $\sigma\in \mathbb{C}$ with $\text{Re}(\sigma)>>1$ and $\text{Im}(\sigma)\geq 0$ bounded. For this a robust understanding of the trapped set of Myers-Perry is required, which will be the contents of part I. The high energy estimates will follow as a corollary of this in part II.
    \item Uniform Fredholm estimates for $\hat{\Box}_g(\sigma)$ down to $\sigma=0$. This will be done in part II as well.
    \item The points 1-4 above are enough to justify the contourshifting and thus give the forward solution (\ref{eq:ForwardWaveSolution}). However for asymptotic decay we will need to refine the regularity of $\hat{\Box}^{-1}_g(\sigma)$ in the spectral parameter near $\sigma=0$. This will be treated in part III, where we ultimately develop a resolvent expansion of $\hat{\Box}(\sigma)$ near $\sigma=0$. Appropriate control of the output of $\hat{\Box}_g^{-1}(\sigma)$ then gives Price's law for (\ref{eq:ForwardWaveSolution}).
\end{enumerate}

\begin{figure}[H]
    \centering
    \tikzset{every picture/.style={line width=0.75pt}} 
 \makebox[\linewidth][c]{%
        \hspace*{8mm}%
        \begin{tikzpicture}[x=0.75pt,y=0.75pt,yscale=-1,xscale=1]

\draw    (155.14,201.29) -- (495.14,200.29) ;
\draw [shift={(497.14,200.29)}, rotate = 179.83] [color={rgb, 255:red, 0; green, 0; blue, 0 }  ][line width=0.75]    (10.93,-3.29) .. controls (6.95,-1.4) and (3.31,-0.3) .. (0,0) .. controls (3.31,0.3) and (6.95,1.4) .. (10.93,3.29)   ;
\draw    (324.14,271.29) -- (321.17,28.29) ;
\draw [shift={(321.14,26.29)}, rotate = 89.3] [color={rgb, 255:red, 0; green, 0; blue, 0 }  ][line width=0.75]    (10.93,-3.29) .. controls (6.95,-1.4) and (3.31,-0.3) .. (0,0) .. controls (3.31,0.3) and (6.95,1.4) .. (10.93,3.29)   ;
\draw [color={rgb, 255:red, 255; green, 0; blue, 0 }  ,draw opacity=1 ][fill={rgb, 255:red, 179; green, 46; blue, 92 }  ,fill opacity=1 ][line width=0.75]    (210.14,77.29) -- (211.14,202.29) ;
\draw [color={rgb, 255:red, 255; green, 0; blue, 0 }  ,draw opacity=1 ][line width=0.75]    (210.14,77.29) -- (434.14,77.29) ;
\draw [color={rgb, 255:red, 255; green, 0; blue, 0 }  ,draw opacity=1 ][line width=0.75]    (434.14,77.29) -- (435.14,201.29) ;
\draw [color={rgb, 255:red, 255; green, 0; blue, 0 }  ,draw opacity=1 ][line width=0.75]    (211.14,201.29) -- (291.56,201.41) ;
\draw [color={rgb, 255:red, 255; green, 0; blue, 0 }  ,draw opacity=1 ][line width=0.75]    (355.74,201.4) -- (435.14,201.29) ;
\draw  [draw opacity=0][line width=0.75]  (291.56,202.41) .. controls (291.56,202.2) and (291.55,202) .. (291.55,201.79) .. controls (291.28,184.85) and (305.44,170.89) .. (323.16,170.62) .. controls (340.89,170.34) and (355.47,183.84) .. (355.74,200.78) .. controls (355.74,200.99) and (355.74,201.19) .. (355.74,201.4) -- (323.64,201.29) -- cycle ; \draw  [color={rgb, 255:red, 255; green, 0; blue, 0 }  ,draw opacity=1 ][line width=0.75]  (291.56,202.41) .. controls (291.56,202.2) and (291.55,202) .. (291.55,201.79) .. controls (291.28,184.85) and (305.44,170.89) .. (323.16,170.62) .. controls (340.89,170.34) and (355.47,183.84) .. (355.74,200.78) .. controls (355.74,200.99) and (355.74,201.19) .. (355.74,201.4) ;  
\draw  [color={rgb, 255:red, 74; green, 144; blue, 226 }  ,draw opacity=1 ][line width=1.5]  (319.5,201.29) .. controls (319.5,199) and (321.35,197.14) .. (323.64,197.14) .. controls (325.93,197.14) and (327.79,199) .. (327.79,201.29) .. controls (327.79,203.57) and (325.93,205.43) .. (323.64,205.43) .. controls (321.35,205.43) and (319.5,203.57) .. (319.5,201.29) -- cycle ;
\draw  [color={rgb, 255:red, 189; green, 16; blue, 224 }  ,draw opacity=1 ][dash pattern={on 0.84pt off 2.51pt}] (165.14,132.29) -- (259.14,132.29) -- (259.14,201.29) ;
\draw  [color={rgb, 255:red, 189; green, 16; blue, 224 }  ,draw opacity=1 ][dash pattern={on 0.84pt off 2.51pt}] (476.14,133.29) -- (382.14,133.29) -- (382.14,202.29) ;
\draw  [draw opacity=0][dash pattern={on 0.84pt off 2.51pt}] (449.18,78.45) .. controls (433.3,91.21) and (380.52,100.82) .. (317.79,101.18) .. controls (248.01,101.59) and (190.4,90.4) .. (182.18,75.58) -- (317.62,71.18) -- cycle ; \draw  [color={rgb, 255:red, 245; green, 166; blue, 35 }  ,draw opacity=1 ][dash pattern={on 0.84pt off 2.51pt}] (449.18,78.45) .. controls (433.3,91.21) and (380.52,100.82) .. (317.79,101.18) .. controls (248.01,101.59) and (190.4,90.4) .. (182.18,75.58) ;  
\draw [color={rgb, 255:red, 0; green, 35; blue, 255 }  ,draw opacity=1 ]   (153.14,86.79) -- (490.14,84.79) ;

\draw (328.93,15.4) node [anchor=north west][inner sep=0.75pt]  [xslant=-0.01]  {$Im\ \sigma $};
\draw (474,213.4) node [anchor=north west][inner sep=0.75pt]    {$Re\ \sigma $};
\draw (438,113.4) node [anchor=north west][inner sep=0.75pt]    {$Im\ \sigma \leq i\nu $};
\draw (331,212.4) node [anchor=north west][inner sep=0.75pt]    {$\sigma =0$};
\draw (430,60) node [anchor=north west][inner sep=0.75pt]    {$\frac{1}{2\pi}\int_{\mathbb{R}+iC}e^{-i\sigma t_*}\hat\Box ( \sigma )^{-1}\hat{f}(\sigma,x)d\sigma$};
\draw (323.64,89.19) node [anchor=north west][inner sep=0.75pt]    {$iC$};
\draw (200,50.4) node [anchor=north west][inner sep=0.75pt]    {$\|\hat{\Box}(\sigma)\|\leq\langle\sigma\rangle^{-N} $};

\end{tikzpicture}
        \hspace*{-8mm}%
    }

    \caption{The contour used in the Payley-Wiener theorem. Red lines prepresent the contour lines which for whose edges are taken to infinity, and around $\sigma=0$ the half circle is contracted. Polynomial bounds (orange) for high energy estimates. Around $\sigma=0$ the resolvent $\hat{\Box}_g(\sigma)^{-1}$ satisfies low energy bounds. The energy estimates on the strips $\text{Im }\sigma\leq i \nu$ (purple) allow for control of tails of the integral, and are a Corollary of Theorem \ref{thm:Theorem1}.}
    
\end{figure}
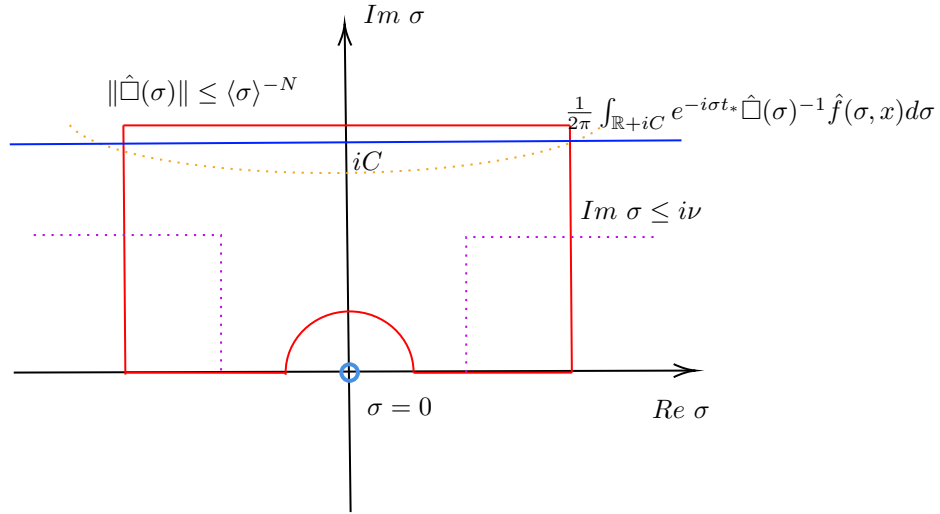

\tikzset{every picture/.style={line width=0.75pt}} 

In part I, we establish that the Hamiltonian flow associated to the Hamiltonian vector field $H_G$ of the dual metric function $g^{-1}$ exhibits normally hyperbolic trapping in the sense of \cite{wunsch2011resolvent},\cite{dyatlov2015asymptotics} \cite{dyatlov2016spectral}, and appendix $E$ of \cite{dyatlov2019mathematical}. Moreover, we will also show that $H_G$ admits a radial set at fibre infinity first shown for Kerr by Vasy \cite{vasy2013microlocal}. For both goals we will follow a more recent exposition given in Chapter $4$ of \cite{hintz2026non}. In short, we will show that the Hamiltonian flow associated to $g^{-1}$ via the vector field $H_G$ admits a trapped set $\Gamma$ inside the dual lightcone $\Sigma_G$, for the full subextremal range $(|\bfa|+|\bfb|)^2<\mu$ which is given as the intersection of two codimension $1$ manifolds inside $\Sigma_G$, in the domain of exterior communication of $M$, i.e. we have that
\begin{align}
    \Gamma = \Gamma^+\cap \Gamma^-.
\end{align}
Here $\Gamma^+$ and $\Gamma_-$ represent the stable and unstable manifolds for the flow. The argument presented in \cite{hintz2026non} and \cite{dyatlov2016spectral} is symplectic in nature for Kerr and it involves finding the correct constraints on $\Sigma_G$ and $\Gamma$. It basically employs the fact that for the exact Kerr metric, the dual metric function on $\Sigma_G$ can be separated as $G=\mathscr{G}+\mathscr{C}$ where $\mathscr{C}$ is the Carter constant restricted to $\Sigma_G$ and $\mathscr{G}$ captures the radial $(r,\xi)$ dynamics associated to $H_G$.\\
For Myers-Perry this argument is considerably more involved as in every part of the argument one has to content with the increased degrees of freedom. Obviously, the wave equation is known to be separable due to \cite{frolov2003particle}, \cite{frolov2007hidden} but instead of $4$ conserved quantities (including $\mathscr{C}$) associated to the flow $H_G$ one now has to deal with $5$ conserved quantities. 
Additionally, the fact that there are two spin momenta $(\bfa,\bfb)$ that need to be considered as well as the fact that the phase space $T^*M$ is $10$-dimensional instead of $8$, makes it hard to force through the same constraints. We will show however that once one finds an appropriate diagonalization of $g$ (following Wu \cite{wu2008separability}) one can uncover additional hidden constraints on $\Sigma_G$ and $\Gamma$ in the case of Myers-Perry allowing to push through analogous arguments to Kerr. For this we will also identify the appropriate formula for $\mathscr{C}$ restricted to $\Sigma_G$ that is analogous to the one used in \cite{dyatlov2015asymptotics}. We will also use $\mathscr{C}$ as a defining function for the horizon radial set for $H_G$ by establishing it represents a well defined positve definite quantity on the hypersurface $\{\sigma=0\}$ at $\mathcal{H}^+\subset T^*M$.\\
It remains to emphasize here that it has proven to be difficult up until now however to regularize the flow associated to $\mathscr{C}$ on $\mathbb{S}^3$. This obstructs the generalization of standard method used on Kerr to establish $r$-normally hyperbolicity for every $r\in (0,1)$ for the case of $1+4$ dimensional Myers-Perry. The higher dimensional nature of that problem has so far not allowed to push through the analogous arguments as in \cite{dyatlov2016spectral}. As a consequence of this, we emphasize that the Price law result in Theorem \ref{thm:MainPriceLaw} is not yet shown to hold for small perturbations of $g$ away from the standard Myers-Perry black hole $(\tilde{M},g)$ as the perturbative microlocal argument in Part II require this as the high energy input for the trapped set. \\
\\
In part II, we will bridge the dynamics established in part $I$ to eventually establish microlocal estimates associated to $\hat{\Box}_g(\sigma)$ as a meromorphic family of operators indexed by $\sigma\in \mathbb{C}$. We will work on the compactification of the spatial manifold $X$ at infinty and will develop $b$- and $sc$- pseudodifferential calculi \cite{melrose1993atiyah}\cite{melrose1994spectral} that are adapted to the boundary. We then translate the dynamics established in part I into semiclassical estimates, the horizon and trapped set estimate coming in from the interior $X^\circ$. For the radial point structure at infinity we will briefly recall the $b-sc$ calculus developed \cite{vasy2021limiting} to transition from semiclassical scattering estimates at infinity to semiclassical $b$-estimates, as the radial structure for $H_G$ in the pure $b$-calculus is not resolved. We will then use the $b-sc$ formalism to get high energy estimates of the type 
\begin{align}
    \|u\|_{\mathcal{X}^{s,r}}\lesssim |\sigma|^{-1+\delta}\|\hat{\Box}(\sigma)u\|_{\mathcal{Y}^{s',r'}}
\end{align}
For $\text{Re }\sigma>>1$ and $\text{Im } \sigma$ finite as well as non-elliptic Fredholm estimates of the type 
\begin{align}
    \|u\|_{\mathcal{X}^{s,r}} \lesssim \|\hat{\Box}(\sigma)u\|_{\mathcal{Y}^{s',r'}}+\|u\|_{\mathcal{Z}^{-N,-N}}
\end{align}
at finite energies $\sigma\in \mathbb{C}\backslash \{0\}$ as well as uniformily in the limit $\sigma \to 0$ where $\mathcal{Z}^{-N,-N}\hookrightarrow \mathcal{X}^{s,r}$ compactly. We will define the relevant Sobolev spaces where these estimates are made later on.\\
\\
In Part III, we will analyze the precise regularity output of $\hat{\Box}_g(\sigma)^{-1}$ in the uniform limit of $\sigma\to 0$ and build a resolvent expansion to eventually give asymptotic decay in $t_*$ for (\ref{eq:ForwardWaveSolution}). The main inspiration for this part is \cite{hintz2022sharp},\cite{hintz2023linear} as well as \cite{hintz2026non}. Other notable works have been based on this method as well, most notably \cite{millet2023optimal} for the Teukolsky equation and \cite{sussman2026complete} for larger polyhomogenous expansions for Sch\"odinger operators with short-ranged potentials.. The main task of this method relies on controlling the regularity of $\hat{\Box}_g(\sigma)^{-1}$ at $\sigma=0$. Higher regularity in this limit eventually implies faster temporal decay in $t_*$ for $\phi$. We note here that as in \cite{hintz2022sharp} the $t_*$ in the Fourier transform of (\ref{eq:ForwardWaveSolution}) is given by an appropriate timelike variable whose level sets are transverse to null infinity. As pointed out in \cite{vasy2021limiting}\cite{vasy2021resolvent} the limiting resolvent $\hat{\Box}(\sigma)^{-1}$ for $\sigma\neq 0$ produces outgoing solutions with respect to this time coordinate $t_*$. Concretely setting $t_*=t-r$ on $1+4$ dimensional Minkowski $\mathbb{M}^{1+4}$ solutions to (\ref{eq:ForwardWaveSolution}) satisfy time dependence $e^{-i\sigma t}$ and leading order spatial dependence $r^{-1}e^{i\sigma r}$, i.e. an overall factor $r^{-1}e^{-i\sigma t_*}$. The outgoing condition for $\hat{\Box}_g(\sigma)^{-1}$ is therefore nonoscillatory $\sigma$-independent $r^{-1}$ at infinity. We will show that the output is in fact conormal, i.e. upon repeated application of $r\partial_r$ and rotational vector fields the decay output remains $r^{-1}$ to leading order.\\
We then proceed to analyze the low energy resolvent $\hat{\Box}_g(\sigma)^{-1}$ following \cite{vasy2021limiting} uniformily in the limit of $\sigma 
\to 0$ on conormal spaces $\mathcal{A}^\alpha$ where $\alpha$ is a decay that is allowed to vary over the size of the indicial gap. The model operator $\Delta_{\mathbb{R}^4}$, which is invertible with a domain of functions decaying faster than $r^0$ and slower than $r^{-2}$. Overall in this case one expect invertibility 
\begin{align}\label{eq:InvertibleRange}
    \hat{\Box}_g(0)^{-1}:\mathcal{A}^{2+\alpha}\to \mathcal{A}^\alpha\;\text{ with } \alpha\in (0,2),
\end{align}
and hence uniform bounds for $\hat{\Box}_g(\sigma)^{-1}:\mathcal{A}^{2+\alpha}\to\mathcal{A}^\alpha$ for small $\sigma$. The basic mechanism underlying the proof of the main theorem follows an analogous algorithm to resolvent expansion laid out in \cite{hintz2022sharp}. In the current setting, up to relevant precision $\hat{\Box}_g(\sigma)$ is equal to
\begin{align}\label{eq:RelevantPrecision}
    L(\sigma)=-2i\sigma r^{-1}\left(r\partial_r+\frac{3}{2}\right)+ \Delta_{\mathbb{R}^4}.
\end{align}
for large $r$. The transport term in the above equation only gains one order of the decay $r^\alpha\mapsto r^{\alpha-1}$, whereas the model Laplacian $L(0)$ gets two. 
We expand the resolvent near zero frequency by writing 
\begin{align}\label{eq:ResExpac1}
    L(\sigma)^{-1}f&=L(0)^{-1}f+(L(\sigma)^{-1}-L(0)^{-1})f\nonumber\\
    &=u_0+\sigma L(\sigma)^{-1}f_1,
\end{align}
where we defined
\begin{align}\label{eq:ResExpac2}
    u_0=L(0)^{-1}f, \; f_1=-\sigma^{-1}(L(\sigma)-L(0))u_0.
\end{align}

In the above only $u_0$ is $\sigma$-independent whilst $L(0)^{-1}$ losses two orders of decay in $r$ and $\sigma^{-1}(L(\sigma)-L(0))$ gains an order, hence $f_1$ has one order of decay less than $f$. We then iterate
\begin{align}\label{eq:ResExpac3}
    u_k=L(0)^{-1}f_k \text{ and }  f_{k+1}=-\sigma^{-1}(L(\sigma)-L(0))u_k
\end{align}
as often as possible in the invertible range (\ref{eq:InvertibleRange}). The resolvent expansion becomes
\begin{align}\label{eq:ResExpac4}
    u_0+\sigma u_1 +...+\sigma^k u_k,
\end{align}
where the iteration stops once $f_k$ has $r^{-2}$ decay or less: the correction term that satisfies
\begin{align}
    \sigma^k L(\sigma)^{-1}f_k, \;\;\; f_k \gtrsim r^{-2},
\end{align}
in the expansion is then no longer bounded as $\sigma\to 0$, as $L(\sigma)^{-1}r^{-2}$ has a logarithmic singularity at $\sigma=0$ with an explicit coefficient, which in this case is $\sigma^k\log(\sigma+i0)$. Upon taking the Fourier transform in (\ref{eq:ForwardWaveSolution}) this term gives rise to a $t_*^{-k-1}$ term. The overal philosophy is the same as it has been for other spacetimes and equations analyzed in this matter: the long time asymptotics are yielded by the low energy $\sigma\to 0$ regime for the resolved operator $\hat{\Box}(\sigma)$, which are by Fourier correspondence in turn dominated by the asymptotically flat end of the spacetime. 
As it will turn out, our algorithm will hit the borderline decay term in $k=2$ steps iterating across the indicial gap giving exactly a $\sigma^2\log(\sigma + i0)$ and corresponding $t_*^{-3}$. On first sight the reader might ask why this is exactly the same as the $t_*^{-3}$ computed for Schwarzschild and Kerr in \cite{hintz2022sharp}. This similarity is deceptive however, as there are four key differences that make the resolvent expansion for $1+4$-dimensional Myers-Perry conceptually different than it was for its $1+3$ dimensional Kerr counterpart that depend fundamentally on the dimensionality of the problem.

\begin{enumerate}
    \item \underline{The increase of the indicial gap \& shift of the roots for the Mellin transformed Laplacian.} The indicial gap, obtained after Mellin transforming (\ref{eq:RelevantPrecision}) from $\rho=r^{-1}$ to $\xi$  and analyzing the first non-zero root in $\xi$ increases depending on dimension. In short they can be read off by the eigenvalues of the spherical Laplacian $\Delta_{\mathbb{S}^{d-1}}$. In the case of $4$ spatial dimensions these are given by $-\ell(\ell+2)$ for $\ell \in \mathbb{Z}$ increasing the gap for $\alpha$ in (\ref{eq:InvertibleRange}) from $\alpha \in (0,1)$ in the case of Kerr to $\alpha \in (0,2)$ in the case of Myers-Perry. Effectively increasing the dimensionality of the problem increase the number of steps one can use to run the resolvent expansion in equations (\ref{eq:ResExpac1})-(\ref{eq:ResExpac3}). The shifting of this roots to the Mellin transformed Laplacian also prevents the resolvent expansion in Myers-Perry also from having occurring a double root when contourshifting and solving for $u_1$ in the first iteration (which would lead to a $\sigma\log(\sigma+i0)$-term otherwise).
    \item \underline{Increased decay on the mass term.} As is clear from $1+d$-dimensional Schwarzschild-Tangherlini, the mass term in $\hat{\Box}_g(\sigma)$ scales as $\frac{2\mfrm}{r^{d-2}}$. For this reason it is also not occuring up to relevant precision in (\ref{eq:RelevantPrecision}) and it does not occur as a criterion for spectral admissibility in e.g. Lemma \ref{lem:AdmissibleWaveEquation}. In many instances of the resolvent expansion the higher decay on the mass term allows it to be absorbed away along with other terms in the asymptotic expansion of $u$.
    \item \underline{Fractional indicial roots in the transport term \& the weak Huygens principle.} The only way the mass term under 3. can still dominate a term in the resolvent expansion in (\ref{eq:ResExpac4}) is when the mass term in a lower decay term cancels out against a higher order decay term thereby shifting the decay by some integer value. In Kerr, where the transport term is given by $-2i\sigma r^{-1}\left(r\partial_r + 1\right)$ this indeed occurs, and the resolvent expansion for that reason can be carried out an extra time then what one expects based on the size of the indicial gap $\alpha \in (0,1)$ in 1. alone. In Myers-Perry the indicial root of the transport term in (\ref{eq:RelevantPrecision}) is $-i\frac{3}{2}$ and because the first order mass outputs of $\hat{\Box}(\sigma)^{-1}$ are always integer valued they can never cancel out against the fractional indicial roots of the transport term. The indicial root is only determined by the asymptotically flat regime (and thus the dimension $d$) of the spacetime, and hence if the mass term is of higher decay than the $\frac{3}{2}$ treshold, the late-time tail is up to first order determined by the Minkowski limit of the wave operator. From the perspective of the resolvent expansion one can interpret this as an instance the weak Huygens principle for $1+d$-dimensional spacetimes where $d$ is even.
\end{enumerate}

Point $1$,$2$ and $3$ essentially allow the resolvent expansion in $1+4$-dimensional spacetimes to be executed exactly over the size of the indicial gap $\alpha\in (0,2)$ with as starting point initial date $f\in \mathcal{A}^{2+\alpha}$. The shift of the poles of the Mellin-transformed Laplacian, combined with a higher decay mass term and fractional indicial roots of the transport term cause no log terms to appear any earlier then when one hits $\alpha =0$ (i.e. $r^{-2}$-decay in the output). Thus for completely seperate reasons the singular term one lands on in this case $\sigma^2\log(\sigma+i0)$ as it is with Kerr. We emphasize that the resulting $t_*^{-3}$ is in fact the decay rate that one would obtain by plainly analyzing $1+4$-dimensional Minkowski spacetime. Our results are therefore completely in agreement with general Price law results obtained with similar methods for odd dimensional asymptotically flat spacetimes, c.f. Proposition $6.4$ in \cite{hintz2023linear}, or earlier work as in e.g. \cite{baskin2015asymptotics}. The generic result for linear waves on odd-dimensional Minkowski dates back to H\"ormander \cite{hormander1997lectures}, Section 6.2. In \cite{luk2024late}, where the authors consider the Price law on all even dimensional Schwarzschild-Tangerherlini black holes based on the interaction of high radiation fields with the strong Huygens principle, the authors remark regarding the weak Huygens principle:\\

"The wave equation in even space-time dimension $d$ falls into the generic case, where already on Minkowski space the $l$-th spherical harmonic generally attains an $\mathcal{O}(t^{-1+\sqrt{1+4\alpha}})=\mathcal{O}(t^{-(d-1+2l)})$ tail."\\
\\
Thus, Theorem \ref{thm:MainPriceLaw} is ultimately a pure reflection of the asymptotic nature of the $1+4$-dimensional Myers-Perry spacetime. However, the only way that the contour shifting of (\ref{eq:ForwardWaveSolution}) is allowed in the first place is by a rigorous global analysis of the trapped set at high energies in the full extremal range. The Price-law decay of $t_*^{-3}$ should therefore be interpreted as a corollary of the fact that Myers-Perry that this trapped set persists for any $(\bfa,\bfb)$ satistfying $(|\bfa|+|\bfb|)^2< 2\mfrm$, the key result of part I.

\subsection*{Acknowledgments}

The author is deeply grateful to his PhD advisor, Peter Hintz, for suggesting Myers-Perry black holes and resolvent expansions as the overarching theme of his doctoral research, for his guidance and advice throughout this project, and for sharing an unpublished manuscript of \cite{hintz2026non}. The author also thanks Bernard Whiting, Rita Teixeira da Costa, and Christoph Kehle for useful discussions on real mode stability, and Pascal Millet for discussions on the intricacies of the resolvent expansion.

\subsection*{AI statement}

The proofs in Parts I–III were developed by the author,  during the last three years of graduate research. The results on full hyperbolic trapping and Price’s law were first presented at GR24 \& Amaldi 16 in Glasgow on 16 July 2025. ChatGPT 6.0 Astra was used during the final two weeks of preparation for proofreading and typographical corrections. The author takes full responsibility for the manuscript.

\section{$1+4$ dimensional Myers-Perry spacetimes}\label{sec:Myers-Perry spacetimes}
\subsection{The Myers-Perry metric in Boyer-Lindquist coordinates.}\label{sec:MyersPerryBoyerLindquist}
We define the Myers-Perry outer region $(\tilde{M},g_{\bfa,\bfb,\mfrm})$ as a solution to the vacuum Einstein equations in $1+4$ dimenions
\begin{equation} \label{eq:EinsteinFieldEquation}
    \ric(g) = 0
\end{equation}
where $\tilde{M}\subseteq \mathbb{R}^5$ is given by 
\begin{equation}
    M = \mathbb{R}_t \times (r_+,\infty)_r \times \mathbb{S}^3_{\phi,\psi,\theta},
\end{equation}
and the metric given by $g_{\bfa,\bfb,\mfrm} \in T^*\tilde{M} \otimes T^*\tilde{M}$ as \cite{horowitz2012black}
\begin{align}\label{eq:MyersOriginalMetric}
    g_{\bfa,\bfb,\mfrm}= -dt^2 +\frac{2\mfrm}{\bfr^2}&(dt-\bfa\sin^2\theta\, d\phi-\bfb\cos^2\theta\, d\psi)^2
    \nonumber\\
    &+\frac{\bfr^2}{\Delta_r}dr^2+\bfr^2d\theta^2+(r^2+\bfa^2)\sin^2\theta\, d\phi^2+(r^2+\bfb^2)\cos^2\theta\,d\phi^2.
\end{align}
This is the form for the $1+4$-dimenional Myers-Perry metric with $(-,+,+,+,+)$ signature in the domain of exterior communication, given in Boyer-Lindquist coordinates. It represents a spinning black hole across two orthogonal axes, with rotational angles $\phi$ and $\psi$. It has Killing vector fields $\partial_\phi,\partial_\psi$ and $\partial_t$, making $(\tilde{M},g_{\bfa,\bfb,\mfrm})$ a stationary solution to (\ref{eq:EinsteinFieldEquation}). The variables $(\bfa,\bfb)\in \mathbb{R}^2$ are spin parameters for each axis of rotation, $\mfrm>0$ is the mass of the black hole, and $\bfr^2, \Delta_r$  are given below\footnote{Throughout this article $\bfr^2$ will assume a similar role to the well known quantity $r^2_{\bfa}=r^2+\bfa^2\cos^2\theta$ for the standard Kerr metric.}:
\begin{align}
    \bfDelta &= \frac{(r^2 + \bfa^2)(r^2 +\bfb^2)}{r^2} - 2 \mfrm \\
    \bfr^2 &= r^2 + \bfa^2 \cos^2\theta + \bfb^2 \sin^2\theta.
\end{align}
We note that $\mathbb{S}^3$ is parametrized using Hopf coordinates for symmetry reasons, stemming from the pullback of the Hopf-fibration $\mathbb{S}^3\hookrightarrow\mathbb{S}^2$. Concretely this comes down to a parametrization in two orthogonal planes parametrized in polar angles $\phi,\psi \in [0,2\pi)$ combined with an azimuthal angle $\theta\in [0,\frac{\pi}{2})$ given as\footnote{When embedding $\mathbb{S}^3$ into $\mathbb{C}^2$ equation (\ref{eq:HopfParametrization}) expresses nicely as $(z_1,z_2)=(e^{i\phi}\sin\theta,e^{i\psi}\cos\theta)$. } 
\begin{equation}\label{eq:HopfParametrization}
    \begin{pmatrix}
        \cos(\phi)\sin(\theta)\\
        \sin(\phi)\sin(\theta)\\
        \cos(\psi)\cos(\theta)\\
        \sin(\psi)\cos(\theta)\\
    \end{pmatrix}
\end{equation}
with volume form $dg_{\mathbb{S}^3}= \sin^2\theta\, d\phi^2+\cos^2\theta\, d\psi^2+d\theta^2$. Using this it is then easy to check that the metric in (\ref{eq:MyersOriginalMetric}) is indeed asymptotically flat, we have 
\begin{align}\label{eq:MinkowskiMetric}
    \lim_{r\to \infty} g_{\bfa,\bfb,\mfrm}:= g_{\mathbb{M}^{1+4}}=-dt^2 + dr^2 +r^2(d\theta^2+ \sin^2\theta\,d\phi^2+\cos^2\theta\,d\psi^2),
\end{align}
which is indeed the Minkowski metric on $\mathbb{R}^{1+4}$ written in Hopf coordinates.
We wish to emphasize here that in (\ref{eq:MyersOriginalMetric}) we used the opposite sign convention for $(\bfa,\bfb)$, with respect to \cite{horowitz2012black} to make our notation consistent with later chapters, i.e. we mapped $(\bfa,\bfb)\mapsto(-\bfa,-\bfb)$ . This is \textit{not} a discrete symmetry of the metric but it does land us in the same class of spacetimes. With regard to symmetries, within Hopf-coordinates a true discrete symmetry of the spacetime $(M,g_{\bfa,\bfb,\mfrm})$ is given by
\begin{align}
    (\bfa,\bfb,\phi,\psi,\theta)\mapsto (\bfb,\bfa,\psi,\phi,\frac{\pi}{2}-\theta).
\end{align}
We mention here as well that $(M,g_{\bfa,\bfb,\mfrm})$ contains several subclasses of blackholes, including the biaxial case where $\bfa=\bfb$, the single-spinning case, where either $\bfa=0$ or $\bfb=0$, the counterrotating case where $\bfa=-\bfb$, and of course the Schwarzschild-Tangherlini case, which is given by $\bfa=\bfb=0$ i.e.
\begin{align}
    g_{0,0,\mfrm}&=\lim_{\bfa,\bfb\to 0} g_{\bfa,\bfb,\mfrm}\\
    &=-\left(1-\frac{2\mfrm}{r^2}\right)dt^2+\left(1-\frac{2\mfrm}{r^2}\right)^{-1}dr^2 +d\theta^2+r^2(\sin^2\theta\,d\phi^2+\cos^2\theta\,d\psi^2).
\end{align}

Boyer-Lindquist coordinates on $1+4$-dimensional Myers-Perry degenerate in two different ways. First of all, at the poles $\theta\in\{0,\frac{\pi}{2}\}$ the Hopf-coordinates for $\mathbb{S}^3$ collapse into two disjoint rings, which can be dealt with by defining local smooth coordinates. More importantly, (\ref{eq:MyersOriginalMetric}) degenerates at the hypersurfaces $\{r=r_+\}$ $\{r=r_-\}$ given by the roots of the quartic polynomial $\Delta_r$ in a way reminiscent of Kerr, we have that
\begin{align}
    r^2\Delta_r = (r^2-r_+^2)(r^2-r_-^2),
\end{align}
with $r_+$ and $r_-$ usually given in quadratic form through 
\begin{align}
    2 r_+^2 &= 2\mfrm - \bfa^2 - \bfb^2 +\sqrt{(2\mfrm - \bfa^2 - \bfb^2)^2 - 4 \bfa^2 \bfb^2} \label{eq:EventHorizonMP} \\
    2 r_-^2 &= 2\mfrm - \bfa^2 - \bfb^2 - \sqrt{(2\mfrm - \bfa^2 - \bfb^2)^2 - 4 \bfa^2 \bfb^2} \label{eq:CauchyHorizonMP},
\end{align}
but it should noted that all $4$ roots to $\Delta_r$ can be recast in the elegant form \cite{an2018five}
\begin{align}
    r=\pm\frac{1}{2}\sqrt{2\mfrm - (\bfa+\bfb)^2}\pm \frac{1}{2}\sqrt{2\mfrm -(\bfa-\bfb)^2}.
\end{align}
We will define the event horizon $r_+$ to be the positive square root of (\ref{eq:EventHorizonMP}) and the Cauchy horizon to be the positive square root of (\ref{eq:CauchyHorizonMP}). To ensure that $r_+$ in (\ref{eq:EventHorizonMP}) and $r_-$ in (\ref{eq:CauchyHorizonMP}) are real valued quantities we have to impose an extremality constraint given by 
\begin{align}
    2\mfrm-\bfa^2-\bfb^2\geq 4\bfa^2\bfb^2,
\end{align}
whence 
\begin{align}
    2\mfrm \geq \bfa^2+\bfb^2+2|\bfa\bfb|=(|\bfa|+|\bfb|)^2\label{eq:ExtremalityMP}.
\end{align}
Note that if (\ref{eq:ExtremalityMP}) is an equality, we have that $r_+=r_-$, i.e. two horizon radii coincide. We know very well from Kerr that when the Cauchy and event horizons overlap the black hole becomes \textit{extremal} exactly when $\bfa=\mfrm$ before transitioning into a spacetime with a naked singularity for $\bfa>\mfrm$. For Myers-Perry this carries over analogously for the double spinning case$\bfa\neq 0$, $\bfb\neq 0$, with a word of caution pertaining to the single spinning cases. Indeed, the surface gravity $\kappa$ and black hole temperature $T_g$ associated to (\ref{eq:MyersOriginalMetric}) are related via
\begin{align}
\kappa= 2\pi T_g,   
\end{align}
 where $\kappa$ is defined through 
 \begin{align}
     \nabla_\chi \chi_+=\kappa\chi_+,
 \end{align}
 with $\chi_+$ the horizon Killing vector field (see (\ref{eq:HorizonKillingField})).
  In the case of $1+4$-dimensional Myers-Perry their explicit expressions are given by \cite{lopez2014black}, \cite{an2018five}, \cite{horowitz2012black}
\begin{align}\label{eq:SurfaceGravityIntro}
    \kappa &= \frac{2r_+^2+\bfa^2+\bfb^2-2\mfrm}{2\mfrm r_+},
\end{align}
which exactly vanishes when $2\mfrm=(|\bfa|+|\bfb|)^2$, i.e. when $r_+^2=|\bfa\bfb|$, which as mentioned previously happens when
\begin{align}
    \Delta_r =\Delta_r' =0.
\end{align}

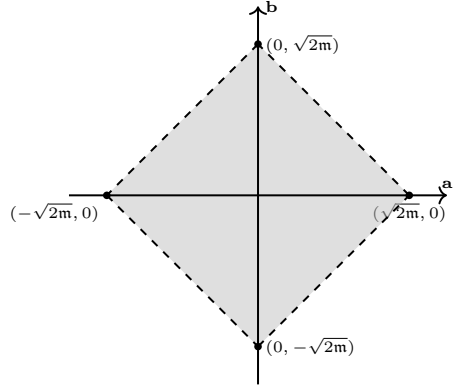
\begin{figure}[H]\label{fig:MyersPerryDiamond}
\begin{center}
\begin{tikzpicture}[scale=0.5]
\node (I)    at (4,0)   {};
\tiny
\path 
   (I) +(90:4)  coordinate[label=0:{$(0,\sqrt{2\mfrm})$}] (Itop)
       +(-90:4) coordinate[label=0:{$(0,-\sqrt{2\mfrm})$}] (Ibot)
       +(180:4) coordinate[label=225:{$(-\sqrt{2\mfrm},0)$}] (Ileft) 
       +(0:4)   coordinate[label=-90:{$(\sqrt{2\mfrm},0)$}] (Iright)
       
       ;
\node (J)    at (4,0)   {}; 

\path 

 (J) +(90:5)  coordinate[label=0:{$\bfb$}] (Jtop)
       +(-90:5) coordinate (Jbot)
       +(180:5) coordinate (Jleft) 
       +(0:5)   coordinate[label=90:{$\bfa$}] (Jright)
;

\path[fill=gray!50,opacity=.5] (Ileft) to (Ibot) to (Iright) to (Itop);

\path [fill=black, fill opacity = 1] (Ileft) circle (0.1);
\path [fill=black, fill opacity = 1] (Iright) circle (0.1);
\path [fill=black, fill opacity = 1] (Itop) circle (0.1);
\path [fill=black, fill opacity = 1] (Ibot) circle (0.1);

\draw[dashed] (Iright) -- (Itop);
\draw [dashed] (Itop) --  (Ileft);
\draw [dashed] (Ileft) -- (Ibot);
\draw[dashed] (Ibot) -- (Iright);

\draw[->] (Jbot) -- (Jtop);
\draw [->] (Jleft) -- (Jright);

\end{tikzpicture}
\caption{The subextremal parameter diamond $D_{\mfrm}\subset \mathbb{R}^2$. Its vertices represent the extremal vanishing horizon limit, whilst the dotted lines represent the regular extremal limit.}
\end{center}
\end{figure}

It should be noted that one can use 
\begin{align}
    \Delta_r|_{r=r_+}=0 \longleftrightarrow \frac{\bfa^2\bfb^2}{r_+^2}=(2\mfrm - r_+^2-\bfa^2-\bfb^2)
\end{align}
to show 
\begin{align}
    \Delta'_r|_{r=r_+} =\frac{2}{r_+}(2r_+^2+\bfa^2+\bfb^2-2\mfrm)
\end{align}
and hence that 
\begin{align}
    \kappa = \frac{\Delta'_r|_{r=r_+}}{4\mfrm}.
\end{align}
When discussing extremality in the single spinning case, without loss of generality $\bfa\neq 0$ $\bfb=0$, the discussion is slightly more nuanced as in this case we have that 
\begin{align}
    r_+^2=2\mfrm-\bfa^2, \qquad r_-^2=0,
\end{align}
and the Myers-Perry black hole does not exhibit a regular Cauchy horizon as it is moved to the singular locus $\{r=0\}$. The equality $r_+^2=r_-^2$ persists at the vertices of the diamond in Figure \ref{fig:MyersPerryDiamond} however and in the same sense the surface gravity still tends to zero in the limits $(\bfa,0)\to (\pm\sqrt{2\mfrm},0),\,(0,\bfb)\to (0,\pm\sqrt{2\mfrm})$ in (\ref{eq:SurfaceGravityIntro}). In this so-called extremal vanishing horizon (EVH) limit limit \cite{di2022black} however, we have as well that $r_+\to 0$, so the event horizon also moves to the singular locus, and $\kappa$ is no longer associated with the surface gravity of a regular event horizon, but should be interpreted as a singular limit with vanishing horizon area.  
The observant reader will also immediately notice that for extremal Myers-Perry $r_+=|\bfa\bfb|$ is given independent of the mass $\mfrm>0$, unlike in Kerr, except of course having an implicit dependence through the equality in (\ref{eq:ExtremalityMP}).
The considerations above lead us to a straightforward notion of subextremality, based on the notion of vanishing surface gravity in case of both the regular extremal limit as well as the EHV limit, see also \cite{emparan2008black}. We say that $(M,g_{\bfa,\bfb,\mfrm})$ is \textit{subextremal} if the spin parameters $(\bfa,\bfb)\in \mathbb{R}^2$ lie within the parameter diamond 
\begin{align}\label{eq:ParameterDiamond}
    D_{\mfrm}=\{(\bfa,\bfb)\in\mathbb{R}^2\:|\:(|\bfa|+|\bfb|)^2 < 2\mfrm\},
\end{align}
i.e. when the inequality (\ref{eq:ExtremalityMP}) is strict. When this holds the ''best'' lower bound for the magnitude of $r_+^2$ is given by its limiting value in the extremal case,
\begin{align}\label{eq:BestPossibleBound}
    r_+^2>\frac{r_+^2+r_-^2}{2}=\mfrm-\frac{1}{2}(\bfa^2+\bfb^2)\geq |\bfa\bfb|,
\end{align}
which is a fact that we will use a lot in the chapters to come. Another important relation that can be gained directly from (\ref{eq:EventHorizonMP}) and (\ref{eq:CauchyHorizonMP}) is that 
\begin{align}\label{eq:VietaRelation}
    r_+^2r_-^2= \bfa^2\bfb^2.
\end{align}

We note that throughout this paper we will only ever assume that $(\bfa,\bfb)\in D_{\mfrm}$, hence our treatment in particular encompasses all the special cases mentioned before (e.g. single-spinning, biaxially symmetric, or Schwarzschild-Tangherlini) within the subextremal diamond.

It is worth giving a remark on the ergo-region following \cite{horowitz2012black} here as well. Within our current signature convention, we say that a vector $X_p\in T_p\tilde M$ for $p\in \tilde M$ is \textit{timelike} when $g_{\bfa,\bfb,\mfrm}(X_p,X_p)<0$ and \textit{spacelike} when $g_{\bfa,\bfb,\mfrm}(X_p,X_p)>0$. The ergo-region $\mathcal{E}_{MP}\subset \tilde M$ is then logically the region where $\partial_t$ fails to be negative, i.e. 
\begin{align}\label{eq:ErgoRegionDef}
    \mathcal{E}_{MP}:=\{(t,r,\phi,\psi,\theta)\in \mathbb{R}_t\times (r_+,\infty)_r\times\mathbb{S}^3_{\phi,\psi,\theta}\:|\: g_{\mfrm,\bfa,\bfb}(\partial_t,\partial_t)\geq0 \}.
\end{align}
By direct evaluation in (\ref{eq:MyersOriginalMetric}) we see that 
\begin{equation}\label{eq:PartialtEval}
    g_{\mfrm,\bfa,\bfb}(\partial_t,\partial_t) = -\left(1-\frac{2\mfrm}{\bfr^2}\right).
\end{equation}
allowing us to rewrite (\ref{eq:ErgoRegionDef}) as
\begin{align}\label{eq:ErgoRegionDef2}
    \mathcal{E}_{MP}:=\{(t,r,\phi,\psi,\theta)\in \mathbb{R}_t\times (r_+,\infty)_r\times\mathbb{S}^3_{\phi,\psi,\theta}\:|\: \bfr^2\leq 2\mfrm\}.
\end{align}
Therefore plotting the outer boundary of the ergoregion as all $(r_\mathcal{E}(\theta),\theta)$ for which $\bfr^2=2\mfrm$ we see that 
\begin{align}\label{eq:ErgoRegionvsEventHorizon}
    r^2_{\mathcal{E}}(\theta) = 2\mfrm-\bfa^2\cos^2\theta-\bfb^2\sin^2\theta \geq 2\mfrm-\bfa^2-\bfb^2\geq r_+^2.
\end{align}
as long as we assume $(\bfa,\bfb)\in D_{\mfrm}$ and used (\ref{eq:EventHorizonMP}). We emphasize here that the first inequality in (\ref{eq:ErgoRegionvsEventHorizon}) is strict when $\bfa$ and $\bfb$ are non-vanishing, in which case the outer boundary of the ergoregion and the event horizon do not touch anywhere, in contrast to Kerr.
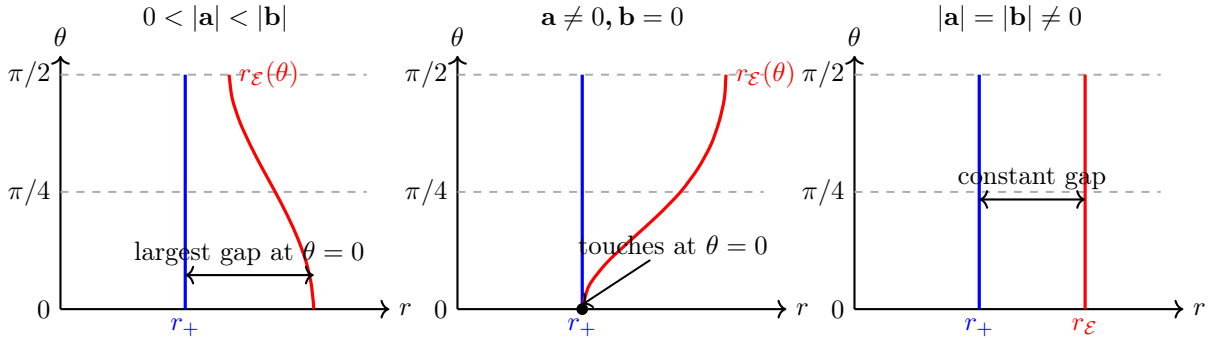
\begin{figure}[H]
\centering
\begin{tikzpicture}[scale=1.0]

\newcommand{\axes}[3]{%
\begin{scope}[shift={(#1,#2)}]
  \node[font=\bfseries,align=center] at (2.05,3.85) {#3};
  \draw[->,thick] (0,0) -- (4.35,0) node[right] {$r$};
  \draw[->,thick] (0,0) -- (0,3.35) node[above] {$\theta$};
  \node[left] at (0,0) {$0$};
  \node[left] at (0,1.55) {$\pi/4$};
  \node[left] at (0,3.10) {$\pi/2$};
  \draw[dashed,gray!65] (0,1.55) -- (4.05,1.55);
  \draw[dashed,gray!65] (0,3.10) -- (4.05,3.10);
\end{scope}}

\axes{0}{0}{$0<|\bfa|<|\bfb|$}
\begin{scope}[shift={(0,0)}]
  \def\rplus{1.65}
  \draw[blue,very thick] (\rplus,0) -- (\rplus,3.10);
  \node[blue,below] at (\rplus,0) {$r_+$};
  \draw[red,very thick,smooth] plot coordinates {
    (3.35,0.000) (3.32,0.310) (3.25,0.620) (3.14,0.930)
    (3.00,1.240) (2.84,1.550) (2.67,1.860) (2.51,2.170)
    (2.37,2.480) (2.27,2.790) (2.23,3.100)
  };
  \node[red,right] at (2.23,3.10) {$r_{\mathcal{E}}(\theta)$};
  \draw[<->,thick] (\rplus,0.45) -- (3.35,0.45);
  \node[above] at (2.50,0.45) {largest gap at $\theta=0$};
\end{scope}

\axes{5.25}{0}{$\bfa\neq0$,\,$\bfb=0$}
\begin{scope}[shift={(5.25,0)}]
  \def\rplus{1.65}
  \draw[blue,very thick] (\rplus,0) -- (\rplus,3.10);
  \node[blue,below] at (\rplus,0) {$r_+$};
  \draw[red,very thick,smooth] plot coordinates {
    (1.65,0.000) (1.74,0.310) (1.98,0.620) (2.31,0.930)
    (2.65,1.240) (2.95,1.550) (3.18,1.860) (3.35,2.170)
    (3.46,2.480) (3.53,2.790) (3.55,3.100)
  };
  \fill[black] (\rplus,0) circle (2.2pt);
  \draw[->] (2.55,0.65) -- (\rplus,0.06);
  \node[align=center] at (2.85,0.85) {touches at $\theta=0$};
  \node[red,right] at (3.55,3.10) {$r_{\mathcal{E}}(\theta)$};
\end{scope}

\axes{10.5}{0}{$|\bfa|=|\bfb|\neq0$}
\begin{scope}[shift={(10.5,0)}]
  \def\rplus{1.65}
  \def\re{3.05}
  \draw[blue,very thick] (\rplus,0) -- (\rplus,3.10);
  \draw[red,very thick] (\re,0) -- (\re,3.10);
  \node[blue,below] at (\rplus,0) {$r_+$};
  \node[red,below] at (\re,0) {$r_\mathcal{E}$};
  \draw[<->,thick] (\rplus,1.45) -- (\re,1.45);
  \node[above] at (2.35,1.45) {constant gap};
\end{scope}

\end{tikzpicture}
\caption{
Blue: event horizon $r=r_+$, independent of $\theta$.\quad
Red: ergosurface $r=r_\mathcal{E}(\theta)$.\quad
For $0<|\bfa|<|\bfb|$, $r_{\mathcal{E}}(\theta)$ decreases monotonically from $\theta=0$ to $\pi/2$, so the largest gap is at $\theta=0$.}
\end{figure}

\subsection{Dual metric, wave equation}

Whilst the form for $g_{\bfa,\bfb,\mfrm}$ given by (\ref{eq:MyersOriginalMetric}) is most commonly used in the literature, it does not represent the metric with an orthonormal frame. In order to diagonalize (\ref{eq:MyersOriginalMetric})  we first define the function $\bfp^2$ through
\begin{align}
    \bfr^2=r^2+\bfp^2, \;\text{ i.e.} \;\;\;\; \bfp^2=\bfa^2\cos^2\theta+\bfb^2\sin^2\theta,
\end{align}
which for $\theta \in [0,\frac{\pi}{2})$, $\bfp$ gives the norm of the first quadrant of an ellipsoid with semi-axes $\bfa, \bfb$ traced out by the function $\theta\mapsto (\bfa\cos\theta,\bfb\sin\theta)$.\\
Using $\bfp(\theta)$ one can then recast (\ref{eq:MyersOriginalMetric}) an extremely elegant form, first proposed by Wu in \cite{wu2008separability}, whom we will follow for the first part of this section; 

\begin{align}
    g_{\bfa,\bfb,\mfrm} &= -\frac{\bfDelta}{\bfr^2}\left(dt-\bfa\sin^2\theta d\phi -\bfb\cos^2\theta d\psi\right)^2 + \frac{\bfr^2}{\bfDelta}dr^2 + \bfr^2 d\theta^2 \nonumber \\
    &+\frac{\sin^2\theta \cos^2\theta}{\bfp^2 \bfr^2}\left[(\bfb^2-\bfa^2)dt + (r^2 + \bfa^2)\bfa d\phi - (r^2 + \bfb^2)\bfb d\psi \right]^2 \nonumber\\
    &+\frac{1}{r^2 \bfp^2}\left(-\bfa \bfb dt + (r^2 +\bfa^2)\bfb\sin^2\theta d\phi + (r^2 + \bfb^2) \bfa \cos^2\theta d\psi \right)^2 \label{eq:5DMPMetric}.
\end{align}

A straightforward but tedious computation shows that (\ref{eq:5DMPMetric}) is indeed equivalent to (\ref{eq:MyersOriginalMetric}) and it will be the form of the metric that we will use for $(M,g_{\bfa,\bfb,\mfrm})$ from now on. We remark here that it is indeed because of \cite{wu2008separability} that we use the opposite convention for the sign of $(\bfa,\bfb)$ as Myers himself uses in \cite{horowitz2012black}. \\
The form of (\ref{eq:5DMPMetric}) is indeed diagonal, with a co-frame spanned by co-vectors 
\begin{align}
    dT&:= dt-\bfa\sin^2\theta d\phi -\bfb\cos^2\theta d\psi\label{eq:CoframedT}\\
    d\Phi&:= (\bfb^2-\bfa^2)dt + (r^2 + \bfa^2)\bfa \: d\phi - (r^2 + \bfb^2)\bfb \: d\psi \label{eq:CoframedPhi}\\
    d\Psi &:= -\bfa \bfb dt + (r^2 +\bfa^2)\bfb\sin^2\theta d\phi + (r^2 + \bfb^2) \bfa \cos^2\theta d\psi \label{eq:CoframedPsi}
\end{align}

It can easily be used to compute the metric determinant, to find 
\begin{align}\label{eq:MetricDeterminant}
    \det (-g)=  r^2 \bfr^4 \sin^2\theta \cos^2\theta,
\end{align}
which is non-degenerate over $r_+$, unlike in the Kerr case. Its associated spacetime density is given by
\begin{align}\label{eq:MPSpacetimeMetric}
    dg_{\bfa,\bfb,\mfrm} = r \bfr^2 dt dr dg_{\mathbb{S}^3} \text{ with } dg_{\mathbb{S}^3}=\cos\theta\sin\theta d\phi d\psi d\theta,
\end{align}
with $dg_{\mathbb{S}^3}$ given by Hopf coordinates on $\mathbb{S}^3$. Furthermore, we with this we also introduce the stationary spatial density $dg_{\bfa,\bfb,\mfrm}|_X$ through
\begin{align}\label{eq:StationarySplitting}
    dg_{\bfa,\bfb,\mfrm}= dt_* dg_{\bfa,\bfb,\mfrm}|_X, \text{ with } dg_{\bfa,\bfb,\mfrm}|_X = r\bfr^2 dr dg_{\mathbb{S}^3}.
\end{align}

We note that in the Schwarzschildt-Tangerherlini limit $(\bfa,\bfb)\to (0,0)$ this again reduces to the metric determinant for $1+4$ Minkowski with Hopf coordinates for $\mathbb{S}^3$,
\begin{align}
     \lim_{(\bfa,\bfb)\to (0,0)}\det (\sqrt{-g})=  r^4 \sin^2\theta \cos^2\theta.
\end{align}
Inverting the frame in (\ref{eq:5DMPMetric}) we then get to a definition of the dual metric tensor $g^{-1}_{\bfa,\bfb,\mfrm}\in TM\otimes TM$, which we will rescale by $\bfr^2$:

\begin{align}
    \bfr^2 \cdot g^{-1}_{\bfa,\bfb,\mfrm} &:= -\frac{1}{r^4\bfDelta}\left[(r^2+\bfa^2)(r^2 +\bfb^2) \partial_t + \bfa(r^2 + \bfb^2) \partial_\phi + \bfb(r^2 + \bfa^2)\partial_\psi \right]^2 + \bfDelta \partial_r^2 + \partial^2_\theta \nonumber\\
    &+ \frac{\sin^2\theta \cos^2\theta}{\bfp^2} \left[ (\bfa^2-\bfb^2)\partial_t + \frac{\bfa}{\sin^2\theta}\partial_\phi - \frac{\bfb}{\cos^2\theta}\partial_\psi \right]^2 + \frac{\bfr^2}{r^2\bfp^2}[\bfa\bfb \partial_t + \bfb\partial_\phi+ \bfa \partial_\psi]^2 \label{eq:DualMPMetric}.  
\end{align}
A similar dual metric tensor was also reported by Frolov and Stojkovic \cite{frolov2003particle} with oppsite spin convention, however these authors did not have access to the orthogonal frame that Wu proposed in \cite{wu2008separability} and computed their dual metric using GRTensor. \\
With these computations we can write down the wave operator, defined as the Laplace-Beltrami operator via
\begin{align}\label{eq:LaplaceBeltrami}
    \Box_{g} u=-\frac{1}{\sqrt{-\det(g)}}\partial_\mu \left(\sqrt{-\det(g)}g^{\mu\nu} \partial_\nu \right)u,\;\;\; u\in C^\infty(M),
\end{align}
explicitly as
\begin{align}\label{eq:FullWaveEquation}
    \bfr^2\Box_{g} u&=-\left(-\frac{(r^2\Delta_r+2\mfrm r^2)^2}{r^4\Delta_r} \left[\partial_t+\frac{\bfa}{r^2+\bfa^2}\partial_\phi +\frac{\bfb}{r^2+\bfb^2}\partial_\psi\right]^2 +\frac{1}{r}\partial_r(r\Delta_r\partial_r)+\frac{\partial_\theta(\cos\theta\sin\theta \partial_\theta)}{\cos\theta\sin\theta} \nonumber\right.\\
    &+\left.\frac{\cos^2\theta\sin^2\theta}{\bfp^2}\left[(\bfa^2-\bfb^2)\partial_t +\frac{\bfa}{\sin^2\theta}\partial_\phi -\frac{\bfb}{\cos^2\theta}\partial_\psi\right]^2+\frac{\bfr^2}{r^2\bfp^2}\left[\bfa\bfb\partial_t+\bfb\partial_\phi+\bfa\partial_\psi\right]^2\right)u.
\end{align}

This equation can be separated, using an Ansatz 
\begin{align}\label{eq:SeparationAnsatz}
u(t,r,\phi,\psi,\theta)=e^{-i\sigma t} e^{i(m\phi+k\psi)}R(r) S(\theta)     
\end{align} 
giving 
\begin{align}
    r^{-1}\Delta_r\partial_rR &+\left[\frac{(r^2+\bfa^2)^2(r^2+\bfb^2)^2}{r^4\Delta_r}\left(\sigma -\frac{m\bfa}{r^2+\bfa^2}-\frac{k\bfb}{r^2+\bfb^2}\right)^2 \right.\nonumber\\
    &\left.-\frac{1}{r^2}(\bfa\bfb\sigma-m\bfb-k\bfa)^2-\lambda_{\ell k m}\right]R(r)=0 \label{eq:RadialSeparatedEquation}\\
    \frac{1}{\cos\theta\sin\theta}\partial_\theta&(\cos\theta\sin\theta\,\partial_\theta S_{\ell k m})- \left[\frac{\cos^2\theta\sin^2\theta}{\bfp^2}\left(\bfa^2-\bfb^2)\sigma-\frac{m\bfa}{\sin^2\theta}+\frac{k\bfb}{\cos^2\theta}\right)\right.\nonumber\\
    &+\left.\frac{1}{\bfp^2}(\bfa\bfb\sigma-m\bfb-k\bfa)^2-\lambda^2_{\ell km}\right]S_{\ell k m}(\theta)=0\label{eq:AngularSeparatedEquation}.
\end{align}
In the above the $\lambda_{\ell m k}$ are eigenvalues belonging to the angular equation for the spherical harmonic (Wigner-D functions) $S_{\ell m k}(\theta)$ where  $l=0,1,2,3...$, $|m|+|k|\leq l$, $\ell -|m|-|k|\in 2\mathbb{Z}$  and 
\begin{align}
    \lambda_{\ell m k}= \ell(\ell+2).
\end{align}

It can be shown that the radial equation (\ref{eq:RadialSeparatedEquation}) is of confluent Heun type \cite{kunduri2005integrability}, \cite{ronveaux1995heun}, \cite{cvetivc2022mode}. The latter reference, \cite{cvetivc2022mode} is particularly notable in this context as they prove mode stability in the upper half plane for $(M,g_{\bfa,\bfb,\mfrm})$. That is they show that there exist no solutions of the Ansatz (\ref{eq:SeparationAnsatz}) with $\text{Im}\, \sigma>0$ that are smooth across the horizon $\mathcal{H}^+$ and satisfying suitable outgoing decay, which we will make more precise later.

\subsection{Extension and compactification}\label{subsec:ExtensionCompactification}

As discussed in section \ref{sec:MyersPerryBoyerLindquist}, in Boyer-Lindquist coordinates, $(M,g_{\bfa,\bfb,\mfrm})$ exhibits coordinate singularities at $r=r_{\pm}$ and $\theta=0,\frac{\pi}{2}$. To deal with the former, as in Kerr there exist $*$-coordinates that allows for a smooth extension of the metric across the event horizon $\mathcal{H}^\pm$. In order to describe late-time asymptotics correctly we will need to work on a slightly larger domain than just $M\subset \mathbb{R}^5$ warranting such extensions. Let us therefore introduce a slightly more flexible notion for the spacetime $(M,g_{\bfa,\bfb,\mfrm})$ that we will consider $\Box_g$ on and simultaneously introduce its radial compactification.
\begin{definition}\label{def:RadiallyCompactifiedR4}
    For $r_X\geq 0$, we define the spatial manifold $X^\circ\subseteq\mathbb{R}^4$ as 
    \begin{align}
        X^\circ = (r_X ,\infty)_r \times \mathbb{S}^3. 
    \end{align}
    Now let $X\subseteq \overline{\mathbb{R}^4}$ be the closure of this manifold under radial compactification, that is 
    \begin{align}
        X= [r_X,\infty]_r\times \mathbb{S}^3 
    \end{align}
    which is well defined as the gluing of 
    \begin{align}
        [r_X,\infty)_r\times\mathbb{S}^3 \sqcup [0,r_X^{-1}]_\rho \times \mathbb{S}^3/\sim
    \end{align}
    along $(0,r_X]_r$ $(0,r_X^{-1})_\rho$ using the identification with the inverted radial coordinate $\rho = r^{-1}$. If $r_X=0$ then $X=\overline{\mathbb{R}^4}$ and the gluing simply takes place along $(0,\infty)_r$ and $(0,\infty)_\rho$. When $r_X>0$ we let $\partial X|_{\infty}$ be the boundary $\partial \overline{\mathbb{R}^4}\cap X$, for which the function $(t,\rho,\phi,\psi,\theta)\mapsto \rho$ is a boundary defining function. The boundary at $\rho=0$ is diffeomorphic to the $3$-sphere, i.e. $\partial X \cap \partial X|_{\infty} \cong \mathbb{S}^3$.
\end{definition}
Now the $*$-coordinates that allow for an extension of $g_{\bfa,\bfb,\mfrm}$ across $r=r_\pm$ are defined given $r_+^2\leq 2\mfrm$ from (\ref{eq:EventHorizonMP}) as (see also \cite{horowitz2012black})
\begin{align}
    t_*&:= t+\int^r_{\sqrt{3\mfrm}}\frac{(\tilde{r}^2+\bfa^2)(\tilde{r}^2+\bfb^2)}{\tilde{r}^2\Delta_{\tilde{r}}}\:d\tilde{r}\label{eq:MPStartcoordinate2}\\
    \phi_* &:= \phi + \int^r_{\sqrt{3\mfrm}} \frac{\bfa(\tilde{r}^2 +\bfb^2)}{\tilde{r}^2\Delta_{\tilde{r}}} \:d\tilde{r},\label{eq:MPStarphicoordinate2} \\
     \psi_* &:= \psi + \int^r_{\sqrt{3\mfrm}} \frac{\bfb(\tilde{r}^2 +\bfa^2)}{\tilde{r}^2\Delta_{\tilde{r}}} \:d\tilde{r}\label{eq:MPStarpsicoordinate2},  
\end{align}
giving rise to the change in the radial vector field $\partial_r$:
\begin{align}
    \partial_r|_{t,\phi,\psi}\mapsto \partial_r|_{t_*,\phi_*,\psi_*} + \frac{(r^2+\bfa^2)(r^2+\bfb^2)}{r^2\Delta_r}\partial_{t_*}+\frac{\bfa (r^2+\bfb^2)}{r^2\Delta_r}\partial_{\phi_*}+\frac{\bfb(r^2+\bfa^2)}{r^2\Delta_r}\partial_{\psi_*}.
\end{align}

Under these transformations, the dual metric $g_{\bfa,\bfb,\mfrm}^{-1}$ transforms in a way very reminiscent of Kerr
\begin{align}
    &\bfr^2  g^{-1}_{\mfrm,\bfa,\bfb} =  \bfDelta \partial_r^2 +\partial_\theta^2   +2 \left[\frac{(r^2+\bfa^2)(r^2+\bfb^2)}{r^2}  \partial_{ t_*} + \frac{\bfa(r^2+\bfb^2)}{r^2}  \partial_{\phi_*} + \frac{\bfb(r^2+\bfa^2)}{r^2} \partial_{\psi_*}\right]\partial_r \nonumber \\
    &+ \frac{\sin^2\theta\cos^2\theta}{\bfp^2}\left[(\bfa^2 - \bfb^2) \partial_{t_*} + \frac{\bfa}{\sin^2\theta} \partial_{\phi_*} -\frac{\bfb}{\cos^2\theta}\partial_{\psi_*}  \right]^2 + \frac{\bfr^2}{r^2 \bfp^2}\left[\bfa\bfb\partial_{t_*} + \bfb \partial_{\phi_*} + \bfa \partial_{\psi_*}\right]^2. \label{eq:MPStarDualMetric2}
\end{align}

And we also see that the metric itself transforms in an analogous way. If we would define the $*$-coframe vectors analogous to (\ref{eq:CoframedT})-(\ref{eq:CoframedPsi}) in the $(t_*,\phi_*,\psi_*)$ coordinates
\begin{align}
    dT_*&:= dt_*-\bfa\sin^2\theta d\phi_* -\bfb\cos^2\theta d\psi_*\label{eq:CoframedT}\\
    d\Phi_*&:= (\bfb^2-\bfa^2)dt_* + (r^2 + \bfa^2)\bfa \: d\phi_* - (r^2 + \bfb^2)\bfb \: d\psi_* \label{eq:CoframedPhi}\\
    d\Psi_* &:= -\bfa \bfb dt_* + (r^2 +\bfa^2)\bfb\sin^2\theta d\phi_* + (r^2 + \bfb^2) \bfa \cos^2\theta_* d\psi \label{eq:CoframedPsi}
\end{align}

we see that 
\begin{align}
    dT\mapsto dT_*-\frac{\bfr^2}{\Delta_r}dr,\;\;d\Phi\mapsto d\Phi_*,\;\; d\Psi\mapsto d\Psi_* 
\end{align}
where the $d\Phi,d\Psi$ transforamtions have the $r$-dependent additions cancel out exactly, giving
\begin{align}
    g_{\bfa,\bfb,\mfrm}&= -\frac{\Delta_r}{\bfr^2}dT_*^2+  2 dT_* dr + \bfr^2d\theta^2 +\frac{\cos^2\theta \sin^2\theta}{\bfp^2\bfr^2} d\Phi^2+\frac{1}{\bfp^2r^2} d\Psi^2
\end{align}

\begin{align}
    g_{\bfa,\bfb,\mfrm} &= -\frac{\Delta_r}{\bfr^2}(dt_*-\bfa\sin^2\theta d\phi_* -\bfb\cos^2\theta d\psi_*)^2 + 2(dt_*-\bfa\sin^2\theta d\phi_* -\bfb\cos^2\theta d\psi_*)dr\nonumber\\
    &+\frac{\cos^2\theta \sin^2\theta}{\bfp^2\bfr^2}((\bfb^2-\bfa^2)dt_* + (r^2 + \bfa^2)\bfa \: d\phi_* - (r^2 + \bfb^2)\bfb \: d\psi_*)^2\nonumber\\
    &+\frac{1}{r^2\bfp^2}(-\bfa \bfb dt_* + (r^2 +\bfa^2)b\sin^2\theta d\phi_* + (r^2 + \bfb^2) \bfa \cos^2\theta_*)^2\label{eq:EddingtonFinkelsteinMPMetric}.
\end{align}
These coordinates allow for an extension into the inner domain of the Myers-Perry spacetime and is valid for all $r\in (0,\infty)$. 
The rescaled wave equation reads in this coordinate system for $u\in C^\infty(\mathbb{R}_{t_*}\times (0,\infty)_r\times \mathbb{S}^3_{\phi_*\psi_*\theta}$):
\begin{align}
    \bfr^2 &\Box_{g_{\bfa,\bfb,\mfrm}} u = \frac{1}{r} \partial_r\left(\Delta_r \partial_r + r\left[\frac{(r^2+\bfa^2)(r^2+\bfb^2)}{r^2}\partial_{t_*} +\frac{\bfa(r^2+\bfb^2)}{r^2}\partial_{\phi_*}+ \frac{\bfb(r^2+\bfa^2)}{r^2}\partial_{\psi_*} \right]\right)u \nonumber \\
    &+\left[\frac{(r^2+\bfa^2)(r^2+\bfb^2)}{r^2}\partial_{t_*} +\frac{\bfa(r^2+\bfb^2)}{r^2}\partial_{\phi_*}+ \frac{\bfb(r^2+\bfa^2)}{r^2}\partial_{\psi_*} \right]\partial_r u+ \frac{1}{\cos\theta\sin\theta}\partial_\theta(\cos\theta\sin\theta\, \partial_\theta)u  \nonumber\\
    &+\frac{\cos^2\theta\sin^2\theta}{\bfp^2}\left[(\bfa^2-\bfb^2)\partial_{t_*}+\frac{\bfa}{\sin^2\theta}\partial_{\phi_*}-\frac{\bfb}{\cos^2\theta}\partial_{\psi_*}\right]^2u + \frac{\bfr^2}{r^2\bfp^2}[\bfa\bfb \partial_{t_*}+\bfb\partial_{\phi_*}+\bfa\partial_{\psi_*}]^2u. 
    \label{eq:FullWaveEquationEddingtonFinkelsteinCoord}
\end{align}

With these Eddington-Finkelstein coordinates adapted for Myers-Perry we are able to quantify the full spacetime that we will solve $\Box_{g}u=0$ on. We will consider $(M,g_{\bfa,\bfb,\mfrm})$ with $M$ given as 
\begin{align}\label{eq:FullBaseManifold}
    M=\mathbb{R}_{t_*} \times X,
\end{align}
where $X$ is the set
\begin{align}
    [r_X,\infty]_r\times \mathbb{S}^3_{\phi_*\psi_*\theta}
\end{align}
with coordinates given by (\ref{eq:MPStartcoordinate2})-(\ref{eq:MPStarpsicoordinate2}) and metric $g_{\bfa,\bfb,\mfrm}$ by (\ref{eq:EddingtonFinkelsteinMPMetric}). The cut-off that we choose $r_X$ with is arbitrary, after which we will stop keeping track of waves, but following \cite{hintz2022sharp} it is conventient to choose it somewhere between the event and Cauchy horizon. An algebraically elegant choice is by setting $r_X = \mathfrak{r}_{\bfa,\bfb,\mfrm}$ where we take $\mathfrak{r}_{\bfa,\bfb,\mfrm}$ from\footnote{Another possible choice would have been $\mathfrak{r}_{\bfa,\bfb,\mfrm}=\sqrt{|\bfa\bfb|}$ although this quantity can be smaller than $r_-$, for example in the single rotating case.} (\ref{eq:BestPossibleBound}) 
\begin{align}
    \mathfrak{r}_{\bfa.\bfb,\mfrm}^2=\frac{1}{2}(r_+^2+r_-^2) = 2\mfrm - \frac{\bfa^2+\bfb^2}{2}.
\end{align}
In addition to the full spatial manifold $X$ we also want to sometimes distinguish to the part that is restricted to the exterior region only:
\begin{align}
    \tilde X = [r_+,\infty]_r \times \mathbb{S}^3_{\phi,\psi,\theta}.
\end{align}
It should be noted that it is also possible to define a tortoise radial coordinate $r_*$ through the equation
\begin{align}
    r_* = \int^r_{\sqrt{3\mfrm}}\frac{(\tilde{r}^2+\bfa^2)(\tilde{r}^2+\bfb^2)}{\tilde{r}^2\Delta_{\tilde{r}}}d\tilde r,
\end{align}
which is useful to describe the Penrose compactification of $M$. 
Specifically the analogous construction with Kerr goes through by setting retarded/forward null coordinates $u$ and $v$ as
\begin{align}
    u &= t_* - 2r_*\\
    v&=t_*.
\end{align}
A structure completely analogous to Kerr then arises for the exterior domain of communication $\tilde{M}$. Whilst we will not use these definitions later on in this paper, we sketch it here for completeness, starting with future null/past null infinity
\begin{align}
    \mathscr{J}^+ &= \cup_{u\in\mathbb{R}}\{ v\to +\infty, u \in \mathbb{R}\}\\
    \mathscr{J}^- &= \cup_{v\in\mathbb{R}}\{ u\to -\infty, v \in \mathbb{R}\}
\end{align}
and (defining with $J^\pm (S)$ the causal future/past of a set $S$ within the Penrose compactification)  we have that
\begin{align}
    \mathcal{H}^+&=\partial J^-(\mathscr{J}^+) = \cup_{v\in\mathbb{R}}\{ u\to +\infty, v \in \mathbb{R}\}\\
    \mathcal{H}^-&=\partial J^+(\mathscr{J}^-)=\cup_{u\in\mathbb{R}}\{ v\to -\infty, u \in \mathbb{R}\}\\
    \mathcal{H}^+ &\cup \mathcal{H}^- = \{r=r_+\},
\end{align}
owing to the fact that $r_*(r)\to+\infty$ as $r\to \infty$. More precisely 
we have that \cite{cvetivc2022mode}
\begin{align}
    r_*(r) &= \frac{(r_+^2+\bfa^2)(r_+^2+\bfb^2)}{2r_+(r_+^2-r_-^2)}\log\left( \frac{r^2-r_+^2}{r_+^2-r_-^2}\right) +\mathcal{O}(1) 
\end{align}
near $\mathcal{H}$.

$\mathcal{H}^\pm$ are in this case integral curves of the null Killing field 
\begin{align}\label{eq:HorizonKillingField}
    \chi_+&= \partial_{t_*}+\frac{\bfa}{r_+^2+\bfa^2} \partial_{\phi_*} +\frac{\bfb}{r^2_++\bfb^2}\partial_{\psi_*}.
\end{align}
The edges of the Penrose diagram, future/past timelike infinity $i^\pm$ and spatial infinity $i^0$ are then definitions that also come naturally, and correspond to
\begin{align}
    i^\pm &= \left\{u,v \to \pm \infty, v-u=2r_*,\; \frac{r_*}{|t|}\to 0\right\}\\
    i^0&=\cup_{t\in\mathbb{R}}\left\{u\to-\infty,v\to+\infty,  u+v= 2t \right\}.
\end{align}

We note that these Kerr-$*$ type coordinates fit within a framework within the maximal analytic extension of $(\tilde{M},g_{\bfa,\bfb,\mfrm})$. In particular, Myers and Perry used the $1+d$-dimensional generalization of the coordinates discussed above to define Kruskal-type coordinates that are completely regular across $\mathcal{H}^+$ and $\mathcal{H}^-$ simultaneously, showing that the overlapping Eddington patches are in fact smoothly connected at the bifurcation sphere $\mathcal{H}^+\cap \mathcal{H}^-=\mathcal{B}\cong \mathbb{S}^{3}$ \cite{myers1986black}.\\ Although we will stop keeping track of waves just within the inide of the black hole interior, the causal structure of the maximal analytic extension $(\mathcal{M}_{max},g_{\bfa,\bfb,\mfrm})$ can be described similarly to the Kerr black hole. Denoting the future domain of dependence by $D^+(\Sigma_0)$ for a smooth closed achronal regions joining the asymptotically bifurcation sphere we denote its future domain of dependence by, $D^+(\Sigma_0)$, and $D^-(\Sigma_0)$ by reversing time orientation. Then the future/past Cauchy horizons are given by 
\begin{align}
    \mathcal{C}\mathcal{H}^\pm(\Sigma_0) :=\overline{D^\pm(\Sigma_0)}\backslash I^{\mp}(D^\pm(\Sigma_0)),
\end{align}
with $I^{\pm}$ denoting the chronical future past. We then identify two cases: \begin{enumerate}
    \item \textbf{Double Spinning case.} $(|\bfa|+|\bfb|)^2<2\mfrm$ and $\bfa\neq 0$ and $\bfb \neq 0$. In this case we retrieve a causal structure for the maximal extension that is similar to Kerr.
    We have 
    \begin{align}
        \mathcal{C}\mathcal{H}^{\pm}\subset \{r=r_-\}
    \end{align}
 with Killing fields 
 \begin{align}
     \chi_- =\partial_{t_*}+\frac{\bfa}{r_-^2+\bfa^2}\partial_{\phi_*}+\frac{\bfb}{r_-^2+\bfb^2}\partial_{\psi_*},
 \end{align}
 where we can cover the inner horizon branch with forward coordinates with 
 \begin{align}
     u=t_*-2r_*,\qquad v=t_*
 \end{align}
 and we have $2$ branches for the future black-hole
 \begin{align}
     \mathcal{C}\mathcal{H}^+_1&=\{u\to-\infty, v\in \mathbb{R}\}\\
     \mathcal{C}\mathcal{H}_2^+&=\{v\to+\infty,u\in \mathbb{R}\}
 \end{align}
 \item \textbf{Single Spinning case.} Let us assume without loss of generality that $\bfb=0$. Then, even when $\bfa\neq 0$, the causal structure of the maximal extension reverts to that of Schwarzschild. Note that we have from (\ref{eq:EventHorizonMP}) and (\ref{eq:CauchyHorizonMP}) that
 \begin{align}
     r_+^2=2\mfrm-\bfa^2,\qquad r_-^2=0.
 \end{align}
 We then find that 
 \begin{align}
     D(\Sigma_0) =\mathcal{M}_{max},\qquad \mathcal{C}\mathcal{H}^+(\Sigma_0)=\mathcal{C}\mathcal{H}^-(\Sigma_0)=\emptyset,
 \end{align}
with spacelike future and past singular boundaries $\mathcal{S}^\pm$.

\end{enumerate}

As a result one can extend the doubly spinning case to a global tower structure  for the maximal extension of the Myers-Perry, similar to Kerr spacetime, by successive analytic extensions across the event and Cauchy horizons. Moreover, similar to Kerr, Myers-Perry allows for an analytic continuation through $x=0$ to $x<0$ in the coordinate $x=r^2$.  The geometric singularity at $\{\bfr^2 =0 \}$ is however of a substantially different nature than the ring singularity at $\{r^2_{\bfa}=0\}$, in the doubly spinning case\footnote{Note that $\bfr^2=r^2 +\bfa^2\cos^2\theta+\bfb^2\sin^2\theta$ is still strictly positive if $r=0$ and $\bfa,\bfb \neq 0$, for all $\theta$}. Finally, the earlier discussed notions regarding extremality discussed in Section \ref{sec:MyersPerryBoyerLindquist} should still be taken into account when discussing the causality structure of Myers-Perry: for $\bfa\neq 0$, $\bfb\neq 0$ and $(|\bfa|+|\bfb|)^2=2\mfrm$ we have that extremal Myers-Perry exhibts a smooth degenerate event horizon, whilst for either $\bfa=0$ or $\bfb =0$ at the vertices of the extremality gives an event horizon that is not even regular.

\begin{figure}[H]
\centering
\resizebox{\linewidth}{!}{%
\begin{tikzpicture}[x=.78cm,y=.78cm,font=\small]

\begin{scope}
\coordinate (itop)   at (0,4);
\coordinate (ibot)   at (0,-4);
\coordinate (bif)    at (-4,0);
\coordinate (iright) at (4,0);

\coordinate (innerleft) at (-8,4);
\coordinate (innertop)  at (-4,8);

\fill[orange!18]
  (bif) -- (itop)
  .. controls (-1.25,2.25) and (-1.25,-2.25) ..
  (ibot) -- cycle;

\fill[blue!10]
  (itop) -- (iright) -- (ibot)
  .. controls (-1.25,-2.25) and (-1.25,2.25) ..
  (itop) -- cycle;

\fill[violet!16]
  (bif) -- (innerleft) -- (innertop) -- (itop) -- cycle;

\draw[thick] (bif) --
  node[midway,below right] {$\mathcal H^+$} (itop);
\draw[thick] (bif) --
  node[midway,below left] {$\mathcal H^-$} (ibot);

\draw[thick] (bif) --
  node[midway,below left] {$r=r_+$} (innerleft);

\draw[thick,dashed] (itop) --
  node[midway,above right] {$\mathscr{J}^+$} (iright);
\draw[thick,dashed] (ibot) --
  node[midway,below right] {$\mathscr{J}^-$} (iright);

\draw[thick,dash dot] (innerleft) --
  node[midway,above left,sloped] {$r=r_-$} (innertop);
\draw[thick,dash dot] (innertop) --
  node[midway,above right,sloped] {$r=r_-$} (itop);

\draw[thick,red,dashed]
  (itop)
  .. controls (-1.25,2.25) and (-1.25,-2.25) ..
  node[pos=.70,left] {$r_{\mathcal E}(\theta)$}
  (ibot);

\node[align=center] at (-2.55,0)
  {ergoregion $\mathcal E$\\$(g_{tt}>0)$};
\node[align=center] at (1.65,0)
  {timelike $\partial_{t_*}$\\$(g_{tt}<0)$};
\node[align=center] at (-4,4) {interior region};

\foreach \p in {itop,ibot,bif,iright,innerleft,innertop}
  \fill (\p) circle (1.7pt);
\node[right=3pt] at (itop) {$i^+$};
\node[below] at (ibot) {$i^-$};
\node[right] at (iright) {$i^0$};
\node[below left=2pt] at (bif) {$\mathcal B$};
\node[align=center] at (-2,-5.15)
  {\textbf{a)} Doubly spinning ($ab\ne0$)};
\end{scope}

\begin{scope}[shift={(14.4,0)}]
\coordinate (itop)   at (0,4);
\coordinate (ibot)   at (0,-4);
\coordinate (bif)    at (-4,0);
\coordinate (iright) at (4,0);

\coordinate (innerleft) at (-8,4);

\fill[orange!18]
  (bif) -- (itop)
  .. controls (-1.25,2.25) and (-1.25,-2.25) ..
  (ibot) -- cycle;

\fill[blue!10]
  (itop) -- (iright) -- (ibot)
  .. controls (-1.25,-2.25) and (-1.25,2.25) ..
  (itop) -- cycle;

\fill[violet!16]
  (bif) -- (innerleft) -- (itop) -- cycle;

\draw[thick] (bif) --
  node[midway,below right] {$\mathcal H^+$} (itop);
\draw[thick] (bif) --
  node[midway,below left] {$\mathcal H^-$} (ibot);

\draw[thick] (bif) --
  node[midway,below left] {$r=r_+$} (innerleft);

\draw[thick,dashed] (itop) --
  node[midway,above right] {$\mathscr{J}^+$} (iright);
\draw[thick,dashed] (ibot) --
  node[midway,below right] {$\mathscr{J}^-$} (iright);

\draw[thick,decorate,
  decoration={snake,amplitude=.7mm,segment length=3mm}]
  (innerleft) -- node[midway,above=6pt] {$r=0$} (itop);

\draw[thick,red,dashed]
  (itop)
  .. controls (-1.25,2.25) and (-1.25,-2.25) ..
  node[pos=.70,left] {$r_{\mathcal E}(\theta)$}
  (ibot);

\node[align=center] at (-2.55,0)
  {ergoregion $\mathcal E$\\$(g_{tt}>0)$};
\node[align=center] at (1.65,0)
  {timelike $\partial_{t_*}$\\$(g_{tt}<0)$};
\node[align=center] at (-4,2.6) {interior region};

\foreach \p in {itop,ibot,bif,iright,innerleft}
  \fill (\p) circle (1.7pt);
\node[right=3pt] at (itop) {$i^+$};
\node[below] at (ibot) {$i^-$};
\node[right] at (iright) {$i^0$};
\node[below left=2pt] at (bif) {$\mathcal B$};
\node[align=center] at (-2,-5.15)
  {\textbf{b)} Singly spinning ($ab=0$)};
\end{scope}

\end{tikzpicture}
}
\caption{Schematic exterior and future-interior Pensrose diagrams for
five-dimensional Myers--Perry black holes.
\textup{(a)} Doubly spinning, $\bfa\bfb\ne0$ and $(|\bfa|+|\bfb|)^2<2\mfrm$:
 the purple interior lies between the outer and inner horizons. Whilst not visible here that for no angle $\theta\in [0,\frac{\pi}{2}]$ does $r_{\mathcal{E}}(\theta)$ pinch off at $\mathcal{H}^\pm$ even at $\theta=0,\theta=\frac{\pi}{2}$.
\textup{(b)} Singly spinning, exactly with either $\bfa$ or $\bfb$ (or both) vanishing and
$(|\bfa|+|\bfb|)^2<2\mfrm$: the purple interior ends at the spacelike singularity
$r=0$. The exterior ergoregions are depicted at a generic angle;
adjacent exteriors, past interiors, and further extensions are omitted.
$\mathcal B$ denotes the outer bifurcation sphere.}
\label{fig:myers-perry-interiors-comparison}
\end{figure}
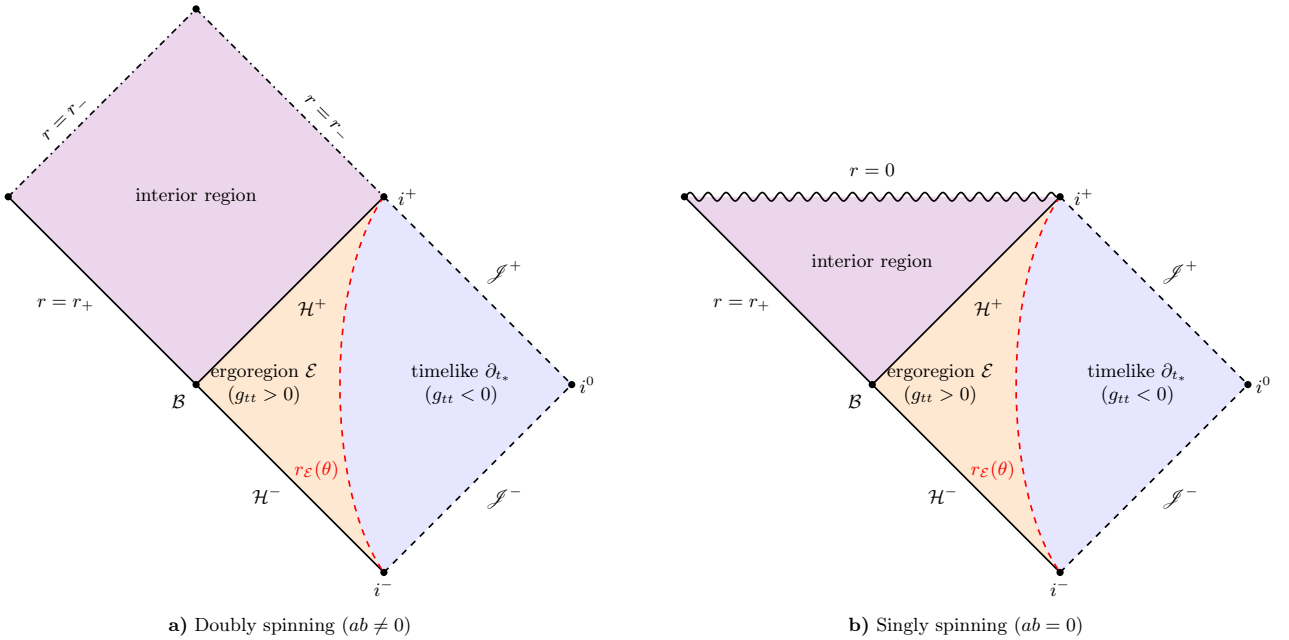

\begin{figure}[H]
\centering
\resizebox{\linewidth}{!}{%
\begin{tikzpicture}[x=1cm,y=1cm,font=\small,
 horizon/.style={thick,dashed},curve/.style={thick},
 every node/.style={inner sep=2pt}]
\begin{scope}
\coordinate (L) at (0,1.4);
\coordinate (T) at (2.8,3.3);
\coordinate (R) at (5,1.1);
\draw[curve] (L)--node[midway,above left] {$r=r_X$}(T);
\draw[horizon] (-.35,-.1)--(T);
\node at (-.3,.7) {$\mathcal H^+$};
\draw[horizon] (T)--(5.7,.4);
\node at (5.55,1.15) {$\mathscr I^+$};
\draw[curve] (L) .. controls (1.4,-.2) and (3.6,-.15) .. (R);
\node at (2.65,.62) {$t_*=t_{*,0}$};
\draw[red!80!black,curve] (1,-.25)--(T);
\node[red!80!black,below] at (1,-.25) {$r=r_0$};
\draw[blue!80!black,curve] (4.15,-.35) .. controls (4.22,1.1) and (3.7,2.5) .. (T);
\node[blue!80!black,below] at (4.15,-.35) {$r=t_*/v$};
\draw[->] (4.95,1.9)--(3.6,3.2) node[above] {$t_*$};
\draw[fill=white] (T) circle (2pt);
\node[above] at (T) {$i^+$};
\node at (2.6,-1.1) {\textbf{a)} Forward region};
\end{scope}
\begin{scope}[shift={(7.6,0)}]
\coordinate (L) at (0,1.4);
\coordinate (A) at (2.6,3.0);
\coordinate (B) at (3.65,3.3);
\coordinate (C) at (5.05,3.3);
\coordinate (R) at (6.9,1.1);
\draw[curve] (L)--node[midway,above left] {$r=r_X$}(A);
\draw[horizon] (.35,.2)--(2.76,2.76);
\node at (.8,1.52) {$\mathcal H^+$};
\draw[horizon] (A) .. controls (2.85,2.35) and (3.55,2.6) .. (B);
\node at (3.12,3.1) {$\mathcal K^+$};
\draw[horizon] (B)--(C);
\node[above] at (4.35,3.3) {$I^+$};
\draw[horizon] (C)--(7.5,.5);
\node at (7.35,1.2) {$\mathscr I^+$};
\draw[curve] (L) .. controls (1.5,-.2) and (4.6,-.2) .. (R);
\node at (3.2,.6) {$t_*=t_{*,0}$};
\draw[red!80!black,curve] (1.45,-.25)--(2.98,2.67);
\node[red!80!black,below] at (1.45,-.25) {$r=r_0$};
\draw[blue!80!black,curve] (5.05,-.35) .. controls (5,1.05) and (4.7,2.5) .. (4.5,3.3);
\node[blue!80!black,below] at (5.05,-.35) {$r=t_*/v$};
\draw[->] (6.8,2.05)--(5.65,3.2) node[above] {$t_*$};
\node at (3.65,-1.1) {\textbf{b)} Resolved future boundary};
\end{scope}
\end{tikzpicture}%
}
\caption{Schematic forward region for five-dimensional Myers--Perry and its resolved future boundary, following the layout of \cite{hintz2022sharp}, Figure 1.1. The solid upper-left boundary is the artificial spacelike cutoff $r=r_X$, with $r_-<r_X<r_+$ (taking $r_-=0$ in the singly spinning case). The dashed event horizon $\mathcal H^+$ lies inside the retained region. The red curve has $r=r_0>r_+$; the blue curve has $t_*/r=v>0$. Here $t_*$ denotes the time coordinate adapted to the future event horizon and future null infinity. In (b), $\mathcal K^+$ and $I^+$ are the successive resolved future faces; $\mathscr I^+$ denotes null infinity. The geometry in $b$ will be used when deriving the global asymptotics in Theorem \ref{thm:GlobalAsymptotics}.}
\label{fig:myers-perry-forward}
\end{figure}
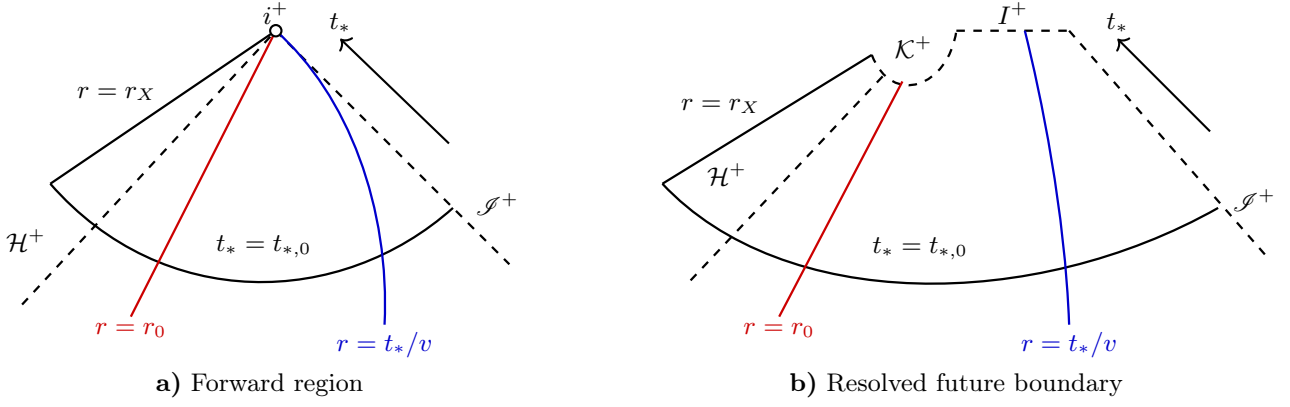

\section{Null geodesic flow and normally hyperbolic trapping}
\subsection{Hamiltonian dynamics of the dual metric function.}\label{subsec:HamiltonianDynamics}
In this first section we will analyze the dynamics associated to the Hamiltonian flow of the dual metric function, and elucidate its source-sink structure on the interior region $M^\circ = \mathbb{R}_{t_*}\times X^\circ$. We recall that given a canonical $1$-form $\upalpha \in T^*M^\circ$, its dual metric function $\tilde{G}$ is given by the contraction 
\begin{align}
    \tilde{G}=g^{-1}_{\bfa,\bfb,\mfrm}(\upalpha,\upalpha).
\end{align}
The object of our interest will be the \textit{rescaled} dual metric function, i.e. 
\begin{align}\label{eq:GeneralRescaledDualMetricFunction}
    G = \bfr^2 \cdot g^{-1}_{\bfa,\bfb,\mfrm}(\upalpha,\upalpha).
\end{align}

In local coordinates of $T^*_U M^\circ$ $(x,\zeta) = (x_i,\zeta_i)$ on a chart $(U,\phi)$ of $M^\circ$ we recall that the associated symplectic form is given by 
\begin{align}
    \omega = d\zeta = \sum_i d\xi_i\wedge dx_i.
\end{align}
The associated Hamiltonian vector field $H_G\in T(T^*M^\circ)$ associated to $G$ is given invariantly by 
\begin{align}
    \iota_{H_G} \omega = -dG
\end{align}
where $\iota_X$ denotes contraction by a vector field $X\in TM^\circ$ in the first component of $\omega$. Locally this reads as the well known equation
\begin{align}
    H_G = \sum_i \frac{\partial G}{\partial \xi_i}\frac{\partial}{\partial x_i}-\frac{\partial G}{\partial x_i}\frac{\partial}{\partial \xi_i}.
\end{align}

Integral curves $\gamma:I \to T^*M^\circ$ of this vector field are referred to in the context of microlocal analysis as \textit{null characteristics} (see also section \ref{subsec:bandscPseudodiffernetialOperators}, Proposition \ref{prop:bscPropagation}), and satisfy the equation 
\begin{align}
    \gamma'= H_G\gamma.
\end{align}
and lie in the characteristic set
\begin{align}
    \gamma(I)\subset\Sigma_G := \Char (G) = \{G=0\},
\end{align}
which by construction is an invariant set under the flow of $H_G$, we have 
\begin{align}
    e^{tH_G}(p)(\Sigma_G)\subseteq \Sigma_G, \; \forall \; (x_0,\xi_0) \in \Sigma_G, t\in I,
\end{align}
where $e^{tH_G}(x_0,\xi_0):I\to T^*M^\circ$ represents the exponential function associated to the bicharacteristic $\gamma:I\to T^*M^\circ$ with $\gamma(0)=(x_0,\xi_0)$.
In what follows we shall choose on the exterior region $\tilde M = \mathbb{R}_t \times \tilde X^\circ$ the canonical $1$-form 
\begin{align}\label{eq:Canonical1FormInteriorRegion}
    \upalpha= \sigma dt + \xi dr + \eta_\phi d\phi+\eta_\psi d\psi + \eta_\theta d\theta,
\end{align}
in the coordinates (\ref{eq:MyersOriginalMetric})
and whenever we will need to analyze the flow near the future event horizon $\mathcal{H}^+$ we will choose (with slight abuse of notation) 
\begin{align}\label{eq:Canonical1FormHorizon}
    \upalpha= \sigma dt_* + \xi dr + \eta_\phi d\phi_* + \eta_{\psi}d\psi_* +\eta_\theta d\theta.
\end{align}
with local smooth Kerr-star-like coordinates (\ref{eq:MPStartcoordinate2})-(\ref{eq:MPStarpsicoordinate2}).
We note that the characteristic set of $G$ locally represents the disjoint union of the forward/backward dual light cone $\Sigma_G\cap \{\pm \sigma >0\}$ of $G$ in $T^*M^\circ$. We also introduce the sets where $G$ is postive resp. negative definite as
\begin{align}
    \Sigma_G^\pm = \{\pm G >0\}
\end{align}
The main goal of this chapter is to prove that $(M,g_{\bfa,\bfb,\mfrm})$ admits a hyperbolically trapped set. 
\begin{definition}(Normally hyperbolic trapping)\label{def:TrappedSet} Let $\Gamma \subset T^*M^\circ$ be a co-dimension $2$ manifold with Lorentzian metric $g$ with dual metric function $G$. We say that $\Gamma$ is a normally hyperbolic trapped set with respect to $H_G$ if:  
    \begin{description}
        \item[($\Gamma_1$)] There exists two codimension $1$ submanifolds $\Gamma^\pm$ inside $\Sigma_G$ (hence codimension $2$ in $T^*M$ which intersect transversally in $\Gamma$:
        \begin{align}
            \Gamma =\Gamma^+ \cap \Gamma^-
        \end{align}
        and for which $H_G$ is tangent to $\Gamma^\pm$.
        \item[($\Gamma_2$)] There exists smooth defining function $\phi_\pm\in C^\infty(\Sigma_G)$ of $\Gamma^\pm$ near $\Gamma$, that is there exists an open neighbourhood $U\subset T^*M$ of $\Gamma$ such that $\Gamma^\pm \cap U = \{\phi_\pm=0\}$ and $d\phi_\pm \neq 0$ on $\Gamma^\pm$ such that 
    \begin{align}
        H_{G}\phi_+ =\nu_+ \phi_+ \text{ and } H_{G}\phi_-=-\nu_- \phi_-.
    \end{align}
    on $\Sigma_G$ near $\Gamma$ and $\nu_+,\nu_- \in C^\infty(U)$ with $\nu_{\min}:=\min(\inf_\Gamma \nu_+,\inf_\Gamma\nu_-)>0$.
    \item[($\Gamma_3$)] We have $\{\phi_-,\phi_+\}>0$ and $d G \neq 0$ on $\Gamma$.
    \end{description}

\end{definition}
The main result of this first part of the paper is then simply described as follows.

\begin{theorem}
    Let $(\tilde{M},g_{\bfa,\bfb,\mfrm})=\mathbb{R}_t\times \tilde {X}^\circ$ denote the domain of exterior communication $\{r>r_+\} \cap M^\circ$ with $g_{\bfa,\bfb,\mfrm}$ given by (\ref{eq:5DMPMetric}). Then for the full subextremal range (\ref{eq:ParameterDiamond}), i.e. $(\bfa,\bfb)\in D_{\mfrm}\subset \mathbb{R}^2$, $(\tilde{M},g_{\bfa,\bfb,\mfrm})$ admits a unique normally hyperbolic trapped set $\Gamma\subset T^*\tilde M^\circ$ satisfying conditions $(\Gamma_1)$-$(\Gamma_3)$.
\end{theorem}
    
We will connect these trapping dynamics associated to $H_G$ to the source structure at of the flow the future event horizon $\mathcal{H}^+$ and in Section \ref{subsec:bandscPseudodiffernetialOperators} to the source/sink structure on the compactified shell at infinity $\mathbb{R}_{t_*}\times\partial X|_\infty$, which will only involve using the asymptotically flat structure (\ref{eq:MinkowskiMetric}) at $r\to \infty$. This will complete the dynamical structure of the flow of $H_G$ on the full compactified spacetime $M = \mathbb{R}_{t_*}\times X$ allowing us to convert these dynamics into semiclassical microlocal estimates for $\bfr^2\Box_{g_{\bfa,\bfb,\mfrm}}$.\\
Thus to complete the preliminary definitions needed for the theory of this chapter we will introduce the concepts of sources/sinks for the Hamiltonian flow using defining functions as well, such that they are compatible with the estimates we will derive for them in Section \ref{subsec:bandscPseudodiffernetialOperators}. We will work on an arbitrary smooth $d$-dimensional manifold $\mathcal{M}$ with general local coordinates $(x,\xi)\in T^*M$. Let us introduce fibrewise-radial coordinates $(r_\xi,\omega_\eta)$ via $r_\xi=\|\xi\|$ and $\omega_\eta \in \mathbb{S}^3$ such that 
\begin{align}
    \xi = r_\xi \cdot \omega_\eta.
\end{align}
Furthermore let $P \in C^\infty (T^* \mathcal{M})$ denote a smooth degree $m$ homogeneous polynomial in the variables $\xi_i$ that is 
\begin{align}
    P(x,\xi) = \sum_{|\alpha|\leq m} a_\alpha(x,\xi) \xi^\alpha,
\end{align}
where $\alpha=(\alpha_1,...,\alpha_d) \in \mathbb{N}_0^d$ is a multi-index and $\xi^\alpha = \xi_1^{\alpha_1}...\xi_d^{\alpha_d}$. 
A radial point $(x,\xi) \in T^*M$ for the flow of $H_p$ of $P$ is then a point where the Hamiltonian vector field is tangent to the generator of dilations in $T^*_x \mathcal{M}$,
i.e. 
\begin{align}
    H_P\sim r_\xi \partial_{r_\xi}.
\end{align}
Now let $M_\lambda: T^*\mathcal{M}\to T^*\mathcal{M}$ denote the fibrewise dilation by $\lambda\in \mathbb{R}$, i.e. $M_\lambda: (x,\xi)\mapsto (x,\lambda\xi)$.
When $P\in T^*\mathcal{M}$ is degree $m$ homogeneous, one can easily show that $M_\lambda^* H_p = \lambda^{m-1}H_p$ where we denote with $*$ the pullback map. This implies that the integral curves of $H_p$ through $(x,\lambda)$ and $(x,\lambda\xi)$ agree up to rescaling. From this perspective the rescaled Hamiltonian vector field $\rho_\xi^{-(m-1)} H_p$ can be projected onto the sphere bundle via $ \pi_{S^*}:T^*\mathcal{M}\to S^*\mathcal{M}$. Radial points of $H_p$ then become points where $d\pi_{S^*}H_p$ vanishes. \\
Another way to interpret radial points and interpolate between $T^*\mathcal{M}$ and $S^*\mathcal{M}$ is to introduce the fibre-radially compactified cotangent space $\overline{T^*}\mathcal{M}$. In principle this space is obtained through the same the radial compactification of $\mathbb{R}^d$ introduced in Definition \ref{def:RadiallyCompactifiedR4} but now the gluing is done fibrewise over the tangent variables.
\begin{definition}\label{def:FibrewiseCompactification}
    Let $\mathcal{M}$ be a $d$-dimensional smooth manifold with cotangent bundle $T^*\mathcal{M}$. For a fibre $T^*_x\mathcal{M}$, let $\xi=(\xi_1,...,\xi_d)$ denote local fibre coordinates and define  $r_\xi$ and $\omega_\xi \in S^*_x\mathcal{M}\cong \mathbb{S}^d$ through $\rho_\xi=\|\xi\|$ and $\xi = r_\xi \cdot \omega_\xi$. We then define 
    \begin{align}
        \overline{T^*_x}\mathcal{M} = T^*_x\mathcal{M}\sqcup [0,r_\xi^{-1})_{\rho_\xi}\times\mathbb{S}^d/\sim
    \end{align}
    where the gluing $\sim$ is done along $\varrho_\xi=r_\xi^{-1}$. Then $\varrho_\xi$ is a boundary defining coordinate for $\partial T^*_x\mathcal{M} \cong S^*_x\mathcal{M}$.
    We define the compactified space space as the union
    \begin{align}
        \overline{T^*}\mathcal{M} = \bigsqcup_{x\in \mathcal{M}}\overline{T}^*_x\mathcal{M}.
    \end{align}
\end{definition}
With this definition one easily checks that $\partial \overline{T^*}\mathcal{M}\cong S^*\mathcal{M}$. One can additionally show that when $P$ is homogeneous of degree $m$, its rescaled Hamiltonian vector field $\varrho_\xi^{m-1}H_p$ extends smoothly to $\partial\overline{T^*}\mathcal{M}$.  We thus arrive at the equivalent notion that a $H_p$ has a radial point at $(x,\xi)$ if it is normal to $S^*\mathcal{M}$, vanishes on projective space $\mathbb{P}T^*\mathcal{M}$ or1 when its rescaling $\rho_\xi^{-(m-1)} H_p$ vanishes at fibre infinity $\partial\overline{T^*}\mathcal{M}$. We will phrase this with slightly more generality where we set 
\begin{align}
    \Char(P) = \{(x,\xi)\in T^*\mathcal{M}\;|\;P(x,\xi)=0\}.
\end{align}

\begin{definition}\label{def:GeneralizedRadialSet} (Generalized radial set at fibre infinity.)
 Let $\Lambda\subset \Char(P)\subset \partial\overline{T^*}\mathcal{M}\cong S^*\mathcal{M}$ be a smooth manifold and $\rho_\infty\in C^\infty(T^*\mathcal{M}\backslash 0)$ be a defining function of $\partial \overline{T^*}\mathcal{M}$ that is homogenous of degree $-1$ and elliptic in a neighbourhood of $\Lambda$.
    Let $P\in C^\infty(T^*\mathcal{M})$ be homogeneous of degree $m$ such that its associated rescaled Hamiltonian vector field $\rho_{\infty}^{m-1}H_p$ extends smoothly to $\overline{T^*}\mathcal{M}$. Then $\Lambda\subset \Char(P)$ is called a radial source/sink ($+/-)$ for $H_p$ if the following conditions hold:
    \begin{description}
        \item[($R_1$)] We have that
        \begin{align}\label{eq:RadialCondition1}
            \rho_\infty^{m-1}H_p \rho_\infty = \pm\alpha_\infty \rho_\infty
        \end{align}
        with $\alpha_{\infty}\in C^\infty(\Lambda)$ homogeneous of degree $0$ and satisfies $\alpha_\infty>0$.
        \item[($R_2$)] There exists a $\mathfrak{q}\in C^\infty(U)$ with $U$ a conic neighbourhood of $\Lambda$ which is a quadratic defining function of $\Lambda\subset \Char (P)$ such that 
        \begin{align}\label{eq:RadialCondition2}
            \Lambda= \{\mathfrak{q}=0 \}\cap \Char(P)
        \end{align}
        and $\mathfrak{q}$ vanishes quadratically at $\Lambda$ such that in $U$ we have 
        \begin{align}\label{eq:RadialCondition3}
            \pm\rho_\infty^{m-1}H_p\mathfrak{q}\geq \alpha_{U}\mathfrak{q}+\mathcal{O}(\mathfrak{q}^{3/2})
        \end{align}
        where $\alpha_U\in C^\infty(U)$ is homogenous of degree $0$ and satisfies $\alpha_U>0$.
    \end{description}
  We define a \textit{generalized} radial set as a set for which conditions ($R_1$) and ($R_2$) hold, but we do not exclude $\rho^{m-1}_\infty H_p$ have drift tangential to $\Lambda$. That is, we only require that $\Lambda$ is invariant under the rescaled Hamiltonian vector field $\rho_\infty^{m-1} H_p$:
    \begin{align}\label{eq:GeneralizedRadialCondition}
        \rho^{m-1}_\infty H_p \in T\Lambda.
    \end{align}
\end{definition}

We note that very often it will be useful to work on $\overline{T^*}\mathcal{M}$ not in radial fibre coordinates but in projective fibre coordinates, whose charts are given by 
\begin{align}
    U_i :=\left\{(x,\hat\xi_i) \in T^*\mathcal{M} \;|\; \hat{\xi}_1= \frac{\xi_1}{|\xi_i|},...,\hat{\xi_i}=|\xi_i|^{-1}, \hat{\xi}_d= \frac{\xi_d}{|\xi_i|} \right\},
\end{align}
used for the open regions $|\xi_i|>>|\xi_j| $ for all $j\neq i$. In these conic sectors this gives a boundary defining function of $\partial\overline{T^*}\mathcal{M}$ given by $\rho_\infty = |\xi_i|^{-1} =|\hat{\xi}_i|$, for which the formulation in Definition \ref{def:GeneralizedRadialSet} goes through unchanged.\\

To conclude this preliminary section on dynamics, we would like to emphasize that because $(M,g_{\bfa,\bfb,\mfrm})$ is stationary, we have that 
\begin{align}\label{eq:LittleInvariance}
    H_G\sigma = -\partial_{t_*}G =0.
\end{align}
Therefore, as the dynamics of $H_G$ are time-translation in variant, one could opt to restrict to the fixed frequency slice
\begin{align}
    \{\sigma \equiv z\}\cap \{t_*=0\}\subset T^*M^\circ,
\end{align}
when studying dynamics on $T^*M^\circ=T^*(\mathbb{R}_{t_*}\times X^\circ)$. We note that because of (\ref{eq:LittleInvariance}), slices $\{\sigma\equiv z\}$ are invariant under the flow of $H_G$ but in general the slice $\{t_*=0\}$ is not. Thus the intersection should be seen as the projection of the flow of $H_G$ inside $\{\sigma\equiv z\}$ \textit{down to } $\{\sigma \equiv z\}\cap \{t_*=0\}$.  In the context of semiclassical microlocal analysis (c.f. Section \ref{subsec:SemiclassicalAnalysis}) there are two hypersurfaces of particular interest. One is the dynamics of $G|_{\sigma=0}$, i.e. $H_G$ restricted to the slice $\{\sigma=0\}$, which we will refer to as the \textit{classical} dynamics of $G$. The other are the dynamics $G|_{\sigma=\pm1}$ associated to $\{\sigma =\pm1\}$ which we will refer to as the \textit{semiclassical} dynamics of $G$. The latter two slices are conjugate to one another in the sense of time reversal; the direction of the associated flow flips sign upon inverting $\sgn\, \sigma$. In former treatments of the dynamical flow on Kerr, e.g. \cite{millet2023optimal}, the analysis of the flow on these two hypersurfaces were seperated for microlocal considerations. Here we opt to do a full treatment of the null geodesic flow first, before specializing to the classical and semiclassical dynamics in part II.

\subsection{Dual metric function and Carter's constant.}

The starting point for $1+4$ Myers-Perry is equation (\ref{eq:DualMPMetric}) on the exterior region $\mathbb{R}_{t}\times \tilde{X}^\circ$.

\begin{align}
    \bfr^2 \cdot g^{-1}_{\bfa,\bfb,\mfrm} &:= -\frac{1}{r^4\bfDelta}\left[(r^2+\bfa^2)(r^2 +\bfb^2) \partial_t + \bfa(r^2 + \bfb^2) \partial_\phi + \bfb(r^2 + \bfa^2)\partial_\psi \right]^2 + \bfDelta \partial_r^2 + \partial^2_\theta \nonumber\\
    &+ \frac{\sin^2\theta \cos^2\theta}{\bfp^2} \left[ (\bfa^2-\bfb^2)\partial_t + \frac{\bfa}{\sin^2\theta}\partial_\phi - \frac{\bfb}{\cos^2\theta}\partial_\psi \right]^2 + \frac{\bfr^2}{r^2\bfp^2}[\bfa\bfb \partial_t + \bfb\partial_\phi+ \bfa \partial_\psi]^2 \nonumber.  
\end{align}
We note that
\begin{align}
    \frac{\bfr^2}{r^2\bfp^2} = \frac{1}{r^2}+\frac{1}{\bfp^2},
\end{align}
which allows us to further separate the dual metric in its radial and $\mathbb{S}^3$-spherical components by writing
\begin{align}
    \bfr^2 \cdot g^{-1}_{\bfa,\bfb,\mfrm} &:= -\frac{1}{r^4\bfDelta}\left[(r^2+\bfa^2)(r^2 +\bfb^2) \partial_t + \bfa(r^2 + \bfb^2) \partial_\phi + \bfb(r^2 + \bfa^2)\partial_\psi \right]^2 + \bfDelta \partial_r^2 + \frac{1}{r^2}[\bfa\bfb \partial_t + \bfb\partial_\phi+ \bfa \partial_\psi]^2 \nonumber\\
    &+\partial_\theta^2 + \frac{\sin^2\theta \cos^2\theta}{\bfp^2} \left[ (\bfa^2-\bfb^2)\partial_t + \frac{\bfa}{\sin^2\theta}\partial_\phi - \frac{\bfb}{\cos^2\theta}\partial_\psi \right]^2 + \frac{1}{\bfp^2}[\bfa\bfb \partial_t + \bfb\partial_\phi+ \bfa \partial_\psi]^2  \label{eq:DualMPMetricSep}.  
\end{align}
We will use the canonical one-form $\upalpha$ on $T^*\tilde M^\circ$ to introduce covariables $(\sigma,\xi,\eta_\phi,\eta_\psi,\eta_\theta)$ as 
\begin{align}
    \upalpha = \sigma dt +\xi dr +\eta_\phi d\phi +\eta_\psi d\psi +\eta_\theta d\theta.
\end{align}
Recall that we use $G$ as the symbol for the rescaled dual metric function and $\tilde{G}$ for the true dual metric function, we have
\begin{align}
    G=\bfr^2g^{-1}_{\bfa,\bfb,\mfrm}(\upalpha,\upalpha) \text{ and } \tilde{G}=g^{-1}_{\bfa,\bfb,\mfrm}(\upalpha,\upalpha).
\end{align}

We note that $G$ is given via (\ref{eq:DualMPMetricSep}) as

\begin{align}
    G&:= -\frac{1}{r^4\bfDelta}\left[(r^2+\bfa^2)(r^2 +\bfb^2) \sigma + \bfa(r^2 + \bfb^2) \eta_\phi + \bfb(r^2 + \bfa^2)\eta_\psi \right]^2 + \bfDelta \xi^2 + \frac{1}{r^2}[\bfa\bfb \sigma + \bfb\eta_\phi+ \bfa \eta_\psi]^2 \nonumber\\
    &+\eta_\theta^2 + \frac{\sin^2\theta \cos^2\theta}{\bfp^2} \left[ (\bfa^2-\bfb^2)\sigma + \frac{\bfa}{\sin^2\theta}\eta_\phi - \frac{\bfb}{\cos^2\theta}\eta_\psi \right]^2 + \frac{1}{\bfp^2}[\bfa\bfb \sigma + \bfb\eta_\phi+ \bfa \eta_\psi]^2  \label{eq:DualMPMetricSep2}.  
\end{align}

Decomposing $G$ in its radial and spherical components, i.e. defining

\begin{align}
    \mathscr{G} &:= -\frac{1}{r^4\bfDelta}\left[(r^2+\bfa^2)(r^2 +\bfb^2) \sigma + \bfa(r^2 + \bfb^2) \eta_\phi + \bfb(r^2 + \bfa^2)\eta_\psi \right]^2 + \bfDelta \xi^2 + \frac{1}{r^2}[\bfa\bfb \sigma + \bfb\eta_\phi+ \bfa \eta_\psi]^2\\
    \mathscr{C} &:= \eta_\theta^2 + \frac{\sin^2\theta \cos^2\theta}{\bfp^2} \left[ (\bfa^2-\bfb^2)\sigma + \frac{\bfa}{\sin^2\theta}\eta_\phi - \frac{\bfb}{\cos^2\theta}\eta_\psi \right]^2 + \frac{1}{\bfp^2}[\bfa\bfb \sigma + \bfb\eta_\phi+ \bfa \eta_\psi]^2, 
\end{align}

we have that 
\begin{align}\label{eq:SymbolSplitting}
    G=\mathscr{G}+\mathscr{C},
\end{align}
where $\mathscr{G}$ is $(\theta,\eta_\theta)$-independent, and $\mathscr{C}$ is $(r,\xi)$-independent. Consequently they are both conserved constants\footnote{One reason for looking at the rescaled flow $\bfr^2 \tilde{G}$ is that otherwise one would not have a fully separated flow in $r$.} of motion and we have 
\begin{align}
    H_G \mathscr{G}=0,\;H_G\mathscr{C}=0\;\text{ and } H_G=H_{\mathscr{G}}+ H_{\mathscr{C}}.
\end{align}

For further simplification let us introduce the following quantities $A,B$ and $C$ that are related to the co-frame $dT,d\Phi$ and $d\Psi$ in (\ref{eq:CoframedT})-(\ref{eq:CoframedPsi}) and as

\begin{align}
    A &= \frac{1}{r^2}[(r^2+\bfa^2)(r^2+\bfb^2)\sigma + \bfa (r^2+\bfb^2)\eta_\phi + \bfb(r^2+\bfa^2)\eta_\psi]\label{eq:MatrixEntry1}\\
    B&=\frac{\cos\theta\sin\theta}{\bfp}\left[(\bfa^2-\bfb^2)\sigma + \frac{\bfa}{\sin^2\theta}\eta_\phi -\frac{\bfb}{\cos^2\theta}\eta_\psi\right]\label{eq:MatrixEntry2}\\
    &=\frac{1}{\bfp}\left[(\bfa^2-\bfb^2)(\cos\theta\sin\theta)\sigma + \bfa\cot\theta\:\eta_\phi -\bfb\tan\theta\:\eta_\psi\right]\label{eq:MatrixEntry2Expanding}\\
    C&=\frac{\sqrt{\Delta_r}}{r}[\bfa\bfb\sigma + \bfb \eta_\phi +\bfa \eta_\psi],\label{eq:MatrixEntry3}
\end{align}

Using these quantities we can recast $\mathscr{G}$ and $\mathscr{C}$ as 
\begin{align}
    \mathscr{G}&=-\Delta_r^{-1}(A^2-C^2) + \Delta_r\xi^2\\
    \mathscr{C}&=\eta_\theta^2+B^2 + \frac{r^2}{\bfp^2}\Delta^{-1}_rC^2,\label{eq:CartersConstantonCharacteristicSet}
\end{align}
and for a given $A,B,C$ one can fix the hypersurface $\{\sigma\equiv 
\text{const} \}$ by inverting (\ref{eq:MatrixEntry1})-(\ref{eq:MatrixEntry3}):

\begin{align}\label{eq:sigmaInverted}
    \sigma=\frac{A}{\bfr^2
    }-\frac{\bfa^2-\bfb^2}{\bfp\bfr^{2}}B-\frac{\bfa\bfb}{r\sqrt{\Delta_r}\bfp^2}C.
\end{align}

Comparison with the relevant literature for the Kerr black hole (e.g. \cite{dyatlov2015asymptotics}) easily shows that $\mathscr{G}$ and $\mathscr{C}$ are equivalent to $G_r$ resp. $G_\theta$ in that setting.
Two main references for Carter's constant in the context of $1+4$-dimensional Myers-Perry is given by Frolov and Stojkovic \cite{frolov2003particle} and by the aforementioned Wu \cite{wu2008separability}. The former reference computed Carter's constant through the programme GRTensor, whilst the latter used a method analagous to \cite{carter1968hamilton} to construct Carter's constant from an orthogonal pentad construction whose frame diagonalizes the dual metric $g^{-1}_{\bfa,\bfb,\mfrm}$ as in (\ref{eq:DualMPMetric}). \\
In particular Wu defined the pentad $\{\omega_i\}$ for $i=1,...,5$ as
\begin{align}
    \omega^{(0)}&=\frac{1}{r^2\sqrt{\bfr^2\Delta_r}}\left((r^2+\bfa^2)(r^2+\bfb^2)\partial_t + \bfa(r^2+\bfb^2)\partial_\phi+\bfb(r^2+\bfa^2)\partial_\psi\right)\\
    \omega^{(1)}&=\sqrt{\frac{\Delta_r}{\bfr^2}}\partial_r\\
    \omega^{(2)}&=\frac{1}{\sqrt{\bfr^2}}\partial_\theta\\
    \omega^{(3)}&=\frac{1}{\bfp\sqrt{\bfr^2}}\left[(\bfa^2-\bfb^2)\cos\theta\sin\theta\:\partial_t+\bfa\cot\theta\:\partial_\phi-\bfb\tan\theta\:\partial_\psi\right]\\
    \omega^{(4)}&=\frac{1}{r\bfp}(\bfa\bfb\:\partial_t+\bfb\partial_\phi +\bfa\partial_\psi).\\
\end{align}
for which 
\begin{align}
    g^{-1}_{\bfa,\bfb,\mfrm}|_{(\omega_{(i)},\omega_{(j)})}=\diag{(-1,+1,+1,+1)}
\end{align}

It comes with a dual coframe 
\begin{align}
    \omega_{(0)}&=\sqrt{\Delta_r\bfr^{-2}}(dt-\bfa\sin^2\theta d\phi-\bfb \cos^2\theta d\psi)=\sqrt{\Delta\bfr^{-2}}dT\\
    \omega_{(1)}&=\sqrt{\bfr^2\Delta^{-1}_r}dr\\
    \omega_{(2)}&=\bfr^2 d\theta \\
    \omega_{(3)}&=\frac{\cos\theta\sin\theta}{\bfp\bfr}((\bfb^2-\bfa^2)dt + \bfa(r^2+\bfa^2)d\phi -\bfb(r^2+\bfb^2)d\psi)=\frac{\cos\theta\sin\theta}{\bfp\bfr}d\Phi\\
    \omega_{(4)}&=\frac{1}{r\bfp}(-\bfa\bfb dt + \bfb(r^2+\bfa^2)\sin^2\theta +\bfa(r^2+\bfb^2)\cos^2d\psi)=\frac{1}{r\bfp}d\Psi.
\end{align}

Wu then uses these frames to derive a St\"ackel-Killing tensor $K^{\mu\nu} \in TM\otimes TM$, which we will introduce here with the opposite sign, via
\begin{align}
    K^{\mu\nu}\partial_\mu\partial_\nu=\bfp^2 (\omega^{(0)})^2-\bfp^2 (\omega^{(1)})^2+r^2(\omega^{(2)})^2+r^2(\omega^{(3)})^2 -(\bfp^2-r^2) (\omega^{(4)})^2.
\end{align}
It can be easily checked that $K^{\mu\nu}$ satisfies the cyclic covariant derivative condition through
\begin{align}\label{eq:CyclicCovariantCondition}
    K_{\mu\nu;\rho}+K_{\nu\rho;\mu}+K_{\rho\mu;\nu}=0.
\end{align}

Let us denote $p_i\in C^\infty(T^*M)$ for the contraction with the canonical $1$-form, 
\begin{align}
    p_{(i)}:=\upalpha (\omega^{(i)})=\omega_{(i)}(\upalpha^\sharp) \text{ and } \upalpha =\sum_{i=1}^5 p_{(i)} \omega_{(i)}. 
\end{align}
We then have that 
\begin{align}
    G=\bfr^2(-p_{(0)}^2+p_{(1)}^2+p_{(2)}^2+p_{(3)}^2+p_{(4)}^2) \text{ and } \tilde G = -p_{(0)}^2+p_{(1)}^2+p_{(2)}^2+p_{(3)}^2+p_{(4)}^2.
\end{align}
If we then define Carter's constant as defined by Wu in the covariables as 
\begin{align}
    \mathscr{C}_{Wu} := K^{\mu\nu}(\upalpha,\upalpha),
\end{align}
and observe that the constant $\mathscr{C}$ defined by (\ref{eq:CartersConstantonCharacteristicSet}) is given by
\begin{align}
    \mathscr{C}=\bfr^2(p_{(2)}^2+p^2_{(3)})+r^2p^2_{(4)}\end{align}

we then see that 
\begin{align}\label{eq:CartersConstantRelation}
    \mathscr{C}_{Wu}&=\underbrace{\bfr^2(p_{(2)}^2+p^2_{(3)})+r^2p^2_{(4)}}_{\mathscr{C}}-\bfp^2\underbrace{(-p_{(0)}^2+p_{(1)}^2+p_{(2)}^2+p_{(3)}^2+p_{(4)}^2)}_{=\tilde{G}}\\
    &=\mathscr{C}-\frac{\bfp^2}{\bfr^2}G.
\end{align}
i.e. on the characteristic set $\Sigma_G =\{G=0\}=\{\tilde G=0\}=\Sigma_{\tilde{G}}$ the two definitions agree.
The other commonly cited reference regarding Carters constant is that of Frolov and Stojkovich \cite{frolov2003particle}. Not having access to Wu's diagonal frame, they computed it straight from the metric elements $g^{\mu\nu}$ with GRTensor which they gave in Boyer-Lindquist coordinates as 
\begin{align}
    g^{tt}&=\frac{1}{\bfr^2}\left[(\bfa^2-\bfb^2)\sin^2\theta-\frac{(r^2+\bfa^2)(r^2\Delta_r +2\mfrm(r^2+\bfb^2)}{r^2\Delta_r} \right]\label{eq:FrolovExpressiontt}\\
    g^{t\phi}&=-2\mfrm\frac{\bfa(r^2+\bfb^2)}{\bfr^2r^2\Delta_r},\;\;\;\;\;g^{t\psi}=-2\mfrm\frac{\bfb(r^2+\bfa^2)}{\bfr^2r^2\Delta_r}\\
    g^{\phi\phi}&=\frac{1}{\bfr^2}\left[\frac{1}{\sin^2\theta}-\frac{(\bfa^2-\bfb^2)(r^2+\bfb^2)+2\mfrm\bfb^2}{r^2\Delta_r}\right]\\
    g^{\psi\psi}&=\frac{1}{\bfr^2}\left[\frac{1}{\cos^2\theta}+\frac{(\bfa^2-\bfb^2)(r^2+\bfa^2)+2\mfrm\bfa^2}{r^2\Delta_r}\right]\\
    g^{\phi\psi}&=-2\mfrm \frac{\bfa\bfb}{\bfr^2r^2\Delta_r},\;\;\;\;\; g^{rr}=\frac{\Delta_r}{\bfr^2},\;\;\;\;\; g^{\theta\theta}=\frac{1}{\bfr^2}\label{eq:FrolovExpressionThetaTheta},
\end{align}
where we note that they used the same convention for the sign of $\bfa,\bfb$ as Myers did in (\cite{horowitz2012black}) which we introduced here with the opposite convention to keep consistency with \cite{wu2008separability}.
They then computed the St\"ackel Killing tensor as 
\begin{align}\label{eq:StackelKillingTensor}
    K^{\mu\nu}=-\bfp^2(g^{\mu\nu}+\delta^\mu_t\delta^\nu_t)+\frac{1}{\sin^2\theta}\delta^\mu_\phi\delta^\nu_\phi+\frac{1}{\cos^2\theta
    }\delta_\psi^\mu\delta^\nu_\psi +\delta^\mu_\theta\delta^\nu_\theta.
\end{align}
which again satisfies (\ref{eq:CyclicCovariantCondition}).
Contracting this with $\zeta$ gives a neat expression $\mathscr{C}_{FS}$ for the Carter's constant as proposed by (\cite{frolov2003particle}), i.e.
\begin{align}
    \mathscr{C}_{FS}:= -\frac{\bfp^2}{\bfr^2}G+ \eta_\theta^2 +\frac{\eta_\phi^2}{\sin^2\theta}+\frac{\eta_\psi^2}{\cos^2\theta}-\sigma^2\bfp^2.
\end{align}

It is clear that $\mathscr{C}_{FS}$ is regular in the limit $(\bfa,\bfb)\to 0$ and, that in this limit we recover the $1+4$ Schwarzschild-Tangherlini \footnote{Our reason for introducing $K^{\mu\nu}$ (\ref{eq:StackelKillingTensor}) with the opposite sign of that in \cite{wu2008separability} is to keep this limit positive-definite.}{angular momentum} as one would expect:
\begin{align}\label{eq:SchwarzschildCarterConstant}
    \lim_{(\bfa,\bfb)\to (0,0)} \mathscr{C}_{FS}= \eta_\theta^2+\frac{\eta_\phi^2}{\sin^2\theta}+\frac{\eta_\psi^2}{\cos^2\theta}.
\end{align}
It is however not so clear that $\mathscr{C}$ that we defined in (\ref{eq:CartersConstantonCharacteristicSet}) on the characteristic set, and therefore Wu's $\mathscr{C}_{Wu}$ are in fact regular in this limit. Fortunately this in fact still holds.

\begin{proposition}\label{prop:RegularCarterConstant}
    Let $\mfrm >0$. We have for all $(\bfa,\bfb)\in D_{\mfrm}$ that 
    \begin{align}\label{eq:RegularCarterConstant}
        \mathscr{C}=\eta_\theta^2 + \frac{\eta_\phi^2}{\sin^2\theta}+\frac{\eta_\psi^2}{\cos^2\theta}+\sigma^2(\bfa^2\sin^2\theta+\bfb^2\cos^2\theta) + 2(\bfa\eta_\phi+\bfb\eta_\psi)\sigma.
    \end{align}
    Hence it is regular in $(\bfa,\bfb)=(0,0)$ and tends to the Schwarzschild-Tangherlini limit (\ref{eq:SchwarzschildCarterConstant}).
\end{proposition}

\begin{proof}
    To prove we start out from 
    \begin{align}
        \mathscr{C}&=\eta_\theta^2 + \frac{\sin^2\theta \cos^2\theta}{\bfp^2} \left[ (\bfa^2-\bfb^2)\sigma + \frac{\bfa}{\sin^2\theta}\eta_\phi - \frac{\bfb}{\cos^2\theta}\eta_\psi \right]^2 + \frac{1}{\bfp^2}[\bfa\bfb \sigma + \bfb\eta_\phi+ \bfa \eta_\psi]^2\\
        &=\eta_\theta^2 +B^2+\tilde{C}^2,
    \end{align}
    where we introduced the renormalized momentum term 
    \begin{align}
        \tilde{C} = \frac{1}{\bfp^2}(\bfa\bfb\sigma+\bfb\eta_\phi+\bfa\eta_\psi) = \frac{\bfp\sqrt{\Delta_r}}{r}C
    \end{align}

Expanding the $B$ term we see that
    \begin{align}
        \bfp^2B^2&=\cos^2\theta\sin^2\left[ (\bfa^2-\bfb^2)\sigma + \frac{\bfa}{\sin^2\theta}\eta_\phi - \frac{\bfb}{\cos^2\theta}\eta_\psi \right]^2\\
        &=\sin^2\theta\cos^2\theta(\bfa^2-\bfb^2)^2\sigma^2+\bfa^2\cot^2\theta\,\eta_\phi^2+\bfb^2\tan^2\theta\,\eta_\psi^2\\
        &+2\bfa\cos^2\theta(\bfa^2-\bfb^2)\sigma\eta_\phi -2\bfb\sin^2\theta(\bfa^2-\bfb^2)\sigma\eta_\psi -2\bfa\bfb\eta_\phi\eta_\psi.
    \end{align}
and for the $C$ term we see that 
\begin{align}
    \frac{r^2}{\bfp^2 \Delta_r}C^2&=(\bfa\bfb\sigma+\bfb\eta_\phi+\bfa\eta_\psi)^2=\bfa^2\bfb^2\sigma^2+\bfb^2\eta_\phi^2+\bfa^2\eta_\psi^2+2\bfa\bfb^2\sigma\eta_\phi+2\bfa^2\bfb\sigma\eta_\psi +2\bfa\bfb\eta_\phi\eta_\psi
\end{align}

Now adding $B^2+\tilde{C}^2=\frac{\bfp^2B^2+\bfp^2\tilde{C}^2}{\bfp^2}$ we see that the $\eta_\phi,\eta_\psi$terms come in as $-2\bfa\bfb\bfp^2\eta_\phi\eta_\psi$ and $2\bfa\bfb\bfp^2\eta_\phi\eta_\psi$ thus cancel. We then have that 
\begin{align}
    \bfp^2(B^2+C^2)&=\left[\sin^2\theta\cos^2\theta(\bfa^2-\bfb^2)^2+\bfa^2\bfb^2\right]\sigma^2\\
    &+\left(\bfa^2\cot^2\theta+\bfb^2\right)\eta_\phi^2+\left(\bfb^2\tan^2\theta +\bfa^2 \right)\eta_\psi^2\\
    &+(2\bfa \cos^2\theta(\bfa^2-\bfb^2)+2\bfa\bfb^2)\sigma\eta_\phi\\
    &+(-2\bfb\sin^2\theta(\bfa^2-\bfb^2)+2\bfa^2\bfb)\sigma\eta_\psi  
\end{align}
Now we use some identities for each of the terms to clean up the expression. For $\eta_\phi^2$ and $\eta_\psi^2$ we find
\begin{align}
(\bfa^2\cot^2\theta+\bfb^2)\eta_\phi^2&=\frac{\bfp^2}{\sin^2\theta}\eta_\phi^2\\
(\bfb^2\tan^2\theta+\bfa^2)\eta_\psi^2&=\frac{\bfp^2}{\cos^2\theta}\eta_\psi^2   
\end{align}
For the $\sigma\eta_\phi$ and $\sigma\eta_\psi,\sigma\eta_\psi$ coefficients we see that
\begin{align}
    2\bfa(\cos^2\theta(\bfa^2-\bfb^2)+\bfb^2)\sigma\eta_\phi&= 2\bfa\bfp^2\sigma\eta_\phi\\
    2\bfb(-\sin^2\theta(\bfa^2-\bfb^2))+\bfa^2)\sigma\eta_\psi&=2\bfb\bfp^2\sigma\eta_\psi
\end{align}
For the $\sigma^2$ coefficient we use $\cos^4\theta+\sin^4\theta=1-2\cos^2\theta\sin^2\theta$ and we ultimately find
\begin{align}
    (\cos^2\theta\sin^2\theta((\bfa^2-\bfb^2)^2+\bfa^2\bfb^2)\sigma^2=\bfp^2(\bfa^2\sin^2\theta+\bfb^2\cos^2\theta).
\end{align}
Putting it all together we see that 
\begin{align}
    \mathscr{C}=\eta_\theta^2+\frac{\eta_\phi^2}{\sin^2\theta}+\frac{\eta_\psi^2}{\cos^2\theta}+\sigma^2(\bfa^2\sin^2\theta+\bfb^2\cos^2\theta)+2 (\bfa\eta_\phi+\bfb\eta_\psi)\sigma
\end{align}
which is regular in $(\bfa,\bfb)\to (0,0)$ and satisfies (\ref{eq:SchwarzschildCarterConstant}) in this limit.
\end{proof}

Finally we are now in a position where we can neatly compare Frolov/Stokjkovic's \cite{frolov2003particle} and Wu's \cite{wu2008separability} Carter constants away from the characteristic set $\Sigma_G$. We see that within our current conventions, 
\begin{align}
    \mathscr{C}_{Wu}-\mathscr{C}_{FS}&=\left[-\frac{\bfp^2}{\bfr^2}G +\mathscr{C}\right]-\left[-\frac{\bfp^2}{\bfr^2}G-\bfp^2\sigma^2 +\frac{\eta_\phi^2}{\sin^2\theta}+\frac{\eta_\psi^2}{\cos^2\theta}+\eta_\theta^2\right]\\
    &=\mathscr{C}+\bfp^2\sigma^2-\frac{\eta_\phi^2}{\sin^2\theta}-\frac{\eta_\psi^2}{\cos^2\theta}-\eta_\theta^2\\
    &=(\bfa^2+\bfb^2)\sigma^2+2(\bfa\eta_\phi +\bfb\eta_\psi)\sigma.
\end{align}
Where we used Proposition \ref{prop:RegularCarterConstant} to expand $\mathscr{C}$ in the final step.
This shift is clearly invariant under the $(r,\xi)$-flow and is exaclty the additive constant that is mentioned in \cite{wu2008separability}. The reason we inverted the sign in Wu's Killing-St\"ackel tensor was exactly to make both constants have the same sign when taking the Schwarzschild limit.
Finally we will show that $\mathscr{C}$ constitutes a positive definite quantity away from the zero section when considered on the classical hypersurface $\{\sigma=0\}$. This observation will be important as we will use $\mathscr{C}$ later on as a quadratic defining function for the source of the classical flow at $\mathcal{H}^+$.

\begin{proposition}
    Let $(t_*,\phi,\psi,\theta,\sigma,\eta_\phi,\eta_\psi,\eta_\theta)\in T^*(\mathbb{R}_{t_*}\times \mathbb{S}^3)\backslash o$. Let $(\bfa,\bfb)\in D_{\mfrm}$. Then $\mathscr{C}\geq0$ and moreover the inequality is strict on the classical hypersurface $\{\sigma\equiv 0\}$. 
\end{proposition}

\begin{proof}  Note first off that (with our slight abuse of notation) the expression for $\mathscr{C}$ in (\ref{eq:CartersConstantonCharacteristicSet}) stays effectively the same in both Boyer-Lindquist coordinates and in $*$-coordinates (\ref{eq:MPStartcoordinate2})-(\ref{eq:MPStarpsicoordinate2}), as only the radial momentum $\xi$ truly gets shifted when changing between these charts.
    Away from the polar rings on $\mathbb{S}^3$, i.e. away from $\theta=0$ and $\theta=\frac{\pi}{2}$ we have that $\mathscr{C}=\eta_\theta^2+B^2+\frac{r^2}{\bfp^2\Delta_r}C^2$ and thus the statement is immediate. For the purposes of this proof however it is more suitable to work in the form of $\mathscr{C}$ (\ref{eq:RegularCarterConstant}) that is manifestly regular in the limit $(\bfa,\bfb)\to (0,0)$. Indeed away from the poles this can also be factorized as
    \begin{align}\label{eq:FactorizedRegularCartersConstant}
        \mathscr{C} = \eta_\theta^2 +\left(\frac{\eta_\phi}{\sin\theta }+\bfa\sigma \sin\theta\right)^2 +\left(\frac{\eta_\psi}{\cos\theta}+\bfb\sigma \cos\theta\right)^2\geq 0.
    \end{align}
    Near the poles, where Hopf coordinates exhibit a degeneracy into two collapsing rings parametrized by $\phi\in [0,2\pi)$ and $\psi\in [0,2\pi)$, we work with regular coordinates. For $\theta=0$ let us set
    \begin{equation}
        y_1=\cos\phi\sin\theta,\;\;\;y_2=\sin\phi\sin\theta,
    \end{equation}
    with smooth dual variables $\eta_1,\eta_2$. We find that
    \begin{align}
        dy_1 &=\cos\phi\cos\theta d\theta -\sin\phi\sin\theta d\phi\\
        dy_2 &=\sin\phi\cos\theta d\theta +\cos\phi\sin\theta d\phi.
    \end{align}
   We then have that equality must hold for the canonical $1$-forms expressed in both coordinate charts: 
    \begin{align}
        \sigma dt_*+\eta_1dy_1+\eta_2dy_2+\eta_\psi d\psi =\sigma dt_*+ \eta_\theta d\theta +\eta_\phi d\phi +\eta_\psi d\psi. 
    \end{align}
    We compute by matching terms 
    \begin{align}
        \eta_\theta&=\cos\theta(\eta_1 \cos\phi +\eta_2\sin\phi)\\
        \eta_\phi&=\sin\theta(-\eta_1\sin\phi+\eta_2\cos\phi).
    \end{align}
    With these coordinates $\mathscr{C}$ transforms as 
    \begin{align}
        \mathscr{C}&=\eta_1^2+\eta_2^2 -(y_1\eta_1 +y_2\eta_2)^2 +2\bfa\sigma (-y_2\eta_1+y_1\eta_2)+\bfa^2\sigma^2(y_1^2+y_2^2)\\
        &+\left(\frac{\eta_\psi}{\sqrt{1-y_1^2-y_2^2}}+\bfb\sigma\sqrt{1-y_1^2-y_2^2}\right)^2
    \end{align}
    At the pole $y_1=y_2=0$ this then reduces as follows
    \begin{align}
        \mathscr{C}|_{\theta=0} = \eta_1^2+\eta_2^2+(\eta_\psi^2+\bfb\sigma)^2.
    \end{align}
    Similarly near $\theta=\frac{\pi}{2}$ we set 
    \begin{align}
        y_3=\cos\psi \cos\theta, \; y_4=\sin\psi\cos\theta,
    \end{align}
    with conjugate variables $\eta_3,\eta_4$
    and we find 
    \begin{align}
        dy_3 &=-\cos\phi\sin\theta d\theta -\sin\phi\cos\theta d\phi\\
        dy_4 &=-\sin\phi\sin\theta d\theta +\cos\phi\sin\theta d\phi.
    \end{align}
    and hence 
    \begin{align}
        \eta_\theta&=-\sin\theta(\eta_3\cos\psi +\eta_4\sin\psi)\\
        \eta_\psi&=\cos\theta(-\eta_3\sin\psi+\eta_4\cos\psi),
    \end{align}
    whence 
    \begin{align}
        \mathscr{C}&=\left(\frac{\eta_\phi}{\sqrt{1-y_3^2-y_4^2}}+\bfa\sigma\sqrt{1-y_3^2-y_4^2}\right)^2\\
        &+\eta_3^2+\eta_4^2 -(y_3\eta_3 +y_4\eta_4)^2 +2\bfb\sigma (-y_4\eta_3+y_3\eta_4)+\bfb^2\sigma^2(y_3^2+y_4^2),
    \end{align}
    and at at $y_3=y_4=0$ we have
    \begin{align}
        \mathscr{C}|_{\theta=\frac{\pi}{2}}=\eta_3^2+\eta_4^2+(\eta_\phi+\bfa\sigma)^2
    \end{align}
\end{proof}

\subsection{Generalized radial source at the horizon.}\label{subsec:GeneralizedRadialSourceatHorizon}

We will now first analyze the flow of $H_G$ at $\mathcal{H}^+$. Starting point here is truly to switch to Eddington-Finkelstein/star-like coordinates, and re-express $G$ accordingly. We will emphasize the nature of this coorinate system by using the star notation everywhere. Note that we still maintain the convention $\Sigma_G^\pm =  \Sigma_G\cap \{\pm \sigma>0\}$.\\ 
In star-like coordinates we have that the dual metric function $G$ is given by
\begin{align}
    G&=  \bfDelta \xi^2 +\eta_\theta^2   +2 \left[\frac{(r^2+\bfa^2)(r^2+\bfb^2)}{r^2}  \sigma + \frac{\bfa(r^2+\bfb^2)}{r^2}  \eta_{\phi_*} + \frac{\bfb(r^2+\bfa^2)}{r^2} \eta_{\psi_*}\right]\xi \nonumber \\
    &+ \frac{\sin^2\theta\cos^2\theta}{\bfp^2}\left[(\bfa^2 - \bfb^2) \sigma + \frac{\bfa}{\sin^2\theta} \eta_{\phi_*} -\frac{\bfb}{\cos^2\theta}\eta_{\psi_*}  \right]^2 + \frac{\bfr^2}{r^2 \bfp^2}\left[\bfa\bfb\sigma + \bfb \eta_{\phi_*} + \bfa \eta_{\psi_*}\right]^2 \label{eq:SemiclassicalFTSymbol2},
\end{align}
which can again be decomposed as 
\begin{align}
    G=\mathscr{G}+\mathscr{C}\nonumber
\end{align}
with  

\begin{align}
    \mathscr{G}&:=\bfDelta \xi^2 +\eta_\theta^2   +2 \left[\frac{(r^2+\bfa^2)(r^2+\bfb^2)}{r^2}  \sigma + \frac{\bfa(r^2+\bfb^2)}{r^2}  \eta_{\phi_*} + \frac{\bfb(r^2+\bfa^2)}{r^2} \eta_{\psi_*}\right]\xi \\
    \mathscr{C}&:= \eta_\theta^2 + \frac{\sin^2\theta\cos^2\theta}{\bfp^2}\left[(\bfa^2-\bfb^2)\sigma +\frac{\bfa}{\sin^2\theta}\eta_{\phi_*}-\frac{\bfb}{\cos^2\theta}\eta_{\psi_*}\right]^2+\frac{1}{\bfp^2}[\bfa\bfb\sigma+\bfb\eta_{\phi_*}+\bfa\eta_{\psi_*}]^2.
\end{align}

We will first turn our attention to null-bicharacteristics of $H_G$ falling \textit{into} the Myers-Perry black hole.

\begin{lemma}\label{lem:InfallingParticles}
    We have $\xi H_Gr\leq0$ on $\Sigma_G$  for $r<r_+$. Moreover, strict inequality holds on the classical hypersurface $\{\sigma=0\}$.
\end{lemma}
\begin{proof}
    Observe that 
    \begin{align}
        H_Gr&=\partial_{\xi}G=2 \Delta_r \xi +2\left[\frac{(r^2+\bfa^2)(r^2+\bfb^2)}{r^2}\sigma + \frac{\bfa(r^2+\bfb^2)}{r^2}\eta_{\phi_*}+\frac{\bfb(r^2+\bfa^2)}{r^2}\eta_{\psi_*}\right].\\
        &= 2\left[\Delta_r\xi + (\Delta_r+2\mfrm)\sigma + \frac{\bfa(r^2+\bfb^2)}{r^2}\eta_{\phi_*}+\frac{\bfb(r^2+\bfa^2)}{r^2}\eta_{\psi_*}\right].
    \end{align}
    Thus on $\Sigma_G$ we have that 
    \begin{align}
        0&=\underbrace{\bfDelta \xi^2   +2 \left[\frac{(r^2+\bfa^2)(r^2+\bfb^2)}{r^2}  \sigma + \frac{\bfa(r^2+\bfb^2)}{r^2}  \eta_{\phi_*} + \frac{\bfb(r^2+\bfa^2)}{r^2} \eta_{\psi_*}\right]\xi}_{:=\xi H_G r-\Delta_r\xi^2} +\eta_\theta^2  \nonumber \\
    &+ \frac{\sin^2\theta\cos^2\theta}{\bfp^2}\left[(\bfa^2 - \bfb^2) \sigma + \frac{\bfa}{\sin^2\theta} \eta_{\phi_*} -\frac{\bfb}{\cos^2\theta}\eta_{\psi_*}  \right]^2 + \frac{\bfr^2}{r^2 \bfp^2}\left[\bfa\bfb\sigma + \bfb \eta_{\phi_*} + \bfa \eta_{\psi_*}\right]^2
    \end{align}
which in turn implies that 

\begin{align}
    \xi H_G r \leq \Delta_r\xi^2 -\mathscr{C},
\end{align}
This is in fact negative, due to $\Delta_r\leq 0$ if $r\leq r_+$ and
$\mathscr{C}$ being positive definite on $\{\sigma=0\}$. 

\end{proof}

We can furthermore state the following more general result.

\begin{lemma}\label{lem:SignFlipInterior}
    For $r_-<r<r_+$ have $\pm H_G r<0$ in $\Sigma_G^\pm$. This holds true at $\{r=r_+\}$ as well and at the complement $N^*\partial \mathcal{H}^+\backslash\{\mathcal{O}\}$ in $\Sigma_G$.
\end{lemma}
\begin{proof}
    Notice that that $dr$ is timelike
    \begin{equation}
        |dr|^2=G(dr,dr)=\Delta_r<0 \text{ on } (r_-,r_+).
    \end{equation}
    Moreover $dr$ is \textit{past causal} for $r\leq r_+$ which follows at $r=r_+$ from
    \begin{align}
        G(dr,-dt_*+c dr)&=-G(dr,dt_*)+c\,G(dr,dr)|_{r=r_+}\\
        &=-\frac{(r_+^2+\bfa^2)(r_+^2+\bfb^2)}{r_+^2}<0,
    \end{align}
and can be extended to $r<r_+$ by continuity. Now let $\varpi \in T_{\partial\mathcal{H}^+}^*M$ we have by $H_Gr|_{\varpi}=\frac{1}{2} G(dr,\varpi)$ and the past light-like character of $dr$ that $G(dr,\varpi)<0$ unless $\varpi = c dr\in N^*\partial\mathcal{H}^+$ with $c\neq 0$. Hence we conclude that $H_Gr|_{\varpi}<0$.

\end{proof}

Hence, as we know that $dr$ is future timelike for $r_-<r<r_+$ we know that 
\begin{align}\label{eq:MontonicityinInterior}
    H_G r < 0 \;\text{ on } \Sigma^+\cap \{r_-<r<r_+\}.
\end{align}

For completeness, we note that the full Hamiltonian vector field is given by 
\begin{align}
    H_G&:=2\left[(\Delta_r+2\mfrm)\xi+\frac{\bfa\bfb}{r^2}[\bfa\bfb\sigma +\bfb\eta_{\phi_*}+\bfa\eta_{\psi_*}]\right]\partial_{t_*}\\
    &+2\left(\Delta_r \xi +\left[\frac{(r^2+\bfa^2)(r^2+\bfb^2)}{r^2}\sigma + \frac{\bfa(r^2+\bfb^2)}{r^2}\eta_{\phi_*}+\frac{\bfb(r^2+\bfa^2)}{r^2}\eta_{\psi_*}\right]\right)\partial_r \nonumber \\
    &-\left(\Delta_r'\xi^2 +2\left[\Delta_r'\sigma - 2\bfb^2 r^{-3}\eta_{\phi_*}-2\bfa^2r^{-3}\eta_{\psi_*}\right]\xi-\frac{2}{r^3}[\bfa\bfb\sigma + \bfb\eta_{\psi_*}+\bfa\eta_{\phi_*}]^2 \right)\partial_\xi\nonumber\\
    &+2\left[\frac{\bfa(r^2+\bfb^2)}{r^2}\xi+\frac{\bfb}{r^2}[\bfa\bfb\sigma +\bfb\eta_{\phi_*} +\bfa\eta_{\psi_*}]\right] \partial_{\phi_*} \\
    &+ 2\left[\frac{\bfb(r^2+\bfa^2)}{r^2}\xi+\frac{\bfa}{r^2}[\bfa\bfb\sigma +\bfb\eta_{\phi_*} +\bfa\eta_{\psi_*}]\right] \partial_{\psi_*}+H_{\mathscr{C}},
\end{align}
where we wrote out $H_{\mathscr{G}}$ explicitly and surpressed the flow on $\mathbb{R}_{t_*}\times\mathbb{S}^3$ using the Hamiltonian flow associated to Carter's constant $H_{\mathscr{C}}$.

\begin{lemma}
    At $\mathcal{H}^+$ we have $H_G r<0$ on $\Sigma^+\cap \{r\leq r_+\}\backslash N^*\mathcal{H}^+$.
\end{lemma}

\begin{proof}
    Observe that the inequality $H_Gr\leq 0$ on $\{r_-\leq r \leq r_+\}$ follows from continuity. At $r=r_+$ we have that $\Delta_{r}|_{r=r_+}=0$ and therefore we have $H_Gr=0$ if and only if 
    \begin{align}
        \frac{(r^2+\bfa^2)(r^2+\bfb^2)}{r^2}\sigma + \frac{\bfa(r^2+\bfb^2)}{r^2}\eta_{\phi_*} +\frac{\bfb (r^2+\bfa^2)}{r^2}\eta_{\psi_*}=0. \label{eq:FirstHomogeneousEquation}
    \end{align}
    Combined with $G=0$ we see that $\mathscr{C}=0$. Due to linear independence this implies then that $\eta_\theta=0$ and combined with (\ref{eq:FirstHomogeneousEquation}) this we have the homogeneous system
    \begin{align}
        \frac{\cos^2\theta\sin^2\theta}{\bfp^2}\left[(\bfa^2-\bfb^2)\sigma +\frac{\bfa}{\sin^2\theta}\eta_{\phi_*}-\frac{\bfb}{\cos^2\theta}\eta_{\psi_*} \right]&=0\\
        \frac{1}{\bfp^2}\left[\bfa\bfb\sigma+\bfb\eta_{\psi_*} +\bfa\eta_{\phi_*} \right]&=0.
    \end{align}
    In this, we recover the the expressions derived from co-frame that diagonalized $g^{-1}_{\bfa,\bfb,\mfrm}$, i.e. equations (\ref{eq:MatrixEntry1})-(\ref{eq:MatrixEntry3}). This implies that the only solution is the trivial one i.e. $(\sigma, \eta_{\phi_*},\eta_{\psi_*})=(0,0,0)$. Combined at $r=r_+$, this restricts us precisely to
    \begin{align}
    N^*\mathcal{H}^+:=\{t_*,r=r_+,\phi_*,\psi_*,\theta,\sigma=0,\xi, \eta_{\phi_*}=0, \eta_{\psi_*}=0, \eta_{\theta}=0\}.
\end{align}
\end{proof}

We now introduce the generalized radial set at the event horizon conform our previous Definition \ref{def:GeneralizedRadialSet} and its associated projective coordinates. We will then show that they will act as a sink to the (renormalized) Hamiltonian flow $H_G$. 

\begin{definition}
    We introduce the inward/outward component of the co-normal bundle at $\mathcal{H}^+$ by 
    \begin{align}
        \mathcal{R}^\pm_{\mathcal{H}^+}&:=N^*\mathcal{H}\cap \{\mp\xi>0\}\cap \chr(G),
    \end{align}
    which is given explicitly by 
    \begin{align}
        \mathcal{R}^{\pm}_{\mathcal{H}^+}&=\{r=r_+,\theta,\phi,\psi,\sigma=0,\eta_\theta=0,\eta_\phi=0, \eta_\psi=0,G=0,\mp\xi>0\}.
    \end{align}
   Moreover we denote with $\partial \mathcal{R}^\pm$ the boundary of this set at fibre infinity. In particular, introducing projective coordinates in the $\xi$-direction, with 
   \begin{align}
       \hat{\xi}:=\frac{1}{\xi},\; \hat{\sigma}:=\frac{\sigma}{\xi},\; \hat{\eta}_{\phi_*} = \frac{\eta_{\phi_*}}{\xi},\;\hat{\eta}_{\psi_*} = \frac{\eta_{\psi_*}}{\xi},\; \hat{\eta}_\theta=\frac{\eta_\theta}{\xi}.
   \end{align}
   then the boundary $\partial\mathcal{R}^\pm$ can be identified with 
   \begin{align}
       \partial\mathcal{R}^\pm&:=\{t_*,r=r_+,\phi_*,\psi_*,\theta,\hat{\sigma}=0,\hat{\xi}=0,\hat{\eta}_{\phi_*}=0,\hat{\eta}_{\psi_*}=0,\hat{\eta}_\theta=0\}.
   \end{align}
\end{definition}

{To analyze} the flow in these coordinates observe that $\partial_\xi$ transforms as
\begin{align}
    \partial_\xi \mapsto -\hat{\xi}^2\partial_{\hat{\xi}}-\hat{\xi}\hat{\sigma}\partial_{\hat\sigma}-\hat{\xi}\hat\eta_{\phi_*}\partial_{\hat\eta_{\phi_*}}-\hat\xi\hat{\eta}_{\psi_*}\partial_{\hat\eta_{\psi_*}}-\hat{\xi}\hat{\eta}_{\theta}\partial_{\hat\eta_{\theta}}.
\end{align}

To this end, let us introduce the fibre-radial vectorfield

\begin{align}
    \Xi:= \hat\xi\partial_{\hat\xi} + \hat\sigma\partial_{\hat\sigma}+\hat{\eta}_{\phi_*}\partial_{\hat\eta_{\phi_*}}+\hat{\eta}_{\psi_*}\partial_{\hat\eta_{\psi_*}}+\hat{\eta}_\theta\partial_{\hat\eta_\theta},
\end{align}
which is the most important quantity to keep track of in the upcoming calculations.\\
We then see that the rescaled Hamiltonian vector field $\hat{H}_G:= \hat{\xi}H_G$ in these coordinates can be expressed as 
\begin{align}\label{eq:RescaledHamiltonianVField}
    \hat{H}_G=\,&2\left[\Delta_r+2\mfrm + \frac{\bfa\bfb}{r^2}[\bfa\bfb\hat{\sigma} +\bfb\hat{\eta}_\phi+\bfa \hat{\eta}_\phi]\right]\partial_{t_*}\nonumber\\
    &+2\left(\Delta_r+\left[\frac{(r^2+\bfa^2)(r^2+\bfb^2)}{r^2}\hat\sigma + \frac{\bfa(r^2+\bfb^2)}{r^2}\hat\eta_\phi +\frac{\bfb(r^2+\bfa^2)}{r^2}\hat{\eta}_\psi\right]\right)\partial_r\nonumber\\
    &+\left(\Delta_r' + 2[\Delta_r' \hat\sigma-2\bfb^2 r^{-3} \hat{\eta}_\phi -2\bfa^2 r^{-3}\hat\eta_\psi]-\frac{2}{r^3}[\bfa\bfb\hat\sigma+\bfb\hat\eta_{\psi_*}+\bfa\hat\eta_{\phi_*}]\hat{\xi}\right)\Xi\nonumber\\
    &+2\left[\frac{\bfa(r^2+\bfb^2)}{r^2}+\frac{\bfb}{r^2}[\bfa\bfb\hat\sigma +\bfb\hat\eta_{\phi_*} +\bfa\hat\eta_{\psi_*}]\right] \partial_{\phi_*}\nonumber \\
    &+ 2\left[\frac{\bfb(r^2+\bfa^2)}{r^2}+\frac{\bfa}{r^2}[\bfa\bfb\hat\sigma +\bfb\hat\eta_{\phi_*} +\bfa\hat\eta_{\psi_*}]\right] \partial_{\psi_*}+H_{\mathscr{\hat C}}.
\end{align}

Now observe that at $\partial\mathcal{R}^+$ we have that
\begin{align}\label{eq:SigmaDefiningQuadratic}
    \hat{H}_G \hat{\sigma}^2|_{\partial \mathcal{R}^+}= 2\Delta_r'|_{r=r_+}\hat\sigma^2+\mathcal{O}_3
\end{align}
and that 
\begin{align}\label{eq:CarterDefiningQuadratic}
    \hat{H}_G (\hat \xi^2 \mathscr{C})&= 2\Delta'_r|_{r=r_+}\hat\xi^2\mathscr{C} + \mathcal{O}_3
\end{align}
where $\mathcal{O}_3$ represents all terms that vanish cubically at the zero section of $T^*_{\mathcal{R}^\pm}X$. Here we used that $\partial_r\mathscr{C}=0$ which implies 
\begin{align}
     \hat{H}_{\mathscr{C}}(\hat{\xi}^2\mathscr{C})=\hat\xi H_{\mathscr{C}}(\hat{\xi}^2\mathscr{C})=\hat{\xi}^2 \mathscr{C} H_{\mathscr{C}}(\hat{\xi})=0 
\end{align}
and the fact that 
\begin{align}
    \hat{\xi}^2\mathscr{C}=\hat{\eta}_\theta^2 + \frac{\cos^2\theta\,\sin^2\theta}{\bfp^2}\left[(\bfa^2-\bfb^2) \hat{\sigma} +\frac{\bfa}{\sin^2\theta} \hat{\eta}_{\phi_*}-\frac{\bfb}{\cos^2\theta}\hat{\eta}_{\psi_*} \right]^2+\frac{1}{\bfp^2}[\bfa\bfb \hat\sigma +\bfb\hat{\eta}_{\phi_*} + \bfa\hat{\eta}_{\psi_*}]^2
\end{align}
is $2$-homogeneous in $(\hat\sigma,\hat{\eta}_\phi,\hat{\eta}_\psi,\hat{\eta}_\theta)$.

\begin{lemma}
    (Quadratic defining function of $\partial \mathcal{R}^\pm$). There exists a neighbourhood $U$ of $\partial \mathcal{R}^{\pm}_{\mathcal{H}^+}$ such that the function 
    \begin{align}
        \mathfrak{q}:=\hat{\sigma}^2 + \hat{\xi}^2\mathscr{C}
    \end{align}
    is smooth on $U$, a quadratic defining function of $\partial \mathcal{R}^+$ inside of $U$ and satisfies $H_G \mathfrak{q}\geq 2(r_+-\bfa^2\bfb^2r^{-3}) \mathfrak{q}$ on $U$.
\end{lemma}

\begin{proof}
    Recall that $\partial_r G = H_G \xi=-\Delta'_r|_{r=r_+} \xi^2=2(r_+-\bfa^2\bfb^2r_+^{-3})\xi^2$ at $N^*\mathcal{H}^+=\{r=r_+,\sigma,\eta_\phi=\eta_\psi=\eta_\theta=0\}$, hence $\partial_r(\hat{\xi^2 G)\neq 0}$ at $\partial\mathcal{R}^\pm_{\mathcal{H}}$. Then by the implicit function theorem we can locally express $r$ as a smooth function $r=r(\hat{\sigma},\phi_*,\psi_*,\theta,\hat{\eta}_{\phi_*},\hat{\eta}_{\psi_*},\hat{\eta}_\theta)$ at $T^*\mathbb{S}^3$ and $r=r_+$ at $\partial \mathcal{R}^\pm_{H^+}$. Now $\mathfrak{q}$ has a positive definite Hessian at $\{\hat\sigma=\hat{\eta}_{\phi_*}=\hat{\eta}_{\psi_*}=\hat{\eta}_{\theta}=0\}$ as $\hat{\xi}^2\mathscr{C}$ vanishes quadratically. The final claim follows directly from (\ref{eq:SigmaDefiningQuadratic}) and (\ref{eq:CarterDefiningQuadratic}), and the fact that $\Delta'_r|_{r=r_+}=2(r_+-\bfa^2\bfb^2r_+^{-3})$. Recall also that $\Delta'_r|_{r=r_+}>0$ as $r_+^2>|\bfa\bfb|$ if we choose the neighbourhood $U$ of $\partial \mathcal{R}^\pm$ small enough.
    
\end{proof}

To see now that $\mathcal{R^\pm}$ are a source of the rescaled Hamiltonian flow, we now only need to still observe that 
\begin{align}
    \hat{\xi}^{-1}\hat{H}_G \hat{\xi}= \Delta_r'|_{r=r_+}=2(r_+-\bfa^2\bfb^2r_+^{-3})>0,
\end{align}
at $\partial\mathcal{R}^{\pm}$. Finally, we compute here as well that 
\begin{align}\label{eq:HGTexpression}
    H_G t_*|_{\mathcal{R}^+} =4\mfrm>0.
\end{align}

We have thus shown 
\begin{theorem}
    Let $\mfrm>0$ and $(\bfa,\bfb)\in D_{\mfrm}$, and $(M^\circ,g_{\bfa,\bfb,\mfrm})$ be a Myers-Perry space time as defined by (\ref{eq:EddingtonFinkelsteinMPMetric}) and (\ref{eq:FullBaseManifold}). Then for $\partial\mathcal{R}^\pm_{\mathcal{H}^+}$ are a generalized radial source in the sense of Definition \ref{def:GeneralizedRadialSet} for the Hamiltonian flow $H_G$ of the dual metric function associated to $\bfr^2g^{-1}_{\bfa,\bfb,\mfrm}$. In particular, the coefficients $\alpha_\infty$, $\alpha_U$ in $(R.1)$ and $(R.3)$ of Definition \ref{def:GeneralizedRadialSet} is given by 
    \begin{align}\label{eq:InftyThresholdCoefficient}
        \alpha_U=2\alpha_\infty = 2\Delta'_r|_{r=r_+}=8\mfrm \kappa,
    \end{align}
    where we recall the surface gravity $\kappa $ from(\ref{eq:SurfaceGravityIntro}).
\end{theorem}

\subsection{Convexity of the flow at $\mathcal{H}^+$ and $i^0$}

We now proceed to show that the flow of $H_G$ satisfies $r$-convexity properties near $\mathcal{H}^+$ and $i^0$. As the trapped set $\Gamma$ where $H_G$ will be hyperbolic will exhibit a critical point for $H_G$, we will thus exclude the possibility that $\Gamma$ lies arbitrarily close to $\mathcal{H}^+$ or on the boundary at infinity $\partial M|_\infty$, i.e. we guarantee that if  $\Gamma$ exists, we have that $\Gamma\subset \tilde{M}$, in the region $\{r_+<r<\infty\}$.\\
As these arguments mostly constrain again to the $(r,\xi)$-flow of $\mathscr{G}$ it is helpful to introduce the following radial potential 
\begin{align}
    V = -\Delta_r^{-1}(A^2-C^2),
\end{align}
and its rescaled $r$-derivative as the quantity
\begin{align}\label{eq:RescaledPotentialrDerivative}
    \Psi:=-\Delta_r^2\,\partial_rV.
\end{align}
such that we have 
\begin{align}
    \mathscr{G}&=\Delta_r\xi^2+V\\
    \mathscr{C}&=\eta_\theta^2+B^2+\frac{1}{\bfp^2}C^2.
\end{align}
As such the Hamiltonian flow to $\mathscr{G}$ is given in Boyer-Lindquist coordinates as 

\begin{align}
    &H_{\mathscr{G}}=\Delta_r'\,\xi^2\partial_\xi+\frac{\partial V}{\partial_r}\partial_\xi -\frac{\partial V}{\partial \sigma}\partial_t - \frac{\partial V}{\partial \eta_\phi}\partial_\phi -\frac{\partial V}{\partial \eta_\psi}\partial_\psi-2\Delta_r \xi \partial_r\nonumber\\
    &=2\left[\frac{A}{r^2\Delta_r}((r^2+\bfa^2)(r^2+\bfb^2)\partial_t+\bfa(r^2+\bfb^2)\partial_\phi+ \bfb(r^2+\bfa^2)\partial_\psi) - \frac{C}{\Delta_r}\frac{\sqrt{\Delta_r}}{r}(\bfa\bfb\partial_t + \bfb \partial_\phi + \bfa\partial_\psi)
    \right]\nonumber\\
    &\;\;\;\;\;\;\;\;\;\;\;\;\;\;\;\;\;\;\;\;\;\;\;\;\;\;\;\;\;\;\;\;\;\;\;\;\;\;\;\;\;\;\;\;\;\;\;\;\;\;\;\;\;\;\;\;\;\;\;\;\;\;\;\;\;\;\;\;\;\;\;\;\;\;\;\;\;\;\;\;\;\;\;\;\;\;\;\;\;\;\;\;\;\;\;-2\Delta_r \xi \partial_r +\left[\Delta'_r \xi^2 -\frac{\Psi}{\Delta_r^2} \right]\partial_\xi.\label{eq:RadialHamiltonianVectorfield}
\end{align}

We then proceed to first show $r$-convexity near $\mathcal{H}^+$.

\begin{lemma}\label{lem:rConvexityonCharacteristicSet}(r-convexity near $\{r=r_+\}$)
   There exists $\epsilon_+>0$ on $\mfrm$ and $(\bfa,\bfb)\in D_{\mfrm}$ such that for $r\in (r_+,r_++\epsilon_+)$, $H_Gr=0$ implies that $H^2_Gr<0$.
\end{lemma}

\begin{proof}
    (Condition for non-trapped $\pi\circ \gamma$) Suppose that $H_G r =-2\Delta_r \xi=0$ then we have that $\xi=0$. Thus, working on $\Sigma_G$, we have 
    \begin{align}
        \Delta_r^{-1}A^2&=\eta_\theta^2+B^2+\Delta_r^{-1}\left(1+\frac{r^2}{\bfp^2}\right)C^2
    \end{align}
    From this we can derive three important inequalities, that are in general valid whenever we are on the set $\Xi:=\{\xi=0,\: G\leq 0\}\subset T^*M$ and on the set $\{G\leq 0\}\cap \{r>r_+\}$\footnote{As $\Delta_r\xi^2\geq 0$ whenever $r>r_+$.}:
    \begin{align}
        &\Delta_r^{-1}A^2\geq B^2 \label{eq:BAInequality}\\
        &A^2\geq C^2 \label{eq:CAInequality1}\\
        &A^2\geq\frac{r^2}{\bfp^2}C^2\label{eq:CAInequality2}.
    \end{align}
    Note that if $A=0$ we therefore must have that $B=0$ and $C=0$ (away from the poles), which would imply $\sigma=\eta_\theta=\eta_\phi=\eta_\psi=0$, following from the invertibility of (\ref{eq:MatrixEntry1})-(\ref{eq:MatrixEntry3}) is matrix equations relating $(A,B,C)$ to $(\sigma,\eta_\phi,\eta_\psi)$ i.e. we are at the $0$-section. Hence we may assume $A\neq 0$. Similarly one quickly sees that near the poles we have that  Now we calculate $H^2_G r$ as 
    \begin{align}
        H^2_G r&=H_G(-2\Delta_r\xi)\nonumber\\
        &=2\xi^2\Delta_r\Delta_r'-2\left(\Delta_r\Delta'_r\xi^2 - \frac{2\Psi}{\Delta_r}\right)\nonumber\\
        &=\frac{2\Psi}{\Delta_r
        }
    \end{align}
    Thus if we can show that $\Psi<0$ holds uniformly under the limit $r\to r_+$ we are done. To deal with uniformity, notice that $G$ is $1$-homogeneous under dilatation under the co-fibres of $T^*M$, so without loss of generality, we may work on the co-sphere bundle 
    \begin{align}
       S^*M&=\left\{(p,\zeta)\in T^*M \;|\; ||\zeta||=1  \right\}.
    \end{align}

    Now we need an explicit expression for $\Psi$. Notice that 
    \begin{align}
    \Psi&:=-\Delta_r^2 \partial_r V =\Delta_r^2 \left(\frac{2\Delta_r(A\partial_r A - C\partial_r C)-(A^2-C^2)\Delta_r'}{\Delta_r^2} \right)
\end{align}
Furthermore we note the following useful $\partial_r$-identities for $A$ and $C$.

\begin{align}\label{eq:PartialAIdentity}
    \partial_r A &= 2 r \sigma-2r^{-3}\left[\bfa^2\bfb^2\sigma + \bfa\bfb^2\eta_\phi + \bfa^2\bfb \eta_\psi\right]\nonumber\\
    &=2r\sigma - 2r^{-3}\bfa\bfb\frac{r}{\sqrt{\Delta_r}}\frac{\sqrt{\Delta_r}}{r}[\bfa\bfb\sigma + \bfb\eta_\phi + \bfa\eta_\psi]\nonumber\\
    &=2 \left(r\sigma-\frac{\bfa\bfb}{r^2\sqrt{\Delta_r}}C\right)\\
    &\nonumber\\
    \partial_r C &= \partial_r\left(\frac{\sqrt{\Delta_r}}{r}\right)[\bfa\bfb\sigma + \bfb\eta_\phi + \bfa\eta_\psi]\nonumber\\
    &=\partial_r\left(\frac{\sqrt{\Delta_r}}{r} \right)\frac{r}{\sqrt{\Delta_r}}\frac{\sqrt{\Delta_r}}{r}[\bfa\bfb\sigma + \bfb\eta_\phi + \bfa\eta_\psi]\nonumber\\
    &=\partial_r\left(\frac{\sqrt{\Delta_r}}{r}\right)\frac{r}{\sqrt{\Delta_r}}C\nonumber\\
    &=\left[\partial_r \log\left(\frac{\sqrt{\Delta_r}}{r}\right)\right]C\nonumber\\
    &=\left(\frac{\Delta_r'}{2r\Delta_r}-\frac{1}{r} \right)C.
\end{align}
We then find by straightforward computation
\begin{align}
A^{-2}\Psi&=2\Delta_r \left[\frac{2r\sigma}{A}-\frac{2\bfa\bfb}{r^2\Delta_r}\frac{C}{A}-\partial_r\left(\frac{\sqrt{\Delta_r}}{r}\right)\frac{r}{\sqrt{\Delta_r}}\frac{C^2}{A^2} \right]-\left(1-\frac{C^2}{A^2}\right)\Delta_r'\label{eq:ExplicitPsiExpression}
\end{align}
Noting that $C\propto \frac{\sqrt{\Delta_r}}{r}$ in the limit $r\to 0$ we have that by ordering (\ref{eq:ExplicitPsiExpression}) in terms of $\Delta_r$ (which vanishes in the limit $r\to r_+$) we see that
\begin{align}
    A^{-2}\Psi\leq \mathcal{O}(\Delta_r^{1/2})-\left(1-\frac{C^2}{A^2} \right)\Delta_r',
\end{align}
where we used (\ref{eq:CAInequality1}). Now given that the we work on $S^*M$, which is compact in its fibres, we conclude that the limit
\begin{align}
    \lim_{r \searrow r_+} A^{-2}\Psi\leq \lim_{r \searrow r_+} \mathcal{O}(\Delta_r^{1/2})-\left(1-\frac{C^2}{A^2} \right)\Delta_r'=-\Delta_r'(r_+)<0
\end{align}
holds uniformily on $S^*M$. Here we also used the fact that $r_+^2>|\bfa\bfb|$ whenever $(\bfa,\bfb)\in D_{\mfrm}$, whence $\Delta'(r_+)>0$.

\end{proof}

\begin{remark}
    We remark that the observation from (\ref{eq:CAInequality1}) that $C^2<A^2$ on $\Sigma_G \cap \{\xi=0\}$ is one of the most important observations in this entire chapter as it allows us to reduce the additional dimensionality of Myers-Perry whenever restricting to the characteristic set with vanishing $\xi$-momentum. As it turns out, the trapped set $\Gamma$ will be a subset of $\Xi$ as well.
\end{remark}

Similarly, we proof $r$-convexity as $r\to\infty$. To this end let $\rho_0=r^{-1}$ be the smooth boundary defining function of $i_0$.
We then have the following lemma.
\begin{lemma}\label{lem:ConvexityNearInfinity}
    ($r$-convexity near $i^0$) There exists $\epsilon>0$ such that for $\varpi \in \chr(G)$ with $\rho_0<\epsilon$ and $H_G\rho_0=0$ we have $H^2_G\rho_0<0$.
\end{lemma}

\begin{proof}
    We first calculate $H_G\rho_0$. Using (\ref{eq:RadialHamiltonianVectorfield}) one quickly calculates that 
    \begin{align}
        H_G\rho_0=-2\Delta_r\xi\partial_r\left(\frac{1}{r}\right)=\frac{2\Delta_r\xi}{r^2}=2\Delta_r\xi\rho_0^2,
    \end{align}
    hence from $H_G\rho_0=0$ we once more find that $\xi=0$.
    Then again taking $H^2_G\rho_0$ and using the constraint $\xi = 0$ we see that 
    \begin{align}
        H_G^2\rho_0&=H_G\left(\frac{2\Delta_r\xi}{r^2}\right)\\
        &=\underbrace{-\left.2\Delta_r\xi\partial_r\left(\frac{2\Delta_r\xi}{r^2}\right)\right|_{\xi=0}}_{=0}+\left.\left(\underbrace{\Delta_r'\xi^2}_{=0} -\frac{\Psi}{\Delta_r^2}\right)\right|_{\xi=0}\partial_\xi\left(\frac{2\Delta_r\xi}{r^2}\right)\\
        &=-\frac{2\Psi}{r^2\Delta_r}.
    \end{align}
    So if we can prove $\Psi>0$ as $r\to \infty$ we are done. Using (\ref{eq:sigmaInverted}) we see that in the long $r$-limit, 
    \begin{align}
        \frac{\sigma}{A}\propto r^{-2}.
    \end{align}
    Moreover we also know that $\frac{C}{A}\propto \mathcal{O}(r^{-1})$ as $r\to \infty$. Thus plugging this into (\ref{eq:ExplicitPsiExpression}) we see that 
    \begin{align}
        A^{-2}\Psi \propto 4 r\Delta_r r^{-2}+\mathcal{O}(r^{-1})-\left(1-\mathcal{O}(r^{-2}\right))\Delta_r'\propto (4r^{-1}\Delta_r-\Delta_r')+\mathcal{O}(r^{-1})>0,\; \text{ as }\; r\to \infty.
    \end{align}
    We conclude that $\lim_{r\to\infty} \Psi >0$, completing the proof.
\end{proof}

\subsection{Phase space dynamics away from the event horizon and spatial infinity.}

We now describe the dynamics of $G$ in the region $\{r+\epsilon_+<r<r_{1}\}$ where $r_1$ is large enough such that $H_G r = 0$, $r>r_1$ implies that  {$H^2_G r>0$} by Lemma \ref{lem:ConvexityNearInfinity}. In this chapter, the observations made regarding equations (\ref{eq:BAInequality})-(\ref{eq:CAInequality2}) on the set $\Sigma_G\cap \{\xi=0\}$ will be paramount for identifying the trapped set $\Gamma$ as defined in Definition \ref{def:TrappedSet}, and reducing the degrees of freedom on it. Additionally the following basic observation,  which is basically a corollary of the GM-AM inequality, will be important to consider certain limits and asymptotics.
\begin{lemma}\label{lem:GeometricMeanDeltar}
    Let $(\bfa,\bfb)\in D_{\mfrm}$. We then have that 
    \begin{align}
        r^3\Delta_r'\geq 2r^2\Delta_r, \text{ for all } r\in(r_+,\infty).
    \end{align}
\end{lemma}
\begin{proof}
    Notice first of all that for the arithmetic mean of $r_+^2$ and $r_-^2$ it holds that 
    \begin{align}
        r_-^2\leq \frac{r_+^2+r_-^2}{2}\leq r_+^2, 
    \end{align}
    and recall that when $(\bfa,\bfb)\in D_{\mfrm}$, we have by equation (\ref{eq:ExtremalityMP}) and (\ref{eq:EventHorizonMP})-(\ref{eq:CauchyHorizonMP}) that
    \begin{align}
        \frac{r_+^2+r_-^2}{2}=\frac{1}{2}(2\mfrm-\bfa^2-\bfb^2) \geq |\bfa\bfb|.
    \end{align}
    Now we claim that 
    \begin{align}\label{eq:GeometricMeanPolynomial}
        (r^2-r_+^2)(r-r_-)^2\leq \left(r^2-\frac{r_+^2+r_-^2}{2}\right)^2.
    \end{align}
    Indeed, this equation reduces readily to 
    \begin{align}
        r_1^2r_2^2 \leq \left(\frac{r_1^2+r_2^2}{2}\right)^2,
    \end{align}
    whence 
    \begin{align}
        0\leq r_1^4+r_2^4-2r_1^2r_2^2 = (r_1^2-r_2^2)^2,
    \end{align}
    which is a perfect square.
    Finally recall that 
    \begin{align}
        2r^2\Delta_r=2(r^2-r_+^2)(r^2-r_-^2),\;\text{ and }\; 2\Delta'_r = 2\left(r-\frac{|\bfa\bfb|^2}{r^3}\right)
    \end{align}
    From our observation about the geometric mean of $r_+^2$ and $r_-^2$ in (\ref{eq:GeometricMeanPolynomial}) we then conclude that
    \begin{align}
        r^3\Delta'_r&= 2(r^4-|\bfa\bfb|^2)\nonumber\\
        &=2(r^2+|\bfa\bfb|)(r^2-|\bfa\bfb|)^2\nonumber\\
        &\geq 2(r^2-|\bfa\bfb|)^2\nonumber\\
        &\geq 2\left(r^2-\frac{r_+^2+r_-^2}{2}\right)^2\nonumber\\
        &\geq 2 (r^2-r_+^2)(r^2-r_-^2) = 2r^2\Delta_r,
    \end{align}
    as desired.
\end{proof}

We already described some of the symmetries of $H_G$ in prior sections. In summary, we have due to $\partial_t,
\partial_\phi,\partial_\psi$ being Killing and $\mathscr{C}$ being independent of $(t,r,\phi,\psi,\xi)$ that 
\begin{align}
    H_G\sigma=0,\; H_G\eta_\phi=0, H_G\eta_\psi=0,H_G\mathscr{C}=0,\;H_{\mathscr{C}}r=0,\; H_{\mathscr{C}}\xi=0,
\end{align}
which holds everywhere on $T^*M$. Due to this separation of variables, and $\sigma,\eta_\phi,\eta_\psi$ being constant along the flow of $H_{G}$ we can describe the dynamics of $H_{\mathscr{G}}$ by projecting to the $(r,\xi)$-variables where its flow can be studied seperately. We introduce the notation $z=(t,\omega)$, $\zeta=(\sigma,\eta_\omega)$ for the remaining variables, with $\omega=(\phi,\psi,\theta)$ and $\eta_\omega=(\eta_\phi,\eta_\psi,\eta_\theta)$, i.e. $(z,\omega)\in T^*(\mathbb{R}_t \times\mathbb{S}^3)$. Let us fix a point $\varpi^0=(r^0,z^0;\xi^0,\zeta^0)\in T^*M\cap\Sigma_G$, where $r^0>r_+$, and let
\begin{align}\label{eq:FixedrxiDynamics}
    \mathscr{G}^0_{(z^0,\zeta^0)}:=\Delta_r\xi^2+\Phi^0(r),\; \; \Phi^0_{(z^0,\zeta^0)}(r):=(V+\mathscr{C})|_{(z^0,\zeta^0)},
\end{align}
where we will omit the superscripts "$0$" when it does not cause confusion. Let $\pi:T^*M\to T^*(r_+,\infty)_r$. Then the integral curve of $H_{\mathscr{G}^0}$ through the point $(r^0,\xi^0)\in T^*(r_+,\infty)_r$ is then equal to $\pi\circ\gamma$, where $\gamma\subset\chr(G)\subset T^*M$ is the integral curve of $H_G$ through through $\varpi^0$. By assumption on $\varpi \in \chr(G)$ we must have that $\mathscr{G}^0(r^0,\xi^0)=0$, hence 
\begin{align}
    \mathscr{G}^0(r,\xi)=0\;\;\;\text{along} \;\;\;\pi\circ\gamma,
\end{align}
hence by (\ref{eq:FixedrxiDynamics}) we must have $\Phi^0(r)\leq 0$ along $\pi \circ \gamma$. Let us study the structure of the zero set $(\mathscr{G}^0)^{-1}(0)\subset T^*(r_+,\infty)_r$.
\begin{lemma}(Bicharacteristics $\pi\circ \gamma$ that are not trapped escape to the horizon $\{r=r_+\}$ or to spatial infinity.)
    Fix $\mfrm>0$ and let $(\bfa,\bfb)\in D_{\mfrm}$. Suppose that $(\Phi^0(r),\partial_r\Phi^0(r))\neq (0,0)$ for all $r>r_+$. Then $(\mathscr{G}^0)^{-1}(0)$ is a smooth $1$-dimensional submanifold, and in both directions along $\pi\circ \gamma$ we have either $r\to r_+$ or $r\to \infty$.
\end{lemma}
\begin{proof}
    We have $d\mathscr{G}^0=(\Delta_r'\xi^2+\partial_r\Phi^0)dr + 2 \Delta_r \xi \,d\xi$. This is non-zero whenever $\xi\neq 0$. Vice versa, if we have that $\xi=0$ for some $(r,\xi)\in \pi\circ \gamma$ then it also holds that $\Phi^0(r)=0$, because $\Phi^0(r)\leq 0$ on all $\pi\circ \gamma\subset T^*(r_+,\infty)_r$ as mentioned before. By assumption then, we should have that $\partial_r\Phi^0(r)\neq 0$ from which follows that $d\mathscr{G}^0\neq 0$. This proves that $(\mathscr{G}^0)^{-1}(0)$ is smooth and $1$-dimensional, and moreover implies that $H_{\mathscr{G}^0}$ never vanishes on $(\mathscr{G}^0)^{-1}(0)$. The proof is complete if we can show that $(\mathscr{G}^0)^{-1}(0)$ has no connected components. Assuming to the contrary that $\Phi^0(r)$ has a local minimum on $(\mathscr{G}^0)^{-1}(0)$, we obtain a contradiction with the lemma below.
\end{proof}

\begin{lemma}\label{lem:NoLocalMinimum} (No local minimum for $\Phi^0(r)$) Let $\mfrm>0$ and $(\bfa,\bfb)\in D_{\mfrm}$. Let $(r,0)\in \pi\circ\gamma\subset  T^*(r_+,\infty)_r$, where $\pi$ and $\gamma$ are given as before. Then, if $\partial_r\Phi^0(r)=0$ we have $\partial_r^2\Phi^0(r)<0$.
    
\end{lemma}

\begin{proof}
    By assumption of $(r,0)\in \pi \circ \gamma\subset T^*(r_+,\infty)_r$ we have $\xi=0$ and $\Phi^0(r)\leq 0$. From this it follows that 
    \begin{align}
        0\geq \Phi^0(r)=-\Delta_r^{-1} (A^2-C^2)+\mathscr{C}= -\Delta_r^{-1}(A^2-C^2)+\eta_\theta^2+B^2 +\frac{\bfp^2}{r^2}C^2.
    \end{align}
    From this it follows that $(r,0) \in \pi\left(\{\xi=0,G\leq 0\}\right)$, and hence that we are allowed to use estimates (\ref{eq:BAInequality})-(\ref{eq:CAInequality2}) from lemma \ref{lem:rConvexityonCharacteristicSet} as we have once again that 
    \begin{align}
        B^2+\frac{\bfp^2}{r^2} C^2\leq \mathscr{C}\leq \Delta_r^{-1}(A^2-C^2).
    \end{align}
    In particular, we note again that we work away from the zero section, i.e. $A\neq 0$.
    Now, $\partial_r\Phi^0(r)=\partial_rV=-\Delta_r\Psi =0$ gives (after some calculations, using (\ref{eq:ExplicitPsiExpression})) 
    \begin{align}\label{eq:SigmaExpression}
        \sigma=\frac{\bfa\bfb}{r^3\sqrt{\Delta_r}}C + \frac{\Delta_r'}{4r\Delta_r}A-\frac{1}{2r^2}\frac{C^2}{A}.
    \end{align}
Furthermore computing $\partial_r^2\Phi^0(r)$ we see that 
\begin{align}
    \partial_r^2\Phi^0(r)&=-\partial_r(\Delta_r^{-2}\Psi)=2\Delta_r'\Delta_r^{-3}\underbrace{\Psi}_{=0}-\Delta_r^{-2}\partial_r\Psi,
\end{align}
and thus that the sign for $\partial_r^2\Phi^0(r)$ is opposite that of $\partial_r\Psi$. Therefore, if we can show that $\partial_r\Psi>0$ we are done. Calculating $\partial_r\Psi$ we find that 
\begin{align}\label{eq:ExplicitPsirDerivative}
    \partial_r\Psi=2\Delta_r\left[4r^2\sigma^2-\frac{8\bfa\bfb}{r\sqrt{\Delta_r}}\sigma C + \frac{4\bfa^2\bfb^2}{r^4\Delta_r}C^2 + 2\sigma A + \frac{6\bfa\bfb}{r^3\sqrt{\Delta_r}}AC + \left(\frac{2\Delta_r'}{r\Delta_r} - \frac{3}{r^2}\right)C^2\right]-\Delta''_rA^2.
\end{align}

After some simple but tedious calculations one can write $\partial_r\Psi$ as a general mixed polynomial in $A$ and $C$ once one plugs back expression (\ref{eq:SigmaExpression}) into expression (\ref{eq:ExplicitPsirDerivative})
\begin{align}
    \partial_r\Psi=S_{A^2}\cdot A^2+S_{C^2}\cdot C^2 + S_{AC}\cdot AC + S_{\frac{C^3}{A}}\cdot \frac{C^3}{A}+ S_{\frac{C^4}{A^2}}\cdot \frac{C^4}{A^2}.
\end{align}
where the terms $S_{A^2},S_{C^2},S_{AC},S_{\frac{C^3}{A}}$ and $S_{\frac{C^4}{A^2}}$ are explicitly given by
\begin{align}
    S_{A^2}&=\frac{1}{2r\Delta_r}\left[\Delta_r'(r\Delta_r' + 2\Delta_r) - 2r\Delta_r\Delta_r'' \right]\label{eq:A2Coefficient}\\
    S_{C^2}&=\frac{2\Delta'_r}{r}-\frac{8\Delta_r}{r^2}\label{eq:C2Coefficient}\\
    S_{AC}&=\frac{16\bfa\bfb\sqrt{\Delta_r}}{r^3}\label{eq:ACCoefficient}\\
    S_{\frac{C^3}{A}}&=0 \;\text{ (full cancellation)}\label{eq:C3ACoefficient}\\
    S_{\frac{C^4}{A^2}}&=\frac{2\Delta_r}{r^2}\label{eq:C4A2Coefficient}.
\end{align}
To prove that $\partial_r\Psi>0$ using expressions (\ref{eq:A2Coefficient})-(\ref{eq:C4A2Coefficient}) we rescale $\partial_r\Psi$ by a factor of $\frac{2r\Delta_r}{A^2}$, and introduce a new variable $x=\frac{C}{A}$. We then obtain a rescaled polynomial $p(x)$ given by 
\begin{align}
    p(x)&=a(r) x^4 + b(r)x^3 + c(r)x^2 + d(r)x +e(r),
\end{align}
with
\begin{align}
    a(r)&=\frac{4\Delta_r^2}{r}\\
    b(r)&=0\\
    c(r)&=4\Delta_r'\Delta_r-16\frac{\Delta_r^2}{r}\\
    d(r)&=\frac{32\bfa\bfb\Delta_r^{3/2}}{r^2}\\
    e(r)&=r(\Delta_r')^2+2\Delta_r'\Delta_r-2r\Delta_r\Delta_r'',
\end{align}
which when considered on the domain of relevance $x\in(-1,1)$ following from $|x|=\left|\frac{C}{A}\right|<1$ (i.e. estimate (\ref{eq:CAInequality1}))  should be of the same signature as $\partial_r\Psi$, by virtue of the fact that $2r\Delta_r>0$ for $r>r_+$.  Upon inspection we see that $a(r)>0$ for all $r>r_+$. We can therefore reduce this polynomial once more by rescaling by $a(r)$ giving 
\begin{align}
 \bar{p}(x):=\frac{p(x)}{a(r)}=x^4+\frac{c(r)}{a(r)}x^2+\frac{d(r)}{a(r)}x+\frac{e(r)}{a(r)}.
\end{align}
This polynomial is clearly a depressed monic quartic, and we will refer to it colloquially as the monic depressed Myers-Perry quartic and we will show that $\bar p(x)>0$ for $x\in (-1,1)$. We first claim that 
\begin{align}
    \frac{e(r)}{a(r)}\geq 1 \text{ and } \frac{c(r)}{a(r)}\geq -2.
\end{align}

To this end, let us first estimate the quotient $\frac{(\Delta_r')^2}{2\Delta_r''\Delta_r}$:
\begin{align}
    \frac{(\Delta_r')^2}{2\Delta_r''\Delta_r}&=\frac{(r^3\Delta'_r)^2}{2r^4\Delta_r'' r^2\Delta_r}\\
    &=\frac{4(r^2-|\bfa\bfb|)^2(r^2+|\bfa\bfb|)^2}{4(r^4+3\bfa^2\bfb^2)(r^2-r_+^2)(r^2-r_-^2)}\label{eq:LittleQuotientEstimate0}\\
    &\geq \frac{(r^2+|\bfa\bfb|)^2}{r^4+3\bfa^2\bfb^2}\label{eq:LittleQuotientEstimate1}\\
    &=\frac{(r^4+2|\bfa\bfb|r^2+\bfa^2\bfb^2)}{r^4+3\bfa^2\bfb^2}\label{eq:LittleQuotientEstimate2}\\
    &>\frac{r^4 +3\bfa^2\bfb^2}{r^4 +3\bfa^2\bfb^2}\geq 1\label{eq:LittleQuotientEstimate3}.
\end{align}
Where we used the fact that $|\bfa\bfb|<\frac{r_+^2+r_-^2}{2}$ for $(\bfa,\bfb)\in D_{\mfrm}$ combined with Lemmas \ref{lem:GeometricMeanDeltar} to estimate (\ref{eq:LittleQuotientEstimate0}) into (\ref{eq:LittleQuotientEstimate1}), and the fact that $r^2>r_+^2>|\bfa\bfb|$ for estimate (\ref{eq:LittleQuotientEstimate2}) and beyond. Then we see that 
\begin{align}
    \frac{e(r)}{a(r)}&=\frac{r\Delta'_r}{2\Delta_r}+\underbrace{\frac{r^2(\Delta'_r)^2-2r^2\Delta_r\Delta_r''}{4\Delta_r^2}}_{> 0 \;\text{by equation (\ref{eq:LittleQuotientEstimate3})}}\\
    &\geq \frac{r\Delta'_r}{2\Delta_r}\\
    &=\frac{r^3\Delta'_r}{2r^2\Delta_r}\\
    &=\frac{r^3(2r^{-3}(r^2-|\bfa\bfb|)(r^2+|\bfa\bfb|))}{2(r^2-r_+^2)(r^2-r_-^2)}\\
    &=\frac{(r^2-|\bfa\bfb|)(r^2+|\bfa\bfb|)}{(r^2-r_+^2)(r^2-r_-^2)}>1,
\end{align}
where we use again Lemma \ref{lem:GeometricMeanDeltar} and $r^2>r_+^2>|\bfa\bfb|$. Similarly we see that
\begin{align}
    \frac{c(r)}{a(r)}&=\frac{4 \Delta_r'\Delta_r- 16r^{-1}\Delta_r^2}{\frac{4\Delta_r^2}{r}}\\
    &=\frac{r\Delta_r'}{\Delta_r}-4\\
    &=\frac{2(r^2-|\bfa\bfb|)(r^2+|\bfa\bfb|)}{(r^2-r_+^2)(r^2-r_-^2)}-4\\
    &>-2\label{eq:SecondLittleQuotientEstimate}
\end{align}

In fact it can be shown that in the long $r$-limit we have
\begin{align}
   \frac{c(r)}{a(r)}\searrow -2,\;\frac{d(r)}{a(r)}\searrow 0,\;\text{ and } \frac{e(r)}{a(r)}\searrow 1.
\end{align}
With this in mind, we complete the proof by showing that $\bar p(x)>\bar q(x)$, where $\bar q(x)$ is the $r$-limiting polynomial of $\bar p(x)$ given by
\begin{align}\label{eq:LimitingPolynomial}
    \bar q(x)=x^4-2x^2+1.
\end{align}
To this end, let us consider the difference $\Omega(x)=\bar p(x) -\bar q(x)$, given by
\begin{align}\label{eq:DifferenceEquation}
    \Omega(x)= \underbrace{\left(\frac{c(r)}{a(r)}+2\right)}_{>0 \; \text{by equation (\ref{eq:SecondLittleQuotientEstimate})}}\cdot x^2+\frac{d(r)}{a(r)}\cdot x + \underbrace{\frac{e(r)}{a(r)}-1}_{>0 \; \text{by equation (\ref{eq:LittleQuotientEstimate3})}}
\end{align}
If we can show that this polynomial is always positive, we know that the graph of $\bar p(x)$ should always lie above its $\bar q(x)$ and hence shows that in particular $\bar p(x)>0$ for $x\in (-1,1)$.
Given the coefficients in (\ref{eq:DifferenceEquation}) this would follow if the associated associated discriminant $D_\Omega$ of $\Omega$ is negative, given the  i.e.
\begin{align}
    D_\Omega &= \frac{d(r)^2}{a(r)^2}-4\left(\frac{c(r)}{a(r)}+2\right)\left(\frac{e(r)}{a(r)}-1\right)<0
\end{align}
Another tedious calculation shows that this is true, by virtue of the fact that after rescaling once more, we have that
\begin{align}\label{eq:MathematicaFile}
    \Delta_r^3 D_\Omega = -8((2\mfrm -\bfa^2-\bfb^2)^2-4\bfa^2\bfb^2)(2r^2-(2\mfrm -\bfa^2-\bfb^2))
\end{align}
which is in fact negative due to the assumption that $(\bfa,\bfb)\in D_{\mfrm}$, the associated extremality condition (\ref{eq:ExtremalityMP}) and the fact that $2r^2>2r_+^2>2\mfrm -\bfa^2-\bfb^2$ by equation (\ref{eq:EventHorizonMP}). Thus we have shown that $\bar p(x)>\bar q(x)>0$ for all $x\in (-1,1)$ proving the statement that $\partial_r\Psi>0$.
\end{proof}
\begin{remark}
    Lemma \ref{lem:NoLocalMinimum} is the central lemma to proving the hyperbolic trapping result. It basically asserts that the sign of the derivative of the rescaled potential, $\partial_r\Psi$, given through $\bar{p}(x)$ can be determined by comparing its value against a limiting polynomial (in this case $\overline{q}$ in (\ref{eq:LimitingPolynomial}) that always bounds both the long range asymptotics $r\to \infty$ using Lemma \ref{lem:GeometricMeanDeltar} as well as the extremal limit $(\bfa,\bfb)\to \partial D_{\mfrm}$ of $\partial_r\Psi$ from below, after multiple rescalings. Originally the function $\overline{q}$ was found by dynamically plotting $\partial_r\Psi$ for different values of $\bfa$ and $\bfb$ in both the limits $r\to r_+$ and $r\to r_-$ in Mathematica 11.3, before a rigorous attempt was made to analyze asymptotics using \ref{lem:GeometricMeanDeltar}, upon (\ref{eq:MathematicaFile}) was found using Mathematica for computing the discriminant of the difference.
\end{remark}

\begin{figure}[H]\label{fig:drPsiPlot}
\centering
\includegraphics[scale=0.8]{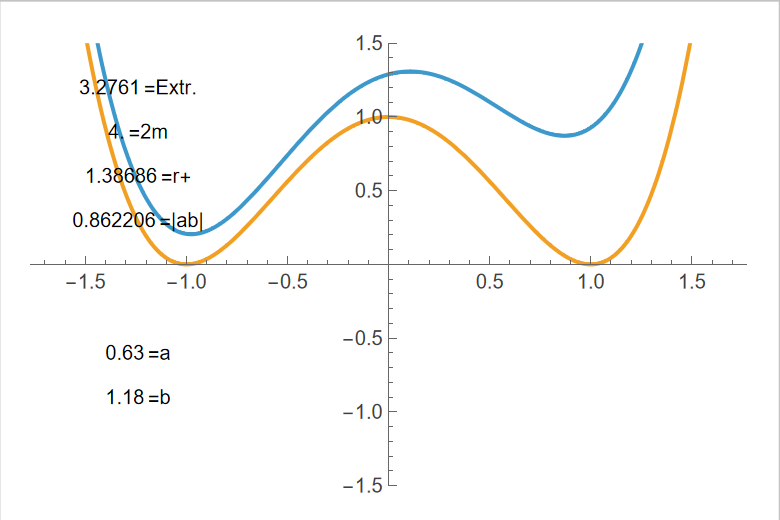}

\caption{Mathematica plot of the rescaled potential $\partial_r\Psi$ (blue) versus $\bar{q}(x)$ (orange) as its limiting polynomial. The Mathematica code can be found in Appendix $B$, values are displayed in the inset.}
\label{fig:myers-perry-forward}
\end{figure}

To conclude the study of the characteristic set we also prove the following property.
\begin{lemma}\label{lem:NondegenerateMax}    (Points $\varpi'\in T^*M$ with $\Phi^0(r')=0$ and $\partial_r\Phi^0(r')=0$ are trapped.)
    Let $\varpi'\in \chr(G)$ and $\pi\circ \gamma \subset \Sigma_G$ the projected integral curve of $H_G$ through $\varpi^0$. If $\Phi^0(r)$ attains a maximum at $r'$
    and $\Phi^0(r')=0$ then $\Phi^0(r)<0$ for all $\varpi \in \gamma\backslash\{\varpi^0\}\subset\Sigma_G.$ 
\end{lemma}
\begin{proof}
    Suppose there exists a $\varpi'$ along $\gamma$ such that 
    \begin{align}
        \Phi^0(r')&=\partial_r\Phi^0(r')=0 \;
        \text{ at } \varpi',
    \end{align}
    then this is equivalent to the statement that $\xi=0$ and $H_G\xi=0$. We note that we are then immediately in the setting of lemma \ref{lem:NoLocalMinimum}, and thus we have that $\partial_r^2\Phi^0(r')<0$. We claim that $(r-r')\partial_r\Phi^0(r)$ for all $(r,\xi)\in \pi\circ \gamma$. Near $r'$ the statement is clear as $\varpi'$ represents a local $r$-maximum for $\Phi^0(r)$. To see that the statement holds globally as well, note that if it would fail in some $r''$ then this must represent a critical point of $\Phi^0(r)$ which must be either a saddle point or a minimum, for which $\partial_r^2\Phi^0(r')\leq 0$. This however leads to another contradicition with lemma \ref{lem:NoLocalMinimum}. We conclude that $\Phi^0(r)<0$ for all $r\neq r'$ on $\gamma$.
\end{proof}

Finally we will prove the following lemma that will allow us to break down the stable and unstable manifolds of the flow of $H_{\mathscr{G}^0}$ depending on the sign of $\sigma$. 
\begin{lemma}\label{lem:SignofSigma} Let $\mfrm>0$, $(\bfa,\bfb)\in D_{\mfrm}$.
    Furthermore $\varpi^0=(r^0,z^0;\xi^0,\zeta^0)\in \chr(G)$ with $\zeta^0=(\sigma^0,\eta_\omega^0)$ and once again let $\gamma$ be the integral curve of $H_G$ through the point $\varpi^0$. Then, for any other point $\varpi =(r,z,\xi,\zeta)$ with $\zeta=(\sigma,\eta_\omega)$ we have that 
    \begin{equation}
        \text{sign}(\sigma^0)=\text{sign}(\sigma).
    \end{equation}
\end{lemma}

\begin{proof}
  We will show that at $\varpi$, $G(dt,\varpi)$ has the same sign as $A$, which in turn has the same sign as $\sigma$. We have that 
  \begin{align}
      G(dt,\varpi)&=\frac{1}{2}\frac{\partial}{\partial\sigma}G
    =\Delta_r^{-1}A\partial_\sigma A-\Delta_r^{-1}C\partial_\sigma C- B\partial_\sigma B -\frac{r^2}{\bfp^2\Delta_r}C\partial_\sigma C.
\end{align}
Plugging back in the identities (\ref{eq:MatrixEntry1})-(\ref{eq:MatrixEntry3}), we then find that 
\begin{align}
    G(dt,\varpi)&=\underbrace{\frac{(r^2+\bfa^2)(r^2+\bfb^2)}{r^2\Delta_r}}_{\frac{\Delta_r+2\mfrm}{\Delta_r}}A-\frac{\bfa\bfb}{r\sqrt{\Delta_r}}C-\frac{(\bfa^2-\bfb^2)\cos^2\theta\sin^2\theta}{\bfp}B-\frac{\bfa\bfb r}{\bfp^2\sqrt{\Delta_r}}C.
\end{align}
Now let us assume without loss of generality that $A>0$ and $|\bfa|\leq|\bfb|$, we then see that
\begin{align}
    \frac{G(dt,\varpi)}{A}&\geq \frac{\Delta_r+2\mfrm}{\Delta_r}-\frac{|\bfa\bfb|}{r\sqrt{\Delta_r}}\left|\frac{C}{A}\right|-\frac{(\bfb^2-\bfa^2)\cos^2\theta\sin^2\theta}{\bfp}\left|\frac{B}{A}\right|-\frac{|\bfa\bfb|r}{\bfp^2\sqrt{\Delta_r}}\left|\frac{C}{A}\right|.
\end{align}
Using again that $\varpi \in \gamma\subset \chr(G)$ we plug in the inequalities (\ref{eq:BAInequality}) and (\ref{eq:CAInequality2}), i.e. $\left|\frac{B}{A}\right|<\Delta_r^{-1/2}$ and $\left|\frac{C}{A}\right|\leq \frac{\bfp}{r}$ and rescaling by $\bfp\Delta_r$ we find that 
\begin{align}
   \bfp\Delta_r \frac{ G(dt,\varpi)}{A}&\geq \bfp\Delta_r-\left[\bfp^2+(|\bfb|^2-|\bfa|^2)\sin^2\theta\cos^2\theta\ +|\bfa\bfb|\right]\Delta_r^{1/2} +2\mfrm \bfp\\
   &\geq  \bfp\Delta_r-\left[\bfp^2+(|\bfb|^2-|\bfa|^2)\sin^2\theta\cos^2\theta\ +|\bfa\bfb|\right]\Delta_r^{1/2} +\bfp(|\bfa|+|\bfb|)^2
\end{align}

Now the statement follows as we claim that for all $(\bfa,\bfb)\in D_\mfrm$
\begin{align}
   M^2:=\left[\bfp^2+(|\bfb|^2-|\bfa|^2)\sin^2\theta\cos^2\theta\ +|\bfa\bfb|\right]^2<4\bfp^2(|\bfa|+|\bfb|)^2.
\end{align}

Now observe that $M$ can be decomposed as 
\begin{align}
    M:= (|\bfa|+|\bfb|)\underbrace{(|\bfa|\cos^4\theta+|\bfb|(1-\cos^4\theta))}_{:=K}
\end{align}

Now $K$ in turn can be written as 
\begin{align}
    K&=(|\bfa|\cos\theta)(\cos\theta)^{3}+(|\bfb|\sqrt{1-\cos^2\theta}(1+\cos^2\theta)\sqrt{1-\cos^2\theta}
\end{align}

Thus by Cauchy-Schwarz, we have that 
\begin{align}
    K^2&\leq (\bfa^2\cos^2\theta+\bfb^2(1-\cos^2\theta))(\cos^6\theta+(1-\cos^2\theta)
    (1+\cos^2\theta)^2).\\
    &=\bfp^2(\cos^6\theta+(1-\cos^2\theta)
    (1+\cos^2\theta)^2).
\end{align}
However
\begin{align}
    (\cos^6\theta+(1-\cos^2\theta)(1+\cos^2\theta)^2)&=1+\cos^2\theta-\cos^4\theta \\
    &=\frac{5}{4}-\left(\cos^2\theta -\frac{1}{2}\right)^2 \leq \frac{5}{4}<2.
\end{align}
Thus we see that 

\begin{align}
    M^2 = (|\bfa|+|\bfb|)^2K^2 \leq \frac{5}{4}(|\bfa|+|\bfb|)^2\bfp^2<2(|\bfa|+|\bfb|)^2\bfp^2,
\end{align}
as required.

\end{proof}

\begin{remark}
    Lemma \ref{lem:SignofSigma} simply proves that the future and past oriented lightcones $\Sigma_G\cap \{\pm \sigma>0\}$ stay conserved under the Hamiltonian flow at the trapped set $\Gamma$, which is also a well known property of the trapped set of Kerr. One could also use a ''soft'' argument to prove this fact, as this property holds clearly on Schwarzschild-Tangerlini in $1+4$ dimensions. One then perturbs from $(\bfa,\bfb)=(0,0)$ to the $(\bfa,\bfb)\in D_{\mfrm}$ under consideration here and shows that the lightcones stay disjoint under null-geodesics along the path connecting these two parameters.
\end{remark}

With lemmas \ref{lem:NondegenerateMax} and \ref{lem:SignofSigma} in hand, we can write write the characteristic set $(\mathscr{G}^0)^{-1}(0)$ as the following union for a fixed $\varpi \in \chr(G)^{\sgn\,\sigma^0}$:

\begin{align}\label{eq:DefinitionofGamma}
    (\mathscr{G}^0)^{-1}(0)&=\Gamma^+_{(z_0,\zeta_0)}\cup \Gamma^-_{,(z_0,\zeta_0)}, \\
    \text{ where } \Gamma_{\pm,(z_0,\zeta_0)}&= \left\{(r,\xi)\in T^*(r_+,\infty)_r\;|\; \xi=\mp\sgn(\sigma^0)\sgn(r-r')\sqrt{-\Delta_r^{-1}\Phi^0(r)} \right\} \subset T^*(r_+,\infty)_r\nonumber,
\end{align}
and we take the positive square root, which is well-defined and real due $\Delta_r>0$ for $r>r_+$ and $\Phi^0(r)\leq0$ due to lemma \ref{lem:NondegenerateMax}. We have that 
\begin{align}
    \Gamma^+_{(z_0,\zeta_0)}\cap \Gamma^-_{(z_0,\zeta_0)}=\{\varpi\},
\end{align}
from which follows as well that at $\varpi$ we should have that $r=r'$ and $\xi=0$.
Finally we note that $\Phi^0(r)$ has a non-degenerate maximum, again by lemma \ref{lem:NondegenerateMax} hence we can conclude that the submanifolds $\Gamma_{\pm,(z_0,\zeta_0)}$ are smooth. \\
We have that $H_{\mathscr{G}^0}|_{\varpi}=0$, so that $H_G$ integral curves from $\varpi$ lie along constant $r=r'$. Hence they do not escape to $r=r_+$  or $r\to \infty$ and are \textit{trapped}. Let us in general define $\chr(G)^\pm=\{(z,\zeta)\in \chr(G)\subset T^*M\;|\; \sigma\gtrless0\}$, then on $\Gamma^+_{(z^0,\zeta^0)}\cap \chr(G)^\pm$ we have that 
\begin{align}
    \pm(r-r')H_Gr=\mp 2(r-r')\Delta_r\xi>0,
\end{align}
whereas on $\Gamma^-_{(z^0,\zeta^0)}\cap \chr(G)^\pm$ we have that  
\begin{align}
    \pm(r-r')H_Gr=\pm 2(r-r')\Delta_r\xi<0.
\end{align}
From this we conclude that future directed integral curves within $\Gamma_{+,(z_0,\zeta_0)}\backslash \{(r',0)\}$ escape to $r\to r_+$ or $r\to \infty$ whilst past directed integral curves tend to $r'$. Vice versa we also see that future directed integral curves on $\Gamma_{-,(z_0,\zeta_0)}\backslash \{(r',0\}$ tend to $r_+$ and past directed integral curves tend to $r\to r_+$ or $r\to \infty$. From this we deduce that $\Gamma_{+,(z_0,\zeta_0)}$ and $\Gamma_{-,(z_0,\zeta_0)}$ define the forward resp. backward trapped sets of the flow of $H_{\mathscr{G}^0}$ on the projection $\chr(\mathscr{G}^0)\subset T^*(r_+,\infty)_r$.\\
Under pullback the above considerations lead us to the following definition of the trapped set. 
\begin{definition}\label{def:rxiTrappedSet}
    The trapped set $\Gamma=\Gamma_+\cup \Gamma_-\subset T^*M\backslash O$, with $\Gamma^{\pm}=\Gamma\cap \chr(G)^\pm$ is defined as 
    \begin{align}\label{eq:TrappedSet}
        \Gamma:=\left\{\varpi\in \chr(G)\cap T^*M\;|\; r(\varpi)>r_+, \; \xi(\varpi)=(H_G\xi)(\varpi)=0\right\}
    \end{align}
\end{definition}
\begin{remark}
    We note equivalent conditions for (\ref{eq:TrappedSet}) as well. Given that we work on the region $r>r_+$ it is easy to see that (\ref{eq:TrappedSet}) is equivalent to 

        \begin{align}
        \Gamma:=\left\{\varpi\in \chr(G)\cap T^*M\;|\; H_G r=0, H_G^2 r=0\right\}.
    \end{align}
    Similarly one can easily show that (\ref{eq:TrappedSet}) is equivalent to 
     
    \begin{align}
        \Gamma:=\left\{\varpi\in \chr(G)\cap T^*M\;|\; \xi(\varpi)=0, (\partial_rG)(\varpi)=0\right\}.
    \end{align}
\end{remark}

With regards to the discussion and the lemmas above we then have the following.

\begin{lemma}\label{lem:TrappedSetCodimension}
    We have $\pm\sigma >0$ on $\Gamma$ and $\Gamma$ is a smooth codimension $2$ submanifold $\chr(G)$ hence a smooth codimension $3$ submanifold of $T^*M$.
\end{lemma}

\begin{proof}
    The sign of $\sigma$ was discussed in lemma \ref{lem:SignofSigma}. Notice now that $\xi=0$ $\partial_rG=0$ and $G=0$ (the definition of $\Gamma$) is equivalent to $(\xi,\Psi,G)=(0,0,0)$. We will show that the differentials $d\xi$, $d\Psi$ and $dG$ are linearly independent on $\{(\xi,\Psi,G)=(0,0,0)\}$. Working away from the zero section as in lemma \ref{lem:rConvexityonCharacteristicSet}, we know that $A\neq 0$. Notice that we then have
   \begin{equation}
    \begin{aligned}
    &\partial_\xi \xi = 1, \qquad
    &\partial_\xi \Psi = 0, \qquad
    &\partial_\sigma G = 0, \\
    &\partial_r \xi = 0, \qquad
    &\partial_r \Psi \neq 0, \qquad
    &\partial_r G = \Delta^{-2}\Psi = 0, \\
    &\partial_\sigma \xi = 0, \qquad
    & \qquad
    &\partial_\sigma G = 2G(dt,-) \neq 0.
\end{aligned}
\end{equation}
    Contrasting this against the fact that by the proof of lemma \ref{lem:NoLocalMinimum}, we have that $\partial_r \Psi\neq 0$, linear independence of $d\xi, d\Psi$ and $dG$ follows regardless of the value of $\partial_\sigma\Psi$.
\end{proof}

Let us consider again the forward/backward trapped sets in $T^*(r_+,\infty)_r$ from definition \ref{def:rxiTrappedSet}. Let us now define the projection of  $
\Pi:T^*M\to T^*(\mathbb{R}_t\times\mathbb{S}^3_{\phi,\psi,\theta})$ via
\begin{align}
    \Pi(r,z;\xi,\zeta) = (z,\zeta).
\end{align}

 of the trapped set, $\Pi(\Gamma)$ is explicitly given by 
\begin{align}
    \tilde\Gamma&=\Pi(\Gamma)=\{(z,\zeta)\in T^*(\mathbb{R}_t\times \mathbb{S}^3_{\phi,\psi,\theta})\;|\; \exists(r;\xi)\in T^*(r_+,\infty) \;\text{s.t.} (r,z;\xi,\zeta)\in \Gamma\}.
\end{align}

and we have explicitly that for the backward $(-)$ and forward $(+)$ trapped sets $\Gamma_\pm$, defined in lemma \ref{lem:TrappedSetCodimension} that 
\begin{align}
    \Gamma^\pm&=\bigcup_{(z;\zeta)\in\tilde{\Gamma}}\{(r,z;\xi,\zeta) \in\chr(G)\cap T^*M\;|\; (r;\xi)\in \Gamma_{\pm,(z,\zeta)}\}\nonumber\\
    &=\left\{(r,z;\xi,\zeta) \in\chr(G)\cap T^*M\;|\;\xi=\mp\sgn(\sigma(\zeta))\sgn(r-r'_{(z,\zeta)})\sqrt{-\Phi^0_{(z,\zeta)}(r)/\Delta_r}\right\}\label{eq:TrappedSetDeffromzAndzeta},
\end{align}
with $\Phi^0_{(z,\zeta)}(r)$ and its unique non degenerate maximum $r'_{(z;\zeta)}$ given by equation (\ref{eq:FixedrxiDynamics}) and lemma \ref{lem:NondegenerateMax}, respectively.
Note that the above is well defined as by construction the sets $\Gamma_{\pm,(z;\zeta)}\subset T^*(r_+,\infty)_r$ do not depend on the choice of $(z,\zeta)\in \Pi\circ \gamma$ for null-bicharacteristics $\gamma$ of the flow of $H_G$, showing that $\Gamma_\pm$ are invariant under the $H_G$-flow. We then have the following lemma.

\begin{lemma}\label{lem:TransversalIntersection}
    The set $\tilde\Gamma$ is a smooth conic codimension $1$ submanifold of $T^*(\mathbb{R}_t\times \mathbb{S}^3)$. The forward$(-)$/backward$(+)$ trapped sets $\Gamma_\pm$ intersect transversally at $\Gamma$ as submanifolds of $\chr(G)$.
\end{lemma}
\begin{proof}
    We have to show here that $\Gamma_\pm$ is smooth, conic and of codimension $1$ within $T^*M$ and that $\Gamma_+$, $\Gamma_-$ intersect transversally at $\Gamma$ as submanifolds of $\chr(G)$. Let us consider $(z;\zeta)\in \tilde{\Gamma}$ and let $\varpi=(r',z;0,\zeta)\in \Gamma$ as the preimage $\Pi^{-1}\{(z;\zeta)\}$. We will show that $T_{\varpi}\Gamma\cap\ker D_\varpi \Pi = \{0\}$. We know by lemma \ref{lem:TrappedSetCodimension} that 
    \begin{align}\label{eq:KernelIntersections}
        T_\varpi\Gamma = \ker d\xi\cap\ker d\Psi\cap \ker G.
    \end{align}
    Thus, for any $v=c_1\partial_r+c_2\partial_\xi\in \ker D_\varpi\Pi$ we have by (\ref{eq:KernelIntersections}) that 
    \begin{equation}
        0=d\xi(v)=c_2,
    \end{equation}
    which forces $c_2=0$, and
    \begin{equation}
        0=d\Psi(v)=c_1\partial_r\Psi.
    \end{equation}
    which forces $c_1=0$ by virtue of $\partial_r\Psi>0$ in lemma \ref{lem:NoLocalMinimum}. This shows that for a small enough open neighbourhood $\mathcal{U}$ of $\varpi$ inside $\Gamma$, $\Pi(\mathcal{U})$ will be a small open neighbourhood of $(z;\zeta)\in T^*(\mathbb{R}\times \mathbb{S}^3)$ of codimension $3-2=1$.\\
    The uniqueness of $r'_{(z;\zeta)}$ and $\varpi$ given $(z;\zeta)$ implies there exists a neighbourhood $\mathcal{W}$ of $\Pi(\varpi)$ inside $\tilde{\Gamma}$ such that $\Pi^{-1}(\mathcal{W})=\mathcal{U}'$ where $\mathcal{U}'\subset \mathcal{U}$ is a small neighbourhood of $\varpi$. Finally transversality follows from $\Phi^0_{(z;\zeta)}(r)$ having a non-degenerate maximum at $r=r'_{(z;\zeta)}$ and the explicit definition of $\Gamma_\pm$ from (\ref{eq:TrappedSetDeffromzAndzeta}).
\end{proof}

This also shows that the map $(z;\xi)\mapsto r'_{(z,\zeta)}$ is smooth. Let us introduce defining functions for $\Gamma_\pm$ as 
\begin{align}\label{eq:DefiningFunctionExtension}
    \phi_\pm(r,z;\xi,\zeta)&=\xi\pm\sgn(\sigma(\zeta))\sgn(r-r'_{(z,\zeta)})\sqrt{-\Phi^0_{(z,\zeta)}(r)/\Delta_r},
\end{align}
for $(z,\zeta)\in \tilde{\Gamma}$. Then $\phi_\pm$ are stationary and homogenous of degree $1$ in the $(\xi,\zeta)$-variables inside $\tilde \Gamma\times T^*(r_+,\infty)_r$. Extending $\phi_\pm$ in an arbitrary, smooth, stationary and homogenous way to an open neighbourhood inside $T^*M$ we have that $H_G\phi_\pm=0$ at $\Gamma_\pm$. As the codimension of $\Gamma_\pm\subset\chr(G)$ is $1$ and $H_G$ is tangent to $\Sigma$, we can write 
\begin{align}
    \sigma^{-1}H_G \phi_\pm|_{\chr(G)} = \mp\nu_\pm \phi_\pm|_{\chr(G)},
\end{align}
where $\nu_\pm$ are smooth on $\chr(G)$ and homogeneous of degree $0$ by definition. We note that $\sigma^{-1}H_G =\pm|\sigma|^{-1} H_G$. We conclude with the following lemma.

\begin{lemma}
    The restrictions of $\nu_\pm|_{\Gamma}\in C^\infty(\Gamma)$ are independent of the choice of extensions of $\phi_\pm$. They are given by 
    \begin{align}\label{eq:DefiningCoefficients}
        \nu_\pm(r'_{(z,\zeta)},z;0,\zeta)=\nu(z,\zeta):=|\sigma|^{-1}\sqrt{-2\Delta(r'_{(z,\zeta)})\partial_r^2\Phi^0_{(z,\zeta)}(r'_{(z,\zeta)})}.
    \end{align}
    In particular, $\inf_\Gamma \nu_\pm>0$. Furthermore, $\{\phi_+,\phi_-\}\in C^\infty(\Gamma_0)$ is everywhere non-zero.
\end{lemma}

\begin{proof}
    We work on $\chr(G)$. Let $\tilde{\phi}_\pm$ is another extension of $\tilde{\phi}$ of (\ref{eq:DefiningFunctionExtension}), then for $a_\pm:=\phi_\pm-\tilde{\phi}_\pm$ we have $a_\pm=0$ and $H_Ga_\pm =0$ on $\Gamma_+\cap\Gamma_-$. By lemma \ref{lem:TransversalIntersection}, $\Gamma_+$ and $\Gamma_-$ intersect transversally, so $a_\pm=\phi_+\phi_-a_\pm'$ and $H_G a_\pm=\phi_+\phi_- a''_\pm$ with $a_\pm',a_\pm''$ smooth. Therefore, we have that 
    \begin{align}
        \frac{H_G\tilde{\phi}_\pm}{\tilde{\phi}_\pm}&=\frac{H_G\phi_\pm+\phi_+\phi_-a''_\pm}{\phi_\pm+\phi_+\phi_-a'_\pm}=\frac{H_G\phi_+}{\phi_\pm(1+\mathcal{O}(\phi_\mp)}+\mathcal{O}(\phi_\mp)\to \frac{H_G\phi_\pm}{\phi_\pm}\;\text{ at }\; \Gamma.
    \end{align}
    Therefore it suffices for the evaluation of $\nu_\pm$ to compute the limit of $(H_G\phi_\pm)/\phi_\pm$ at $\Gamma$ along $\Gamma_\mp$. At $\Gamma_{\mp}$ we have
    \begin{align}
        H_G\phi_\pm=H_G(\phi_\pm+\phi_\mp)=H_G(2\xi)=-\frac{2\Psi}{\Delta_r^2}+2\Delta_r'\xi^2.
    \end{align}
    We know that $\xi$ vanishes quadratically at $\Gamma$ hence $\frac{\xi^2}{\phi_\pm}\to 0$ at $\Gamma$. On the other hand we have that 
    \begin{align}
        \mp\phi_\pm=\mp(\phi_\pm+\phi_\mp)=\mp2\xi =-2\sgn(\sigma)\sgn(r-r'_{(z,\zeta))})\sqrt{\frac{-\Phi^0_{(z,\zeta)}(r)}{\Delta_r}}.
    \end{align}
    Since $\partial_r^2\Phi^0_{(z,\zeta)}(r'_{(z,\zeta)})<0$ by lemma \ref{lem:NoLocalMinimum}, we see that by expanding around $r=r'_{(z,\zeta)}$
    \begin{align}\label{eq:TaylorExpansionsatRadialMax}
        \sqrt{-\frac{\Phi^0_{(z,\zeta)}(r)}{\Delta_r}}&=\left.\sqrt{-\frac{\partial_r^2\Phi^0_{(z,\zeta)}(r)}{2\Delta_r}}\right|_{r=r'_{(z,\zeta)}}(r-r'_{(z,\zeta)})\sgn(r-r'_{(z,\zeta)})+\mathcal{O}((r-r'_{(z,\zeta)})^2)\\
        \Delta_r^{-2}\Psi&=\underbrace{\Delta_r^{-2}\Psi(r'_{(z,\zeta)})}_{=0}+\Delta_r^{-2}\partial_r\Psi(r'_{(z,\zeta)})(r-r'_{(z,\zeta)})+\mathcal{O}((r-r'_{(z,\zeta)})^2),\nonumber
    \end{align}
    where we used again that $\Psi=0$ on $\Gamma$. Then, using 
    \begin{align}
        \Delta_r^{-2}\partial_r\Psi=\partial_r(\Delta_r^{-2}\Psi)=-\partial_r^2V=\partial_r^2\Phi^0_{(z,\zeta)}
    \end{align}
    at $r'_{(z,\zeta)}$, we verify the formulas for $\nu_\pm$.\\
    Finally we see that for the Poisson bracket $\{\phi_+,\phi_-\}$ we have at $r=r'_{(z,\zeta)}$ that 
    \begin{align}
        \{\phi_+,\phi_-\}=H_{\phi_+}\phi_-=H_{\phi_++\phi_-}
        \phi_-=2H_\xi \phi_- = 2\partial_r\phi_- = 2\sgn(\sigma)\left.\sqrt{-\frac{\partial_r^2\Phi^0_{(r,\zeta)}}{2\Delta_r}}\right|_{r=r'_{(z,\zeta)}},
    \end{align}
    which is non-zero.
    \end{proof}

\begin{figure}[H]
    \centering
    \begin{tikzpicture}[
  x=1.15cm,y=1.15cm,
  every node/.style={font=\large},
  axis/.style={black,line width=.55pt,-{Stealth[length=2.1mm]}},
  trapped/.style={black,line width=.65pt},
  nearby/.style={black!32,line width=.65pt},
  flow/.style 2 args={
    postaction={decorate},
    decoration={markings,
      mark=at position #1 with {
        \arrow[#2]{Stealth[length=2mm,width=1.4mm]}
      }
    }
  }
]
  \draw[axis] (-.5,0)--(7,0) node[right] {$r$};
  \draw[axis] (0,-2.5)--(0,2.5) node[above left=-1pt] {$\xi$};

  \draw[nearby,flow={.38}{black!32}]
    (.5,2)
    .. controls (1.5,2) and (2.5,1.5) .. (3.5,1.5)
    .. controls (4.5,1.5) and (5.5,2) .. (6.5,2);
  \draw[nearby,flow={.38}{black!32}]
    (6.5,-2)
    .. controls (5.5,-2) and (4.5,-1.5) .. (3.5,-1.5)
    .. controls (2.5,-1.5) and (1.5,-2) .. (.5,-2);

  \draw[nearby,flow={.38}{black!32}]
    (.5,1)
    .. controls (1.833333,1) and (1.833333,-1) .. (.5,-1);
  \draw[nearby,flow={.38}{black!32}]
    (6.5,-1)
    .. controls (5.166667,-1) and (5.166667,1) .. (6.5,1);

  \draw[trapped,flow={.53}{black}]
    (.5,1.5)
    .. controls (1.5,1.5) and (2.4,.9) .. (3.5,0);
  \draw[trapped,flow={.53}{black}]
    (6.5,-1.5)
    .. controls (5.5,-1.5) and (4.6,-.9) .. (3.5,0);

  \draw[trapped,flow={.53}{black}]
    (3.5,0)
    .. controls (4.6,.9) and (5.5,1.5) .. (6.5,1.5);
  \draw[trapped,flow={.53}{black}]
    (3.5,0)
    .. controls (2.4,-.9) and (1.5,-1.5) .. (.5,-1.5);

  \fill (3.5,0) circle[radius=2.2pt];
  \node[below=7pt] at (3.5,0) {$r'$};
  \node[anchor=west,inner sep=2pt] at (6.58,1.5)
    {$\Gamma^{-}_{(z^0,\zeta^0)}$};
  \node[anchor=west,inner sep=2pt] at (6.58,-1.5)
    {$\Gamma^{+}_{(z^0,\zeta^0)}$};
\end{tikzpicture}
    \caption{Picture adapted from \cite{hintz2026non} figure 4.2 of the normally hyperbolic trapped null geodesic flow projected to the $(r,\xi)$ plane of $T^*M$. $\Gamma^+_{(z_0,\zeta_0)}$ represents the projected stable manifold of the flow and $\Gamma^-_{(z_0,\zeta_0)}$ the projected unstable manifold of the flow. The trapping persists throughout the full subextremal range $(|\bfa|+|\bfb|)^2<2\mfrm$.}
    \label{fig:NormallyHyperbolictrapped}
\end{figure}
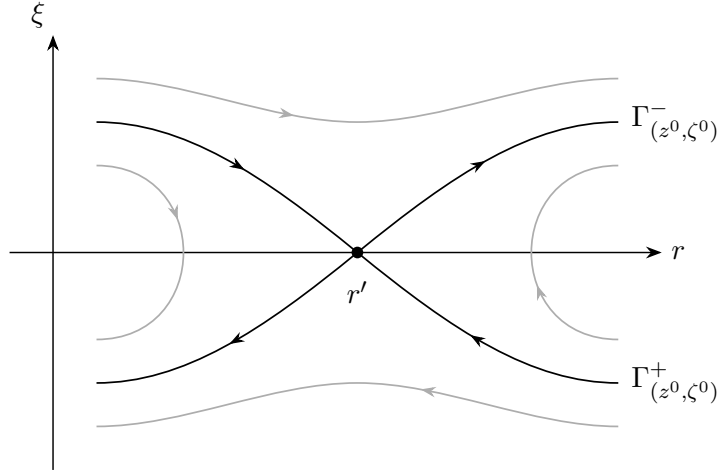

\section{Part II: Microlocal Analysis.}

To analyze the forced wave problem (\ref{eq:ForwardWaveSolution}) on the Fourier side, we essentially make use of stationary of the problem and analyze the spectral family 
\begin{align}\label{eq:SpectralFamily}
    \hat{\Box}_{g_{\bfa,\bfb,\mfrm}}(\sigma):= e^{i\sigma t_*}\circ \Box_{g_{\bfa,\bfb,\mfrm}}\circ e^{-i\sigma t_*}, \text{ for } \sigma \in \mathbb{C},
\end{align}
on function spaces over $X$ as defined in Definition \ref{def:RadiallyCompactifiedR4}, intuitively by replacing $\partial_{t_*}$ by $-i\sigma$. To this end we develop the geometric $b$- and $sc$- structures first introduced by Melrose \cite{melrose1993atiyah}\cite{melrose1994spectral}, necessary to obtain uniform Fredholm estimates for $\sigma \in \overline{\mathbb{H}}$.

\subsection{b- and scattering structures.}

Let $X$ be a compact $d$-dimensional manifold with boundary\footnote{We can and will in some cases below allow $X$ to have corners} $\partial X \neq \emptyset$ and denote with $\rho\in C^\infty(X)$ its boundary defining fucntion, that is $\partial X =\rho^{-1}(0)$ and $d\rho|_{\partial X}\neq 0$. We then define the Lie algebras of $b$-vector fields and scattering vector fields through
\begin{align}
    \mathcal{V}_b(X)=\{V\in \mathcal{V}(X)\;|\; V \text{ is tangent to } X\},\;\;\;\;\;\; \mathcal{V}_{sc}(X)=\rho \mathcal{V}_b(X),
\end{align}
where as usual with being tangent to $\partial X$ we mean that for all $p\in \partial X$ we have $V_p\in T_p\partial X$, i.e. $V\rho = \rho C^\infty$, implying $d\rho(V_p)=V_p(\rho)=0$. The $b$ vector fields are by definition the vector fields that are tangent to the boundary, and the scattering field are the vector fields that commute up to the boundary $\partial X$ as we have 
\begin{align}
    [\mathcal{V}_{sc}(X),\mathcal{V}_{sc}(X)]=[\rho\mathcal{V}_b(X),\rho\mathcal{V}_b(X)]\subset \rho^2\mathcal{V}_b(X) =\rho\mathcal{V}_{sc}(X).
\end{align}
In local adapted coordinates near the boundary $\partial X$, given by $\mathbb{R}_x^+\times \mathbb{R}^{d-i}_{y_j}$ for $j=1,2,...,d-1$, with $x$ direction of the boundary, we have that elements of $\mathcal{V}_b(X)$ are given as $a(x,y)x\partial_x+\sum^{d-1}_{i=1}b_i(x,y)\partial_{y_i}$ and elements of $\mathcal{V}_{sc}(X)$ are given as $a(x,y)x^2\partial_x+x\sum_{i=1}^{d-1}b_i(x,y)\partial_{y_i}$ where $a(x,y),b_i(x,y)\in C^\infty(X)$. Thus there natural vector bundles of $X$ (including its boundary) given by 
\begin{align} 
    {}^bTX\to X \;\;\text{ and }\;\;{}^{sc}TX\to X,
\end{align}
with local frames given by $\text{Span}_{C^\infty(X)}\{x\partial_x,\partial_{y_i}\}$ resp. $\text{Span}_{C^\infty(X)}\{x^2\partial_x,x\partial_{y_i}\}$ such that $\mathcal{V}_b(X)=\Gamma(X;{}^{b}TX)$ and $\mathcal{V}_{sc}(X) = \Gamma(X;{}^{sc}TX)$.\\ 
The manifold $X$ that is of primary interest to us in this case would of course be the spatial manifold $X=\overline{\mathbb{R}^4}$ defined in Definition \ref{def:RadiallyCompactifiedR4}, which is why we abuse notation on purpose here. In this specific case we often simply use $\rho=r^{-1}$ as our boundary defining coordinate $x$ and boundary defining function simultaneously (i.e. setting $a(\rho,y)=1$, working on a collar neighbourhood $[0,1)_\rho\times \mathbb{S}^3$ of $\partial X\cong \mathbb{S}^3$. In this particular case we can choose $3$ sections $\Omega_j\in \mathcal{V}(\mathbb{S}^3)$ to span $T\mathbb{S}^3$ (using paralellizability of $\mathbb{S}^3$), to write 
$b-$ and $sc-$ vector fields as 
\begin{align}
   \underbrace{-a(\rho,y)\rho\partial_\rho}_{a(\rho,y)r\partial_r} + \sum_{j=1}^3 b_j(\rho,y) \rho \Omega_j\; \text{ and }\; \underbrace{-a(\rho,y)\rho^2\partial_\rho}_{a(\rho,y)\partial_r} + \sum_{j=1}^3 b_j(\rho,y) \rho \Omega_j.
\end{align}
Moreover, choosing Cartesian coordinates $(x_1,...,x_4)$ for $X^\circ = \mathbb{R}^4$ we see that that the $\Omega_j$ away from the boundary can be given explicitly as 
\begin{align}
    \Omega_j= r(r\partial_{x_j} -\omega_j\partial_r),
\end{align}
in spherical coordinates with $x=r\omega$ and $\omega\in \mathbb{S}^3$,
from which follows that 
\begin{align}
    \partial_{x_j}=\omega_j\partial_r+r^{-1}\Omega_j= -\omega_j\rho^2\partial_\rho +\rho\Omega_j.
\end{align}
showing that on $X^\circ$ there is a natural identification of ${}^{sc}TX^\circ$ with $TX^\circ$ and by extension also with ${}^bTX^\circ$, as $\rho\neq 0$  on $X^\circ$ and in particular the vector fields $\partial_{x_i}$ extend to scattering vector fields on $\overline{\mathbb{R}^4}$ Indeed, in general it follows that from the natural inclusions 
$V_b(X)\hookrightarrow\mathcal{V}(X)$ and $\mathcal{V}_{sc}(X)$ that there are natural bundle maps 
\begin{align}
    \iota_b:{}^bTX\to TX\; \text{ and } \iota_{sc}:{}^{sc}TX\to TX,
\end{align}
that are injective when restricted to $TX^\circ$, satisfying 
\begin{align}
    \iota_b(\rho\partial_\rho)= &-r\cdot \partial_r, \;\;\;\; \iota_b(\partial_{y_i})=\partial_{y_i}\\
    & \iota_{sc} \circ m_\rho=\rho\cdot \iota_b
\end{align}
where $m_\rho:{}^bTX\to {}^{sc}TX$ is multiplication by $\rho$ identifying ${}^{sc}TX=\rho\cdot{}^bTX$ (i.e. we have $\iota_{sc}(\rho^2\partial_\rho)=-\partial_r$). Away from the boundary, we have that the maps $\iota_{\bullet}:{}^\bullet TX^\circ \to TX^\circ$ are fibrewise isomorphisms, however at $\partial X$ this picture changes. We have that $\ker (\iota_{sc}|_{\partial X})={}^{sc}T_{\partial X}X$, but, more interestingly we have that the kernel of $\iota_b$ is $1$-dimnensional, defining the $b$-normal bundle at $\partial X$ 
\begin{align}
{}^bN\partial X=\ker(\iota_{b}|_{\partial X})\subset {}^bTX, 
\end{align}
which is locally spanned over $C^\infty(X)$ by $\rho\partial_\rho$. Note that $(\rho\partial_\rho)|_{\partial X}\in C^\infty(\partial X;{}^bTX)$ independent of the chart/collar neighbourhood chosen at $\partial X$ and hence gives a canonical trivialisation of ${}^b N\partial X$ and a natural splitting of ${}^bT_{\partial X}X$ through the sequence 
\begin{align} \label{eq:TangentSpaceExactSeq}
    0\to{}^bN\partial X\to {}^bT_{\partial X}X \xrightarrow{\iota_b} T_{\partial X}X\to N\partial X\to 0,
\end{align} 
where we set $N\partial X:=T_{\partial X}X/T\partial X$.
From the above, it makes sence to introduce the dualizations ${}^bT^*X$ and ${}^{sc}T^*X$, which in local adapted coordinates is given by 
\begin{align}
    {}^bT^*X=\text{span}_{C^\infty(X)}\left\{\frac{dx}{x},dy_1,...,dy_d\right\} \text{ and } {}^{sc}T^*X=\text{span}_{C^\infty(X)}\left\{\frac{dx}{x^2},\frac{dy_1}{x},...,\frac{dy_{d-1}}{x}\right\}, \label{eq:LocalSpanscoTangentspaces}
\end{align}
where we induce bundle maps $\iota_b^*:{}^bT^*X\to T^*X$, $\iota_{sc}:{}^{sc}T^*X\to T^*X$ from the dualization of ${}^bTX$ and ${}^{sc}TX$ as well. We have again that $\ker(\iota_{sc}^*|_{\partial X})=0$, however for $p\in \partial X$ we have
\begin{align}
    \ker(\iota_b^*)=\{\alpha\in T^*X\;|\;\alpha(T_p\partial X)=0\text{ for all } p \in \partial X\}=N^*\partial X,
\end{align}
i.e. its kernel is exactly the ordinary normal bundle of $\partial X$, $N^*\partial X \subset T^*X$, which can also be found from dualizing (\ref{eq:TangentSpaceExactSeq})
\begin{align}
    0\to N^*\partial X\to T^*_{\partial X}X \xrightarrow{\iota_b^*} {}^bT^*_{\partial X}X \to({}^bN\partial X)^*\to 0.
\end{align}
When one defines the $b$-conormal bundel of a submanifold $Z\subset X$, as 
\begin{align}
    {}^bN^*Z=\{\alpha \in {}^bT^*_ZX\;|\;\alpha(V)=0, \forall V\in {}^bT_ZM, V|_Z\in \mathcal{V}_b(Z)\}.
\end{align}
and by a boundary defining function $\rho \in C^\infty(X)$ one has that
\begin{align}
    {}^bN^*\partial X = \text{Span}\left(\frac{d\rho}{\rho}\right).
\end{align}
Indeed in that case we have the isomorphism
\begin{align}
    {}^bN^*\partial X \cong ({}^bN\partial X)^*,
\end{align}
albeit not canonical as it depends on the choice of $\rho$.\\
The existence ${}^\bullet T^*X$ also gives rise to so-called $b$- and $sc$-metrics on $X$.In general a b- or scattering metric $g_{\bullet}\in C^\infty(X;S^2{}^{\bullet}T^*X)$ are given as a non-degenerate quadratic form on each fibre ${}^\bullet T^*_pX$. The standard $b$- and $sc$- metrics are usually just given by 
\begin{align}
    g_{sc}=\frac{dx^2}{x^4}\wedge \frac{dy_1^2}{x^2}\wedge...\wedge \frac{dy_{d-1}}{x^2} \;\text{ and } g_b=\frac{dx^2}{x^2}\wedge dy_1\wedge...\wedge dy_{d-1},
\end{align}
and they come with associated scattering densities, i.e. positive definite sections $\omega_\bullet\in {}^{\bullet}\mathcal{D}(X):=C^\infty(X;|\Lambda^d({}^\bullet T^*X)|)$ given by 
\begin{align}
    |\omega_{sc}|=a(x,y)\frac{dx dy}{x^{d+1}} \text{ and }  |\omega_{b}|=a(x,y)\frac{dx dy}{x}, \text{ with } a(x,y)>0.
\end{align}

On the model case of $X=\overline{\mathbb{R}^4}$ we will fix a standard representation for these metrics. Letting $x=\rho =r^{-1}$ and choosing a collar neighbourhood of $\partial X$ as $[0,\epsilon)_{\rho}\times \partial X \cong [0,\epsilon)_\rho\times \mathbb{S}^3$ as

\begin{align}\label{eq:bscmetrics}
    g_b =\frac{d\rho^2}{\rho^2}+h(\rho, \omega_i,d\omega_i) \text{ and } g_{sc}=\frac{d\rho^2}{\rho^4}+\frac{h(\rho,\omega_i,d\omega_i)}{\rho^2} 
\end{align}
where $\lim_{\rho\to 0} h(\rho,\omega_i,d\omega_i)=S^2(T^*\partial X)=S^2(T^*\partial\mathbb{S}^3)$ is the standard metric on $\partial X\cong \mathbb{S}^3$.
The associated standard $b$-densities and $sc$-densities on $X=\overline{\mathbb{R}^4}$ are then given as
\begin{align}
    dg_b &:= \frac{d\rho}{\rho} dg_{\mathbb{S}^3}\label{eq:bDensity} \\
    dg_{sc}&:= \frac{d\rho}{\rho^5}dg_{\mathbb{S}^3}.\label{eq:scDensity},
\end{align}
where in Hopf coordinates we give the standard volume density and metric on $\mathbb{S}^3$ through (\ref{eq:MinkowskiMetric}) by
\begin{align}
    dg_{\mathbb{S}^3}= \cos\theta\sin\theta  d\psi d\phi d\theta, \text{ and } g_{\mathbb{S}^3}=d\theta^2+\sin^2\theta d\phi^2+\cos^2\theta d\psi^2.
\end{align}
The scattering densitiy and the stationary the spacetime denisity $dg_{\bfa,\bfb,\mfrm}|_X$ (\ref{eq:StationarySplitting})differ only by only a smooth positive factor extending to the compactification of $X$.

It should be noted as well that in a fashion similar to the extension of the Euclidean vector fields $\partial_{x_i}$ to scattering vector fields on the boundary so too does the Eucleadean metric $dx_1^2\wedge...\wedge dx_4^2$ extend to the standard scattering metric $g_{sc}\in C^\infty(\overline{\mathbb{R}^4};S^2{}^{sc}T^*\overline{\mathbb{R}^4})$. Indeed writing out (\ref{eq:scDensity}) in standard spherical coordinates we see that
\begin{align}
     \frac{d\rho}{\rho^5}dg_{\mathbb{S}^3} = r^3dr dg_{\mathbb{S}^3}.
\end{align}
By contrast, the $b$-metric allows to identify the spherical compactification $\overline{\mathbb{R}}^4$ as a compactification of a Riemannian manifold with cylindrical end metric. This is obtained by choosing a chart of logarithimic coordinates near the boundary by 
\begin{align}\label{eq:CylindricalEndTransform}
    z=-\log \rho,
\end{align}
which transforms the $b$-metric to a cylindrical end metric on the collar neighbourhood $[0,\infty)_z\times \mathbb{S}^3$ of $\partial \overline{\mathbb{R}^4}$ as 
\begin{align}
    g_b= dz^2 + h(e^{-z},\omega_i,d\omega_i).
\end{align}
In general one can introduce $b$- and $sc$- differential operators as one would do with standard notation. 
\begin{definition}\label{def:bscDifferentialOperators}
    Let $X$ be a $d$-dimensional compact manifold with boundary and boundary defining function $\rho$. We define $b$- and $sc$- operators as the filtered algebras 
    \begin{align}
        \text{Diff}_\bullet(X)=\text{Span}_{C^\infty(X)}\{V_1\circ...\circ V_k\;|\; V_j\in \mathcal{V}_\bullet(X), k\leq m\}
    \end{align}
    where $\bullet \in \{b,sc\}$ and \textit{weighted} $b$- and $sc$- differential operators for $l\in \mathbb{R}$ through
    \begin{align}
        \text{Diff}_\bullet^{m,l}(X)=\rho^{-l}\text{Diff}_\bullet^m(X) \;\text{ for } \bullet\in \{b,sc\}.
    \end{align}
\end{definition}

\begin{remark}
    When working only on the interior $X^\circ$ we can simply define $\text{Diff}^m(X^\circ)$ analogously to Definition \ref{def:bscDifferentialOperators} using $V_k\in\mathcal{V}(X)$.
\end{remark}

Using Definition \ref{def:bscDifferentialOperators} gives us a first rough characterization of the spectral operator $\hat{\Box}_{g_{\bfa,\bfb,\mfrm}}(\sigma)$ for $\sigma\in \mathbb{C}$ as defined by (\ref{eq:SpectralFamily}) on $X$.

\begin{proposition}\label{prop:RoughCharacterization}
    We have that $\hat{\Box}_{g_{\bfa,\bfb,\mfrm}}(\sigma)\in \text{Diff}_{sc}^2(X)\cap \rho\text{Diff}_b^2(X)$ for all $\sigma\in \mathbb{C}$
\end{proposition}
\begin{proof}
    This can be checked by performing the conjugation (\ref{eq:FullWaveEquationEddingtonFinkelsteinCoord}) directly by performing the transformation $\partial_{t_*}\mapsto \sigma$ and then changing variables to $r=\rho^{-1}$ and anlyzing the resultant operator as a $b$- and $sc$- operator. It is more insightful however to appeal to asymptotical flatness by (\ref{eq:MinkowskiMetric}) from which follows that in the limit $r\to \infty$, $\rho\to 0$ we must have that $g_{\bfa,\bfb,\mfrm}\to g_{\mathbb{M}^{1+4}}$. It is well known that in Eddington-Finkelstein coordinates we have that 
    \begin{align}
        g_{\mathbb{M}^{1+4}}=-dt_*^2-2dt_*dr +r^2\cancel{g}_{\mathbb{S}^3},
    \end{align}
    and hence that the so-called "model operator" at $\partial\overline{\mathbb{R}^4}$ (see also \cite{vasy2013microlocal}, \cite{vasy2018minicourse}) must be given by 
    \begin{align}\label{eq:MinkowskiModel}
        \hat{\Box}_{\mathbb{M}^{1+4}}(\sigma)=\rho^2(-(\rho\partial_\rho)^2+2\rho\partial_\rho+\Delta_{\mathbb{S}^3})+\underbrace{2i\sigma\rho\left(\rho\partial\rho -\frac{3}{2}\right)}_{\text{transport term}}
    \end{align}
    From this it is immediate that that $\hat{\Box}_{\mathbb{M}^{1+4}}(\sigma) \in \text{Diff}^2_{sc}(X)\cap \rho\text{Diff}^2_b(X)$ for all $\sigma\in \mathbb{C}$ using the pedantic fact that $\rho\text{Diff}^1_b(X)\subset \rho\text{Diff}^2_b(X)$ for the transport term.
\end{proof}
\begin{remark}
    Observe here that as $\sigma\to 0$ we have that $\hat{\Box}_{g_{\mathbb{M}^{1+4}}}(\sigma)$ in equation (\ref{eq:MinkowskiModel}) shifts from being a $\rho\text{Diff}_b^2(X)$ operator to a $\rho^2\text{Diff}^2_b(X)$ operator. More precesily it holds that $\hat{\Box}_{g_{\mathbb{M}^{1+4}}}(\sigma) \in \rho^2\Diff_b^2(X)+\rho\sigma \Diff^1_b(X)$.
\end{remark}

With the $b$- and $sc$- differential operators in hand and their densities in hand we can now define the function spaces that $\hat{\Box}_{g_{\bfa,\bfb,\mfrm}}(\sigma)$ naturally acts on.
\begin{definition}\label{def:bscSobolevSpaces}
    Dentote by $\mathscr{D}'(X)$ the space of distrubitions on $X$. We introduce $b-$ and $sc-$ $L^2$ spaces $L^2_\bullet(X)$ using \textit{only} the scattering density $dg_{sc}$, that is 
    \begin{align}
        L^2_\bullet(X):=L^2(X,|dg_{sc}|)=\left\{u\in \mathscr{D}'(X)\;| \int_X |u|^2\, dg_{sc} <\infty\right\} ,
    \end{align}
with $|dg_\bullet|=d\text{vol}_{g_\bullet}$ with $\bullet \in \{b,sc\}$. For $s\in \mathbb{N}_0$ we then introduce $b$- and $sc$- Sobolev spaces as 
\begin{align}
    H^s_\bullet(X) :=\left\{u\in \mathscr{D}'(X)\;|\; Au\in L^2_\bullet(X),\forall A\in \text{Diff}_{\bullet}^s(X) \right\},
\end{align}
and equip these space with the Hilbert norm 
\begin{align}
    \|u\|^2_{H^s_\bullet}:=\sum^s_{j=0}\sum^j_{k=0}\|A_{jk}u\|^2_{L^2_\bullet(X)}
\end{align}
where the $A_{jk}$ are chosen such that $\text{span}\{A_{jk}\}_{k\leq j} = \text{Diff}_\bullet^j(X)$, for all $j\leq s$. For $s<0$ we extend this definition by duality and for $s\in \mathbb{R}\backslash \mathbb{Z}$ we extend it by interpolation.\\
Additionally, for $\ell \in \mathbb{R}$ we introduce the weighted spaces $H^{s,\ell}_\bullet(X)$ as 
\begin{align}
    H^{s,\ell}_\bullet(X)=\rho^\ell H^s_\bullet(X),
\end{align}
equipped with norms 
\begin{align}
    \|u\|_{H^{s,\ell}_\bullet(X)}:=\|\rho^{-\ell}u\|_{H^s_\bullet(X)}.
\end{align}

\end{definition}
\begin{remark}\label{rem:RestrictiontoOpenInterior}
    When working on the interior $X^\circ$ we define the standard Sobolev spaces $H^s(X^\circ)$ over the interior of $X$ exactly as above using $A_{jk}\in \text{Diff}^s(X^\circ)$. It should be noted that in particular we have that 
    \begin{align}\label{eq:EquivalentNormsonInterior}
        \|\cdot\|_{H^{s,\ell}_\bullet (X^\circ)} \sim \|\cdot\|_{H^s(X^\circ)},
    \end{align}
    as on the interior $X^\circ$ the boundary defining function $\rho \in S^0(T^*X^\circ)$ defines a harmless elliptic factor and gives equivalence of norms in (\ref{eq:EquivalentNormsonInterior}) . Hence we have for $s\in\mathbb{R}$ 
    \begin{align}
        H^{s,\ell}_\bullet(X^\circ)=H^s(X^\circ),
    \end{align}
    for any $\ell \in \mathbb{R}$ and $\bullet \in \{b,sc\}$. In terms of extendible distirubtions we also have 
    \begin{align}
        \|u\|_{H^s(X^\circ)}\asymp \|u\|_{\overline{H}^{s,r}_\bullet(X)}\asymp \|u\|_{\dot{H}^{s,r}_\bullet(X)}, \text{ for } \text{supp u}\subset K \Subset X^\circ.
    \end{align}
\end{remark}

Naturally then an operator $A\in \text{Diff}^{m,l}(X)$ for some $m\in \mathbb{Z}$, $l\in \mathbb{R}$ acts as a continuous operator
\begin{align}
    A:H^{s,\ell}_\bullet(X)\to H^{s-m,\ell-l}_\bullet(X).
\end{align}
Next to these spaces we also introduce the conventions 
\begin{align}
    H_\bullet^{\infty,\ell}(X) =\bigcap_{s\in \mathbb{R}} H^{s,\ell}(X), \; \text{ and } H^{-\infty,\ell}_\bullet(X)=\bigcup_{s\in \mathbb{R}}H^{s,\ell}_\bullet(X).
\end{align}
and the convention for approximate decay as
\begin{align}
    H^{s,\ell+}_\bullet(X) =\bigcup_{\epsilon>0}H^{s,l+\epsilon}_\bullet(X) \text{ and } H^{s,\ell-}_\bullet(X)=\bigcap_{\epsilon>0}H_\bullet^{s,\ell-\epsilon}(X),
\end{align}
for $s\in \mathbb{R}\cup\{\pm \infty\}$.
Our eventual Price-law result will be formulated for conormal functions, which are basically functions that are "smooth" with respect to $b$-differentiation.

\begin{definition}
    Let us define $\mathcal{A}^0(X)$ through 
    \begin{align}
        \mathcal{A}^0(X) = \{u\in L_b^\infty(X)\;|\; Au\in L^\infty_\bullet(X)\;\text{for all } A\in \text{Diff}^1_b(X)\}
    \end{align}
    and we define the space of conormal functions of weight $\ell\in \mathbb{R}$ as
    \begin{align}
        \mathcal{A}^\ell(X) :=\rho^\ell\mathcal{A}^0(X).
    \end{align}
    Similarly we will set 
    \begin{align}
        \mathcal{A}^{\ell+}(X)=\bigcup_{\epsilon>0}\mathcal{A}^{\ell+\epsilon}(X) \text{ and } \mathcal{A}^{\ell-}(X)=\bigcap_{\epsilon>0}\mathcal{A}^{\ell-\epsilon}(X).
    \end{align}
\end{definition}

With respect to the $b$-conormal spaces we will make use of the following sobolev embedding.

\begin{lemma} ($b$-Sobolev embedding) \label{lem:SobolevEmbedding} For $X=\overline{\mathbb{R}^4}$ we have the following inclusions
\begin{align}\label{eq:SobolevEmbedding}
    H_b^{\infty,\ell}(X) \hookrightarrow \mathcal{A}^{\ell+2}(X)\hookrightarrow H^{\infty,\ell-}_b(X). 
\end{align}
\end{lemma}
\begin{proof}
    This is a direct consequence of the Sobolev embeddings on $\mathbb{R}^4$. Let $u\in H^{\infty,\ell}_{b}(X)$, define  
    \begin{equation}
        w:=\rho^{-(\ell+2)}u,
    \end{equation}
    then 
    \begin{equation}
        \|w\|^2_{L^2(d\rho d\omega)}=\int |\rho^{-(\ell+2)}u|^2d\rho d\omega=\int|u|^2\rho^{-2\ell-4}d\rho d\omega,
    \end{equation}
    now using (\ref{eq:bDensity}), 
    \begin{equation}
        \|w\|_{L^2(d\rho d\omega)}=\|\rho^{-\ell}u\|_{L^2(|dg_X|)} =\|u\|_{L^{2,\ell}_b}.
    \end{equation}
    By definition of $w$, we have (writing out in the basis $\rho\partial_\rho$ and $\partial_{\omega_i}$ for $i=1,2,3$) that 
    \begin{equation}
        u\in H_b^{N,\ell}((0,\epsilon)_\rho\times \mathbb{S}^3) \implies w \in H^N((0,\epsilon)_\rho \times \mathbb{S}^3)\; \text{ for } N>2.
    \end{equation}
    Now by the ordinary Sobolev embedding we have \begin{equation}
    H^N((0,\epsilon)_\rho\times\mathbb{S}^3)\hookrightarrow L^\infty((0,\epsilon)_\rho\times\mathbb{S}^3)
    \end{equation}
    for $N>2$. 
    Thus the statement $w\in \mathcal{A}^0(X)$ i.e. $(\rho\partial_\rho)^k \partial^\beta_{\omega} w\in L^\infty((0,\epsilon)_\rho\times\mathbb{S}^3)$ for arbitrary $k,\beta$ then implies by translating back that $(\rho\partial_\rho)^k \partial^\beta_{\omega} u\in \rho^2 L^\infty((0,\epsilon)_\rho\times\mathbb{S}^3)$, i.e. $u\in \rho^2\mathcal{A}^{\ell+2}(X)$, which is the first inclusion.
    To show the second inclusion assume that $u \in \mathcal{A}^\alpha (X)$ then $\rho^{-\alpha}u:=a$ is conormal and thus $a \in H^N((0,\epsilon)_\rho \times \mathbb{S}^3)$. However for the $b$-density (\ref{eq:bDensity}) we have an extra $\rho^{-4}$ weight. Therefore to ensure integrability for every $b$-derivative of $a$ at the endpoint $\rho=0$ we need to weaken the $\rho$-weight slightly to allow for the second inclusion.  
\end{proof}

Finally, to end our discussion on the relevant spaces we will now slightly adjust from the model case $X=\overline{\mathbb{R}^4}$ that we need to the case that $X=\overline{\mathbb{R}^4\backslash B(0,r_X)}$ as defined in Definition \ref{def:RadiallyCompactifiedR4}. Suppose therefore that $X'$ is a compact manifold with boundary and let $X\subset X'$ be a submanifold with boundary, and suppose that its boundary decomposes into two non-empty sets 
\begin{align}
    \partial X=\partial_-X \sqcup \partial_+X,\;\;\;\;\partial_-X = \partial X', \;\partial_+X=\partial X'
\end{align}
and let $\partial_+X$ be the boundary 'at infinity' whilst $\partial_-X$ is the interior boundary. We then define the $b$ and $sc$- vectorfields with respect to the larger set $X'$ as 
\begin{align}
    \mathcal{V}_b(X):=V|_X\;|\; V\in \mathcal{V}_b(X')\},\;\;\;\mathcal{V}_{sc}(X):=\{V|_X\;|\; V\in \mathcal{V}_{sc}(X')\},
\end{align}
which are $b$- $sc$- vector fields at infinity but unrestricted at $\partial_-X$.  In our specific case of $X\subsetneq \overline{\mathbb{R}^4}:=X'$ as in Definition \ref{def:RadiallyCompactifiedR4} with $r_X>0$ we have that $\partial_+X=\partial\overline{\mathbb{R}^4}$ and $\partial_-X=\partial B(0,r_X)$.

We then define the two natural classes of Sobolev spaces, the \textit{extendible distrubitions}, given by 
\begin{align}
    \overline{H}^{s,\ell}_\bullet(X):=\{u|_{X^\circ}\;|\; u \in H^{s,\ell}_\bullet (X')\},\;\bullet \in \{b,sc\}
\end{align}
and the \textit{supported distributions} given by 
\begin{align}
    \dot{H}^{s,\ell}_{\bullet}(X):=\{u \in \mathscr{D}'(X) \;|\; u\in H^{s,\ell}_\bullet(X').\;\text{supp }u\subset X\}.
\end{align}

We note that localizing away from $\partial_-X$ elements of these spaces automatically lie in $H_{\bullet}^{s,\ell}(X)$, thus the subspaces of $\overline H^{s,\ell}_\bullet(X)$ and $\dot{H}^{s,\ell}_\bullet(X)$ consisting of elements that are conormal at $\partial_+X$ are well-defined. With respect to taking formal adjoints on Myers-Perry we will use the following definition

\begin{definition}\label{def:FormalAdjointPairing}
Let $(M,g_{\bfa,\bfb,\mfrm})$ be a Myers-Perry spacetime. We define pairings and formal adjoints on $X$ with respect to the stationary spatial density
$dg_{\bfa,\bfb,\mfrm}|_X$ defined in \eqref{eq:StationarySplitting}.
More precisely, we write
\begin{equation}
  \langle u,v\rangle
  := \langle u,v\rangle_{L^2(X;dg_{\bfa,\bfb,\mfrm}|_X)}
  = \int_X u\,\overline{v}\,dg_{\bfa,\bfb,\mfrm}|_X,
  \qquad u,v\in C_c^\infty(X^\circ).
\end{equation}
We use the same notation whenever this integral converges and for
its continuous extensions by duality. Thus the formal adjoint
$P^*$ of a differential operator $P$ on $X$ is characterized by
\begin{align}
  \langle Pu,v\rangle=\langle u,P^*v\rangle,
  \qquad u,v\in C_c^\infty(X^\circ).
\end{align}
\end{definition}

The idea of the proof laid out in the proof of Proposition \ref{prop:RoughCharacterization} can be made even more rigorous with the concept of the normal operator.

\begin{definition}\label{def:bscNormalOperator}
   Let $X$ be a compact smooth manifold with boundary with boundary defining function $\rho\in C^\infty(X)$. Let $P\in \Diff^{m,l}_\bullet(X)$ and let $l = 0$. We define its $b$- and $sc$- normal operators $N^{\bullet}_{\partial X}(P) \in \text{Diff}^{m,0}_{\bullet}(X)$ through ''freezing'' coefficients at $\rho=0$. Specifically in adapted coordinates $(\rho,y)$ at a collar neighbourhood $[0,\epsilon)_\rho\times \partial X$ of $\partial X$ we set for $P\in \text{Diff}^{m,0}_{sc}(X) $ we have that 
    \begin{align}
        N^{sc}_{\partial X}(P) = \sum_{j+|\alpha|\leq m} a(\rho,y)|_{\rho=0}(\rho^2 D_\rho)^j(\rho D_y)^\alpha.
    \end{align}
    and for $P\in \text{Diff}^{m,0}_{b}(X) $ we have that 
    \begin{align}
        N^{b}_{\partial X}(P) = \sum_{j+|\alpha|\leq m} a(\rho,y)|_{\rho=0}(\rho D_\rho)^j(D_y)^\alpha,
    \end{align}
    where we introduced the convention 
    \begin{align}
        D_x:= \frac{1}{i}\partial_x.
    \end{align}
    We reserve the subscript notation $\partial \mathcal{M}$ to indicate what kind of boundary we freeze coefficients at if needed, but it will almost always be 'at infinity' $\partial_+X$. When $P\in \text{Diff}_\bullet^{m,\ell}(X)$ we define the normal operator by rescaling 
    \begin{align}
        N^{\bullet,(l)}_{\partial X}(P)= N^{\bullet}_{\partial X}(\rho^{-l}P).
    \end{align}
\end{definition}

We will almost never use the $b$-normal operator as defined in Definition \ref{def:bscNormalOperator} but instead use the indicial operator, which is related to the $b$ normal operator through the cylindrical end transform (\ref{eq:CylindricalEndTransform}).

\begin{definition}
    Let $X$ be a compact manifold with boundary $\rho$ and let $u\in H^{s,\ell}_b(X)$. We then define the Mellin transform $\mathscr{M}u$ through 
    \begin{align}
        \mathscr{M}_{\rho\to\lambda}u = \hat{u}(\lambda) :=\int^\infty_0 \rho^{-i\lambda} u(\rho)\frac{d\rho}{\rho}.
    \end{align}
    For an operator $P\in \text{Diff}^{m,l}(X)$ we then define the indicial operator $\hat{N}^b(P,\lambda):H^{s}(\partial X)\to H^{s-m}(\partial X)$ through
    \begin{align}
        \hat{N}(P,\lambda)=\mathscr{M}_{\rho\to\lambda}\circ N^{b,(l)}_{\partial X}(P) \circ \mathscr{M}^{-1}_{\lambda\to\rho} = \sum_{j+|\alpha|\leq m}a(\rho,y)|_{\rho=0}(-\lambda)^j (D_y)^\alpha
    \end{align}
    
\end{definition}
Note that
\begin{align}
    \mathscr{M}_{\rho\to\lambda}(\rho\partial_\rho u)=i\lambda \mathscr{M}_{\rho\to\lambda}u \text{ or } \mathscr{M}_{\rho\to\lambda}(\rho D_\rho u)= \lambda.
\end{align}
and thus that if $P\in \text{Diff}^{m,l}_b(X)$ the indicial operator defines a continuous map 
\begin{align}
    \hat{N}(P,\lambda):H^{s}(\partial X)\to H^{s-m}(\partial X)
\end{align}
where for the case $X=\overline{\mathbb{R}^4}$ we have $\partial X\cong \mathbb{S}^3$. Recall as well that the Mellin transform of a $b$-operator is related to the Fourier transform of $X$ with a cylindrical end metric under the logarithmic change of coordinates $z=-\log\rho$ i.e.
\begin{align}
    \mathscr{M}_{\rho\to\lambda}u=\hat{u}(\lambda) =\int^\infty_{-\infty}e^{i\lambda z}u(e^{-z})dz =\mathscr{F}_{z\to \lambda}u(e^{-z}).
\end{align}
With respect to the Mellin transform we end this section with an introduction to indicial roots to the model equation (\ref{eq:MinkowskiModel}).
\begin{definition}
    Let $X$ be a compact smooth manifold with boundary and denote its boundary defining function by $\rho$. Denote for $P\in \Diff^{m,\ell}_b(X)$ its indicial operator by $\hat{N}^{b,(l)}_{\partial X}(P,\lambda)$. An \textit{indicial root} of $P$ is a complex number $\lambda\in \mathbb{C}$ such that
    \begin{align}
        \ker \hat{N}^{b,(l)}_{\partial X}(P,\lambda)\neq \{0\},
    \end{align}
    or equivalently, there exists a $w\in C^\infty(\mathbb{S}^3)$ such that $N^{b,(\ell)}_{\partial X|_{\infty}}P(\rho^{i\lambda} w) =0$. The \textit{indicial} gap of $P$ is an interval $(\beta_-,\beta_+)$ is an interval such that $\text{Im } \lambda \notin (\beta_-,\beta_+)$. In terms of the the spectral family $\hat{\Box}_g(\sigma)$ and its Mellin traasform they are found as the roots $\lambda \in \mathbb{C}$ for $\hat{N}^{b,(l)}_{\partial X}(\hat{\Box}_g(\sigma),\lambda)$ and $\hat{N}^{b,(l)}_{\partial X}(\hat{\Box}_g(0),\lambda)$ respectively.
\end{definition}

The following proposition is deeply interconnected with the emergence of the weak Huygens principle in the resolvent expansion in Section \ref{sec:PriceLaw}, which is represented by the fact that the indical roots to (\ref{eq:MinkowskiModel}) are fractional.

\begin{proposition}\label{prop:IndicialDecayRootsMinkowski}
    For the model operator (\ref{eq:MinkowskiModel}) on $\mathbb{M}^{1+4}$ we have that the indicial roots are given by $\lambda =0, \lambda=-i\frac{3}{2}$. Its indicial gap is given by $(\beta_-,\beta_+)=(-2,0)$. 
\end{proposition}
\begin{proof}
    Follows from immediate calculation where we observe that the term $-2i\sigma \left(\rho\partial_\rho-\frac{3}{2}\right)$ dominates when $\sigma\neq 0$. The indicial gap is actually found from setting $\sigma=0$ in which case $\hat{\Box}(0)\in \text{Diff}^{2,2}_b(X)$. The indicial $b$-operator is then found as 
    \begin{align}
        \hat{N}^{b,(2)}_{\partial X}(\hat{\Box}_g(0))= (i\lambda)^2+2i\lambda +\Delta_{\mathbb{S}^3}.
    \end{align}
    Decomposing the eigenfunctions of $\Delta_{\mathbb{S}^3}$ in Wigner $D$-functions $S_{\kappa l m}$ we see that after separating variables, the $\lambda \in \mathbb{C}$ satisfy the characteristic equation
    \begin{align}
        ((i\lambda)^2+2i\lambda +l(l+2))S_{kml}=0,
    \end{align}
    for $l \in \mathbb{Z}$, from which it follows that $\lambda=0$ and $\lambda = -2i$ are the first 2 roots.
\end{proof}

\subsection{Pseudodifferential operators and semiclassical analysis.}\label{subsec:bandscPseudodiffernetialOperators}

We now push the notions of the previous chapter further and embed them into a pseudodifferential $b$- and $sc$- framework.\\

We first note that when we assume $X$ to be a smooth manifold with boundary with $n-1$ corners, one can easily extend the Sobolev spaces and conormal spaces defined in Definition \ref{def:bscSobolevSpaces} to the $n$ piecewise smooth boundary segments. That is, if $X$ has $n\in \mathbb{N}_0$ boundary faces one can define 
\begin{align}
    H_\bullet^{s,\mathfrak{l}}(X) \text{ and } \mathcal{A}^{\mathfrak{l}}(X), \text{ with } \mathfrak{l}\in \mathbb{R}^n, \bullet \in \{b,sc\},
\end{align}
with respect to adapted coordinates 
\begin{align}
    [0,\infty)^n_{x_i}\times \mathbb{R}^{d-n}_{y_j},
\end{align}

where we work with boundary weights 
\begin{align}
    \rho_1^{\mathfrak{l}_1}\cdot \rho_2^{\mathfrak{l_2}}\cdot...\cdot \rho_n^{\mathfrak{l}_n}.
\end{align}
The most relevant example in our case are obtained when $X$ is a compact manifold with boundary that is then also fibrewise compactified in ${}^\bullet T^*X$ as already introduced in Definition \ref{def:FibrewiseCompactification}. In this case we write ${}^\bullet\overline{T} X$ and we will denote the fibre-radial boundary defining function by $\varrho$.

The corner of ${}^\bullet \overline{T}X$ is then given by the set $\partial {}^\bullet\overline{T^*}_{\partial X}X=\{\rho=0\}\cap \{\varrho=0\}\subset{}^\bullet \overline{T^*}X$.\\ 

We then define the $b$- and $sc$- symbol spaces $S^{m,l}_\bullet({}^\bullet TX)$ as 
\begin{align}
    S^{m,l}_\bullet({}^\bullet T^*X)=\varrho^{-m}_\bullet\rho^{-l}\underbrace{\mathcal{A}^{0,0}({}^\bullet\overline{T}^*X)}_{:=S^{0,0}_\bullet({}^\bullet T^*X)},
\end{align}
i.e. we measure differentiability order simply as fibre-radial decay w.r.t $\varrho_\bullet$ on ${}^\bullet\overline{T^*}X$. \\
We can then generate the classes of $b$- and $sc$- pseudo-differential operators $\Psi^{m,l}_\bullet(X)$ for $m\in \mathbb{R}$ $l\in \mathbb{R}$ through quantization of these symbol-classes i.e. 
\begin{align}
    \Psi^{m,l}_\bullet(X)=\{A:H_\bullet^{s,\ell}(X)\to H^{s-m,\ell-l}_\bullet(X)\;|\; A=\text{Op}_\bullet(a) \text{ with } a\in S^{m,l}({}^\bullet T^*X) \},
\end{align}
which define bounded operators between $H^{s,\ell}_\bullet(X)$ and $H^{s-m,\ell-l}_\bullet(X)$ in the case that $X$ is compact\footnote{In general properly supported pseudodifferential operators $P \in \Psi^{m,\ell}_\bullet (X)$ only define operators $P: H^{s,r}_{\bullet,c}(X) \to H^{s-m,r-\ell}_{\bullet,loc}(X)$, i.e. between Fr\'echet spaces. In our case we are mostly interested in the case where $X$ is compact however, where they reduce to proper norms.}.

Where for $A=\text{Op}_b(a)\in \Psi_b^{m,l}(X)$ we define the quantization for any $u\in \mathscr{D}'(X)$, $a \in S_b^{m,l}({}^\bullet T^*X)$ in adapted coordinates $(\rho,y)$ via
\begin{align}\label{eq:QuantizationbT*M}
    &Au(\rho,y):=\text{Op}_b(a)u(\rho,y)\\
    &=(2\pi)^{-d}\int_0^\infty\int_{\mathbb{R}^{d-1}}\int_{\mathbb{R}^d} e^{i \log(\rho/\rho')\cdot{}^b\xi+(y-y')\cdot {}^b\eta} a(\rho,y,{}^{b}\xi,{}^b\eta) u(\rho',y') \frac{d\rho'}{\rho'} dy'_id{}^b\xi d{}^b\eta.
\end{align}
Likewise $A=\text{Op}_{sc}(a) \in \Psi^{m,l}_{sc}(X)$ for $a\in S^{m,l}_{sc}({}^{sc} T^*X)$ and $u\in \mathscr{D}'(X)$  are quantized by using the scattering density (\ref{eq:scDensity}) and give
\begin{align}
    &Au(\rho,y) := \text{Op}_{sc}(a)u(\rho,y)\label{eq:QuantizationscT*M}\\
    &= (2\pi)^{-d}\int_0^\infty\int_{\mathbb{R}^{d-1}}\int_{\mathbb{R}^d} e^{i(\rho'^{-1}-(\rho)^{-1})\cdot{}^{sc}\xi+\frac{(y-y')\cdot{}^{sc}\eta}{\rho}} a(\rho,y,{}^{sc}\xi,{}^{sc}\eta)u(\rho',y') \frac{d\rho'}{\rho'^2}\frac{dy_i}{\rho'^{3}}d{}^{sc}\xi d{}^{sc}\eta_i.\nonumber 
\end{align}

In the above we adopted the convention for the Liouville $1$ form ${}^\bullet\zeta\in \Gamma ({}^\bullet T^*X)$ for any bundle in adapted boundary coordinates $(\rho,y)$ to be of the form 
\begin{align}
    {}^b\upalpha= {}^b \xi \frac{d\rho}{\rho} + \sum^3_{i=1} {}^b\eta_i dy_i \text{ and } {}^{sc}\upalpha= {}^{sc} \xi \frac{d\rho}{\rho^2} + \sum^3_{i=1} \frac{{}^b\eta_i}{\rho} dy_i \label{eq:Canonicalbsc1Forms}
\end{align}

To identify kernels of these maps, and also for computations to come we recall the definition of a blow-up space.

\begin{definition}[Real blow-up]
Let $X$ be a manifold with corners, with adapted local coordinates
\begin{align}
 (x,y)\in[0,\infty)^k_{x_i}\times\mathbb R^{n-k}_{y_i},
\end{align}
where the \(x_i\) are boundary-normal coordinates. A closed subset
\(Y\subset X\) is a \(p\)-submanifold if, near each of its points, such
coordinates can be chosen so that
\begin{align}
 Y=\{x_i=0,\; y_j=0\}_{i\in I,j\in J}
\end{align}
for subsets $I\subset\{1,\ldots,k\}$ and
$J\subset\{1,\ldots,n-k\}$.

The real blow-up $[X;Y]$ replaces $Y$ by the inward-pointing
spherical directions of its normal variables $(x_I,y_J)$. Locally
its coordinates are
\begin{align}
    (r_Y,\omega,x_{I^c},y_{J^c})
\end{align}

\begin{align}
     r_Y=\left(\sum_{i\in I}|x_i|^2+\sum_{j\in J}|y_j|^2\right)^{1/2},\qquad
 \omega=\frac{(x_i,y_j)}{r_Y}\in \mathbb{S}_+^{I+J-1}.
\end{align}

where \(\omega_{x_i}\geq0\) for \(i\in I\). The blow-down map
\(\beta:[X;Y]\to X\) is given by
\begin{align}
    \beta(r_Y,\omega,x_{I^c},y_{J^c})=(r_Y\cdot\omega_x,x_{I^c},r_Y\cdot \omega_y,y_{I^c}),
\end{align}
and is a diffeomorphism away from $Y$. We will refer to the hypersurface that is the lift of $Y$ inside $[X;Y]$ as the front face of $Y$, given by
\begin{align}
    \text{ff}_Y =\overline{\beta^{-1}(Y)}=\{r_Y=0\}.
\end{align}
\end{definition}

\begin{remark}[Projective coordinates for blow-ups]
It will be often be useful here as well to work on projective coordinates $\beta^{-1}(Y)$ in $[X;Y]$. Letting 
\begin{align}
    z=(x_I,y_J), 
\end{align}
we pick a chart for which $z_a\gtrsim z_b$ and work on 
\begin{align}
     \hat{z}_a :=|z_a| \text{ and } \hat{z}_b = \frac{z_b}{|z_a|},
\end{align}
whilst leaving the coordinates $(x_{I^c},y_{J^c})$ unchanged. The obvious boundary defining function for $\text{ff}_Y$ in this chart is then obviously 
\begin{align}
    \rho_{\text{ff}_Y} =|x_a|.
\end{align}
Note that this convention is also consistent with the projective coordinates used in the (fibre)-radial compactifications we introduced in Section \ref{sec:MyersPerryBoyerLindquist} and Definition \ref{def:RadiallyCompactifiedR4}. Indeed we have that 
\begin{align}
    \overline{\mathbb{R}^d} \backslash \{0\}\cong [\mathbb{R}^d;\{0\}]
\end{align}
through radial inversion $r=\rho^{-1}$. Indeed in the case of fibre-infinity as given in Section \ref{sec:MyersPerryBoyerLindquist}  the equality
\begin{align}
    \overline{T}^*X^\circ\backslash O = [T^*X^\circ;O]
\end{align}
is made explicit through 
\begin{align}
    \frac{\xi_b}{|\xi_a|}=\hat{z}_b, \text{ and } \frac{1}{|\xi_a|}=\hat{z}_a \cdot\left(1+\sum_{b\neq a} \hat{z}_b^2\right)
\end{align}

\end{remark}

Within the framework of blow-ups, the Schwartz kernels $K_b$ and $K_{sc}$ for equations (\ref{eq:QuantizationbT*M}) and (\ref{eq:QuantizationscT*M}) can be stated on the double spaces
\begin{align}
    X_b^2 = [X\times X;\partial X\times \partial X] \text{ and } X^2_{sc}= [X_b^2;\partial \Delta_b],
\end{align}
with blowdown maps $\beta_b:X^2_b\to X\times X$ and $\beta_{sc}:X^2_{sc}\to X^2_b$
where $\partial \Delta_b$ stands for the boundary of the lifted diagonal $\beta^{-1}_{sc}(\Delta_b)$ of $\Delta_b=$ inside $X_b^2$.
Explictly then $K_b\in \mathscr{D}'(X_b^2)$ and $K_{sc}\in \mathscr{D}'(X_{sc}^2)$ are locally given by 
\begin{align}
    K_b(\rho,\rho',y,y'):=\frac{\chi(\rho,y,\rho',y')}{(2\pi)^d}\int_{\mathbb{R}^d}e^{i \log(\rho/\rho')\cdot{}^b\xi+(y-y')\cdot {}^b\eta} a_b(\rho,y,{}^{b}\xi,{}^b\eta)\; d\bxi d\beeta
\end{align}
and 
\begin{align}
    K_{sc}(\rho,y,\rho',y'):=
    \frac{\chi(\rho,\rho',y,y')}{(2\pi)^d}\int_{\mathbb{R}^d} e^{i(\rho^{-1}-(\rho')^{-1})\cdot{}^{sc}\xi+\frac{(y-y')\cdot{}^{sc}\eta}{\rho}} a_{sc}(\rho,y,{}^{sc}\xi,{}^{sc}\eta)d{}^{sc}\xi d{}^{sc}\eta_i\nonumber,
\end{align}
for cut-offs $\chi \in C_c^\infty(X_\bullet^2)$ and symbols $a_\bullet \in S^{m,l}_\bullet({}^\bullet T^*X)$, which are smooth maps away from the diagonals in $X_\bullet^2$. We say that $K_b$ and $K_{sc}$ are properly supported when their projections on each product are proper maps.

Note that we have 
\begin{align}
    S^{m,l}_\bullet({}^\bullet T^*X)=\rho^{-l} S^{m}_\bullet({}^\bullet T^*X) \text{ and } \Psi_\bullet^{m,l}(X) = \rho^{-l} \Psi^m_\bullet(X).
\end{align}

We note note that the measure used here are the $b$-measure (\ref{eq:bDensity}) and $sc$-measure (\ref{eq:scDensity}). We also remark that (\ref{eq:QuantizationbT*M}) can be functionally be transformed into the classical quantization $\text{Op}(a)u(e^{-x})$ through the logarithmic change of coordinates $x=-\log\rho$. We recall the classical operators $A\in\Psi^m(X^\circ)$ for $m\in \mathbb{R}$ for $X^\circ=\mathbb{R}^d$ 
\begin{align}
    Au(x)=\text{Op}(a) u(x) = (2\pi)^{-d} \int_{\mathbb{R}^d} e^{i x\cdot\xi} a(x,\xi) \hat{u}(\xi)  d\xi, \text{ for } u\in \mathscr{D}'(X)
\end{align}
which can be generalized to arbitrary manifolds $X$ by pulling back from charts $\phi:U\to \mathbb{R}^d$ as in \cite{hintz2025introduction}. On $\overline{\mathbb{R}^4}$ the associated classical kernel
\begin{align}
    K(x,y)= (2\pi)^{-d}\int_{\mathbb{R}^{2d}} e^{i(x-y)\cdot\xi}a(x,\xi)d\xi,
\end{align}
defines a scattering kernel globally as well after extending to $\partial\overline{\mathbb{R}^4}$.
\begin{remark}\label{rem:RestrictionPDOStoOpenInterior}
    Note that also in the context of Remark \ref{rem:RestrictiontoOpenInterior} when $X^\circ$ is not compact, pseudodifferential operators $A\in \Psi^m(X^\circ)$ generally only define continuous operators between Fr\'echet spaces $A:H^s_{loc}(X^\circ)\to H^{s-m}_c(X^\circ)$ with 
    \begin{align}
        H^s_c(X^\circ) =\{u\in H^s_{loc}(X^\circ) \;|\; \text{supp}(u)\subset X^\circ \text{ is compact.}\}
    \end{align}
\end{remark}

For $A\in \Psi^{m,l}_\bullet(X)$ are bounded operators $A:H^{s,\ell}_{\bullet}(X)\to H^{s-m,\ell-l}_{\bullet}(X)$. 
\begin{align}
    \Psi^{m,l}_\bullet(X)\circ \Psi^{m',l'}_\bullet(X)\subseteq \Psi^{m+m',l+l'}_\bullet(X), 
\end{align}
through which we also define residual operators as
\begin{align}
    \Psi^{-\infty,-\infty}_{sc}(X) = \bigcap_{\substack{m\in\mathbb{R} \\ l\in\mathbb{R}}}\Psi_{sc}^{m,l}(X) \text{ and } \Psi_b^{-\infty,l}(X) =\bigcap_{m\in\mathbb{R}}\Psi_b^{m,l}(X).
\end{align}
which are generated by smoothing oeprators
\begin{align}
    S^{-\infty,-\infty}_{sc}({}^{sc}T^*X)=\bigcap_{m,l\in \mathbb{R}} S^{m,l}_{sc}({}^{sc}T^*X)\; \text{ and } \;  S^{-\infty,l}_{b}({}^{sc}T^*X)=\bigcap_{m\in \mathbb{R}} S^{m,l}_{b}({}^{sc}T^*X).
\end{align}

Moreover we have for $m \in \mathbb{N}_0$ the special case that 
\begin{align}
    \Diff^{m,l}_\bullet (X)\subsetneq \Psi^{m,l}_\bullet(X),
\end{align}
with kernel that is smoothly supported on the diagonals of $X_\bullet^2$.\\
Now for an arbitrary symbol $p_\bullet \in S^{m,\ell}_\bullet({}^\bullet T^*X)$ with $P =\text{Op
}_\bullet(p)$ we can associate a canonical Hamiltonian vector field $H_p$ and associate its corresponding dynamics. Given the Liouville $1$-forms (\ref{eq:Canonicalbsc1Forms}) we can define a canonical symplectic form through 
\begin{align}
    {}^\bullet\upomega := d{}^\bullet\upalpha,
\end{align}
with respect to which ${}^\bullet H_p\in \mathcal{V}_\bullet({}^\bullet T^*X)$ is canonically defined as 
\begin{align}
    \iota_{{}^\bullet H_p}{}^\bullet\upomega = -dp_\bullet.
\end{align}
Its integral curves $\gamma\in {}^\bullet T^*X$ satisfying 
\begin{equation}
    \gamma' = {}^\bullet H_p\gamma,
\end{equation}
we will refer to again as bicharacteristic curves.
Concretely for $\bullet\in \{b,sc\}$ we have 
\begin{align}
    {}^b\upomega &= \rho^{-1}(d{}^b\xi\wedge d\rho) +\sum_{i} d{}^b\eta_i \wedge dy_i\\  {}^{sc}\upomega &= \rho^{-2}(d{}^b\xi\wedge d\rho)+ \rho^{-1}\sum_i[ d{}^{sc}\eta_i\wedge dy_i -\rho^{-1} {}^{sc}\eta_i(d\rho\wedge dy_i)],
\end{align}
and we compute using the conventions in (\ref{eq:Canonicalbsc1Forms}) that 
\begin{align}
    {}^bH_{p}&=  \frac{\partial p_b}{\partial\bxi}\rho\partial_\rho +\sum_i \frac{\partial p_b}{\partial\beeta_i}\partial_{y_i}- \rho \frac{\partial p_b}{\partial \rho}\partial_{\bxi} -\sum_i \frac{\partial p_b}{\partial y_i} \partial_{\beeta_i}\label{eq:bHamiltonVFieldConvention}\\
    {}^{sc}H_p&=  \frac{\partial p_{sc}}{\partial\scxi}\rho^2\partial_\rho + \sum_i \frac{\partial p_{sc}}{\partial \eta_i}\rho\partial_{y_i}\nonumber\\
    &- \left(\rho^2 \frac{\partial p_{sc}}{\partial \rho}+\sum_i \eta_i\rho\frac{\partial p_{sc}}{\partial\sceta_i}\right)\partial_{\scxi}-\rho\sum_i\left(\frac{\partial p_{sc}}{\partial y_i}-\sceta_i\frac{\partial p_{sc}}{\partial \scxi}\right)\partial_{\eta_i}\label{eq:scHamiltonVFieldConvention}
\end{align}

It should be mentioned once again that on the interior $X^\circ$ the flow ${}^\bullet H_{p}$ is conjugate to the standard Hamiltonian flow $H_p$ when considering $p_\bullet \in S^0(T^*X^\circ)$, which can be computed easily using the isomorphisms $\iota_\bullet :{}^\bullet T^*X^\circ \to T^*X^\circ$. Thus on $X^\circ$ we often simply work with the Liouville form 
\begin{align}
    \upalpha = \sum_i \xi dx_i, \text{ and } \upomega_i = \sum_i d\xi_i\wedge dx_i
\end{align}
giving the well-known form
\begin{align}
    H_p = \sum_i \frac{\partial p}{\partial\xi_i}\partial_{x_i}-\frac{\partial p}{\partial x_i}\partial_{\xi_i},
\end{align}
which thus separates the dynamics of an operator $P=\text{Op}_\bullet (p) \in \Psi^{s,l}_\bullet (X)$ to the boundary $\partial X$ using the full $b$ and $sc$-structures and the interior $X^\circ$ we we can study the classical Hamiltonian vector field.
In general for $p_\bullet \in {}^\bullet T^*X$ the only dynamics that are invariant under coordinate transformations are those that come from the so-called principal symbol $\upsigma_\bullet(p)$, which come from the \textit{principal symbol map} 
\begin{align}
    {}^\bullet\upsigma: \Psi^{m,l}_\bullet(X) \to S^{m,l}({}^\bullet T^*X),
\end{align}
that is induced via the short exact sequences 
\begin{align}\label{eq:bSES}
    0\to \Psi_b^{m-1,l}(X) \to \Psi_b^{m,l}(X)  \overset{{}^{b}\sigma_{m,l}}{\to} S^{m,l}({}^b T^*X)/S^{m-1,l}({}^{b} T^*X) \to 0.
\end{align}

and 

\begin{align}\label{eq:ScatteringSES}
    0\to \Psi_{sc}^{m-1,l-1}(X) \to \Psi_{sc}^{m,l}(X)  \overset{{}^{sc}\sigma_{m,l}}{\to} S^{m,l}({}^{sc} T^*X)/S^{m-1,l-1}({}^{sc} T^*X) \to 0.
\end{align}

For a representative $a_\bullet \in S_\bullet^{m,l}({}^\bullet\overline{T^*}X)$ of the principal symbol, we define its normalized symbol by 
\begin{align}
    \tilde{a}_\bullet:=\varrho_\bullet^m \rho^l a_\bullet, \;\;\bullet \in \{b,sc\},
\end{align}
if $a'_\bullet \in S^{m,l}_\bullet({}^\bullet \overline{T}^*X)$ is another symbol satisfying 
\begin{align}
{}^{\bullet}\upsigma_{m,l} (a) = {}^\bullet\upsigma_{m,l}(a') 
\end{align}
then the quotients in (\ref{eq:bSES}) and (\ref{eq:ScatteringSES}) give 
\begin{align}
    \tilde a_b -\tilde a'_b &\in \varrho_b \mathcal{A}^{0,0}({}^{b}\overline{T^*}X)\\
    \tilde a_{sc}-\tilde a_{sc}' &\in \varrho_{sc}\rho \mathcal{A}^{0,0}({}^{sc}\overline{T^*}X).
\end{align}
Accordingly, the normalized lower-order remainders tend to zero on the respective loci of microlocalization $\mathcal{Z}_\bullet \subset {}^\bullet\overline{T^*}X$, given by 
\begin{align}
    \mathcal{Z}_b&:=\{\varrho_b=0\}={}^bS^*X\label{eq:bLocus}\\
    \mathcal{Z}_{sc}&:=\{\varrho_{sc}=0\}\cup\{\rho=0\}=\{\varrho_{sc}\rho=0\}=\partial{}^{sc}\overline{T^*}X={}^{sc}S^*X \,\sqcup\, {}^{sc}\overline{T^*}_{\partial X}X\label{eq:ScatteringLocus}.
\end{align}

For $a_\bullet \in S^{m',l'}_\bullet ({}^\bullet T^*X)$, s.t. $A=\Op_\bullet(a_\bullet)\in \Psi^{m',l'}_\bullet(X)$ the relation between the Hamiltonian vector field of $p\in S^{m,l}_\bullet ({}^\bullet T^*X) $ with $P=\Op_\bullet(p_\bullet)\in \Psi^{m,l}_\bullet(X)$ can also be given through the commutator relations
\begin{align}\label{eq:CommutatorHamiltonians}
    {}^b\upsigma_{m+m'-1,l+l'}\left(i[P,A]\right)={}^bH_pa \text{ and } {}^{sc}\upsigma_{m+m'-1,l+l'-1}\left(i[P,A]\right)={}^{sc}H_pa.
\end{align}

We note that on compact sets of the interior $X^\circ$ the decay factors $\rho^\ell \asymp \rho^{\ell'}$ for any $\ell,\ell'$ hence both of these sequences only decay in the fibre-radial direction into account and reduce to the well known classical sequence

\begin{align}
    0\to \Psi^{m-1}(X^\circ)\to \Psi^{m}(X^\circ)  \overset{\sigma_{m,}}\to S^{m}(T^*X^\circ)/S^{m-1}(T^*X^\circ)\to 0\label{eq:ClassicalSES},
\end{align}
so that we have for $A\in \Psi^{m,\ell}_\bullet (X)$ that
\begin{align}
    {}^\bullet\sigma_{m,\ell}(A)|_{X^\circ}\sim \sigma^m(A)|_{X^\circ}\in S^m(T^*X^\circ) \text{ and } {}^\bullet H_a\in \mathcal{V}_\bullet ({}^\bullet\overline{T^*}X^\circ) \sim H_a\in \mathcal{V}(\overline{T^*}X^\circ),
\end{align}
given the fact that the Hamiltonian vector field expressions themselves are conjugate with respect to the same symbol, as observed earlier.
In contrast, the sequences (\ref{eq:bSES}) and (\ref{eq:ScatteringSES}) hold globally up to the boundary $\partial X$. When we restrict (\ref{eq:ClassicalSES}) to $\Diff^m(X^\circ)\subsetneq \Psi^m(X^\circ)$ for $m\in \mathbb{N}_0$ we get that the associated principal map maps into homogeneous symbol spaces
\begin{align}
    0\to \Diff^{m-1}(X^\circ)\to \Diff^m(X^\circ)\to S_{hom}(T^*X^\circ)/S_{hom}(T^*X^\circ)\to 0
\end{align}
where we denote $S^m_{hom}(T^*X^\circ)$ to be the space of homogenous symbols on $T^*X^\circ$, i.e. we have 
\begin{align}\label{eq:ClassicalHomogeneity}
    \upsigma_m(P)(x,\lambda\xi) = \lambda^m \upsigma_m(P)(x,\xi),
\end{align}
which is a fact that we used a lot in Section \ref{subsec:HamiltonianDynamics} to extend the Hamiltonian vector field of a differential operator to fibre infinity after rescaling. For the $b$- and $sc$- calculus on the other hand we have in general 
\begin{align}
    0\to\Diff_b^{m-1,l}(X) \to \Diff_b^{m,l}(X)  \overset{{}^{b}\sigma_{m,l}}{\to} \mathcal{P}_b^{m,l}({}^b T^*X)/\mathcal{P}_b^{m-1,l}({}^{b} T^*X) \to 0.
\end{align}
and 
\begin{align}
    0\to \text{Diff}_{sc}^{m-1,l-1}(X) \to \text{Diff}_{sc}^{m,l}(X)  \overset{{}^{sc}\sigma_{m,l}}{\to} \mathcal{P}_{sc}^{m,l}({}^{sc} T^*X)/\mathcal{P}_{sc}^{m-1,l-1}({}^{sc} T^*X) \to 0,
\end{align}
where $\mathcal{P}^{m,l}_\bullet({}^\bullet T^*X)$ denotes merely polynomial symbols in $b$ or $sc$ cotangent variables. In effect therefore the associated Hamiltonian flow of an operator $A\in \Psi_\bullet^{m,l}(X)$ with $m\in \mathbb{N}_0$ can therefore be interpreted as a classical flow at fibre infinity $S^*X^\circ$ on the interior $X^\circ$ due to extending by (\ref{eq:ClassicalHomogeneity}) and as a flow at possibly finite momenta at the boundary $\partial X$ as can be gleaned from the model operator (\ref{eq:MinkowskiModel}).
\begin{proposition}\label{prop:SymbolCalculuation}
    Consider the model operator $P(\sigma)=\hat{\Box}_{\mathbb{M}^{1+4}}(\sigma)$ from equation (\ref{eq:MinkowskiModel}). Then we have at $\partial X$ that
    \begin{align}
        {}^b\upsigma_{2,-1} (P(\sigma))|_{\partial X}= \rho^2(\bxi^2 + 
    \|\beeta_i\|_{g^{-1}_{\mathbb{S}^3}}^2) 
    \end{align}
    which is homogeneous but
    \begin{align}\label{eq:ScatteringMinkowskiSymbol}
        {}^{sc}\upsigma_{2,0}(P(\sigma))|_{\partial X} = \scxi^2+\|\sceta^2_i\|_{g^{-1}_{\mathbb{S}^3}}-2\sigma\scxi,
    \end{align}
    which is not.
    In the above we denote 
    \begin{align}
        \|{}^\bullet\eta\|_{g^{-1}_{\mathbb{S}^3}}={}^\bullet\eta_\theta^2+\frac{{}^\bullet\eta_\phi^2}{\sin^2\theta}+\frac{{}^\bullet\eta_\psi^2}{\cos^2\theta},
    \end{align}
    in Hopf coordinates.
\end{proposition}
\begin{proof}
    This comes immediately from rewriting and identifying symbols. It is convenient to quantize  (\ref{eq:MinkowskiModel}) as a $b$-operator through $D_\rho=i^{-1}\partial_\rho$
    \begin{align}\label{eq:MinkowskiScatterOperator}
        \hat{\Box}_\mathbb{M}^{1+4}(\sigma)=\rho^2((\rho D_\rho)^2+2i\rho D_\rho +\Delta_{\mathbb{S}^3})-2i\sigma\rho\left(\rho D_\rho +\frac{3}{2}\right),
    \end{align}
    giving principal symbol (using $2i\rho D_\rho,2i\sigma\rho\left(\rho D_\rho +\frac{3}{2}\right) \in \Diff^{1,-1}_b(X) $, 
    \begin{align}
        {}^b\upsigma_{(2,-1)}(\hat{\Box}(\sigma))=\rho^2(\bxi^2 + 
    \|\beeta_i\|_{g^{-1}_{\mathbb{S}^3}}^2). 
    \end{align}
    To write (\ref{eq:MinkowskiModel}) in scattering form note that 
    \begin{align}
        \hat{\Box}_{\mathbb{M}^{1+4}(\sigma)}=-(\rho^2\partial_\rho)^2 + \rho^2\Delta_{\mathbb{S}^3}+2i\sigma (\rho^2\partial_\rho)-3\rho(\rho^2\partial_\rho-i\sigma)
    \end{align}
    or 
    \begin{align}
        \hat{\Box}_{\mathbb{M}^{1+4}(\sigma)}=(\rho^2D_\rho)^2 + \rho^2\Delta_{\mathbb{S}^3}-
        2\sigma (\rho^2D_\rho)+3i\rho(\rho^2D_\rho-\sigma).
    \end{align}
    We note that ${}^{sc}\upsigma_{(2,0)} (2\sigma \rho^2D_\rho)=2\sigma\scxi\neq 0 $ and hence we find 
    \begin{align}
        {}^{sc}\upsigma_{(2,0)}(P(\sigma))=\scxi^2+\|\sceta_i\|^2_{g^{-1}_{\mathbb{S}^3}} - 2\sigma\scxi.
    \end{align}
\end{proof}

Finally, to study the dynamics of pseudo-differential operators $A\in \Psi_\bullet^{m,l}(X)$ and convert these dynamics into corresponding estimateswe need to define their operator wavefront sets $\text{WF}_\bullet'(A)$ and characteristic/elliptic sets $\chr_\bullet(A)$ and $\Ell_\bullet(A)$. Assume that 
\begin{align}
    A = \text{Op}_\bullet (A_\bullet) \text{ with } a_\bullet \in S^{m,l}_\bullet ({}^\bullet \overline{T}^* X). 
\end{align}

with corresponding renormalized symbol,

\begin{align}
    \tilde{a}_\bullet(x,{}^{\bullet}\xi)= \varrho^{m}_\bullet\rho^{l} \cdot a_\bullet(x,{}^{\bullet}\xi).
\end{align}\footnote{We implicitly identify the base defining function $\rho$ with its pullback ${}^\bullet\overline{T^*}X$ under the base projection $\pi_\bullet:{}^\bullet\overline{T^*}\to X$. Thus at a phase-point $q$ the factor $\rho$ means $\rho \pi_\bullet(q)$ whereas $\varrho_\bullet$ means $\varrho_\bullet(q)$. We surpress the pullback notation and arguments of these boundary defining functions throughout.}
For both calculi we then define $\Ell_\bullet(A)$, $\chr_\bullet(A)$ and $\WF'_\bullet(A)$ on the respective loci of microlocalization $\mathcal{Z}_\bullet\subset \overline{T^*}X$, defined in (\ref{eq:bLocus}) and (\ref{eq:ScatteringLocus}).  
Starting with the the scattering elliptic set $\Ell^{m,l}_{sc}(A)\subset {}^{sc}\overline{T^*}X$ as 
\begin{align}
    \Ell^{m,l}_{sc}(A):=\left\{q\in \mathcal{Z}_{sc}\;|\;\exists c>0, U\ni q \text{ open in } {}^{sc}\overline{T^*}X, \text{ s.t. } |\tilde a_{sc}|\geq c \text{ on } U\cap {}^{sc}T^*X^\circ \right\}.
\end{align}
and we also define the characteristic set as the complement 
\begin{align}
    \Char^{m,l}_{sc}(A)=\mathcal{Z}_{sc}\backslash  \Ell^{m,l}_{sc}(A).
\end{align}

We note again that these sets can be interpreted as lying over ${}^{sc}S^*X^\circ\cup {}^{sc}\overline{T^*}_{\partial X}X$, i.e. at fibre infinity over the interior $X^\circ$ whilst possible at finite momenta at the boundary $\partial X$ where decay is weighted for the computations of the principal symbol ${}^{sc}\upsigma_{m,l} (A)$. These two sets can meet in the corner $\partial^2 {}^{sc}\overline{T^*}X$.

For an operator $A$ we also define the complement of the essential support of its symbol as 
\begin{align}
    \mathcal{Z}_{sc}\backslash\text{esssup}_{sc}(a) &:=\left\{q\in {}^{sc}\overline{T^*}X\;|\; \right.\\
    &\;\;\;\;\;\;\;\;\;\left.  \exists \chi\in C^\infty({}^{sc}T^*X), \chi\equiv 1 \text{ near } (x,{}^{sc}\xi), \text{ for which } \chi a_{sc}\in S^{-\infty,-\infty}_{sc} \right\}
\end{align}

for which we define the operator wavefront set as 
\begin{align}
    \text{WF}_{sc}'(A):=\text{esssup}_{sc}(a)\subset \mathcal{Z}_{sc}.
\end{align}

For a $b$ operator $A=\text{Op}_b(a)\in \Psi_{b}^{m,l}(X)$ we define almost by analogy, but we do not quotient out with respect to base decay order based on (\ref{eq:bLocus}). We emphasize that even at $\partial_+X$, we have that $\mathcal{Z}_b|_{\partial_+X}={}^bS^*_{\partial X}X$. The elliptic and characteristic sets are then defined as
\begin{align}
    \Ell_b^{m,l}(A)&:=\{q\in \mathcal{Z}_b\; | \;\exists c>0, U\ni q \text{ open in } {}^{b}\overline{T^*}X, \text{ s.t. } |\tilde a_b|\geq c \text{ on } U\cap {}^bT^*X^\circ \} \\\text{ and  } &\Char_b(A)=\mathcal{Z}_b \backslash \Ell_b^{m,l}(A),
\end{align}

We define the operator wavefront set for $b$-calcui in almost the same way as $sc$-operator, but with no smoothing requirement in the spatial decay. I.e. for $a\in S_b^{m,\ell}({}^bT^*X)$
\begin{align}
    \mathcal{Z}_b\backslash\text{esssup}_{b}(a)&:=\left\{q\in {}^{b}\overline{T^*}X\;|\; \right.\\
    &\;\;\;\;\;\;\;\;\;\left.  \exists \chi\in C^\infty({}^{b}T^*X), \chi\equiv 1 \text{ near } (x,{}^{b}\xi), \text{ for which } \chi a_b\in S^{-\infty,l}_{b} \right\}
\end{align}
such that 
\begin{align}
    \WF'_b(A) :=\esssup_b{(a)}\subset \mathcal{Z}_b.
\end{align}

\begin{remark}\label{rem:ContinousExtenstions}
    A few remarks are in order with respect to the above definitions. First of all, intersections $\text{Ell}_\bullet^{m,l}(X)\cap X^\circ\subset {}^\bullet S^*X^\circ$ can be interpreted as conic subsets of ${}^\bullet \overline{T^*}X$ by taking the inverse image under fibre radial projections 
    \begin{align}
        \pi^{-1}_{rad,\bullet}:\overline{{}^\bullet T^*}X\backslash \mathcal{O}\to {}^\bullet S^* X.
    \end{align}
    then we can give the conic counterparts to the ellitpic, resp. characteristic sets as 
    \begin{align}\pi_{rad,\bullet}^{-1}(\Ell^{m,l}_\bullet (X) \cap {}^\bullet S^* X^\circ)\subset {}^\bullet \overline{T^*}X^\circ
    \end{align}
    
    and
    
     \begin{align}\pi_{rad,\bullet}^{-1}(\chr^{m,l}_\bullet (X) \cap {}^\bullet S^* X^\circ)\subset {}^\bullet \overline{T^*}X^\circ.
    \end{align} 
    Moreover, when the amplitudes $\tilde{a}(\rho,y,{}^\bullet\xi,{}^\bullet\eta)$ extend continuously to the boundary $\partial{}^\bullet \overline{T^*}X$, the expressions of the characteristic and elliptic sets simplify, i.e. we have  
    \begin{align}
         \Ell_\bullet^{m,l}(A)&=\{q\in \mathcal{Z}_\bullet \; | \; \tilde a_\bullet(q) \neq 0 \}  \text{ and } \chr_\bullet^{m,l}(A)=\{q\in \mathcal{Z}_\bullet \; | \; \tilde a_\bullet(q) =0 \}.
    \end{align}
    An example of a normalized symbol which does not extend continuously to $\partial_+X$ is 
    \begin{align}
        a(\rho)=2+\sin(\log \rho),
    \end{align}
    whose $b$-quantization belongs to $\Psi^{0,0}_b(X)$, and it is in fact elliptic on all of ${}^bS^*X$.
\end{remark}

We finally note that the distributional wavefront sets are for a distribution $u\in \mathscr{D}'(X)$ as 
\begin{align}
    \WF^{s,\ell}_\bullet(u)^c &=\{q\in {}^\bullet\overline{T^*}X\;|\; \exists A\in \Psi^{0,0}_\bullet(x), q\in \Ell^{0,0}_{\bullet}(A), Au\in H^{s,\ell}_\bullet(X)\}\\
    &=\bigcap_{A\in \Psi^{0,0}_\bullet, Au \in H^{s,\ell}_\bullet } \chr_{\bullet}^{0,0}(A).
\end{align}

Indeed we could have defined our Sobelev spaces through
\begin{align}
    H^{s,\ell}_\bullet (X):=\{u\in \mathscr{D}'(X)\;|\;\exists A \in \Psi^{s,\ell}_\bullet(X) \text{ elliptic s.t. } Au\in L^2_\bullet(X)\} 
\end{align}
with norm 
\begin{align}
    \|u\|_{H^{s,r}_\bullet(X)}=\|\Lambda ^s(\rho^{-r}u)\|_{L^2_\bullet(X)},
\end{align}
with $\Lambda^s = \text{Op}(\langle \xi\rangle^s)\in \Psi^{s,r}_\bullet(X)$\footnote{One could also take the graph norm for any elliptic $A\in \Psi^{m,\ell}_\bullet(X)$.}. Indeed, we can and will use this definition to generalize to variable order Sobolev spaces $H^{\mathfrak{s},\mathfrak{r}}_\bullet(X)$ where $\mathfrak{s},\tilde{\mathfrak{\ell}}\in C^\infty({}^{\bullet}T^*X)$ are allowed to be smooth functions. Analogously one can also define pseudodifferential operators $\Psi^{\mathfrak{m},\mathfrak{l}}_\bullet(X)$ for $\mathfrak{m},\mathfrak{l}\in C^\infty({}^{\bullet}T^*X)$.

With these constructs we can do several estimates. First of all, we have the elliptic estimate.

\begin{proposition} (Elliptic estimate)
    Let $P\in \Psi^{m,l}_\bullet(X)$, $K\subset X$ compact. Suppose $B,G\in \Psi^{0,0}_\bullet(X)$ are such that 
    \begin{align}
        \WF_\bullet'(B)\subset \Ell^{0,0}_\bullet(G)\cap\Ell^{m,l}_\bullet(P).
    \end{align}
    Let $K\subset X$ compact and assume the Schwarz kernels ${}^\bullet K_{B}$ and ${}^\bullet K_{G}$ are supported in $K\times K\subset X\times X$ and let $\tilde K$ be a compact set containing $K$ such that $\supp\, {}^\bullet K_P \cap K\times X \subset K\times\tilde K$. Let $\chi \in C_c^\infty(X)$ be identically $1$ near $\tilde K$, then we have for all $u\in \mathscr{D}'(X)$ 
    \begin{align}\label{eq:bEllipticEstimate}
        \|Bu\|_{H^{s,\ell}_{b}}\leq C(\|GP u\|_{H_b^{s-m,\ell-l}} +\|\chi u\|_{H_b^{-N,l}}),
    \end{align}
    if $\bullet = b$, and
    \begin{align}\label{eq:scEllipticEstimate}
                \|Bu\|_{H^{s,r}_{sc}}\leq C(\|GP u\|_{H_{sc}^{s-m,\ell-l}} +\|\chi u\|_{H_{sc}^{-N,-N}}).
    \end{align}
    if $\bullet = sc$.
\end{proposition}

\begin{remark}
    If $X$ is non-compact, $\|\cdot\|_{H^{s,\ell}_\bullet}$ should be interpreted in the context of seminorms mentioned before, i.e. as defining a norm on $\dot H^{s,\ell}(K')$ where $K'$ is the compact set $K':=\supp \,\chi$. 
\end{remark}

We can glue elliptic estimates together using microlocal partitions of unity.
\begin{proposition}
    (Microlocal Partition of Unity) Let $K\subset {}^\bullet S^*X$ be compact. Let $\{U\}_{1\leq j \leq N}$ be compact. Then there exist compactly suppported operators $\{A_j\}_{1\leq j \leq N}\in \Psi^{0,0}_\bullet (X)$ such that $\WF'(A_j)\subset U_j$ for all $j$ and $\WF'(Id - \sum_j A_j)\cap K =\emptyset$.
\end{proposition}

The key feature of any microlocal approach to get global estimates on $X$ is to supplement the elliptic estimates of an operator $P\in \Psi^{m,l}_\bullet (X)$ on $\Ell_\bullet(P)$ with propagation of regularity on the characteristic set $\Char(P)$ along the Hamiltonian vector field ${}^\bullet H_p$. Characteristic curves $\gamma:I \to {}^\bullet T^*X$ satisfying 
\begin{align}
    \gamma' = {}^\bullet H_p \gamma 
\end{align}
are referred to as null-bicharacteristics.
For our purposes it will be sufficient to assume that $P\in \Psi_\bullet^{m,\ell}(X)$ admits a real representative $p$ such that on a neighbourhood $U\subset {}^\bullet T^*X$ under consideration,
\begin{align}
    \tilde{p} = \varrho_\bullet^m\rho^l p \in C^\infty (U), 
\end{align}
and to consider the order rescaled vector fields 
\begin{align}
    {}^b\hat{H}_p :=\varrho_b^{m-1}\rho^l H_p \text{ and } {}^{sc}\hat{H}_p :=\varrho_{sc}^{m-1}\rho^{l-1}H_p.\label{eq:RescaledbscHfields}
\end{align}
Under these assumptions, the vector fields (\ref{eq:RescaledbscHfields}) extend smoothly to the boundary, and are tangent to it i.e. 
\begin{align}
    {}^\bullet\hat{H}_p \in \mathcal{V}_b({}^\bullet \overline{T^*}X).
\end{align}

In this setting we thus study 
\begin{align}
    \gamma'(t) = {}^\bullet \hat{H}_p(\gamma(t)),\; \gamma(t)\in \chr_\bullet (P).
\end{align}

\begin{proposition}\label{prop:bscPropagation}
    Suppose $P\in \Psi^{m,l}_\bullet(X)$ is properly supported and its principal-symbol class admits a real-valued representative $p$ whose normalized symbol $\tilde{p}=\varrho_\bullet^m\rho^\ell p$ is smooth in a neighbourhood of the bicharacteristic segment under consideration. Let $\gamma:[0,s_0]\to \chr_\bullet(P)\subset \mathcal{Z}_\bullet$ be a null-bicharacteristic of the rescaled vector field ${}^\bullet \hat{H}$, with ${}^\bullet H_p(\gamma(t))\neq 0, \forall0\leq t \leq s_0$. Let $u\in \mathscr{D}'(X)$, and $K\subset \tilde K\subset X$ be compact subset so that $\pi_\bullet(\gamma([0,s_0])\subset K^\circ \cap X^\circ$ and $\supp ({}^\bullet K_P)\cap (K\times M) \subset K\times \tilde K$. Let $\chi \in C^\infty_c(X)$ satisfy $\chi\equiv 1$ near $\tilde K$. Let $U_0,U\subset {}^\bullet \overline{T^*}X$ be open neighborhoods of $\gamma(0)$ and $\gamma([0,s_0])$ respectively. Then there exist $B,G,E\in \Psi^{0,0}_\bullet(X)$ with 
    \begin{enumerate}
        \item $\gamma([0,s_0]) \subset\Ell_\bullet(B)$,\\
        \item $\WF'(B)\subset \Ell_\bullet(G)$, $\text{WF}'(G) \subset U$\\
        \item $\WF'(E)\subset U_0$ and $\gamma(0)\in \Ell_\bullet(E)$,
    \end{enumerate}
    such that for arbitrary $N\in \mathbb{R}$ and a constant $C>0$:
    \begin{align}
        \|Bu\|_{H^{s,\ell}_b}\leq C(\|Eu\|_{H^{s,\ell}_b}+\|GPu\|_{H^{s-m+1,\ell-l}_b}+\|\chi u\|_{H_b^{-N,l}}).
    \end{align}
if $\bullet = b$, and 

\begin{align}
    \|Bu\|_{H^{s,\ell}_{sc}}\leq C(\|Eu\|_{H^{s,\ell}_{sc}}+\|GPu\|_{H^{s-m+1,\ell-l+1}_{sc}}+\|\chi u\|_{H^{-N,-N}_{sc}}).
\end{align}
\end{proposition}

The primary operator $P(\sigma)$ of interest is the one given as 
\begin{align}\label{eq:SpectralFamily2}
    P(\sigma)=\hat{\Box}_g(\sigma) := e^{i\sigma t_*} \circ\Box_g \circ e^{-i\sigma t_*}:H^{s,r}_{sc}(X)\to H^{s-2,r}_{sc}(X)
\end{align}
For fixed $\sigma \in \mathbb{C}\backslash \{0\}$ with $\text{Im}(\sigma)$ we regard the spectral parameter on as a constant. If we were to compute the standard symbol of $P=\Box_{g_{\bfa,\bfb,\mfrm}}$ as a standard 2nd order operator on the full spacetime manifold $M$ we would end up with an interior symbol 
\begin{align}
    \upsigma_2(P) = g^{-1}(\upalpha,\upalpha ) \in S^2(T^*M^\circ).
\end{align}
To get from the null geodesic flow to the classical flow the fundamental identification on $X^\circ$ is
\begin{align}\label{eq:FundamentalClassicalProjection}
\pi_0:(\overline{T^*}M \cap \overline{\{\sigma = 0\}})\to (\overline{T^*}M \cap \overline{\{\sigma = 0\}})/\mathbb{R}_{t_*}\cong T^*X^\circ.
\end{align}
Stationarity gives us $H_G\sigma =0$, and 
\begin{align}
    H_G|_{\sigma=0}=\partial_\sigma G|_{\sigma=0}\partial_{t_*} +H_{G_0},
\end{align}
where 
\begin{align}
    G_0:=G(\sigma dt_*+\upalpha_X)|_{\sigma =0}.
\end{align}

Hence after quotienting by time translations the vector field simply restricts to $H_{G_0}$ on $T^*X^\circ$. This translates identically to the dynamics one obtains when conjugates $\Box_g$ by $e^{i\sigma t_*}$ and freezes $\sigma \in \mathbb{C}$ as a constant, we get with a slight abuse of notation\footnote{Note that we denote both the $t_*$ conjugate momentum in part I and the spectral parameter with $\sigma$. Writing the latter as $\sigma_{spec}$ our convention $\sigma=-\sigma_{spec}$ over the reals.}
\begin{align}
    {}^{sc}\upsigma_{2,0}(P(\sigma)) = g^{-1}_{\bfa,\bfb,\mfrm}(-\sigma dt_* +\upalpha_X) |_{\sigma =0}=\tilde{G}_0, 
\end{align}
where we interpret $\tilde{G}_0$ simply as a quadratic function on $T^*X^\circ$, with slight abuse of notation. We will make this correspondence even more explicit when treating semiclassical analysis in Section \ref{subsec:SemiclassicalAnalysis}.

Thus analyzing the $b$- or $sc$- flow on $X^\circ$ is equivalent to restricting oneselves to the hypersurface $\{\sigma = 0,t_*=0\}$ inside $T^*M^\circ$. The characteristic set of $P$ would therefore be the intersection of the future/backwards directed lightcone intersected with fibre infinity after quotienting out fibre radial rays and time translation. As it turns out its projection onto the base is exactly the ergo-region. The following proof is analogous to Kerr (see e.g. also \cite{millet2023optimal}).

\begin{proposition}\label{eq:ClassicalCharacteristicSet}
    Let $P(\sigma )=\bfr^2(e^{i\sigma t_*} \Box_g(\sigma) e^{-i\sigma t_*})$. Then for $\chr(P(\sigma)=\chr (G) \subset S^* X^\circ$ we have that 
    \begin{align}
        \pi_{S^*X^\circ}(\chr(G)) = \{ r_X\leq r \leq r_{\mathcal{E}}\},
    \end{align}
    where $\pi_{S^*X^\circ}:S^*X^\circ \to X^\circ$ is the projection onto the base.
    Moreover $G|_{\sigma=0}>0$ on $\{\bfr>\sqrt{2\mfrm}+\epsilon\} \cap \mathring{K}$, where $K\subset \tilde X$ is any compact subset and is hence elliptic there.
\end{proposition}

\begin{proof}
    As noted before the $b$ and $sc$ classical symbols can be taken to be the principal symbol $\upsigma_2(P)$ on the interior $X^\circ$. 
  
    We note that in starlike coordinates, $G|_{\sigma =0}$ is given by
    \begin{align}
    G|_{\sigma=0}&=\bfDelta \xi^2+\eta_\theta^2   +2 \left[\frac{\bfa(r^2+\bfb^2)}{r^2}  \eta_{\phi} + \frac{\bfb(r^2+\bfa^2)}{r^2} \eta_\psi\right]\xi \nonumber \\
    &+ \frac{\sin^2\theta\cos^2\theta}{\bfp^2}\left[\frac{\bfa}{\sin^2\theta} \eta_\phi -\frac{\bfb}{\cos^2\theta}\eta_\psi  \right]^2 + \frac{\bfr^2}{r^2 \bfp^2}\left[\bfb \eta_\phi + \bfa \eta_\psi\right]^2\\
    &=\bfDelta \xi^2+\eta_\theta^2   +2 \left[\frac{\bfa(r^2+\bfb^2)}{r^2}  \eta_{\phi} + \frac{\bfb(r^2+\bfa^2)}{r^2} \eta_\psi\right]\xi \nonumber \\
    &+ \frac{1}{\bfp^2}\left[\bfa\cot\theta\, \eta_\phi -\bfb\tan\theta\,\eta_\psi  \right]^2 + \left(\frac{1}{r^2}+\frac{1}{\bfp^2}\right)\left[\bfb \eta_\phi + \bfa \eta_\psi\right]^2.
\end{align}
We can simplify the above even further by noting that 
\begin{align}
    \frac{1}{\bfp^2}\left(\left[\bfa\cot\theta\eta_\phi-\bfb\tan\theta\eta_\psi\right]^2 + [\bfb\eta_\phi+\bfa\eta_\psi]^2\right)=\frac{1}{\bfp^2}\left(\left[\bfa^2\cot^2\theta+\bfb^2\right]\eta_\phi^2+\left[\bfb^2\tan^2\theta+\bfa^2\right]\eta_\psi^2\right),
\end{align}
which by multiplying out above and below by $\sin^2\theta\cos^2\theta$ we see that 

\begin{align}
   \frac{1}{\bfp^2}\left(\left[\bfa^2\cot^2\theta+\bfb^2\right]\eta_\phi^2+\left[\bfb^2\tan^2\theta+\bfa^2\right]\eta_\psi^2\right)&=\frac{1}{\bfp^2\cos^2\theta\sin^2\theta}\left(\bfp^2\eta_\phi^2\cos^2\theta+\bfp^2\sin^2\theta \eta_\psi^2\right)\\
   &=\frac{\eta_\phi^2}{\sin^2\theta}+\frac{\eta_\psi^2}{\cos^2\theta}.
\end{align}
From this we can simplify $G|_{\sigma=0}$ even further:
\begin{align}
    G|_{\sigma=0}&=\bfDelta \xi^2+\eta_\theta^2   +2 \left[\frac{\bfa(r^2+\bfb^2)}{r^2}  \eta_{\phi} + \frac{\bfb(r^2+\bfa^2)}{r^2} \eta_\psi\right]\xi \nonumber \\
    &+ \frac{\eta_\phi^2}{\sin^2\theta}+\frac{\eta_\psi^2}{\cos^2\theta}+\frac{1}{r^2}\left[\bfb\eta_\phi+\bfa\eta_\psi\right]^2
\end{align}

This can then ultimately be factorized,
\begin{align}\label{eq:ClassicalFactorization1}
    G|_{\sigma=0}&=(\bfr^2-2\mfrm)\xi^2+ \left(\bfa\sin\theta\xi +\frac{\eta_\phi}{\sin\theta}\right)^2 +\left(\bfb\cos\theta\xi +\frac{\eta_\psi}{\cos\theta}\right)^2 \nonumber\\
    &+\frac{1}{r^2}(\bfa\bfb\xi +\bfb\eta_\phi+\bfa\eta_\psi)^2 +\eta_\theta^2.
\end{align}
which is a valid expression away from the poles. At the the pole $\theta=0$ we again set 
\begin{align}
    y_1=\cos\phi\sin\theta\;\text{ and } y_2=\sin\phi\sin\theta,
\end{align}
then we get 
\begin{align}
    \eta_\theta^2+\frac{\eta_\phi^2}{\sin^2\theta}\mapsto     (1-y^2_1)\eta_1^2-2y_1y_2\eta_1\eta_2+(1-y_2^2)\eta_2^2
\end{align}
and hence 
\begin{align}
    G|_{\sigma=0}&=(\bfr^2-2\mfrm)\xi^2+ 2\bfa\xi\underbrace{(-y_2\eta_1+y_1\eta_2)}_{=0} +\left(\bfb\cos\theta\xi +\frac{\eta_\psi}{\cos\theta}\right)^2 \nonumber\\
    &+\frac{1}{r^2}(\bfa\bfb\xi +\bfb\underbrace{(-y_2\eta_1+y_1\eta_2)}_{=0}+\bfa\eta_\psi) +\eta_1^2+\eta_2^2\nonumber\\
    &=(\bfr^2-2\mfrm)\xi^2+ \left(\bfb\xi +\eta_\psi\right)^2 +\frac{1}{r^2}(\bfa\bfb\xi +\bfa\eta_\psi)^2 +\eta_1^2+\eta_2^2.\label{eq:ClassicalFactorization2}
\end{align}
Exactly analogously near $\theta=\frac{\pi}{2}$ we find with
\begin{align}
    y_3&=\cos\psi\cos\theta,\; \text{ and } y_4=\sin\psi\cos\theta,
\end{align}
that 
\begin{align}\label{eq:ClassicalFactorization3}
    G|_{\sigma=0}&=(\bfr^2-2\mfrm)\xi^2+ \left(\bfa\xi +\eta_\phi\right)^2 +\frac{1}{r^2}(\bfa\bfb\xi +\bfb\eta_\phi)^2 +\eta_3^2+\eta_4^2.
\end{align}
The second claim thus follows immediately for $\bfr^2>2\mfrm$ as we are only summing positive squares. For $\bfr^2<2\mfrm$ we can construct a characteristic direction, i.e. in (\ref{eq:ClassicalFactorization1}) we can take $\xi =1, \eta_\phi = -\bfa\sin^2\theta, \eta_\psi=-\bfb\cos^2\theta$ and $\eta_\theta = \sqrt{2\mfrm -\bfr^2}$. Similar choices can be made in the other charts for (\ref{eq:ClassicalFactorization2}) and (\ref{eq:ClassicalFactorization3}).
\end{proof}

Apart from that, in the scattering calculus we still have the characteristic set at $\partial X$.  In this case the set $\chr (P(\sigma))$ contains finite scattering covectors $({}^{sc}\xi,\sceta)$. We have that 
\begin{align}
   p_\sigma(\rho,y_i,{}^{sc}\xi,{}^{sc}\eta)|_{\partial X}&= {}^{sc}\upsigma_{2,0}(P(\sigma))|_{\partial X}\\ 
   &= {}^{sc}\upsigma_{2,0}({}^{sc}N (\hat{\Box}_{g_{\bfa,\bfb,\mfrm}}(\sigma)))|_{\partial X}\\
   &={}^{sc}\upsigma_{2,0}(\hat{\Box}_{g_{\mathbb{M}^{1+4}}}(\sigma))|_{\partial X}\\
   &= {}^{sc}\xi^2+\|{}^{sc}\eta_i\|_{g^{-1}_{\mathbb{S}^3}}-2\sigma {}^{sc}\xi.
\end{align}

Thus, we have for finite $\sigma\in \mathbb{C}\backslash \{0\}$ with $\text{Im }\sigma\geq 0$, that 
\begin{align}
    \Char_{sc}(P(\sigma)) \cap {}^{sc}\overline{T^*}_{\partial X}X = \left\{\rho=0,\;\|{}^{sc}\eta\|_{g^{-1}_{\mathbb{S}^3}}\leq |\sigma|, {}^{sc}\xi=\sigma\pm \sqrt{\sigma^2-\|{}^{sc}\eta\|^2_{g^{-1}_{\mathbb{S}^3}}} \right\}.
\end{align}

It then remains to characterize the radial points of the Hamiltonian flow. Restricting ourselves for now on the scattering calculus they come in two flavors: on the interior $X^\circ$ radial points are characterized exactly at fibre infinity. For this we can exactly use Definition \ref{def:GeneralizedRadialSet} once we restrict ourselves to (\ref{eq:FundamentalClassicalProjection}). More precisely we have the following.
\begin{definition}\label{def:ClassicalRadialSet}
   A smooth set 
   \begin{align}
       \Lambda_X\subset \{q\in S^*X\;|\;(\varrho^2G_0)(q)=0\}
   \end{align}
   is a smooth source or sink if it satisfies Definition \ref{def:GeneralizedRadialSet} for the homogeneous Hamiltonian $\varrho H_{G_0}$ on $\overline{T^*}X^\circ $.
\end{definition}

We see then that on $T^*X^\circ$ these sets are given as
\begin{align}
    \partial\mathcal{R}^{\pm}_{\mathcal{H}^+,X^\circ}:=\pi_0(\mathcal{R}^\pm_{\mathcal{H}^+})=\{r=r_+,\phi,\psi,\hat\eta_\theta=\hat\eta_\phi=\hat \eta_\psi=0,\hat\xi=0\}\cap \Sigma_{G_0}.  
\end{align}

Assuming that the symbol $p=\upsigma_m(P)\in S^m(T^*X^\circ)$ is real-valued, we have that $\upsigma_m (P-P^*)=0$ i.e. $P-P^*\in \Psi^{m-1}(X)$. The threshold regularity of the radial set $\Lambda$ in Definition \ref{def:GeneralizedRadialSet} then depends on the subprincipal symbol, defined through a function $\alpha_\Lambda \in C^\infty(\Lambda)$ such that 
\begin{align}
    \varrho^{m-1}\left.\upsigma_{m-1}\left(\frac{1}{2i}(P(\sigma)-P^*(\sigma)\right)\right|_{\Lambda}=\pm \alpha_\Lambda \alpha_\infty
\end{align}\footnote{Note that here $\alpha_\infty = \frac{4\mfrm\kappa}{\bfr^2|_{r=r_+}}$ in contrast to equation (\ref{eq:InftyThresholdCoefficient}) as we are looking at the non-rescaled operator.}

We phrase the following ingoing and outgoing estimates on classical semi-norms as decay in the differential parameter is all that is relevant for interior propagation. A good reference for the statement and the subsequent threshold calculation of $\alpha_\Lambda$ is \cite{stucker2026quasinormal}, Propositions 2.15-2.16 and Section 3.5.

\begin{proposition}\label{prop:ClassicalHorizonEstimate} (Horizon radial point estimate on the interior)
    Let $P\in \Psi^m(X^\circ)$ be properly supported and have real valued homogeneous symbol with extendible renormalized symbol $\tilde{p}=\varrho^{m}p$. Let $\upsigma_{m-1}(P-P^*)$ be homogeneous as well. Let $\Lambda \subset \overline{T^*}X$ be a radial source or radial sink of $P$. Let $G\in \Psi^0(X^\circ)$ be compactly supported with $\Lambda\subset \Ell(G)$. If $s'$ is of high regularity, that is
    \begin{align}
        s\geq s'> \frac{m-1}{2}+\alpha_\Lambda,
    \end{align}
    then for any $s,N\in \mathbb{R}$ with $s\geq s'$ there exist $C>0$, $\chi\in C_c^\infty(X^\circ)$ and $B\in \Psi^0(X^\circ)$ compactly supported with $\Lambda \subset \Ell(B)$ such that for all $u\in H^{-N}_{loc}(X^\circ)$ with $Gu \in H^{s'}_c(X^\circ)$ and $GPu \in H_c^{s-m+1}(X^\circ)$ we have $Bu \in H_c^s(X^\circ)$ and 
    \begin{align}
        \|Bu\|_{H^s}\leq C(\|GPu\|_{H^{s-m+1}}+\|\chi u\|_{H^{-N}}).
    \end{align}
    Conversely if $s'$ is of low regularity, that is if
    \begin{align}
        s< \frac{m-1}{2}+\alpha_\Lambda,
    \end{align}
    there exists $C>0$ and $\chi\in C_c^\infty(X^\circ)$ and additionally an $E\in \Psi^0(X^\circ)$, such that for all $u\in H^{-N}_{loc}(X^\circ)$ with $GPu \in H^{s-m+1}_c(X^\circ)$, $Eu \in H^s_c(X^\circ)$ we have that $Bu\in H_c^s(X^\circ)$ and the following estimate holds,
    \begin{align}
        \|Bu\|_{H^s}\leq C(\|GPu\|_{H^{s-m+1}}+\|Eu\|_{H^s}+\|\chi u\|_{H^{-N}}).
    \end{align}
\end{proposition}

The statement for $P(\sigma)=\hat{\Box}_{g_{\bfa,\bfb,\mfrm}}(\sigma)\in \Psi^2(X^\circ)$ holds in particular for the case $m=2$ as its interior principal symbol $\upsigma_2(P(\sigma))=G|_{\sigma=0}$ does not depend on $\sigma$ and is hence real-valued. Moreover, in the case of the scalar wave equation, it remains for us to determine $\alpha_\Lambda$. For $P(\sigma)= \hat{\Box}_{g_{\bfa,\bfb,\mfrm}}(\sigma)$ let us define the rescaled operator 
\begin{align}
    L(\sigma) = \bfr^2 P(\sigma),
\end{align}
such that we can associated flow of the principal symbol can be directly associated with the flow of the rescaled dual metric function $G|_{\sigma=0} =\bfr^2 g^{-1}(0dt_*+\upalpha_X)$. When we taking are taking formal adjoints with respect to the stationary density $dg_{\bfa,\bfb,\mfrm}|_X=r \bfr^2 dr d\cancel{g}_{\mathbb{S}^3}$, defined in (\ref{eq:StationarySplitting}) we see that $P(\sigma)^*=P(\overline{\sigma})$, and thus
\begin{align}
    L(\sigma)-L(\sigma)^*= \bfr^2(P(\sigma)-P(\overline{\sigma})) + [\bfr^2,P(\overline{\sigma})].
\end{align}
Now recalling that 
\begin{align}
    \mathcal{R}^\pm_{\mathcal{H^+}}=\{r=r_+, \mp\xi>0\}, 
\end{align}
is invariant under the flow, the $r$ and $\theta$ components vanish. The angular drift therefore does not cause an additional term and we have 
\begin{align}
\upsigma_1([\bfr^2,P(\overline{\sigma})])|_{\mathcal{R}^\pm_{\mathcal{H}^+}}=0
\end{align}

Now let us write
\begin{align}
    L(\sigma) = L(0) + \sigma L_1+ \sigma^2L_2, \text{ with } L_1\in \Diff^1(X^\circ), L_2\in C^\infty(X^\circ).
\end{align}
Differentiating the spectral conjugation gives 
\begin{align}
    L_1=-i[\bfr^2\Box_{g_{\bfa,\bfb,\mfrm}},t_*],
\end{align}
from which follows that 
\begin{align}
    \upsigma_1(L_1) = -H_Gt_*
\end{align}

and hence 
\begin{align}
    \upsigma_1\left(\frac{L(\sigma)-L(\sigma)^*}{2i}\right)= -\text{Im}\,\sigma H_G t_*. 
\end{align}
We already computed $\hat{H}_Gt_*$ in (\ref{eq:HGTexpression}) and hence $H_Gt_* = \xi \hat{H}_Gt_*=4\mfrm\xi$. 

And thus that 
\begin{align}
    \left.\upsigma_1 \left(\frac{L(\sigma)-L(\sigma)^*}{2i}\right)\right|_{\mathcal{R^\pm_{\mathcal{H}}}}=-4\mfrm \text{Im} (\sigma)\, \xi.
\end{align}
which is a source for $\xi>0$ and a sink for $\xi<0$.
Using that $\alpha_\infty =\Delta'_r|_{r=r_+}=4\mfrm \kappa$ for the rescaled flow and $|\xi| =\varrho_\infty^{-1}$, we find that $\alpha_{\Lambda}$ depends on the imaginary part of $\sigma$
\begin{align}\label{eq:ThresholdSurfacePressure}
    \alpha_\Lambda(\sigma)=-\beta\text{ Im }\sigma, \text{ where } \beta = \inf_\mathcal{R^{\pm}_{\mathcal{H}}} \frac{\hat{H}_G t_*}{\alpha_\infty}=\frac{4\mfrm}{\Delta_r'|_{r=r_+}}=\frac{1}{\kappa}.
\end{align}
where we recall the surface gravity for Myers-Perry from equation (\ref{eq:SurfaceGravityIntro}) and $\alpha_\infty$ from condition $(R1)$ in (\ref{eq:RadialCondition1}).

Albeit crude, we conclude that uniformly for $\sigma \in \mathbb{C}\backslash \{0\}$, with $\text{Im }\sigma \geq 0$ the uniform differntial order for 
\begin{align}
    s>\frac{1}{2},
\end{align}
is sufficient for horizon propagation.
The only other characteristic region that is relevant are the radial points at $\partial X$. For now we will restrict ourselves only to the scattering calculus and the case that $\sigma \in \mathbb{C}\backslash \{0\}$, $\text{Im} (\sigma)>0$. We know that in this case, the characteristic set, and thus also the radial set is allowed to lie at finite covectors. We will wield the following definition.

\begin{definition} (Radial points at base infinity)\label{def:scRadialPoints}
    Let ${}^{sc}\hat{H}_{p} = \rho^{-1} ({}^{sc}H_p)$ denote the rescaled Hamiltonian vector field by $\rho$. We say that submanifold $\mathcal{R}_{\partial X}^\pm\subset \chr_{sc}(P(\sigma)\cap \partial X$ is a scattering radial source/sink for ${}^{sc}\hat{H}_p$ if
    \begin{enumerate}
        \item It is invariant/stationary under the flow of ${}^{sc}\hat{H}_p$ inside $\chr_{sc}(P(\sigma))\cap \partial X$ i.e.
        \begin{align}
            {}^{sc}\hat{H}_p|_{\mathcal{R}^\pm_{\partial X}} =0.
        \end{align}
        \item There exists a quadratic defining function $\mathfrak{q}$ restricted to $\partial X$ satisfying 
        $\mathfrak{q}^{-1}(0) = \mathcal{R}_{\partial X}^\pm$ such that 
        \begin{align}
            \pm{}^{sc}\hat{H}\mathfrak{q} = \beta \mathfrak{q}+\mathcal{O}(\mathfrak{q}^{3/2}).
        \end{align}
        with $\beta\in C^\infty(U)$ a nonnegative smooth function on a neighbourhood of $\mathcal{R}^\pm$.
        \item There exists positive $\beta_0\in C^\infty(\mathcal{M})$ such that 
        \begin{align}\label{eq:BoundaryRadialPointConditiion}
            \pm{}^{sc}\hat{H}\rho|_{\mathcal{R}^\pm_{\partial X}}=\beta_0\rho +\mathcal{O}(\rho^2),
        \end{align}
        near $\partial X$.
    \end{enumerate}
\end{definition}

In the case of the operator of interest (\ref{eq:SpectralFamily2}) we claim that the the radial source/sink for $\sigma>0$ are given by 
\begin{align}\label{eq:ScatteringRadialSource}
    \mathcal{R}_{\sigma,\partial X}^+ &=\{(\rho,y_i,\scxi, \sceta_i)\in {}^{sc}T^*_{\partial X}X\cap \Sigma_G\;|\; \rho=0, {}^{sc}\xi=2\text{Re } \sigma, {}^{sc}\eta_i=0\}\\
    \mathcal{R}_{\partial X}^-&=\{(\rho,y_i,\scxi, \sceta_i)\in {}^{sc}T^*_{\partial X}X\cap \Sigma_G\;|\; \rho=0, {}^{sc}\xi=0, {}^{sc}\eta_i=0\},\label{eq:ScatteringRadialSink}
\end{align}
with the roles reversed for $\sigma <0$.
Indeed, for $\mathcal{R}^+_{\sigma,\partial X}X$ one can take the defining function 
\begin{align}
    \mathfrak{q}= (\scxi-2\sigma)^2+\|{}^{sc}\eta_i\|^2,
\end{align}
and for $\mathcal{R}^-_{\partial X}X$ one takes 
\begin{align}
    \mathfrak{q}= \scxi^2+\|{}^{sc}\eta_i\|^2,
\end{align}

whence
\begin{align}
    {}^{sc}\hat{H}\mathfrak{q} = \pm 4 \sigma \mathfrak{q} +\mathcal{O}(\mathfrak{q}^{3/2}).
\end{align}
Furthermore, (\ref{eq:BoundaryRadialPointConditiion}) also holds with $\beta_0\equiv 2\sigma$.
For complex flow the dynamical picture away from $\sigma =0$ can still be resolved with this time radial source and radial sink (\ref{eq:ScatteringRadialSource})-(\ref{eq:ScatteringRadialSink}) are found by replacing $\sigma$ by $\text{Re } \sigma$.
We omit the analysis of the $b$-radial points at $\partial X$ as the spectral parameter $\sigma$ does not persist in the $b$-principal symbol in the model operator (\ref{eq:MinkowskiModel}). This can again be resolved with $b-sc$ calculus or by introducing a semiclasical parameter $h = |\sigma|^{-1}$ in the next chapter.\\
At $\mathcal{R}^+_{\sigma,\partial X}$ and $\mathcal{R}^-_{\partial X}$ we can get the following scattering radial point estimates which hold for any wave operator on a stationary, asymptotically flat spacetime $(M,g)$ that satisfies at $\partial X$
\begin{align}
    P_0|_{\partial X}=-2\partial_{t_*}\rho \left(\rho\partial_{\rho}-\frac{d-1}{2}-S\right) +\hat{\Box}_{\mathbb{M}^{1+d}}(0) + Q\partial_{t_*}-g^{00}\partial_{t_*}^2, \label{eq:CrudeWaveEquationindDimensions}
\end{align}
where $S\in C^\infty(X)$ and $Q\in \rho^3\text{Diff}^1_b(X)$, see Proposition 9.10 in \cite{hintz2026non}. We will verify in slightly more detail in Lemma \ref{lem:AdmissibleWaveEquation} that $\Box_{g_{\bfa,\bfb,\mfrm}}$ on $T^*M$ indeed satisfies this structural requirement for $S=0$ and $d=4$.

\begin{proposition}\label{prop:RadialScatteringSourceProposition}
    (Radial point estimate near $\mathcal{R}^+_{\sigma_0,\partial X}$) Let $0\neq \sigma_0\in\mathbb{R}$ and let $\nu>0$. Let $s,N \in \mathbb{R}$. Define 
    \begin{align}
        \theta_{in}=\sup_{\partial X}\left(\upsigma^1_b(\text{Im }(\hat{P}(0))|_{-2 \frac{dr}{r}}- 2 \text{Re } S\right).
    \end{align}
    Let $\mathcal{U}\subset \scT^*_{\partial X} X$ be a compact neighbourhood of $\mathcal{R}_{\sigma_0,\partial X}^+$, $K$ a compact neighbourhood of $\partial X$. Let $\chi \in C_c^\infty(X)$ with $\chi\equiv 1$ near $K$. 
    \begin{enumerate}
        \item (Direct estimate.) Let $\ell>\ell_0>\frac{1}{2}(-1+\theta_{in})$. There exists operators $B,G\in \Psi^{0,0}_{sc}(X)$ with Schwartz kernels in $K\times X$ and operator wavefront sets in $U$, with $B$ elliptic at $\mathcal{R}^+_{\sigma_0,\partial X}$ such that 
        \begin{align}
            \|Bu\|_{H^{s,\ell}_{sc}}\leq C\left(\|GP(\sigma)u\|_{H^{s-2,\ell+1}_{sc}} + \|u \|_{H_{sc}^{-N,\ell_0}}\right), \; \sigma \in \sigma_0+i[0,\nu],
        \end{align}
        holds in the strong sense.
        \item Let $r>\frac{1}{2}(-1+\theta_{in})$. There exists operators $B,E,G\in \Psi^{0,0}_{sc}(X)$ with Schwartz kernels supported in $K\times K$ and operator wave front sets in $\mathcal{U}$ and with $B$ elliptic at $\mathcal{R}^+_{\upsigma_0,\partial X}$ and $\WF'(E) \cap \mathcal{R}^+_{\sigma,\partial X}$ such that 
        \begin{align}
            \|Bv\|_{H^{-s+2,-\ell-1}_{sc}}\leq C \left(\|GP(\sigma)^*u\|_{H^{-s,-\ell}_{sc}}+\|Eu\|_{H_{sc}^{-s+2,-\ell-1}}+\|u\|_{H^{-N,-N}_{sc}}\right),\; \sigma \in \sigma_0 + i [0,\nu]
        \end{align}
    \end{enumerate}
    
\end{proposition}

For the outgoing estimate let us set the infimum to the spectrum of $S$ at $\partial X$ as
\begin{align}\label{eq:Sinfimum}
    \underline{S}= \inf_{p\in \partial X,\lambda \in \text{spec}(S)} \text{Re }\lambda, 
\end{align}
which for our case of the scalar wave equation is again equal to $0$.

\begin{proposition}\label{prop:RadialScatteringSinkProposition}
    (Radial point estimate near $\mathcal{R}^-_{\partial X}$.). Let $U,K$ and $\chi$ be as in Proposition \ref{prop:RadialScatteringSourceProposition}. Let $\mathcal D$ be a compact subset of $\{\sigma \in \mathbb{C}\;|\; \sigma \neq 0, \text{Im }\sigma \geq 0\}.$ We then have the following estimates.
    \begin{enumerate}
        \item Let $\ell<-\frac{1}{2}+\underline{S}$. There exist operators $B,E,G \in \Psi^{0,0}_{sc}(X)$ with Schwartz kernels supported in $K\times K$ and operator wave fornts sets in $U$, with $B$ elliptic at $\mathcal{R}^-_{\partial X}$. and $\WF'(E)\cap \mathcal{R}^-_{\partial X}$ such that 
        \begin{align}
            \|Bu\|_{H^{s,r}_{sc}}\leq C \left(\|G P(\sigma)u\|_{H_{sc}^{s-2,\ell+1}} +\|Eu\|_{H^{s,\ell}_{sc}}+\|\chi u\|_{H^{-N,\ell_0}_{sc}}\right),\; \sigma \in \mathcal{D}.
        \end{align}
        \item (Adjoint estimate) Let $\ell<\ell_0 < -\frac{1}{2}+\underline{S}$. There exists operators $B,G \in \Psi^{0,0}_{sc}(X)$ with Schwartz kernels supported in $K\times K$ and operator wave front sets in $U$ and $B$ elliptic at $\mathcal{R}^-_{\partial X}$ such that 
        \begin{align}
            \|Bu\|_{H^{-s+2,-\ell-1}_{sc}}\leq C\left(\|G\hat{P}_0(\sigma)u\|_{H_{sc}^{-s,-\ell}}+\|\chi u \|_{H_{sc}^{-N,-\ell_0-1}}\right), \sigma \in \mathcal {D}.
        \end{align}
    \end{enumerate}
\end{proposition}

All these estimates combine into a final preliminary Fredholm result for $\hat{\Box}_{g_{\bfa,\bfb,\mfrm}}(\sigma)$ in the scattering calculus. First of all recall the following abstract result from functional analysis that characterizes Fredholmness \cite{vasy2018minicourse}.

\begin{theorem}
Let $I\Subset(-\infty,0)\sqcup (0,\infty)$ be a compact interval and let $\nu>0$, set 
\begin{align}
    D=I+i[0,\nu].
\end{align}
    Fix $s>\frac{1}{2}$. Define the variable order Sobolev spaces 
    \begin{align}
       \mathcal{Y}_\sigma =\overline{H}{}^{s-1,\tilde\ell+1}_{sc}(X)\;\text{ and }\;\mathcal{X}_\sigma = \{u\in \overline{H}{}_{sc}^{s,\tilde{\ell}}(X)\;|\; P(\sigma)u\in \mathcal{Y}_\sigma\},
    \end{align}
    where we define $\mathcal{X}_\sigma$ with respect to the graph norm and $\tilde{\ell}\in C^\infty({}^{sc}T^*X)$ to be a variable order parameter that is non-increasing on geodesics from $\mathcal{R}^-_{\partial X}$ to $\mathcal{R}^+_{\sigma,\partial X}$ that satisfies 
    \begin{align}
       \tilde\ell>-\frac{1}{2} \text{ at } \mathcal{R}^+_{\sigma,\partial X} \text{ and } \tilde\ell<-\frac{1}{2} \text{ at } \mathcal{R}^-_{\partial X} 
    \end{align}
    and assume 
    \begin{align}
        s>\frac{1}{2}
    \end{align}
    Then for fixed $\sigma\in \mathbb{C}$ the operator 
    \begin{align}
       P(\sigma):= \hat{\Box}_{g_{\bfa,\bfb,\mfrm}}(\sigma): \mathcal{X}_\sigma\to \mathcal{Y}_\sigma
    \end{align}
satisfies the estimate
\begin{align}
    \|u\|_{\overline{H}{}^{s,\tilde\ell}_{sc}}\leq C \left(\|P(\sigma) u\|_{\overline{H}{}^{s-1,\tilde \ell+1}_{sc}} + \|u\|_{\overline{H}_{sc}^{-N,-N}}\right).
\end{align}
and its adjoint estimate 
\begin{align}
    \|v\|_{\dot{H}^{1-s,-\tilde\ell-1}_{sc}}\leq C\left(\|P(\sigma)^* v\|_{\dot{H}{}^{-s,-\tilde\ell}_{sc}}+\|v\|_{\dot{H}_{sc}^{-N,-N}}\right).
\end{align}

In particular, $P(\sigma)$ is Fredholm of index $0$.

\end{theorem}

\subsection{$sc-b$ calculus}\label{subsec:scbcalculus}

Whilst we managed to obtain preliminary Fredholm estimates in the scattering calculus for $\hat{\Box}_{g_{\bfa,\bfb,\mfrm}}(\sigma)$ the approach sketched out there is not sufficient for uniform $\sigma\in \mathbb{C}$ with $\text{Im }\sigma >0$. In fact the source/sink structure at $\partial X$ breaks down when we let $\sigma \to 0$. Moreover, the argument breaks on $b$-spaces where the radial source sink structure breaks down even for $\sigma\neq 0$ as $\mathcal{R}^-_{\partial X}$ lies on the double corner $\partial^2({}^b\overline{T^*}X)$. This is problematic as we eventually hope to prove an invertibility statement for $\hat{\Box}_{g_{\bfa,\bfb,\mfrm}}(\sigma)$ on conormal functions, hence we wish to establish a Fredholm result for $b$-Sobolev spaces with appropriate uniform control of $\sigma \to 0$. In this section we will summarily discuss the unresolved $sc,b$-calculus introduced in \cite{vasy2021limiting} that one can use to establish Fredholm estimates on $b$ spaces for finite $\sigma$ by resolving the double corners of ${}^b\overline{T}^*X$. In this construction we are also no longer forced to work in variable order Sobolev spaces as $b-$ and $sc-$ decay are separated, allowing us to satisfy the ingoing and outgoing conditions on the relevant radial sets seperately.\\ 

We let ${}^{sc}\overline{T}^*X$ be the fibre-radially compactified cotangent bundle. We define the $sc-b$ cotangent bundle as
\begin{align}
    {}^{sc,b}\overline{T}^*X:=[{}^{sc}\overline{T}^*X; o_{\partial_+ X}]
\end{align}
as the resulting space obtained from blowing up ${}^{sc}\overline{T}^*X$ along the zero section $o_{\partial X}=\{\scxi=0,\sceta=0\}$ along the boundary $\partial X=\{\rho=0\}$. We denote with 
\begin{align}
    \beta_{sc}:[{}^{sc}\overline{T}^* X;o_{\partial_+ X}]\to {}^{sc}\overline{T}^*X,
\end{align}
the corresponding blow-down map.
Equivalently it is the same as blowing up the co-dimension two corner of ${}^b\overline{T}^*X$, i.e. 
\begin{align}
    {}^{sc,b}\overline{T}^*X\cong [{}^b\overline{T}^*X;\partial^2({}^b\overline{T}^* X)]=[{}^b\overline{T}^*X;{}^bS^*_{\partial_+X}X],
\end{align}
with corresponding blow-down map $\beta_b:[{}^b\overline{T^*}X;\partial^2({}^b\overline{T^*}X)]\to {}^b\overline{T^*}X $,
where we defined
\begin{align}
   \partial^2({}^b\overline{T}^* X)= {}^b\overline{T}^*_{\partial_+ X} X\cap {}^bS^*X,
\end{align}
with the identification of ${}^bS^*X$ with fibre infinity. This new space has $3$ different boundary faces, with defining functions the fibre infinity order $\varrho_{\text{fi}}$, the scattering face $\rho_{\text{sf}}$ which is the lift of the scattering base fase ${}^{sc}T^*_{\partial X} X$ to $[{}^{sc}T^*X; o_{\partial X}]$ and $\rho_{\text{ff}}$ the defining function of the front face of the blow-up of $o_{\partial X}$.
In suitable local coordinates we can formulate local defining functions of each of the faces as 
\begin{align}
    \rho_{\text{ff}}&= (\rho^2 +\scxi^2+\|\sceta_i\|_{g^{-1}_{\mathbb{S}^3}}^2)^{1/2} \label{eq:FrontFaceDefiningFunction}\\
    \rho_{\text{sf}}&= \frac{\rho}{\sqrt{\rho^2+\scxi^2+\|\sceta_i\|_{g^{-1}_{\mathbb{S}^3}}^2}}\\
    \varrho_{\text{fi}}&=\frac{1}{\sqrt{\scxi^2+\|\sceta_i\|_{g^{-1}_{\mathbb{S}^3}}^2}},
\end{align}
and we have that the fibre infinity face and the blow-up front face are disjoint, connected by the scattering face, i.e.
\begin{align}
    \text{fi}\cap \text{ff}=\emptyset,\;\text{fi}\cap \text{sf}\neq \emptyset, \text{ and } \text{sf}\cap \text{ff}\neq\emptyset.
\end{align}
We can define a new symbol class over this space simply as 
\begin{align}
    S^{m,l_{sc},l_b}({}^{b,sc}\overline{T}^*X):= \varrho^{-m}_{\text{fi}}\rho^{-l_{sc}}_{\text{sf}}\rho_{\text{ff}}^{-l_{b}}\mathcal{A}^{0,0,0}({}^{sc,b}\overline{T}^*X),
\end{align}
where $m\in \mathbb{R}$ signifies differential order, $l_{sc}\in \mathbb{R}$ signifies scattering decay order and $l_b \in \mathbb{R}$ represents $b$-decay order. In local adapted coordinates near $\partial X$, one can then generate the $sc-b$ class of pseudodifferential operators $\Psi^{m,l_{sc},l_b}_{sc,b}(X)$ for $a_{sc,b}\in S^{m,l_{sc},l_b}_{sc,b}({}^{b,sc}\overline{T}^*X)$ by expressing $a_{sc,b}$ in $b$ covariables over $X^\circ$ and quantizing them using (\ref{eq:QuantizationbT*M}), i.e. one quantizes $a_b:= a_{sc,b}\circ \left(\beta_b|_{T^*X^\circ}\right)^{-1}$ using $\text{Op}_{sc,b}(a_{sc,b}):=\text{Op}_b(a_b)$. The $sc-b$ algebras once again satisfy the filtration condition
\begin{align}
\Psi_{sc,b}^{m,l_{sc},l_b}(X)\circ \Psi_{sc,b}^{m',l'_{sc},l'_b}(X)\subset \Psi_{sc,b}^{m+m',l_{sc}+l_{sc}',l_b+l_b'}(X).
\end{align}

The short sequence associated to the $sc-b$ symbol behaves like the combination of the $b$ and $sc$- sequences (\ref{eq:bSES}) and (\ref{eq:ScatteringSES}) in their respective decay parameters, which is also what is taken into account by the $sc-b$ principal symbol map
\begin{align}
    0\to \Psi^{m-1,l_{sc}-1,l_b}_{sc,b}(X) \to \Psi^{m,l_{sc},l_b}_{sc,b}(X)&\overset{{}^{sc,b}\sigma_{m,l_{sc},l_b}}{\to}\\ &S^{m,l_{sc},l_b}_{sc,b}({}^{sc,b}\overline{T^*}X)/S^{m-1,l_{sc}-1,l_b}_{sc,b}({}^{sc,b}\overline{T^*}X)\to 0.
\end{align}
Based on the sequence above one concludes that the subprincipal parts of an operator $A\in \Psi^{m,l_{sc},l_b}_{sc,b}(X)$ vanish on the locus 
\begin{align}
    \mathcal{Z}_{sc,b}:=\{\varrho_{\text{fi}} \rho_{\text{sf}}=0\}=\text{fi}\cup \text{sf}.
\end{align}

With this we can again define when an operator $A\in \Psi_{sc,b}^{m,l_{sc},l_b}(X)$ is elliptic. We define the normalized symbol of $A=\text{Op}_{sc,b}(a)$ where $a \in S_{sc,b}^{m,l_{sc},l_b}({}^{sc,b}\overline{T^*}X)$ via
\begin{align}
    \tilde a := \varrho_{\text{fi}}^m \rho_{\text{sf}}^{l_{sc}}\rho_{\text{ff}}^{l_b}a.
\end{align}

and define the corresponding elliptic set $\Ell_{sc,b}^{m,l_{sc},l_b}(A)\subset \mathcal{Z}_{sc,b}\subset {}^{sc,b}\overline{T^*}X$ as 

\begin{align}
    \Ell_{sc,b}^{m,l_{sc},l_b}(A):=\{q\in \mathcal{Z}_{sc,b}\;|\; \exists c>0, U \ni q, \text{ open  in } {}^{sc,b}\overline{T^*}X \text{ s.t } |\tilde{a}|\geq c>0 \text{ on } U\cap T^*X^\circ\}
\end{align}

and once again $\chr^{m,l_{sc},l_b}_{sc,b}(A)$ as 
\begin{align}
    \chr^{m,l_{sc},l_b}_{sc,b}(A):=\mathcal{Z}_{sc,b}\backslash  \Ell_{sc,b}^{m,l_{sc},l_b}(A).
\end{align}
The continuous boundary value argument of Remark \ref{rem:ContinousExtenstions}, applies to the elliptic and characteristic sets as well.

We also define the operator wavefront set $\WF_{sc,b}'(A)\subset \mathcal
{Z}_{sc,b}\subset {}^{sc,b}{\overline{T^*}}X $ as 
\begin{align}
    \WF'(A):=\esssup a_{sc,b}\subset \mathcal{Z}_{sc,b}, 
\end{align}
where 
\begin{align}
   \mathcal{Z}_{sc,b}\backslash \esssup(a_{sc,b}):=\{q\in \mathcal{Z}_{sc,b}\;|\; \exists \chi\in C^\infty({}^{sc,b}\overline{T^*}X, \chi \equiv 1 \text{ near } q, \chi a\in S^{-\infty,-\infty,\ell_b}({}^{sc,b}\overline{T^*}X)\}
\end{align}
where we identified the smoothing symbols through 
\begin{align}
    S^{-\infty,-\infty,\ell_b}_{sc,b}({}^{sc,b}\overline{T^*}X)=\bigcap_{m,l_{sc}\in\mathbb{R}} S_{sc,b}^{m,l_{sc},l_b}({}^{sc,b}\overline{T^*}X)
\end{align}
whose $sc-b$ quantization defines 
\begin{align}
    \Psi^{-\infty,-\infty,\ell_b}_{sc,b}(X)=\bigcap_{m,\ell_{sc}\in \mathbb{R}} \Psi^{m,\ell_{sc},\ell_b}_{sc,b}(X).
\end{align}

\begin{definition}
    Let $A\in \Psi^{s,\ell_{sc}-\ell_b,0}(X)$ be elliptic with respect to the differential order and scattering decay (we can take $A= \Op_{sc,b}(\rho_{\text{fi}}^{-s}\rho_{\text{sf}}^{-(\ell_{sc}-\ell_b)})$)\footnote{Note that locally $\rho^{-\ell_{b}}=\rho_{\text{sf}}^{-\ell_b}\rho_{\text{ff}}^{-\ell_b}$ so that for our choice of $A$ we have $\rho^{-\ell_b}A\in \Psi^{s,\ell_{sc},\ell_b}_{sc,b}(X)$.}. We then define the $sc-b$ Sobolev spaces as 
    \begin{align}
        H^{s,\ell_{sc},\ell_b}_{sc,b}(X) =\left\{u \in H^{\tilde s,\ell_b}_b(X)\;|\; Au\in \rho^{\ell_b} L^2_{sc}(X)\right\},
    \end{align}
    where $\tilde s$ is sufficiently small,
    \begin{align}
        \tilde s \leq \min (s,r_{sc}-r_b).
    \end{align}
    They can be equipped with the norm 
    \begin{align}
        \|u\|^2_{H^{s,r_{sc},r_b}_{b,sc}}=\|\rho^{-r_b}Au\|^2_{L^2_{sc}} +\|u\|^2_{H_b^{\tilde s,r_b}}.
    \end{align}
\end{definition}

Naturally then if $X$ is compact we have that $A\in \Psi^{m,l_{sc},l_b}_{sc,b}(X)$ define bounded maps
\begin{align}
    A:H^{s,\ell_{sc},\ell_b}_{sc,b}(X)\to H^{s-m,\ell_{sc}-l_{sc},\ell_b-l_b}_{sc,b}(X).
\end{align}

It is these Sobolev spaces that yield the translation towards estimates on $b$-spaces.
\begin{lemma}\label{lem:btoscTranslationLemma}
    For any $s,\ell_{sc},\ell_b\in \mathbb{R}$ we have the following identification of $b$ and $sc-b$ spaces.    
\begin{align}\label{eq:TranslationtobSpaces}
    H^{s,\ell_b}_{b}(X)= H^{s,s+\ell_b,\ell_b}_{sc,b}(X).
\end{align}    
\end{lemma}

We have the following propagation and elliptic estimates, consistent with the microlocal locus $\mathcal{Z}_{sc,b}$.

\begin{proposition}
    Let 
    \begin{align}
    P\in \Psi^{m,l_{sc},l_b}_{sc,b}(X),\;\;\;\;B,G\in \Psi^{0,0,0}_{sc,b}(X),
    \end{align}
    and suppose 
    \begin{align}
        \WF_{sc,b}'(B)\subset \Ell_{sc,b}(P)\cap \Ell_{sc,b}(G).
    \end{align}
   Then for $s,r_{sc},r_b\in \mathbb{R}$ and sufficiently large $N$, we have 
   \begin{align}
       \|Bu\|_{H_{sc,b}^{s,\ell_{sc},\ell_b}}\leq C\left(\|GPu\|_{H_{sc,b}^{s-m,\ell_{sc}-l_{sc},\ell_b-l_b}} +\|\chi u \|_{H_{sc,b}^{-N,-N,l_b}}\right).
   \end{align}
    
\end{proposition}

\begin{proposition}\label{prop:scbPropagation}
    Assume $P$ has a real principal-symbol representative $p$, meaning 
    \begin{align}
        P-\text{Op}_{sc,b}(p)\in \Psi^{m-1,l_{sc}-1,l_b}_{sc,b}(X),
    \end{align}
    and that its renormalized symbol 
    \begin{align}
        \tilde p =\varrho_{\text{fi}}^m\rho_{\text{sf}}^{l_{sc}}\rho_{\text{ff}}^{l_b}p \in C^\infty({}^{sc,b}\overline{T^*}X)
    \end{align}
    and consider its rescaled vector field 
    \begin{align}
        {}^{sc,b}\hat{H}_p=\varrho_{\text{fi}}^{m-1}\rho_{\text{sf}}^{l_{sc}-1}\rho_{\text{ff}}^{l_b}H_p.
    \end{align}
    Let $B,E,G \in \Psi^{0,0,0}_{sc,b}(X)$ with 
    \begin{align}
        \WF'_{sc,b}(B) \subset \Ell_{sc,b}(G),
    \end{align}
    Assume that at every point of 
    \begin{align}
        \WF'_{sc,b}(B)\cap \chr_{sc,b}(P),
    \end{align}
    can be reached from $\Ell_{sc,b}(E)$ along a finite forward bicharacteristic segment $\gamma$ of ${}^{sc,b}\hat{H}_p$entirely contained in $\Ell_{sc,b}(G)$.
    Then we have that 
    \begin{align}
        \|Bu\|_{H_{sc,b}^{s,\ell_{sc},l_b}}\leq C\Big(\|GPu\|_{H_{sc,b}^{s-m+1,\ell_{sc}-l_{sc}+1,\ell_b-l_b}} +\|Eu\|_{H_{sc,b}^{s,\ell_{sc},\ell_b}}
         +\|\chi u\|_{H^{-N,-N,\ell_b}} \Big)
    \end{align}

\end{proposition}

The operator of interest is once again $P(\sigma) =\hat{\Box}_{g_{\bfa,\bfb,\mfrm}}(\sigma)$ seen as an operator on the sc-b Sobolev spaces. Recall that we have at $\partial_+X$ that 
\begin{align}\label{eq:MoreRefinedOperatoratInfinity}
    P(\sigma) \in (\rho^2\text{Diff}_b^2(X)+ \sigma\rho \text{Diff}_b^1(X) +\sigma^2\rho^2C^\infty(X)) \cap \Diff^2_{sc}(X),
\end{align}
see also Lemma \ref{lem:AdmissibleWaveEquation}, and so we have 
\begin{align}
    P(\sigma) \in \Psi^{2,0,-1}_{sc,b}(X).
\end{align}

Now let $\Char_{sc}^{2,0}(P(\sigma)) :=\Sigma_\sigma\subset {}^{sc}\overline{T}^*X$ and let $\chr^{2,0,-1}_{sc,b}(P(\sigma)):=\tilde{\Sigma}_\sigma\subset {}^{sc,b}\overline{T^*}X$. For real $\sigma \neq 0$, $\tilde{\Sigma}_\sigma$ is the proper  lift of $\Sigma_\sigma$ to $ {}^{sc,b}\overline{T^*}X\overset{\beta}{\to}{}^{sc}\overline{T^*}X$ i.e. 
\begin{align}
    \tilde{\Sigma}_\sigma=\overline{\beta^{-1}(\Sigma_\sigma\backslash o_{\partial X})}.
\end{align}
We note that at fibre infinity over $\partial_+X$, the operator $P(\sigma)$ is elliptic. Consequently, $\tilde{\Sigma}_\sigma \cap \text{fi}\cap \text{sf}=\emptyset$, and hence the characteristic dynamics of the lifted principal symbol $\beta^*p_\sigma$ take place away from fibre infinity on $\text{sf}$, inlcuding the corner $\text{sf}\cap \text{ff}$.

In local projective coordinates near the front face $\{\rho_{\text{ff}}=0\}$ at the blow up, we have 
\begin{align}
    {}^{sc}\hat\xi= \frac{{}^{sc}\xi}{\rho_{\text{ff}}}=\rho_{\text{sf}}{}^b\xi,\;\;\; {}^{sc}\hat{\eta}_i=\frac{{}^{sc}\eta_i}{\rho_{\text{ff}}}=\rho_{\text{sf}}{}^b\eta_i
\end{align}
Note thus that we have that 
\begin{align}
    \text{ff}^\circ \cong \bT^*_{\partial X}X \text{ and } \text{ff}\cap \text{sf}\cong {}^bS^*_{\partial X}X.
\end{align}

Under the assumptions on the rescaled symbol $\tilde p$ in Proposition \ref{prop:scbPropagation}, we get for $P(\sigma) \in \Psi^{2,0,-1}_{sc,b}(X)$ that its $sc,b$-Hamiltonian vector field can be written as 
 \begin{align}
{}^{sc,b}\hat{H}_{p_\sigma}:=\varrho_{\text{fi}}\rho_{\text{sf}}^{-1}\rho_{\text{ff}}^{-1}H_{\beta_{sc}^*p}=\varrho_{\text{fi}}\rho^{-1}H_{\beta_{sc}^*p},\label{eq:ActualscbHField}
\end{align}
where we used $\rho=\rho_{\text{sf}}\rho_{\text{ff}}$.
Let us therefore consider the rescaled lifted $sc-b$ symbol which is given here at $\text{ff}$ by 
\begin{align}\label{eq:SymbolNearff}
    q_\sigma:=\rho^{-1}\beta_{sc}^*p_\sigma := \rho(\bxi^2+|{}^b\eta|^2)-2\sigma\bxi,
\end{align}
and we also define its renormalized symbol though 
\begin{equation}
    \hat{q}_\sigma :=\rho_{\text{sf}}q_\sigma = \rho_{\text{ff}}^{-1}\beta_{sc}^*p_\sigma=\rho_\text{ff}(1-\rho_{\text{sf}}^2)-2\sigma {}^{sc}\hat{\xi}=\rho_{\text{ff}}^{-1}\beta_{sc}^*p_\sigma.\label{eq:RenormalizedRenormalizedSymbol}
\end{equation}
Away from $\text{fi}$, near $\text{ff}\cap \text{sf}$ these are the correct symbols to consider on $\tilde \Sigma_\sigma$ as we have 
\begin{align}
    H_{q_\sigma}|_{\tilde \Sigma_\sigma}=\rho^{-1}H_{\beta^*_{sc}p_\sigma}|_{\tilde\Sigma_\sigma}
\end{align}
so by (\ref{eq:ActualscbHField}) it agrees with ${}^{sc,b}\hat{H}$ up to a harmless factor of $\varrho_{\text{fi}}$ as we have $\text{fi}\cap \text{ff}=\emptyset$. Similarly for $H_{\hat{q}_\sigma}$ we have 
\begin{align}
    H_{q_\sigma} = \rho_{\text{sf}}^{-1}H_{\hat{q}_\sigma}
\end{align}
allowing us to restrict to $\text{sf}$.
When doing so (\ref{eq:RenormalizedRenormalizedSymbol}) simplifies to
\begin{align}
    \hat{q}_\sigma|_{\text{sf}}=\rho_{\text{ff}}-2\sigma {}^{sc}\hat{\xi},
\end{align}
yielding the characteristic relation for $\sigma \neq 0$ that
\begin{align}
   \tilde{\Sigma}_\sigma \cap \text{sf}= \{\rho_{\text{sf}}=0, \rho_{\text{ff}}=2\sigma{}^{sc}\hat{\xi}\}
\end{align}
Restricting then to $\text{ff}$ we see that
\begin{align}\label{eq:SuspectedRadialffSet}
    \tilde{\Sigma}_\sigma\cap \text{ff} \cap \text{sf} =\{\rho_{\text{sf}}=0,\rho_{\text{ff}}=0, {}^{sc}\hat{\xi}=0,\|{}^{sc}\hat{\eta}_i\|_{g^{-1}_{\mathbb{S}^3}}^2=1\},
\end{align}
where $\|{}^{sc}\hat{\eta}_i\|_{g^{-1}_{\mathbb{S}^3}}^2=1$ follows from (\ref{eq:FrontFaceDefiningFunction}), i.e. 
\begin{align}
    1=\rho_{\text{sf}}^2+({}^{sc}\hat{\xi})^2+\|{}^{sc}\hat{\eta}_i\|_{g^{-1}_{\mathbb{S}^3}}^2
\end{align}
and where we used 
\begin{align}
    \rho=\rho_{\text{ff}}\rho_{\text{sf}}.
\end{align}

Now we turn our attention to the radial point that is problematic in the standard $b$-calculus, $\mathcal{R}^-_{\partial X}$. On ${}^{sc,b}\overline{T^*}X$ we consider 
\begin{align}
    \tilde{{R}}^-_{\sigma,\partial_+ X}:= \tilde{\Sigma}_\sigma \cap \text{ff}.
\end{align}

which is indeed the lift of $\mathcal{R}^-_{\partial X}$ intersected with the $sc-b$ characteristic set $\tilde \Sigma_\sigma$,
\begin{align}
    \tilde{R}^-_{\partial X}:=\overline{\beta^{-1}(\mathcal{R}^-_{\partial X})}\cap \tilde{\Sigma}_\sigma.
\end{align}

\begin{definition}\label{def:RadialscbSet}
    We say that a radial set $\tilde{\mathcal{R}}^-_{\partial{X}}\subset \tilde{\Sigma}_\sigma\cap \text{ff}\cap \text{sf} $ is a radial set, resp. source/sink at the front face $\text{ff}\subset {}^{sc,b}\overline{T}X$ if it is invariant under the flow of $H_{q(\sigma)}$ and the following properties hold.
    \begin{enumerate}
        \item $\tilde{\mathcal{R}}^-_{\partial{X}}\subset \tilde{\Sigma}\cap \text{ff} = \{q_\sigma=0,\rho_{\text{ff}}=0\}$, and moreover we have $d\hat{q}_\sigma|_{\text{sf}} \neq 0$ on $\tilde{\mathcal{R}}^-_{\partial X}$.
        \item There exists an open neighbourhood $U$ of $\tilde{\mathcal{R}}^-_{\partial{X}}$ inside $\text{ff}$ and a positive function $\beta_0\in C^\infty(U)$ with $\beta_0>0$ such that 
        \begin{align}
            \left.\rho_{\text{ff}}^{-1}H_{q_\sigma}\rho_{\text{ff}}\right|_{\tilde{\mathcal{R}}^-_{\sigma,\partial_+X}}= \mp\beta_0
        \end{align}
        \item There exists a non-degenerate quadratic defining function $\mathfrak{q}$ of $\tilde{\mathcal{R}}^-_{\partial{X}}$ such that 
        \begin{align}
            \mp H_{q_\sigma}\mathfrak{q}=2\beta_0\mathfrak{q}+\mathcal{O}(\mathfrak{q}^{3/2}).
        \end{align}
    \end{enumerate}
\end{definition}

In  local coordinates we see that the corresponding Hamiltonian vector field $H_{q_\sigma} \in \mathcal{V}_b({}^{sc,b}\overline{T}^*X)$ is given by the $b$-Hamilton convention in equation (\ref{eq:bHamiltonVFieldConvention}), the lift of which is a smooth $b$ vector field by restricting to $\text{sf}\cap \text{ff}$ away from $\text{fi}$ where more differnetial rescaling is needed.
\begin{align}
    H_{q_\sigma}=(2\rho {}^b\xi-2\sigma)\rho\partial_\rho-&\rho ({}^b\xi^2+|{}^b\eta|^2)\partial_{\bxi}+2\rho \sum_{j,k} (g^{-1}_{\mathbb{S}^3})^{jk}({}^b\eta)_k\partial_{y_j}\\&-\rho \sum_{j,k,l}(\partial_{y_j}g^{-1}_{\mathbb{S}^3})^{kl}\cdot{}^b\eta_k{}^b\eta_l \partial_{{}^b\eta_j}
\end{align}
where we note that when working in Hopf coordinates the inverse boundary metric $g_{\partial X}^{-1}= g^{-1}_{\mathbb{S}^3}$ is given by $\text{diag}(1,\sin^{-2}\theta,\cos^{-2}\theta)$.\\
We now claim that $\tilde{\mathcal{R}}^-_{\sigma,\partial X}=\tilde{\Sigma}_\sigma\cap \text{ff}$ is indeed radial in the sense of Definition \ref{def:RadialscbSet}. By moving to the local variables near $\text{sf}\cap \text{ff}$ we immediately check that 
\begin{align}
    H_{q_\sigma}\rho_{\text{ff}}&= \rho_{\text{ff}}\left[-2\sigma +\rho_{\text{ff}}{}^{sc}\hat{\xi}(1+\rho_{\text{sf}}^2)\right]\label{eq:scbRadialSetCheck1}\\
    H_{q_\sigma}\rho_{\text{sf}}&= \rho_{\text{ff}}\rho_{\text{sf}}{}^{sc}\hat{\xi}(1-\rho_{\text{sf}}^2)\label{eq:scbRadialSetCheck2}\\
    H_{q_\sigma}({}^{sc}\hat{\xi})&=-\rho_{\text{ff}}\left[1-({}^{sc}\hat{\xi})^2\right](1-\rho_{\text{sf}}^2), \label{eq:scbRadialSetCheck3}
\end{align}
with the remaining angular components also carrying a factor $\rho_{\text{ff}}$. Note that $H_{q_\sigma}$ agrees with the smooth extension of $\rho_{\text{sf}}^{-1}H_{\hat{q}_\sigma}$ on characteristic hypersurfaces near $\text{sf}\cap \text{ff}$. These identities immediately vanish on the set (\ref{eq:SuspectedRadialffSet}) and $H_{q_\sigma}$ is tangent to this set. Indeed we have 
\begin{align}
    H_{q_\sigma}\hat{q}_\sigma =\rho_{\text{ff}}{}^{sc}\hat{\xi}(1-\rho_{\text{sf}}^2)\hat{q}_\sigma,
\end{align}
which restricted to $\text{sf}$, $\{\rho_{\text{sf}}=0\}$ reduces to
\begin{align}
    H_{q_\sigma}\hat{q}_\sigma =\rho_{\text{ff}}{}^{sc}\hat{\xi}\hat{q}_\sigma,
\end{align}

so the characteristic hyperurface is invariant.
To check the non-degeneracy of $d\hat{q}_\sigma$ simply observe that 
\begin{align}
    \partial_{{}^{sc}\hat{\xi}}\left(\hat{q}_\sigma|_{\text{sf}}\right)=-2\sigma \neq 0
\end{align}
under the assumption that $\sigma \in \mathbb{R}\backslash\{0\}$. Using the characteristic relations $\rho_{\text{sf}}=0$ and ${}^{sc}\hat{\xi}=\frac{\rho_{\text{ff}}}{2\sigma}$, we compute using the quadratic defining function 
\begin{align}
    \mathfrak{q}=\rho_{\text{ff}}^2,
\end{align}
for $\tilde{\mathcal{R}}^-_{\sigma,\partial_+ X}$ near $\text{sf}\cap \text{ff}$, aand we find 
\begin{align}\label{eq:FixedSigmaQuadraticCheck}
    H_{q_\sigma}\mathfrak{q}&= 2\rho_{\text{ff}}H_{q_\sigma}\rho_{\text{ff}}\\
    &=-4\sigma \rho_{\text{ff}}^2+\frac{\rho_{\text{ff}}^4}{\sigma}\\
    &=-4\sigma \mathfrak{q}+\frac{\mathfrak{q}^2}{\sigma}
\end{align}
on $\tilde{\Sigma}_\sigma$ near the corner $\text{sf}\cap \text{ff}$.

The second radial set we have effectively already analyzed in the scattering picture. Naturally it would be the lift of $\mathcal{R}^+_{\sigma,\partial X}$, to ${}^{sc,b}\overline{T}^*X$, i.e.
\begin{align}
    \tilde{R}^+_{\sigma,\partial X}:= \overline{\beta^{-1}(\mathcal{R}^+_{\sigma,\partial X})},
\end{align}

and lies fully on the interior of the scattering face, $\text{sf}^\circ$, away from the corner $\text{ff}\cap \text{sf}$, We know that $\sigma \neq 0$, $\mathcal{R}^+_{\sigma,\partial X}$ lies away from the blown up zero section $o_{\partial X}$ hence $\beta:[\scT^*X;o_{\partial X}]\to \scT^*X$ is locally a diffeomorphism between $ \tilde{R}^+_{\sigma,\partial X}$ and  $R^+_{\sigma,\partial X}$ and scattering coordinates $(\rho,y,\scxi, \sceta)$ remain completely regular after lifting to ${}^{sc,b}\overline{T}^* X$ in this neighbourhood. The radial set verification stays essentially the same as we have already used before. Taking local coordinates $(\rho,y,\scxi,\sceta_i)$ on $\text{sf}^\circ\subset {}^{sc,b}\overline{T}^*X$ one immediately concludes that 
\begin{align}
    \tilde{\mathcal{R}}^+_{\sigma,\partial X}&=\text{sf}^\circ \cap \beta^{-1}(\mathcal{R}^+_{\sigma,\partial X})\\
    &=\left\{(\rho,y,\scxi,\sceta)\in {}^{sc,b}\overline{T}^*X\;|\; \rho_{\text{sf}}=0,{}^{sc}\xi=2\text{ Re }\sigma,\sceta=0\right\}.
\end{align}
We can now rephrase this in semiclassical radial set estimates for 
\begin{align}
    P(\sigma):\left\{u\in H^{s,\ell_{sc},\ell_b}_{sc,b}(X)\;|\; P(\sigma) u \in H_{sc,b}^{s-1,\ell_{sc}+1,\ell_b+1}(X)\right\}\to H^{s-1,r_{sc}+1,\ell_b+1}_{sc,b}(X),
\end{align}
where we can take $s,r_{sc},r_b\in \mathbb{R}$ to be constant in contrast to the base scattering estimate, as the seperate $r_{sc}$ and $r_b$ decay paramter guarantee threshold conditions on both radial sets. 
With the resolved dynamics for real $\sigma\in \mathbb{R}\backslash \{0\}$ we can then state $sc-b$ Fredholm estimates for fixed $\sigma\in \mathbb{H}\backslash \{0\}\subset \mathbb{C}$ where for $\text{Im}\,\sigma>0$ one can analyze the dynamics for $H_{p_{{}_{\text{Re} \sigma}}}$ between $\tilde{\mathcal{R}}^+_{\text{Re }\sigma, \partial_+ X}$ and $\tilde{\mathcal{R}}^-_{\partial_+X}$.

\begin{proposition}[Incoming radial sc-b estimates]
\label{prop:scbFredHolmEstimateincoming}
Let
\begin{align}
  \mathcal K\Subset
  \{\sigma\in\mathbb C:
    \operatorname{Im}\sigma\geq0,\ 
    \operatorname{Re}\sigma\neq0\},
  \qquad
  \ell_{\mathrm{sc}}>\ell_0>\frac{1}{2}(-1+\theta_{in}),
  \qquad s,\ell_b\in\mathbb R,
\end{align}
with $\theta_{in}$ defined as in Proposition \ref{prop:RadialScatteringSourceProposition}.
For $\sigma\in\mathcal K$, let
\begin{align}
  \widetilde{\mathcal R}^{+}_{\operatorname{Re}\sigma,\partial_+X}
  =
  \{\rho_{\mathrm{sf}}=0,\
    {}^{\mathrm{sc}}\xi=2\operatorname{Re}\sigma,\
    {}^{\mathrm{sc}}\eta=0\}
  \subset\mathrm{sf}^{\circ}.
\end{align}\footnote{For nonreal $\sigma$, this is a reference set for the incoming
estimate, not an assertion that it belongs to the characteristic
set of the full complex symbol.}

Choose $K\Subset U$ as above, containing these reference sets, with
\begin{align}
  \overline U\cap(\mathrm{fi}\cup\mathrm{ff})=\varnothing.
\end{align}
The cutoffs below can be chosen as scattering cutoffs whose full
symbols vanish near the blown-up zero section; in particular, they
have trivial normal operator at $\mathrm{ff}$.

For every sufficiently large $N$, the following estimates hold
uniformly for $\sigma\in\mathcal K$.

\begin{enumerate}
\item \emph{Direct estimate.}
There exist cutoffs $B,G\in\Psi_{sc,b}^{0,0,0}(X)$, with $B$ elliptic at
$\widetilde{\mathcal R}^{+}_{\operatorname{Re}\sigma,\partial_+X}$
and
\[
  \operatorname{WF}'_{\mathrm{sc},b}(B)
  \subset
  \operatorname{Ell}_{\mathrm{sc},b}^{0,0,0}(G),
\]
such that
\begin{equation}
\label{eq:scb-incoming-direct}
  \|Bu\|_{H_{\mathrm{sc},b}^{s,\ell_{\mathrm{sc}},\ell_b}}
  \leq C\left(
    \|GP(\sigma)u\|_{
      H_{\mathrm{sc},b}^{s-2,\ell_{\mathrm{sc}}+1,\ell_b+1}}
    +
    \|\chi u\|_{H_{\mathrm{sc},b}^{-N,\ell_0,\ell_b}}
  \right).
\end{equation}

\item \emph{Adjoint estimate.}
There exist cutoffs $B,E,G\in \Psi^{0,0,0}_{sc,b}(X)$, with $B$ elliptic at
$\widetilde{\mathcal R}^{+}_{\operatorname{Re}\sigma,\partial_+X}$,
\begin{align}
  \operatorname{WF}'_{\mathrm{sc},b}(E)
  \cap
  \widetilde{\mathcal R}^{+}_{\operatorname{Re}\sigma,\partial_+X}
  =\varnothing,
\end{align}
and
\begin{align}
  \operatorname{WF}'_{\mathrm{sc},b}(B)
  \cup
  \operatorname{WF}'_{\mathrm{sc},b}(E)
  \subset
  \operatorname{Ell}_{\mathrm{sc},b}^{0,0,0}(G),
\end{align}
such that
\begin{equation}
\label{eq:scb-incoming-adjoint}
\begin{split}
  \|Bv\|_{
    H_{\mathrm{sc},b}^{2-s,-\ell_{\mathrm{sc}}-1,-\ell_b-1}}
  \leq C\bigl(&
    \|GP(\sigma)^*v\|_{
      H_{\mathrm{sc},b}^{-s,-\ell_{\mathrm{sc}},-\ell_b}}
    +
    \|Ev\|_{
      H_{\mathrm{sc},b}^{2-s,-\ell_{\mathrm{sc}}-1,-\ell_b-1}}
    \\
    &+
    \|\chi v\|_{H_{\mathrm{sc},b}^{-N,-N,-N}}
  \bigr).
\end{split}
\end{equation}
Here $E$ is chosen elliptic on the surrounding characteristic region
providing the propagation input, and $G$ is elliptic along the
connecting bicharacteristic segments for real $\sigma$.
\end{enumerate}
\end{proposition}
\begin{remark}   
Both estimates (\ref{eq:scb-incoming-direct}) and (\ref{eq:scb-incoming-adjoint}) hold in the strong sense: finiteness of the displayed
right-hand side implies finiteness of the left-hand side and the
stated bound. Moreover for $P(\sigma)=\hat{\Box}_{g_{\bfa,\bfb,\mfrm}}(\sigma)$ as chosen we can set $\theta_{in}=0$ yielding the propagation threshold
\begin{align}
    \ell_{sc}>\ell_0>-\frac{1}{2}.
\end{align}
\end{remark}

\begin{proposition}[Outgoing second microlocal radial estimates]
\label{prop:ScbOutgoing}
Let
\begin{align}
  \mathcal K\Subset
  \{\sigma\in\mathbb C:
    \operatorname{Im}\sigma\geq0,\ \sigma\neq0\},
  \qquad
  \ell_b<-\frac{1}{2}+\underline{S},
  \qquad s,\ell_{\mathrm{sc}}\in\mathbb R,
\end{align}
with $\underline{S}$ defined as in equation (\ref{eq:Sinfimum}).

Let
\begin{align}
  \widetilde{\mathcal R}^{-}_{\partial_+X}
  =
  \{\rho_{\mathrm{sf}}=0,\ \rho_{\mathrm{ff}}=0,\
    {}^{\mathrm{sc}}\widehat\xi=0\}
  \subset\mathrm{sf}\cap\mathrm{ff}.
\end{align}
Choose $K\Subset U$ as above around this set, with
\begin{align}
  \overline U\cap\mathrm{fi}=\varnothing.
\end{align}
For every sufficiently large $N$, the following estimates hold
uniformly for $\sigma\in\mathcal K$.

\begin{enumerate}
\item \emph{Direct estimate.}
There exist cutoffs $B,E,G\in\Psi^{0,0,0}_{sc,b}(X)$, with $B$ elliptic at
$\widetilde{\mathcal R}^{-}_{\partial_+X}$,
\begin{align}
  \operatorname{WF}'_{\mathrm{sc},b}(E)
  \cap
  \widetilde{\mathcal R}^{-}_{\partial_+X}
  =\varnothing,
\end{align}
and
\begin{align}
  \operatorname{WF}'_{\mathrm{sc},b}(B)
  \cup
  \operatorname{WF}'_{\mathrm{sc},b}(E)
  \subset
  \operatorname{Ell}_{\mathrm{sc},b}^{0,0,0}(G),
\end{align}
such that
\begin{equation}
\label{eq:scb-outgoing-direct}
\begin{split}
  \|Bu\|_{H_{\mathrm{sc},b}^{s,\ell_{\mathrm{sc}},\ell_b}}
  \leq C\bigl(&
    \|GP(\sigma)u\|_{
      H_{\mathrm{sc},b}^{s-2,\ell_{\mathrm{sc}}+1,\ell_b+1}}
    +
    \|Eu\|_{H_{\mathrm{sc},b}^{s,\ell_{\mathrm{sc}},\ell_b}}
    \\
    &+
    \|\chi u\|_{H_{\mathrm{sc},b}^{-N,-N,r_b}}
  \bigr).
\end{split}
\end{equation}
Here $E$ is chosen elliptic on the surrounding characteristic region
providing the propagation input, and $G$ is elliptic along the
connecting bicharacteristic segments for real $\sigma$.

\item \emph{Adjoint estimate.}
There exist cutoffs $B,G\in \Psi^{0,0,0}_{sc,b}(X)$, with $B$ elliptic at
$\widetilde{\mathcal R}^{-}_{\partial_+X}$ and
\begin{align}
  \operatorname{WF}'_{\mathrm{sc},b}(B)
  \subset
  \operatorname{Ell}_{\mathrm{sc},b}^{0,0,0}(G),
\end{align}
such that
\begin{equation}
\label{eq:scb-outgoing-adjoint}
\begin{split}
  \|Bv\|_{
    H_{\mathrm{sc},b}^{2-s,-\ell_{\mathrm{sc}}-1,-\ell_b-1}}
  \leq C\bigl(&
    \|GP(\sigma)^*v\|_{
      H_{\mathrm{sc},b}^{-s,-\ell_{\mathrm{sc}},-\ell_b}}
    +
    \|\chi v\|_{H_{\mathrm{sc},b}^{-N,-N,-\ell_b-1}}
  \bigr).
\end{split}
\end{equation}
\end{enumerate}

\end{proposition}
\begin{remark}
Both estimates hold in the strong sense at the fixed $b$-weight
displayed: finiteness of the right-hand side implies finiteness of
the left-hand side and the stated bound.
No improvement of the $b$-weight is asserted. Moreover for $P(\sigma)=\hat{\Box}_{g_{\bfa,\bfb,\mfrm}}(\sigma)$ we have again that $\underline{S}=0$ hence giving the propagation threshold
\begin{align}
    \ell_b<-\frac{1}{2}.
\end{align}
\end{remark}

Analogous to the previous $sc$-estimate it is possible to once again obtain a sc-b Fredholm estimate. For the wave equation $P(\sigma)=\hat{\Box}_{g_{\bfa,\bfb,\mfrm}}(\sigma)$.

\begin{proposition}(Uniform $sc-b$ Fredholm estimate on compact spectral data)
    Let $P(\sigma)=\hat{\Box}_{g_{\bfa,\bfb,\mfrm}}(\sigma)$, and
    let $\mathcal{K}\subset\mathbb{C}$ be given such that
    \begin{align}
        \mathcal{K}\Subset\{\sigma \in \mathbb{C}\;|\; \text{Im }\sigma\geq 0, \sigma\neq 0\}.
    \end{align}
    Set $s,r_{sc},r_b\in \mathbb{R}$ such that
    \begin{align}
        s>\frac{1}{2}- \frac{\inf_{\sigma\in \mathcal{K}}\text{Im }\sigma}{\kappa},\; \ell_{sc} >-\frac{1}{2},\text{ and } \ell_b<-\frac{1}{2}. 
    \end{align}
Then for sufficiently small $\delta>0$ and sufficiently large $N$ the following direct estimate hold uniformly for $\sigma \in \mathcal{K}$ $u\in \mathscr{D}'(X)$.
\begin{align}
    \|u\|_{\overline{H}^{s,r_{sc},r_b}_{sc,b}}\leq C\left( \|P(\sigma)u||_{\overline{H}_{sc,b}^{s-1,r_{sc},r_b+1}(X)} + \|u\|_{\overline{H}_{sc,b}^{-N,-N,r_b-\delta}}\right)
\end{align}
and the adjoint estimate 
\begin{align}
    \|v\|_{\dot{H}^{1-s,-r_{sc}-1,-r_b-1}_{sc,b}}\leq C\left(\|P(\sigma)^*v\|_{\dot{H}^{-s,-r_{sc},-r_b}_{sc,b}}+\|v\|_{\dot{H}_{sc,b}^{-N,-N,-1-\delta}}\right).
\end{align}
We conclude that the operator 
\begin{align}
    P(\sigma):\mathcal{X_\sigma:=}\{u\in \overline{H}^{s,r_{sc},r_b}_{sc,b}(X)\;|\; P(\sigma)u\in H^{s-1,r_{sc}+1,r_b+1}_{sc,b}(X)\}\to H^{s-1,r_{sc}+1,r_b+1}_{sc,b}(X))
\end{align}
    is Fredholm when we equip $\mathcal{X}_\sigma$ with its graph norm.
\end{proposition}
\begin{proof}(Sketch)
    The threshold conditions for the orders for $r_{sc}$ and $r_b$ are coming incoming and outgoing conditions on $\tilde{\mathcal{R}}^+_{\sigma,\partial X}$ and $\tilde{\mathcal{R}}^-_{\partial X}$ come from the incoming and outgoing radial $sc-b$ estimates Proposition\ref{prop:scbFredHolmEstimateincoming} and Proposition \ref{prop:ScbOutgoing}. For the horizon radial point estimate at $\mathcal{R}^{\pm}_{\mathcal{H^+}}$ the threshold conditions comes from Proposition \ref{prop:ClassicalHorizonEstimate}, with the surface gravity $\kappa$ coming in via the leading $\alpha_\Lambda(\sigma)=\frac{\text{Im }\sigma}{\kappa}$ from the calculation at (\ref{eq:ThresholdSurfacePressure}). Notice as well that the horizon radial point  estimate of Proposition  \ref{prop:ClassicalHorizonEstimate} forces a differential decay of $s\mapsto s-1$ instead of $s\mapsto s-2$ as allowed by Propositions \ref{prop:scbFredHolmEstimateincoming} and \ref{prop:ScbOutgoing}.
    Finally, standard remainder norm can be improved see \cite{vasy2021limiting} Proposition 4.10, Lemma 4.13, Propositions 4.16 and 4.20. 
    \begin{align}
        \|u\|_{H_{sc,b}^{-N,-N,r_b+1}}\lesssim \|u\|_{H_{sc,b}^{-N,-N,r_b+1-\delta}}.
    \end{align}   
    using a commutator argument for the normal operator 
    \begin{align}
        {}^bN(\rho^{-1}P(\sigma)) = 2i\sigma \left(\rho\partial_\rho - \frac{3}{2}\right),
    \end{align}
    that can be computed from (\ref{eq:MoreRefinedOperatoratInfinity})
    The Fredholm statement then follows from the compact embedding
    \begin{align}\|u\|_{\overline{H}^{s,r_{sc},r_b}_{sc,b}}\hookrightarrow\|u\|_{\overline{H}_{sc,b}^{-N,-N,r_b-\delta}}.
    \end{align}
\end{proof}

Now using Lemma \ref{lem:btoscTranslationLemma} this immediately translates into a Fredholm statement on $b$-spaces.
\begin{corollary}\label{cor:FredholmFiniteSigmabStatement}
Let $\sigma\in \mathcal{K}\Subset \{\sigma \in \mathbb{C}\;|\;\text{Im }\sigma \geq 0, \sigma \neq 0\}$, let $\ell<-\frac{1}{2}$, $s+\ell>-\frac{1}{2}$, and $s>\sup_{\sigma\in\mathcal{K}} \left(\frac{1}{2}-\frac{\text{Im }\sigma}{\kappa}\right)$.
We have
\begin{align}\label{eq:FredholmbStatementFiniteSigma}
    \|u\|_{\overline{H}^{s,\ell}}\leq C_{\mathcal{K}}\left(\|P(\sigma) u\|_{\overline{H}^{s,\ell+1}_b}+\|u\|_{\overline{H}^{-N,\ell-\delta}_b}\right)
\end{align}
and 
\begin{align}\label{eq:FredholmbStatementFiniteSigmaAdjoint}
    \|v\|_{\dot{H}^{-s,-\ell-1}_b} \leq C_{\mathcal{K}}\left(\|P(\sigma)^* v\|_{\dot{H}^{-s,-\ell}_b}+\|v\|_{\dot{H}^{-N,-\ell-1-\delta}}\right).
\end{align}

\end{corollary}

\subsection{Semiclassical analysis}\label{subsec:SemiclassicalAnalysis}

We will now push the micro-local approach one step further and also introduce semiclassical scaling with a parameter $h\in [0,h_0]$ to the problem. This allows us to obtain high energy estimates in a strips with $\text{Im }(\sigma)$ with $\text{Im }\sigma \leq C$ and $\text{Re} (\sigma)>>1$. The operator of interest remains $P(\sigma) =\hat{\Box}_{g_{\bfa,\bfb,\mfrm}}(\sigma)=e^{i\sigma t_*}\circ \Box_{g_{\bfa,\bfb,\mfrm}}\circ e^{-i\sigma t_*}$ on the spatial manifold $X$. A major difference with the the classical propagation we have considered thus far, is that the characteristic set and elliptic sets of $P(\sigma)$ in the resulting semiclassical algebra will also lie at finite semiclassical covectors over the interior $X^\circ$, as well as that the rescaled spectral parameter will play a role at the principal symbol level. As a consequence, the full dynamics of $G$ elucidated in Part I, in particular the trapped set $\Gamma$ will play a role in semiclassical propagation in the interior. \\
\\
Let $h\in [0,h_0)$ with $h_0>0$ small enough and arbitrary we simply define semiclassical vector fields for any $\mathcal{V}_{\bullet}(X)$ through 
\begin{align}
    \mathcal{V}_{\bullet,\hbar}(X):=h C^\infty([0,h_0];\mathcal{V}_\bullet(X)), \; \text{ with } \bullet\in \{b,sc\}.
\end{align}
In local adapted coordinates near the boundary we have for $\bullet\in \{b,sc\}$ that 
\begin{align}
    \mathcal{V}_{sc,\hbar}(X):=\text{span}_{C^\infty([0,h_0]\times X)}\{h\rho^2\partial_\rho,h\rho\partial_{y_i}\}&\text{ and }  \mathcal{V}_{b,\hbar}(X):=\text{span}_{C^\infty([0,h_0]\times X)}\{h\rho\partial_\rho,h\partial_{y_i}\},\\
    &\text{ for } \bullet\in\{sc,b\}
\end{align}
whilst on the interior we simply have 
\begin{align}
    \mathcal{V}_\hbar(X^\circ) = \text{span}_{C^\infty([0,h_0]\times X^\circ)}\{h\partial_{x_i}\}.
\end{align}

We then define the spaces of differential operators through 
\begin{align}
    \text{Diff}^m_{\bullet,\hbar}(X)=\text{span}_{C^\infty([0,h_0]\times X)}\{V_1\circ...\circ V_k\;|\; V_j \in \mathcal{V}_{\bullet,\hbar}(X), 0\leq k \leq m\},\;\bullet\in\{sc,b\}
\end{align}
and of course 
\begin{align}
    \text{Diff}^{m,l}_{\bullet,\hbar}(X)=\rho^{-l} \text{Diff}^m_{\bullet,\hbar}(X),\; \bullet\in\{sc,b\} 
\end{align}
The associated phase-spaces under consideration then become products of the already established ${}^{\bullet}\overline{T^*}X$ with the parameter interval $[0,h_0)_h$, where for most of the rest of the chapter we allow $\bullet\in \{sc,b,(sc,b)\}$ unless otherwise mentioned,
\begin{align}
    {}^{\bullet,\hbar}\overline{T^*}X:={}^\bullet\overline{T^*}X\times [0,h_0)_h, \;\bullet\in\{sc,b,(sc,b)\}.
\end{align}
In the classical approach, one concerns oneself with the limit $h\to 0$, and we define the semiclassical ellipitc, characteristic and operator wavefront sets on the \textit{semiclassical} face $\{h=0\}$\footnote{Not to be confused with the semiclassical hypersurface in Part I.}, which we identify with the corresponding compactified phase space,
\begin{align}
    \left.{}^{\bullet,\hbar}\overline{T^*}X\right|_{h=0}\cong {}^\bullet \overline{T^*}X, \text{ and } \left.{}^{sc,b,\hbar}\overline{T^*}X\right|_{h=0}\cong {}^{sc,b}\overline{T^*}X.
\end{align}

We define the symbol spaces as 
\begin{align}
    S_{\bullet,\hbar}^{m,l}({}^{\bullet,\hbar} \overline{T^*}X)=C^\infty([0,h_0)_h;S^{m,l}_\bullet({}^{\bullet} \overline{T^*}X)).
\end{align}
as well as 
\begin{align}
    S^{m,l_{sc},l_b}_{sc,b,\hbar}({}^{sc,b,\hbar}\overline{T^*}X) =C^\infty([0,h_0)_h;S^{m,l_{sc},l_b}_{sc,b}({}^{sc,b}\overline{T^*}X)),
\end{align}
where we emphasize that we deal with smoothness in $h$ to define the symbols, as before the symbol
\begin{align}
    a(\rho,\xi;h)=2+\sin(\log(\langle \xi\rangle)), 
\end{align}
is smooth in $h$, and bounded towards the boundary $\partial {}^\bullet\overline{T^*}X^\circ$.

We can then define semiclassical pseudodifferential operators $\Psi^{m,\ell}_{\bullet,\hbar}(X)$ through their semiclassical quantizations with cut-off functions $\chi \in C_c^\infty(X)$
\begin{align}
    &\Op_{b,\hbar}(a)u(\rho,y)=(2\pi h)^{-d}\int e^{\frac{i}{h}\log(\rho/\rho'){}^b\xi+(y-y'){}^b\eta}\chi \cdot a(\rho,y,{}^b\xi,{}^b\eta;h) u(\rho',y') \frac{d\rho'}{\rho'}dy' d{}^b\xi d{}^b\eta\\
    &\Op_{sc,h}(a)u(\rho,y)=\\&\;\;\;\;\;\;\;\;\;\;\;\;\;\;\;\;(2\pi h)^{-d}\int e^{\left(\frac{i}{h}\frac{\rho-\rho'}{\rho\rho'}{}^{sc}\xi + \frac{(y-y')\cdot {}^{sc}\eta}{\rho}\right)}\chi \cdot a(\rho,y,{}^{sc}\xi,{}^{sc}\eta;h) u(\rho',y')\frac{d\rho'}{(\rho')^5}dy' d{}^{sc}\xi d{}^{sc}\eta.\\ &\label{eq:SemiclasscalbTXQuantization}
\end{align}

and similarly the quantization for the semiclassical $sc-b$ operators in $\Psi^{m,l_{sc},l_b}_{sc,b,h}(X)$  again through identification of $a_{sc,b,h}\in S^{m,l_{sc},l_b}_{sc,b,h}({}^{sc,b}\overline{T^*}X)$ with $a_{b,h}\in S_b^{m,l_b}({}^b\overline{T^*}X^\circ)$ through composition $a_{b,h}:= a_{sc,b,h}\circ \left(\beta_b|_{T^*X^\circ}\right)^{-1}$ and setting $\text{Op}_{sc,b,h}(a_{sc,b,h}):=\text{Op}_{b,h}(a_{b,h})$ using the quantization in (\ref{eq:SemiclasscalbTXQuantization}). We also note the composition conditions
\begin{align}
    \Psi^{m,l}_{\bullet \hbar}(X)\circ  \Psi^{m',l'}_{\bullet \hbar}(X)&\subset  \Psi^{m+m',l+l'}_{\bullet \hbar}(X)\\
      \Psi^{m,l_{sc},l_b}_{sc,b \hbar}(X)\circ  \Psi^{m',l'_{sc},l'_b}_{sc,b \hbar}(X)&\subset  \Psi^{m+m',l_{sc}+l'_{sc},l_b+l_b'}_{sc,b,\hbar}(X)
\end{align}
as well as semiclassical filtration
\begin{align}
    &\Psi^{m,l}_{\bullet,\hbar}(X)\subseteq \Psi^{m',l'}_{\bullet,\hbar}(X), \text{ if } m\leq m', l\leq l'\\
    &\text{ and } h^{N+1}\Psi^{m,\ell}_{\bullet,\hbar}(X) \subseteq h^N \Psi^{m,\ell}_{\bullet,\hbar}(X),
\end{align}
and similarly for the $sc-b$ calculus, where the filtration in all cases is \textit{descreasing} with respect to powers of $h$.

We then get the following sequences induced by the semiclassical principal symbol map
\begin{align}
    &0\to h\Psi_{b,\hbar}^{m-1,\ell}(X)\to \Psi^{m,\ell}_{b,\hbar}(X)\overset{{}^{b,\hbar}\upsigma_{m,\ell}}{\to} S^{m,\ell}_{b,\hbar}({}^b \overline{T^*}X)/hS^{m-1,\ell}_{b,\hbar}({}^b \overline{T^*}X)\to 0\\
    &0\to h\Psi_{sc,\hbar}^{m-1,\ell-1}(X)\to \Psi^{m,\ell}_{sc,\hbar}(X)\overset{{}^{sc,\hbar}\upsigma_{m,\ell}}{\to} S^{m,\ell}_{sc,\hbar}({}^b \overline{T^*}X)/hS^{m-1,\ell}_{sc,\hbar}({}^b \overline{T^*}X)\to 0\\
    &0\to h\Psi^{m-1,\ell_{sc}-1,\ell_b}_{sc,b,\hbar}(X) \to \Psi^{m,\ell_{sc},\ell_b}_{sc,b,\hbar}(X) \\
    &\;\;\;\;\;\;\;\;\;\;\;\;\;\;\;\;\;\;\;\;\;\;\;\;\;\;\;\overset{{}^{sc,b,\hbar}\upsigma_{m,\ell_{sc},\ell_b}}{\to} S^{m,\ell_{sc},\ell_b}_{sc,b,\hbar}({}^{b,sc}\overline{T^*}X)/hS^{m-1,\ell_{sc}-1,\ell_b}_{sc,b,\hbar}({}^{b,sc}\overline{T^*}X)\to 0,
\end{align}
with the semiclassical Hamiltonian vector field understood as
\begin{align}
    {}^{sc,b,\hbar}\upsigma_{m+m'-1,l+l'-1,l_b+l'_b}\left(\frac{i}{h}[P,A]\right)=H_pa,
\end{align}
for $P=\text{Op}_{sc,b,\hbar}(p) \in S^{m,l_{sc},l_b}_{sc,b,\hbar}({}^{sc,b,h}\overline{T^*}X)$ and $A=\text{Op}_{sc,b,\hbar}(a) \in S^{m',l'_{sc},l'_b}_{sc,b,\hbar}({}^{sc,b,h}\overline{T^*}X)$

Applying principal-order normalization to the remainder classes allows us to identify the following joint semiclassical loci $\mathcal{Z}_{\bullet,\hbar} \subset{}^{\bullet,\hbar}\overline{T^*}X$, which is exactly where lower order remainders disappear upon applying the semiclassical principal symbol 

\begin{align}\label{eq:bandscatteringSemiclassicalLoci}
\mathcal{Z}_{
\bullet,\hbar}&=(\mathcal{Z}_\bullet \times [0,h_0)_h)\,\cup\,({}^\bullet \overline{T^*}X\times \{0\}_h), \text{ for } \bullet\in\{sc,b,(sc,b)\}
\end{align}
and in the $sc-b$ case we have in particular 
\begin{align}\label{eq:bscSemiclassicalLocus}
    \mathcal{Z}_{sc,b,\hbar}=((\text{fi}\,\cup \, \text{sf})\times [0,h_0))\,\cup ({}^{sc,b,\hbar}\overline{T^*}X \times \{0\}_h).
\end{align}

Now given a semiclassical symbol $a_{\bullet,\hbar} \in S^{m,l}_{\bullet,\hbar}({}^{\bullet,\hbar}\overline{T^*}X)$ we define again its normalized symbol $\tilde a_{\bullet,\hbar}$ w.r.t. to the loci (\ref{eq:bandscatteringSemiclassicalLoci}) as
\begin{align}
    C^\infty([0,h_0); \mathcal{A}^{0,0}({}^{\bullet}\overline{T^*}X))\ni \tilde a_\bullet(\rho,y_i,{}^\bullet\xi,{}^\bullet\eta_i;h):= \varrho^m_\bullet\rho^l a(\rho,y_i,{}^\bullet\xi,{}^\bullet\eta_i;h).
\end{align}

and similarly for $a_{sc,b,\hbar}\in S^{m,l_{sc},l_b}_{sc,b,\hbar}({}^{sc,b,\hbar}\overline{T^*}X)$ we define 
\begin{align}
    C^\infty([0,h_0); \mathcal{A}^{0,0,0}({}^{sc,b,\hbar}\overline{T^*}X))\ni \tilde a_{sc,b,\hbar}= \varrho_{\text{fi}}^m\rho^{l_{sc}}_{\text{sf}}\rho_{\text{ff}}^{l_b}a_{sc,b,\hbar}.
\end{align}

In particular when symbols $a_{b,\hbar}, a_{sc,\hbar}$ and $a_{sc,b,\hbar}$ have alternative representatives of the same semiclassical symbol, $a'_{b,\hbar},a_{sc,\hbar}'$ and $a'_{sc,b,\hbar}$ we have that their normalized symbols satisfy 
\begin{align}
    \tilde a_{b,\hbar}-\tilde a'_{b,\hbar}&\in h\varrho_bC^\infty([0,h_0)_h;\mathcal{A}^{0,0}({}^b\overline{T^*}X))\\
     \tilde a_{sc,\hbar}-\tilde a'_{sc,\hbar}&\in h\varrho_b\rho C^\infty([0,h_0)_h;\mathcal{A}^{0,0}({}^b\overline{T^*}X))\\
     \tilde{a}_{sc,b,h}-\tilde a_{sc,b,h}'&\in h\varrho_{\text{fi}}\rho_{\text{sf}}\rho_{\text{ff}}C^\infty([0,h_0)_h;\mathcal{A}^{0,0,0}({}^{sc,b}\overline{T^*}X))
\end{align}

thus giving the expression of the loci $\mathcal{Z}_{\bullet,\hbar}$ and $\mathcal{Z}_{sc,b,\hbar}$ as 
\begin{align}
    \mathcal{Z}_{b,\hbar}&=\{h\varrho_b =0\}\\
    \mathcal{Z}_{sc,\hbar}&=\{h\varrho_{sc}\rho=0\}\\
    \mathcal{Z}_{sc,b,\hbar}&=\{h\varrho_{\text{fi}}\rho_{\text{sf}}=0\}.
\end{align}

Whilst we can define semiclassical elliptic, characteristic and operator wavefront sets on the loci given above, we are only interested in the limit $h\to 0$ so we define them on the intersection of the joint semiclassical locus and the hypersurface $\{h=0\}$, which can be identified with the original bundle ${}^\bullet\overline{T^*}X$ as 
\begin{align}
    \mathcal{Z}_{\bullet,\hbar}\;\cap\;\{h=0\}={}^\bullet \overline{T^*}X|_{h=0}\cong {}^\bullet \overline{T^*}X,\; \bullet\in\{sc,b,(sc,b)\}.
\end{align}

This has as a consequence that on the interior elliptic/characteristic sets of semiclassical operators in general can no longer be viewed as conic subsets or subsets over fibre infinity ${}^\bullet S^*X^\circ$ but can instead run over the entirety of phase space ${}^\bullet \overline{T^*}X$. Observe that away from the semiclassical face we have that the joint locus instead reduces to the classical one 
\begin{align}
    \mathcal{Z}_{\bullet,\hbar}\,\cap\, \{h = c\}= \mathcal{Z}_\bullet \times\{c\}\cong \mathcal{Z}_\bullet,\; c>0.
\end{align}
Moreover we obviously also still have that the classical locus $\mathcal{Z}_\bullet$ is a subset of the semiclassical locus, i.e.
\begin{align}\label{eq:ClassicalLocusInclusion}
    \mathcal{Z}_{\bullet}\cong \mathcal{Z}_{\bullet}\times \{0\}_h \subset \mathcal{Z}_{\bullet,\hbar}\,\cap\, \{h=0\}
\end{align}
and similarily for the sc-b locus. We then have the following expressions for elliptic, characteristic sets, and semiclassical wavefront sets for $\bullet\in \{sc,b,(sc,b)\}$, let $A=(A_h)_{h\in[0,h_0)}$ be a semiclassical operator family of order 
\begin{align}
    A\in \begin{cases}
        \Psi^{m,l}_{\bullet,h}(X), &\;\bullet\in\{b,sc\}\\
        \Psi^{m,l_{sc},l_b}_{sc,b,h}(X), &\;\bullet\in (sc,b),
    \end{cases}
\end{align}
choose a full symbol 
\begin{align}
    a_{\bullet,h}\in \begin{cases}
        S^{m,l}_{\bullet,\hbar}({}^{\bullet,h}\overline{T^*}X)&\;\bullet\in \{b,sc\}\\
        S^{m,l_{sc},l_b}_{\bullet,\hbar}({}^{sc,b,h}\overline{T^*}X)&\;\bullet \in\{(sc,b)\}
    \end{cases}
\end{align}
such that $A=\text{Op}_{\bullet,\hbar}(a)$ and denote its corresponding normalized symbol by 
\begin{align}
    \tilde{a}_{\bullet,\hbar}=\begin{cases}
        \varrho^m_\bullet\rho^l a_{\bullet,h},&\;\bullet\in \{sc,b\}\\
        \varrho^m_\bullet\rho^{l_{sc}}\rho^{l_b}a_{sc,b,h}, &\;\bullet\in (sc,b).
    \end{cases}
\end{align}

We then define the semiclassical elliptic set for $\bullet\in\{sc,b,(sc,b)\}$ as 
\begin{align}
    \Ell_{\bullet,\hbar}(A):=\{q\in& {}^\bullet\overline{T^*}X\;|\; \exists U\ni q \text{ open in } {}^{\bullet}\overline{T^*}X\\ 
    &\exists c>0,\; 0<h_1<h_0,\, |\tilde a_{\bullet,\hbar}(q';h)|\geq c, \forall (q',h)\in U\cap {}\overline{T}^*X^\circ \times (0,h_1)\}
\end{align}\footnote{Note that we can interpret $q\in {}^\bullet\overline{T^*}X$ as $(q,0) \in {}^\bullet\overline{T^*}X|_{h=0}$, but we supress this notation. We also supress order superscripts.}
The characteristic set is then defined as the complement of the elliptic set inside ${}^\bullet \overline{T^*}X$:
\begin{align}
    \chr_{\bullet,\hbar}(A) = {}^{\bullet}\overline{T^*}X\backslash \Ell_{\bullet,\hbar}(A),
\end{align}
which without supressing the identification with the semiclassical face $\{h=0\}$ should be seen as 
\begin{align}
    \chr_{\bullet,\hbar}(X) = ({}^{\bullet,\hbar}\overline{T^*}X|_{h=0})\backslash \Ell_{\bullet,\hbar}(A).
\end{align}
Before defining the wavefront sets we identify the spatially residual symbol and operator classes, and their semi classically residual subclasses, which additionally decay rapidly in $h$. 
\begin{align}
    h^\infty\Psi_{b,\hbar}^{-\infty,l_b}(X)&= \{A_h=\text{Op}_{b,\hbar}(a_h)\;|\; a_h\in h^\infty S^{-\infty,l_b}_{b,\hbar}({}^{b,h}\overline{T^*}X)\}\\
    h^\infty\Psi_{sc,\hbar}^{-\infty,-\infty}(X)&= \{A_h=\text{Op}_{sc,\hbar}(a_h)\;|\; a_h\in h^\infty S^{-\infty,-\infty}_{sc,\hbar}({}^{sc,h}\overline{T^*}X)\}\\
     h^\infty\Psi_{sc,b,\hbar}^{-\infty,-\infty,l_b}(X)&= \{A_h=\text{Op}_{sc,b,\hbar}(a_h)\;|\; a_h\in h^\infty S^{-\infty,-\infty,l_b}_{sc,b,\hbar}({}^{sc,b,h}\overline{T^*}X)\}    
\end{align}
where we identified residual symbols 
\begin{align}
     h^\infty S^{-\infty,l_b}_{b,\hbar}({}^{b,h}\overline{T^*}X)&:=\bigcap_{\substack{N\in\mathbb{N}_0\\m\in\mathbb{R}}}h^NS^{m,l_b}_{b,\hbar}({}^{b,h}\overline{T^*}X)\\
     h^\infty S^{-\infty,-\infty}_{sc,\hbar}({}^{sc,h}\overline{T^*}X)&:=\bigcap_{\substack{N\in\mathbb{N}_0\\m\in\mathbb{R}\\\ell_{sc}\in\mathbb{R}}}h^NS^{m,l_{sc}}_{sc,\hbar}({}^{sc,h}\overline{T^*}X)\\
     h^\infty S^{-\infty,-\infty,\ell_b}_{sc,b,h}({}^{sc,b,h}\overline{T^*}X)&:=\bigcap_{\substack{N\in\mathbb{N}_0\\m\in\mathbb{R}\\\ell_{sc}\in\mathbb{R}}}h^NS^{m,l_{sc},l_b}_{sc,b,\hbar}({}^{sc,b,h}\overline{T^*}X).
\end{align}

We we then define then the complement semiclassical essential support of a symbol as the set 
\begin{align}
    &{}^{\bullet}\overline{T^*}X\backslash \esssup_{\bullet,\hbar}a\\ 
    &= \{q\in {}^{\bullet}\overline{T^*}X\;|\; \exists 0<h_1<h_0, \chi \in C^\infty({}^\bullet\overline{T^*}X), \chi \equiv 1 \text{ near } q \text{ s.t. } \chi a_{\bullet,h} \text{ is residual for } h\in [0,h_1) \} 
\end{align}
where $\chi a_{\bullet,h} $ being residual should be interpreted as 
\begin{align}
    \chi a_{\bullet,h}\in\begin{cases}
        h^\infty S_{b,\hbar}^{-\infty,l}({}^{b,h}\overline{T^*}X),&\;\bullet=b\\
        h^\infty S_{sc,\hbar}^{-\infty,-\infty}({}^{sc,h}\overline{T^*}X),&\;\bullet=sc\\
        h^\infty S_{sc,b,h}^{-\infty,-\infty,\ell_b}({}^{sc,b,h}\overline{T^*}X),&\;\bullet=sc,b.
    \end{cases}
\end{align}

We then define the operator wavefront set of a pseudodiffernetial operator $A_h=\text{Op}_{\bullet,\hbar}(a_{\bullet,\hbar})$ through 
\begin{align}
    \WF_{\bullet,\hbar}'(A_h)=\operatorname{ess\,sup}_{\bullet,\hbar}(a).
\end{align}

\begin{definition}[Semiclassical Sobolev spaces]
    Let $\bullet\in \{b,sc\}$. For $s\in \mathbb{N}_0$ let 
    \begin{align}
        H^s_{\bullet,h}(X) :=H^s_\bullet(X),
    \end{align}
    equipped with the norm 
    \begin{align}
        \|u\|^2_{H^s_{\bullet,h}}:= \sum_{j=0}^sh^{2j}\sum_{k=0}^{j}\|A_{jk} u\|_{L_{sc}^2(X)} \text{ where } \span \{A_{jk}\}_{k\leq j} = \Diff^j_{\bullet}(X),\forall j\leq s,
    \end{align}
    And then define weighted semiclassical spaces as 
    \begin{align}
        H^{s,\ell}_{\bullet,h}:=\rho^\ell H^{s}_{\bullet,h}, \;\; \|u\|_{H^{s,\ell}_{\bullet,}} =\|\rho^{-\ell} u\|_{H^s_{\bullet,h}}
    \end{align}
    Equivalently for $s\geq 0$ we can choose $\Lambda^s_{\bullet,h}\in \Psi^{s,0}_{\bullet,h}(X)$ and define an equivalent graph norm 
    \begin{align}
        \|u\|_{H^{s,\ell}_{\bullet,h}}&=\|\Lambda^s_{\bullet,h}(\rho^{-\ell} u)\|^2_{L^2_{sc}(X)} +\|\rho^{-\ell}u\|_{L^2_{sc}(X)}, \text{ for } s\geq 0\\
        \|u\|_{H^{s,\ell}_{\bullet,h}}&=\|\Lambda^s_{\bullet,h}(\rho^{-\ell} u)\|^2_{L^2_{sc}(X)}, \text{ for } s<0.
    \end{align}
    
    The latter method to $sc-b$ Sobolev spaces as well. Let 
    \begin{align}
        \tilde{s}\leq \min \{s,\ell_{sc}-\ell_b\},
    \end{align}
    with $A_h \in \Psi^{s,\ell_{sc}-\ell_b,0}_{sc,b,\hbar}(X)$ and let $A_h\in \Psi^{s,\ell_{sc}-\ell_b,0}_{sc,b,h}(X)$ be semiclassically elliptic on the elliptic face. We define  
    \begin{align}
        H^{s,\ell_{sc},\ell_b}_{sc,b,h}(X):=\{u\in H^{\tilde s,\ell_b}_{b,h}(X)\;|\; A_hu \in \rho^{\ell_b}L^2_{sc}(X)\},
    \end{align}
    equipped with the norm 
    \begin{align}
        \|u\|_{H^{s,\ell_{sc},\ell_b}_{sc,b,h}}:=\|\rho^{-\ell_b}A_h u\|_{L^2_{sc}}^2+\|u\|^2_{H^{\tilde s,r_b}_{b,h}}.
    \end{align}
\end{definition}
\begin{remark}
    Analogously this construction can be generalized for extendible and supported semiclassical distributions as well $\overline{H}^{s,r}_{\bullet,h}(X)$ and $\dot{H}^{s,\ell}_{\bullet,h}(X)$.\\
    We say a family of distributions $(u_h)_{0<h<h_0}$ is bounded in $h^k H^{s,\ell}_{\bullet,h}(X)$ if 
    \begin{align}
        \sup_{0<h<h_0} h^{-k}\|u_h\|_{H^{s,\ell}_{\bullet,h}(X)}<\infty.
    \end{align}
\end{remark}

We have indeed that given $A_h\in \Psi^{m,l}_{\bullet,\hbar}(X)$ that 
\begin{align}
    A_h:H^{s,\ell}_{\bullet,h}\to H^{s-m,\ell-l}_{\bullet,h}(X)
\end{align}
defines a continuous map when $X$ is compact, uniformly for sufficiently small $h>0$, and a similar statement holds for $A_h\in \Psi^{m,l_{sc},l_b}_{sc,b,\hbar}(X)$ for which 
\begin{align}
    A_h:H^{s,\ell_{sc},\ell_b}_{sc,b,h}\to H^{s-m,\ell_{sc}-l_{sc},\ell_b-l_b}_{sc,b,h}(X).
\end{align}
We also have the equivalent of (\ref{eq:TranslationtobSpaces}), 
\begin{align}
    H_{b,h}^{s,\ell_b}(X)=H^{s,s+\ell_b,\ell_b}_{sc,b,h}(X), \label{eq:TranslationtoSemiClasscialbSpaces}
\end{align}
so we may expect semiclassical Fredholm estimates on sc-b spaces to translate over into semiclassical Fredholm estimates on $b$-spaces.
When it comes to the actual dynamics we will separate dynamics on the interior $X^\circ$, and on the boundary $\partial_+X$. On $X^\circ$ we work on standard semiclassical Sobolev spaces $H^{s}(X^\circ)$ which can be analogously defined as in our earlier approaches using Remarks \ref{rem:RestrictiontoOpenInterior} and \ref{rem:RestrictionPDOStoOpenInterior}. On the boundary $\partial_+X$, where we will work immediately in the semiclassical sc-b formalism. 
When it comes to semiclassical dynamics on $T^*X^\circ|_{h=0}$, let
\begin{align}
    P\in \Psi^m_\hbar(X^\circ) \text{ and } A\in \Psi^{m'}_\hbar(X^\circ),
\end{align},  
with real valued representatives 
\begin{align}
    p= p(x,\xi;h)\in S^m_\hbar({}^h\overline{T^*}X^\circ), \text{ and } a =a(x,\xi ;h)\in S^{m'}_\hbar({}^h\overline{T^*}X^\circ),
\end{align}
meaning
\begin{align}
    P-\text{Op}_{\hbar}(p)\in h\Psi_\hbar^{m-1}(X^\circ) \text{ and } A-\text{Op}_\hbar(a)\in h\Psi^{m'-1}_\hbar(X^\circ)
\end{align}
We then have the following commutation relation, 
\begin{align}
    \frac{i}{h}[P,A]\in \Psi_\hbar^{m+m'-1}(X^\circ),
\end{align}
with principal symbol
\begin{align}
    {}^\hbar\upsigma_{m+m'-1}\left(\frac{i}{h}[P,A]\right)=H_{p(x,\xi;h)}a(x,\xi;h)
\end{align}
where for fixed $h>0$,
\begin{align}
    H_{p(x,\xi;h)}a(x,\xi;h)=\sum_j \frac{\partial p}{\partial\xi_j}\frac{\partial a}{\partial x_j}-\frac{\partial a}{\partial\xi_j}\frac{\partial p}{\partial x_j}.   
\end{align}

We have that for each $h\in [0,h_0)$, $H_{p_h}$ (we surpress $(x,\xi)$ notation) defines characteristic curves $\gamma_h \in T^*X\times \{h\}$ associated to $p$ statsifying 
\begin{align}
 \dot{\gamma}_h(t) = H_{p_{h}}(\gamma_h(t)),\label{eq:SemiclassicalBicharacteristics}
\end{align}
 and indicated before we are interested in the semiclassical face where the semiclassical microlocal locus is defined i.e. taking the limit $h \searrow 0$. For the purposes of the interior semiclassical dynamics we once again assume that for $P\in \Psi_{\hbar}^m(X^\circ)$ its normalized representative  $\tilde p(x,\xi;h) := \varrho^m p(x,\xi;h)$ satisfies $\tilde p \in C^\infty([0,h_0)\times {}\overline{T^*}X^\circ)$
 such that 
 \begin{align}
     \hat{H}_{p_h}:=\varrho^{m-1}H_{p_h},
 \end{align}
extends extends smoothly to ${}^hS^*X$ and in particular to fibre infinity at the semiclassical face ${}^hS^*X|_{h=0}$, such that for $p_0=p_h|_{h=0}$
\begin{align}
    \hat{H}_{p_h}|_{h=0}=\hat{H}_{p_0}=\varrho^{m-1}H_{p_0}.
\end{align}

In this case we have again that
\begin{align}
    \chr_\hbar(P) =\{q\in \overline{T^*}X^\circ\;|\; \tilde p (q;0)=0\}, \Ell_\hbar(P)=\{q\in \overline{T^*}X^\circ\;|\; \tilde p (q;0)\neq0\}
\end{align}
and the null-bicharacteristics $\gamma_h$ of $\hat{H}_{p_h}$ satisfying (\ref{eq:SemiclassicalBicharacteristics}) and $p_h\circ \gamma_h=0$.

We can now set up elliptic and propagation estimates on the interior $X^\circ$, $P\in \Psi^m_\hbar(X^\circ)$ are uniformly properly supported. Assume that all cutuffs $B_h,G_h, E_\hbar\in \Psi^0_\hbar(X^\circ)$ have properly supported Schwarz-kernels supported in a fixed compact $K\times K\subset X^\circ \times X^\circ$. Choose a compact $\tilde K\supset K$ such that 
\begin{align}
    \supp K_{P_h}\cap (K\times X^\circ) \subset K\times \tilde K,
\end{align}
and all microlocal sets are viewed on 
\begin{align}
    {}^\hbar\overline{T^*}X^\circ|_{h=0}\cong \overline{T^*}X^\circ.
\end{align}

We then have an elliptic estimate followed by a propagation estimate.
\begin{proposition}[Semiclassical Elliptic Estimate]
    Suppose $B,G\in \Psi_\hbar^0(X^\circ)$ satisfy 
    \begin{align}
        \WF'_\hbar(B)\subset \Ell_\hbar(P) \cap\Ell_\hbar(G).
    \end{align}
    Then for every $s\in \mathbb{R}$ and $N>0$, there exist $C,h_1>0$ such that 
    \begin{align}
        \|B_h u\|_{H^s_h(X)}\leq C\left(\|G_hP_hu\|_{H^{s-m}_h(X)}+h^N\|\chi u\|_{H_h^{-N}}\right),\; 0<h<h_1.
    \end{align}
\end{proposition}

\begin{proposition}[Semiclassical Interior Propagation]
    Suppose the principal-symbol class of $P$ satisfies a real valued representative $p_h$ satisfying 
    \begin{align}
        P_h -\text{Op}_h(p_h)\in h\Psi^{m-1}_\hbar(X^\circ),\; \tilde p_h =\varrho^m p_h\in C^\infty([0,h_0)\times \overline{T^*}X^\circ),
    \end{align}
and set 
\begin{align}
    p_0=p_h|_{h=0}, \;\hat{H}_{p_0}=\varrho^{m-1}H_{p_0}.
\end{align}
Let $B,G,E\in \Psi^0_\hbar(X^\circ)$ satisfy 
\begin{align}
    \WF'_\hbar(B)\subset \Ell_\hbar(G).
\end{align}
Assume every integral curve $\gamma$ of $\hat{H}_{p_0}$ starting in
\begin{align}
    \WF'(B)\cap \chr_\hbar(P)
\end{align}
enters $\Ell_\hbar(E)$ in finite time while remaining in $\Ell_\hbar(G)$. Then for every $s\in \mathbb{R}$, $N>0$ there exist $C, h_1>0$ such that 
\begin{align}
\|B_hu\|_{H^s_h}\leq C\left(\|E_hu\|_{H^s_h}+h^{-1}\|G_hP_hu\|_{H^{s-m+1}_h}+ h^N\|\chi u\|_{H_h^{-N}}\right),\;0<h<h_1.
\end{align}

\end{proposition}
\begin{remark}
    In both propositions we have that they hold for $u\in \mathscr{D}'(X^\circ)$ in the strong sense, finiteness of the right-hand side implies finiteness of the left-hand side.
\end{remark}

When considering the wave type equations in the semiclassical setting, i.e. for $P(\sigma) =\hat{\Box}_g(\sigma)$ with $g$ some spacetime metric, we set for $\sigma \neq 0$ and introduce the semiclassical family through
\begin{align}
    h:=|\sigma|^{-1}, \text{ and }  z=h\sigma =\frac{\sigma}{|\sigma|},
\end{align}
so that 
\begin{align}
    \sigma = h^{-1}z \text{ and } |z|=1,
\end{align}
and we define the semiclassical wave operator as
\begin{align}
    P_{h,z}:= h^2P(h^{-1}z).
\end{align}
More precisely if 
\begin{align}
    \Box_g= \sum_{j+|\beta|\leq2} a_{j\beta}(x) D^j_{t_*}D^\beta_x
\end{align}
then 
\begin{align}
    P_{h,z}=\sum_{j+|\beta|\leq 2} h^{2-j-|\beta|}c_{j\beta}(x)(-z)^j(hD_x)
\end{align}
using the translation $D_{t_*}=i^{-1}\partial_{t_*}\mapsto -\sigma\mapsto -h^{-1}z$. 
For a forced wave equation we have 
\begin{align}
    P(\sigma)v =f \iff P_{h,z}v=h^2f.
\end{align}
Let us illustrate this first on the Minkowski spacetime, using $g=\diag\{-1,+1,+1,+1,+1\}$ and using (\ref{eq:LaplaceBeltrami}) we have
\begin{align}
    \Box_{g_{\mathbb{M}^{1+4}}}=\partial_t^2 +\Delta_{\mathbb{R}^4},
\end{align}
so 
\begin{align}
    \hat{\Box}_{g_{\mathbb{M}^{1+4}}}(\sigma) = \Delta_{\mathbb{R}^4}-\sigma^2
\end{align}
and hence in this case
\begin{align}
    P_{h,z}=h^2\Delta_{\mathbb{R}^4}-z^2, p_{h,z}(x,\xi)=\|\xi_i\|_{g^{-1}_{\mathbb{R}^4}}^2-z^2.
\end{align}
We immediately observe that in this case, for real $z\in \mathbb{S}^1$ (i.e. $z=\pm 1$) the semiclassical characteristic set of $P_{h,z}$ extends through the interior of $T^*\mathbb{R}^4$ at finite covectors of radius $1$.
\begin{align}
    \chr_\hbar (P_{z,h})=\{(x,\xi)\in \overline{T^*}\mathbb{R}^4\;|\|\xi_i\|_{g^{-1}_{\mathbb{R}^4}}=1\}
\end{align}
Notice in particular that near fibre-infinity, the operator is semiclassically elliptic as $\|\xi_i\|_{g^{-1}_{\mathbb{R}^4}}>>z$, or more formally, choosing $\varrho=\hat{\xi_i}=\frac{1}{|\xi_i|}$ in a projective chart $\{|\xi_i|\gtrsim|\xi_j|, \forall j\neq i$ we have that 
\begin{align}
    \tilde p_{z_0} :=\varrho^2 p_{z_0}= 1+\sum_{j\neq i}|\hat{\xi}_j|^2-\varrho^2z^2,
\end{align}
thus 
\begin{align}\label{eq:SemiClassicalFibreInfinityLimit}
    \tilde{p}_z|_{\varrho=0}=1+\sum_{j\neq i}\hat{\xi}_j^2\geq 1.
\end{align}

Equivalently, in forward coordinates $t_*=t-r$, $g_{\mathbb{M}^{1+4}}=-dt_*^2-2dt_*dr + r^2g_{\mathbb{S}^3}$ and inverting in $\rho=r^{-1}$ in example (\ref{eq:MinkowskiModel}) or (\ref{eq:MinkowskiScatterOperator}) we can write semiclassical scattering quantizations as 
\begin{align}
     P_{b,h,z}&=\rho^2((h\rho D_\rho)^2+2ih( h\rho D_\rho)+h^2\Delta_{\mathbb{S}^3})-2z\rho\left((h\rho D_\rho) +\frac{3ih}{2}\right)\\
     p_{b,h,z}&=\rho^2({}^b\xi^2+\|\beeta_i\|^2_{g^{-1}_{\mathbb{S}^3}}) -2z\rho
\end{align}
and in its scattering form as 
\begin{align}
    P_{h,z}&= (h\rho^2 D_\rho)^2+h^2\Delta_{\mathbb{S}^3}-2z(h\rho^2D_\rho)+3ih\rho\left((h\rho^2 D_\rho) -z\right)\\
    p_{h,z,sc}&=\scxi^2+\|\sceta_i\|^2_{g^{-1}_{\mathbb{S}^3}}-2z\scxi \label{eq:SemiclassicalScatteringSymbol}
\end{align}
where we immediately note the analogy between (\ref{eq:SemiclassicalScatteringSymbol}) and (\ref{eq:ScatteringMinkowskiSymbol}) for $z=1$.
We want to eventually analyze $P(\sigma)=\hat{\Box}_{g_{\bfa,\bfb,\mfrm}}(\sigma)$ in high energy strips with 
\begin{align}
    |\text{Re }\sigma|\to \infty, \text{ and } \text{Im }\sigma \in I\Subset\mathbb{R},
\end{align}
with $h=|\sigma|^{-1}$, $z=h\sigma$, and $\nu=\text{Im }\sigma$ restricted to a fixed compact interval we write 
\begin{align}
    z=z_0+hz_1, \;\;z_1:=\frac{z-z_0}{h}=i\nu+\mathcal{O}(h).
\end{align}

To this end we need to be able to reduce our discussion to using a real valued semiclassical representative $p_{z_0,h}$ of our operator $P_{z,h}$. Let us therefore define
\begin{align}
    z_0 :=\sgn (\text{Re }\sigma) \in \{-1,1\} \text{ and  } \nu=\text{Im }\sigma
\end{align}
and let $\nu$ be bounded. 
We then have that 
\begin{align}
    z=z_0\sqrt{1-h^2\nu^2}+ih\nu=z_0+ih\nu +\mathcal{O}(h^2).
\end{align}

Thus a straightforward computation easily shows that for any of the previous representatitive in the examples given we have that 
\begin{align}
    p_{z,h} - p_{z_0,h}\in h S^{1}(\overline{T^*}X^\circ),
\end{align}
and thus that for $P_{h,z_0}:=\text{Op}_h(p_{z_0})$ we have that
\begin{align}
    P_{h,z}-P_{h,z_0}\in h\Psi^1_h(X^\circ).
\end{align}

To make the notions of the above fully rigorous, let $(M^\circ,g)$ be be a stationary space-time satisfying the Einstein vacuum equations and $X^\circ$ be its associated spatial manifold. To provide a fully rigorous link between the classical, semiclassical flow on $\overline{T^*}X^\circ$ and the null geodesic flow on $\overline{T^*}M$ that we considered in Part I, consider the embedding fixing $t_*=0$ and $\sigma=c$,
\begin{align}
    \iota_c:\overline{T^*}X^\circ \to \overline{T^*}M^\circ,\; \iota_c(x,\xi)=(0,x;c,\xi),
\end{align}
i.e. it identifies the semiclassical face $T^*X^\circ \cong {}^\hbar T^*X^\circ|_{h=0}$ at $\{t_*=0,\sigma=c\}$. It holds 
\begin{align}
    j:= \iota_c|_{S^*X^\circ}, j(S^*X^\circ)= S^*M^\circ \cap \{t_*=0,\hat\sigma=0\}
\end{align}
defines the same boundary embedding $\iota_c(S^*X^\circ)$, for any $c\in \mathbb{R}$.
Then write
\begin{align}
    G_c:= \iota_c^* G= G(cdt_*+\upalpha_X),\;\;\tilde{G}_c:= \iota_c^* \tilde{G}= \tilde{G}(cdt_*+\upalpha_X)
\end{align}
A more tedious computation then shows in general that for $(M,g_{\bfa,\bfb,\mfrm})$ with spatial manifold $X$, a wave-type equatoin
\begin{equation}
    P= \Box_g =-\frac{1}{\sqrt{-\det(g)}}\partial_\mu (\sqrt{-\det g}g^{\mu\nu}\partial_\nu)
\end{equation}
has a real representative of its semiclassical principal symbol on the interior given by pullback to $T^*M^\circ$
\begin{align}
    {}^\hbar\upsigma_2(P_{h,z}) = \tilde G_{-{z_0}}, \;z_0\in \{-1,+1\}
\end{align}
on $T^*X^\circ$,
or for the rescaled $L_{h,z_0}=\bfr^2 P_{h,z_0}$ we have 
\begin{align}
    {}^\hbar\upsigma_2(L_{h,z_0})=G_{-z_0},\; z_0\in \{-1,+1\}.
\end{align}
Hence analyzing the Hamiltonian dynamics of $H_{p_{h,z_0}}$ on the characteristic set at finite covectors, 
\begin{align}
    \Sigma_{-z_0}:=\chr(P_{h,z_0}) \cap T^*X^\circ
\end{align}
coincides with analyzing the dynamics of $H_{G}$ (or $H_{\tilde G}$ for the unscaled problem) on the semiclassical slice $\{\sigma=-z_0\}$ at $t_*=0$
\begin{align}
\iota_{-z_0}(\Sigma_{-z_0})=\Sigma_{G}\cap \{t_*=0,\sigma=-z_0\}\subset T^*M^\circ \text{ with } z_0\in \{-1,+1\}.
\end{align}
where we have 
\begin{align}
    (\iota_{c})_* (H_{p_{h,z_0}}|_{\Sigma_{-z_0}})= H_G-(\partial_\sigma G)\partial_{t_*}|_{\{t_*=0,\sigma=-z_0\}}
\end{align}
Note that because the metric $g$ is stationary, we have that $H_G\sigma = -\partial_tG=0$ and thus we could have $\iota_c$ embed into any other constant $t_*$- hypersurface inside $T^*M$, the choice $t_*=0$ is arbitrary, however $H_G$ is not tangent to $\{\sigma=c,t_*=c'\}$ in general.
When extending $p_{z_0,h}$ to fibre infinity $\partial \overline{T}^*X^\circ=S^*X^\circ$ the terms of the symbol $p_{h,z_0}$ containing $z_0$ vanish in the limit $\varrho\to 0$ similar to the computation done in the simple case (\ref{eq:SemiClassicalFibreInfinityLimit}), i.e.
defining the renormalized symbol as $\tilde{p}_{h,z_0}:=\varrho^2 p_{h,z_0}$, we have
\begin{align}
    \tilde p_{h,z_0}|_{S^*X^\circ}&= \lim_{\varrho\to 0}(\varrho^2\tilde G_{-z_0})= \lim_{\varrho\to 0}(\varrho^2\tilde G(-z_0dt_*+\upalpha_X)\\
    &=\lim_{\varrho\to 0}(\tilde G(-z_0\varrho dt_*+ \varrho\upalpha_X)\\
    &=\tilde{G}(0 dt_* +\hat\upalpha_X)= (\varrho^2\tilde G_0)|_{S^*X^\circ}.
\end{align}
where we interpret $\tilde G$ again as a function quadratic functin $T^*X^\circ$ and introduced $\hat\upalpha_X$ as the normalized Liouville $1$-form on the spatial manifold $\overline{T^*}X$ near fibre infinity, from which we conclude that
\begin{align}
    {}^\hbar\upsigma_2(P_{h,z_0})|_{\partial\overline{T^*}X^\circ}=\tilde G_0.
\end{align}
The classical elliptic and characteristic sets agree with the semiclassical ones 
\begin{align}
    \chr_\hbar(P_{h,z_0})\cap S^*X^\circ = \chr(P(\sigma)) \cap S^*X^\circ, \;     \Ell_\hbar (P_{h,z_0})\cap S^*X^\circ = \Ell(P(\sigma)) \cap S^*X^\circ.
\end{align}
In particular then, analyzing the semiclassical flow over $S^*X^\circ$ is equivalent to analyzing the semiclassical flow of $P_{0,h}$ at $S^*X^\circ$. We have that

and hence
\begin{align}
    j(\chr_\hbar(P_{h,z_0})\cap S^*X^\circ) = j(\chr_\hbar (P_{h,0})\cap S^*X^\circ) = \overline{\Sigma_G}\cap S^*M^\circ \cap \{t_*=0,\hat\sigma =0\},
\end{align}

Working then with the familiar rescaled Hamiltonian vector field $H_{p_{z_0}}$ we have 
\begin{align}
    \varrho H_{p_{z_0}}|_{S^*X^\circ}=\varrho H_{\tilde G_0}|_{S^* X^\circ},
\end{align}
or equivalently
\begin{align}
    \varrho H_{G_{-z_0}}|_{S^*X^\circ}=\varrho H_{G_0}|_{S^*X^\circ}.
\end{align}

Thus for the full semiclassical dynamics on $\overline{T^*}X^\circ$ and the corresponding semiclassical estimates we identify few key regions of importance; we have
\begin{enumerate}
        \item \textbf{The characteristic set at fibre infinity.} Write $\overline{\Sigma}_{-z_0}$ for the extension of $\Sigma_{-z_0}$ to $\overline{T^*}X^\circ$. Then
        \begin{align}
            \overline{\Sigma}_{-z_0}\cap \partial \overline{T^*}X^\circ:=\chr_\hbar(P_{h,z_0})\cap  \partial \overline{T^*}X^\circ=\{q\in  \partial \overline{T^*}X^\circ\;|\; (\varrho^2 \tilde G_0)(q)=0\},
        \end{align}
         where $\chr_\hbar(P_{h,z_0})\cap  \partial \overline{T^*}X^\circ=\chr(P(\sigma))$ which we derived in Proposition \ref{eq:ClassicalCharacteristicSet} to be 
        \begin{align}
            \overline{\Sigma}_{-z_0}\cap \partial \overline{T^*}X^\circ \subset\overline{T^*}X^\circ \cap \{r_X\leq r \leq r_{\mathcal{E}}\}.
        \end{align}
        \item \textbf{The horizon radial sets at fibre infinity.} Within $\overline{\Sigma}_{-z_0}\cap \partial\overline{T^*}X^\circ$ we need to take into account the reduced horizon radial sets 
        \begin{align}\label{eq:ReducedHorizonRadialSet}
            \partial\mathcal{R}^\pm_{\mathcal{H}^+,X^\circ} = j^{-1}(\partial \mathcal{R}^\pm_{\mathcal{H}^+})= \{r=r_+, \hat\eta_\theta=\hat\eta_\phi=\hat\eta_\psi=0, G=0, \hat\xi=0\}\subset \overline{\Sigma}_{-z_0},
        \end{align}
        \item \textbf{The normally hyperbolic trapped set at finite covectors.} At finite covectors on the interior, $T^*X^\circ$, we now have to deal with the reduced normally hyperbolic sets 
        \begin{align}
            \Gamma^\pm_{-z_0}:=\iota_{-z_0}^{-1}(\Gamma^\pm) \subset \Sigma_{-z_0}\cap \{r_+<r<\infty\}.
        \end{align}
        \item \textbf{The retained black hole interior.} Finally in the interior region $X_{int}:=X^\circ\cap\{r_X<r<r_+\}=(r_X,r_+)\times\mathbb{S}^3$. We have here that the characteristic region is 
        \begin{align}
            \chr_\hbar(P_{h,z_0}) \cap \pi^{-1}(X_{int}),
        \end{align}
        where by Lemma \ref{lem:SignFlipInterior}
        \begin{align}
            g^{-1}(dr,dr)=\frac{\Delta_r}{\bfr^2}<0,
        \end{align}
        so the $r$-level hypersurfaces are spacelike there. This will provide an additional semiclassical energy estimate near the artificial inner boundary on collars 
        \begin{align}
            r_X\leq r \leq r_1 <r_+.
        \end{align}
\end{enumerate}

We first turn to the horizon semiclassical radial point estimate at $\partial \overline{T^*}X$. In light of the semiclassical correspondence of $\varrho H_{P_{z_0,h}}$ with the classical dynamics of $\varrho H_{p_\sigma}$ on $\overline{T^*}X^\circ$ we relate these back geodesic flow on $\overline{T^*}M^\circ$ through Definition \ref{def:GeneralizedRadialSet} via Definition \ref{def:ClassicalRadialSet}.  

\begin{definition}[Generalized semiclassical radial sets at fibre infinity]\label{def:SemiclassicalRadialSets}
    Let $P_{h,z_0}= h^2\hat{\Box}_{g_{\bfa,\bfb,\mfrm}}(z_0/h)$ have real symbol class representative 
    \begin{align}
        p_{z_0}=\tilde G_{-z_0}
    \end{align}
    and let 
    \begin{align}
        \hat{H}_{p_{z_0}}:=\varrho H_{p_{z_0}}.
    \end{align}
    A smooth submanifold is a generalized semiclassical radial source or sink at fibre infinity if it is a generalized radial source or sink for the classical flow in the sense of Definition \ref{def:ClassicalRadialSet}.
\end{definition}

The classical and semiclassical characteristic sets agree at fibre infinity, and 
\begin{align}
    \varrho H_{p_{z_0}}|_{S^*X^\circ}=\varrho H_{\tilde G_0},\; (\varrho H_{p_{z_0}}-\varrho H_{\tilde G_0})\varrho =\mathcal{O}(\rho^2)
\end{align}
so the classical and semiclassical boundary dynamics agree. Moreover $\pi_0$ in (\ref{eq:FundamentalClassicalProjection}) and $j$ satisfy
\begin{align}
    \pi_0\circ j = \text{id}_{S^*X^\circ}
\end{align}
and we have
\begin{align}
    j^{-1}(\partial \mathcal{R}^\pm_{\mathcal{H}^+})=\pi_0(\partial \mathcal{R}^\pm_{\mathcal{H}^+}).
\end{align}

and we conclude that (\ref{eq:ReducedHorizonRadialSet}) is indeed again a radial set for the semiclassical flow.
Now in general, for the equivalent formulation of Proposition \ref{prop:ClassicalHorizonEstimate}, let us assume $P_{h,z_0} \in \Psi^m_{\hbar}(X^\circ)$ has a real-principal symbol representative $p_{z_0,h}$ so that 
\begin{align}
    P_{h,z}-P_{h,z}^*\in h \Psi^{m-1}_\hbar(X^\circ).
\end{align}
Let $\Lambda\subset \chr_\hbar(P_{h,z}) \cap S^*X^\circ$ be a generalized radial source or sink, with radial coefficient $\alpha_\infty>0$. Define the semiclassical threshold coefficient $\alpha_{\Lambda,h}(z)\in C^\infty(\Lambda)$ by 
\begin{align}
    \varrho^{m-1}\left.\upsigma_{\hbar,m-1}\left(\frac{P_{h,z}-P^*_{h,z}}{2ih}\right)\right|_{\Lambda} =\pm \alpha_\infty\alpha_{\Lambda,h}(z)
\end{align}
We then have the following semiclassical horizon estimate.

\begin{proposition}(Semiclassical radial point estimates) Let $\Lambda\subset \Char_h(P)\cap \partial \overline{T^*}X^\circ$ be a smooth submanifold that is a radial source/sink for the flow of $P_{h,z_0}\in \Psi^m_{\hbar}(X^\circ)$ in the sense of Definition \ref{def:SemiclassicalRadialSets}. Assume $d(\varrho^mp_{z_0,h})\neq 0$ near $\Lambda$. Let $K\subset X^\circ$ be a compact set such that $\pi(\Lambda)\subset K^\circ$. Let $\tilde K$ be a compact set such that the Schwarz kernels $K_{p_{z_0,h}}$ of $P_{z_0,h}$ satisfy $\supp K_{P_{z_0,h}} \cap (K\times X^\circ)\subset K\times \tilde K$. Fix $\chi\in C_c^\infty(X^\circ)$ with $\chi \equiv 1$ near $\tilde K$. Let $u\in \mathscr{D}'(X^\circ)$ and 
\begin{enumerate}
    \item Let $B,G\in \Psi^0_\hbar(X^\circ)$ and $E\in \Psi^0_\hbar(X^\circ)$ have compactly supported Schwarz kernels in $K\times K$.
    \item $\WF'_\hbar(B)\subset \Ell_\hbar(G)$.
    \item $\Lambda\subset \Ell_{\hbar}(G)$ (source); $\Lambda \subset \Ell_{\hbar}(B)$ (sink)
    \item All integral curves $\gamma:I\to \overline{T^*}X^\circ$ of $\mp\varrho^{m-1} H_{p_{h,z_0}}$ (source/sink) starting from in $\WF'_\hbar(B) \cap \Char_\hbar(P_{h,z_0})$ tend to $\Lambda$ asymptotically and satisfy $\gamma(I)\subset \Ell_\hbar (G)$.
\end{enumerate}
Then if $\Lambda$ is a source, for all $s, s_0,N\in \mathbb{R}$ with 
\begin{align}
    s >s_0>\frac{m-1}{2}+\alpha_{\Lambda,h}(z)
\end{align}
then 
\begin{align}
    \|Bu\|_{H^s_h}\leq C\left(h^{-1} \|GP_{h,z}u\|_{H^{s-m+1}_h} +  h^N\|Gu\|_{H^{s_0}_h} + h^N \|\chi u \|_{H_h^{-N}}\right).
\end{align}
If $\Lambda$ is a sink, we have for all $s,N\in \mathbb{R}$ with 
\begin{align}
    s<\frac{m-1}{2}+\alpha_{\Lambda,h}(z)
\end{align}
that 
\begin{align}
    \|Bu\|_{H^s_\hbar}\leq C\left(h^{-1}\|GP_{h,z}u\|_{H^{s-m+1}_h} +\|Eu\|_{H^s_h}+ h^N\|\chi u\|_{H_h^{-N}}\right).
\end{align}
\end{proposition}

A computation similar to the following Proposition \ref{prop:ClassicalHorizonEstimate} shows the semiclassical equivalent to (\ref{eq:ThresholdSurfacePressure}), and gives 
\begin{align}
    \alpha_\infty =\frac{4\mfrm \kappa}{\bfr^2|_{r=r_+}} \text{ and } \alpha_{\Lambda,h}(z) =-\frac{\nu}{\kappa}
\end{align}
giving the treshold condition as 
\begin{equation}
 s_0 >\frac{1}{2}-\frac{\nu}{\kappa}.
\end{equation}

Now we turn our attention to trapping for $P_{h,z_0}$ for $z_0 \in \{-1,+1\}$. We gave its definition on $T^* M^\circ$ already in Definition \ref{def:TrappedSet}. To reduce the analysis to $T^*X^\circ$ we only need to consider 
\begin{align}
    \Gamma^\pm_{-z_0} :=\iota^{-1}_{-z_0}(\Gamma^\pm),
\end{align}
which intersect transversally as codimension one submanifolds of $\Sigma_{-z_0}$ as
\begin{align}
    \Gamma_{-z_0} :=\Gamma^+_{-z_0}\cap \Gamma^-_{-z_0},
\end{align}
To simplify

Correspondingly we also project their defining functions
\begin{align}
    \phi_{\pm,-z_0}:=\iota^*_{-z_0}\phi_\pm,
\end{align}
and hence we know there exists an open neighbourhood $U\subset T^*X^\circ$ such that 
\begin{align}
    \Gamma_{-z_0}^\pm \cap U=\{ \phi_{\pm,-z_0}=0\}
\end{align}
and $d\phi_{\pm,-z_0} \neq 0$ on $\Gamma^\pm_{-z_0}$.
Furthermore we have
\begin{align}
    H_{G_{-z_0}}\phi_{\pm,-z_0} = \pm \nu_{\pm,-z_0}\phi_{\pm,-z_0},
\end{align}
with $\nu_{\pm,-z_0}\in C^\infty(U)$ and we set $\nu_{min} =\text{min}(\inf_{\Gamma_{-z_0}} \nu_{+,-z_0},\inf_{\Gamma_{-z_0}} \nu_{-,-z_0})$.
We have $\{\phi_{-,-z_0},\phi_{+,-z_0}\}\neq 0$ and $dG_{-z_0}\neq 0$ on $\Gamma_{-z_0}$.

In what follows, we will make the explicit choice $z_0= 1$ and surpress the notation i.e. we work again with 
\begin{align}
    \Gamma^\pm := \Gamma^\pm_{-1},\; \Gamma:=\Gamma_{-1},\;\phi_\pm :=\phi_{\pm,-1},\; \nu_{\pm 1,-1}=-1,
\end{align}
the dynamics for the negative frequence $z_0=-1$ can be obtained conjugation, 
\begin{align}
    P_{h,-\overline{z}}\overline u = \overline{P_{h,z}u},
\end{align}
and hence by analyzing $-H_{G_1}$ instead of $H_{G_{-1}}$.

Introducing $\beta_0 =\min_{(x,\xi)\in \Gamma} g^{-1}_{\bfa,\bfb,\mfrm}(-dt +\xi,dt)>0$, and the semiclassical subprincipal symbol at $z_0 =1$ as $p_1:= \upsigma_1\left(\frac{P_{h,1}-P_{h,1}^*}{2ih}\right)$, and the pair 
    \begin{align}\label{eq:DefinitionoflittleGammas}
        \gamma_+:=\frac{\sup_\Gamma p_1}{2\beta_0} \text{ and } \gamma_-:=\gamma_+-\frac{\nu_{min}}{4\beta_0}.
    \end{align}
    Then for $\eta\geq \gamma_-$ and $\epsilon \geq 0$ the semiclassical loss function satisfies $\delta_\Gamma(\eta)\in [0,1]$ where $\delta_\Gamma$ is given by
    \begin{align}\label{eq:Semiclassicallossfunction}
        \delta_{\Gamma}(\eta):=\max\left\{\frac{\gamma_+-\eta}{\gamma_+-\gamma-},0 \right\}.
    \end{align}

\begin{theorem}\label{thm:NormallyHyperbolicTrappingEstimate}
    (Estimate at normally hyperbolic trapping). Let $\beta>\gamma_-$. Let $\nu>0$ and $N\in \mathbb{R}$. Let $U\in T^*X$ be an open neighbourhood of $\Gamma$ and $K\subset X$ be a compact neighbourhood of the baseprojection of $\Gamma$. Let $\chi\in C^\infty_c(X)$ be a bump function with $\chi \equiv 1$ near $K$. Furthermore, fix 
    \begin{align}
        \beta>\gamma_-, \qquad  \nu>\max\{\beta,0\},
    \end{align}
    and define 
    \begin{align}
        \Omega^\pm_{h,\beta,\nu}:=\{z\in\mathbb{C} \;|\; |z|=1, \pm\text{Re }z>0, \;h\beta \leq \text{Im z} \leq h\nu.
    \end{align}
    
    The the following estimates hold uniformily for 
    \begin{align}
        z\in \Omega^+_{h,\beta,\nu}\cup \Omega^-_{h,\beta,\nu}, \qquad 0<h<h_0
    \end{align}
    with the dynamical hypotheses interpreted using $z_0 H_{G_{z_0}}$, $z_0 =\sgn (\text{Re }z)$.
    
    \begin{enumerate}
        \item (Direct Estimate.) There exist semiclassical ps.d.o.s $B,E,G\in \Psi^0_{\hbar}(X) $ with Schwartz kernels supported in $K\times K$ and wavefront sets in $U$ with $B$ semiclassically eliptic at $\Gamma$ and $WF'_{\hbar}(E)\cap \Gamma^u=\emptyset$, and $C>0$ s.t. 
        \begin{align}
            \|Bu\|_{H^s_h} \leq C_{\epsilon}\left(h^{-1-\delta_{\Gamma}(\text{Im}(h^{-1}z)-\epsilon)}\|GP_{h,z}u\|_{H^{s-1}_h} + h^{-\delta_\Gamma(\text{Im}(h^{-1}z-\epsilon)}\|Eu\|_{H^s_h}+h^N \|\chi u\|_{H_h^{-N}} \right).
        \end{align}
        This holds provided $\Gamma \cup WF_{\hbar}'(B)\subset \text{Ell}_\hbar(G)$ and all null-bicharacteristics starting at a point in $WF'_\hbar \backslash \Gamma^u$ reach $\text{Ell}_\hbar(E)$ in finite time whilst remaining in $\text{Ell}_\hbar(G)$.
        \item (Adjoint Estimate.) There exist semiclassical ps.d.o.s $B,E,G\in \Psi^0_{\hbar}(X) $ with Schwartz kernels supported in $K\times K$ and wavefront sets in $U$ with $B$ semiclassically eliptic at $\Gamma$ and $WF'_{\hbar}(E)\cap \Gamma^s=\emptyset$, and $C>0$ s.t.
        \begin{align}
            \|Bv\|_{H^s_h} \leq C_{\epsilon}\left(h^{-1-\delta_{\Gamma}(\text{Im}(h^{-1}z)-\epsilon)}\|GP^*_{h,z}v\|_{H^{s-1}_h} + h^{-\delta_\Gamma(\text{Im}(h^{-1}z-\epsilon)}\|Ev\|_{H^s}+h^N \|\chi v\|_{H_h^{-N}} \right).
        \end{align}
        This holds provided $\Gamma \cup WF_{\hbar}'(B)\subset \text{Ell}_\hbar(G)$ and all null-bicharacteristics starting at a point in $WF'_\hbar \backslash \Gamma^s$ reach $\text{Ell}_\hbar(E)$ in finite time whilst remaining in $\text{Ell}_\hbar(G)$.
        \end{enumerate}
\end{theorem}
Finally at $\partial_-X$ we can provide the following high energy estimate.
\begin{proposition}\label{prop:HighEnergyEstimate}
    (Semiclassical energy estimate near the inner boundary $\partial_-X$) Let $s,\nu\in \mathbb{R}$, and let
    \begin{align}
        r_-<r_X<r_0<r_1<r_+,
    \end{align}
    and set 
    \begin{align}
        \Omega=r^{-1}([r_X,r_1]),\quad \Omega_{out}=r^{-1}([r_0,r_1]), \quad \Omega_{in}=r^{-1}[r_X,r_0].
    \end{align}
    Let $I\Subset \mathbb{R}$ be a compact interval. Then exists a constant $C>0$, $h_0>0$ such that for all $z\in \mathbb{C}$ with $|z|\in [\frac{1}{2},2]$ and $h^{-1}\text{Im }z\in I$ we have the direct estimate 
    \begin{align}
        \|u\|_{H^s_h(\Omega)}\leq C \left(h^{-1}\|P_{h,z} u\|_{H^{s-1}_h(\Omega)}+\|u\|_{H^{s}_h(\Omega_{out})}\right)
    \end{align}
    and the adjoint estimate 
    \begin{align}
         \|v\|_{H^s_h(\Omega)}\leq C \left(h^{-1}\|P^*_{h,z} u\|_{\dot{H}^{s-1}_h(\Omega)}+\|u\|_{H^{s}_h(\Omega_{in})}\right)
    \end{align}
\end{proposition}

Now we turn our attention to $\partial_+X$. Taking $P_{h,z_0}$ as before we see immediately its semiclassical scattering symbol at $\partial_+ X$ is given  
\begin{align}
    p_{z_0}:={}^{sc,\hbar}\upsigma_{2,0}(P_{h,z_0})={}^{sc}\xi^2+\|\sceta_i\|^2_{g^{-1}_{\mathbb{S}^3}}-2z_0\scxi,
\end{align}
with semiclassical characteristic set 
\begin{align}
    \Sigma_{z_0}:=\chr(P_{h,z_0})= \{p_{z_0,h}=0\} 
\end{align}

So we might expect most of our previous analysis in Section \ref{subsec:scbcalculus} to go through unchanged under the substitution $\sigma = z_0$. Indeed, directly wielding Definition \ref{def:scRadialPoints}, we would find the sets

\begin{align}\label{eq:SemiclassicalScatteringRadialSource}
    \mathcal{R}_{z_0,\partial X}^+ &=\{(\rho,y_i,\scxi, \sceta_i)\in {}^{sc}T^*_{\partial X}X\cap \Sigma_G\;|\; \rho=0, {}^{sc}\xi=2z_0, {}^{sc}\eta_i=0\}\\
    \mathcal{R}_{\partial X}^-&=\{(\rho,y_i,\scxi, \sceta_i)\in {}^{sc}T^*_{\partial X}X\cap \Sigma_G\;|\; \rho=0, {}^{sc}\xi=0, {}^{sc}\eta_i=0\},\label{eq:SemiScatteringRadialSink}
\end{align}
and the analysis would go through unchanged. To resolve $  \mathcal{R}_{\partial X}^-$ we then work on 
\begin{align}
    {}^{sc,b,\hbar}\overline{T^*}X:=[{}^{sc,\hbar}\overline{T^*}X;o_{\partial X}].
\end{align}
for which we have 
\begin{align}
    P_{h,z}\in \Psi^{2,0,1}_{sc,b,\hbar}(X).
\end{align}

The lifted semiclassical characteristic set $\tilde{\Sigma}_{z_0}:=\overline{\beta_{sc}^{-1}(\Sigma_{z_0}\backslash o_{\partial X})}$ then needs to be analyzed again near $\text{ff}$. The equivalent of (\ref{eq:SymbolNearff}) becomes 
\begin{align}
    q_{z_0}=\rho^{-1}\beta_{sc}^* p_{z_0}= \rho (\bxi^2 +\|\beeta_i\|_{g^{-1}_{\mathbb{S}^3}}^2)-2z_0\bxi,
\end{align}
which we again rescale as 
\begin{align}\label{eq:RescaledSemiclassicalq}
    \hat{q}_{z_0}:=\rho_{\text{sf}}q_{z_0}=\rho_{\text{ff}}(1-\rho_{\text{sf}}^2)-2z_0{}^{sc}\hat{\xi}
\end{align}
in projective coordinates near $\text{ff}$.\\
The main difference now is that because the whole resolved phase space ${}^{sc,b}\overline{T^*}X$ belongs to the semiclassical locus $\mathcal{Z}_{sc,b,h}$ at $h=0$, we need to consider the entirety of $\text{ff}$, even when $\rho_{\text{sf}}\neq 0$\footnote{Basically this is due to the fact that the semiclassical characteristic set of $P_{h,z}$ on the interior ${}^{sc,b,\hbar}\overline{T^*}X^\circ$ can contain finite coordinates and hence need to connect to the characteristic set at  ${}^{sc,b,\hbar}\overline{T^*}_{\partial X}X$ via the front face $\text{ff}$} as 
\begin{align}\label{eq:RescaledSemiClassicalqonff}
    \hat{q}_{z_0}|_{\text{ff}}=-2z_0{}^{sc}\hat{\xi}.
\end{align}

hence we see that 
\begin{align}
    \tilde\Sigma_{z_0} \cap \text{ff} =\{\rho_{\text{ff}}=0,{}^{sc}\hat{\xi}=0\}. 
\end{align}
We claim that the semiclassical radial outgoing set is now given by 
\begin{align}
    \tilde {\mathcal{R}}_{z_0,\partial_+X}:=\tilde{\Sigma}_{z_0}\cap \text{ff}=\{\rho_{\text{ff}}=0,{}^{sc}\hat\xi=0\}.
\end{align}

We perform the same checks of Definition \ref{def:RadialscbSet}. 
From (\ref{eq:RescaledSemiClassicalqonff}) one instantly checks 
\begin{align}
    d(\hat{q}_{z_0}|_{\text{ff}})=-2z_0 d({}^{sc}\hat{\xi})\neq 0 \text{ on } \tilde {\mathcal{R}}_{z_0,\partial_+X}.
\end{align}

Moreover (\ref{eq:RescaledSemiclassicalq}) gives us on $\tilde{\Sigma}_{z_0}$, 
\begin{align}
    {}^{sc}\hat{\xi}=\frac{\rho_{\text{ff}}(1-\rho_{\text{sf}}^2)}{2z_0}
\end{align}

Now the remaining conditions in Definition \ref{def:RadialscbSet}, we have that 
\begin{align}
 \left.\rho_{\text{ff}}^{-1}H_{q_{z_0}}\rho_{\text{ff}}\right|_{\mathcal{R}^-_{z_0,\partial_+X}}=-2z_0
\end{align}
and taking the quadratic defining function $\mathfrak{q}=\rho_{\text{ff}}^2$ we get 
\begin{align}\label{eq:Fixedz0QuadraticCheck}
    H_{q_{z_0}}\mathfrak{q}=-4z_0\mathfrak{q}+\frac{1-\rho_{\text{sf}}^2}{z_0}.
\end{align}

where we note that (\ref{eq:FixedSigmaQuadraticCheck}) is retrieved from (\ref{eq:Fixedz0QuadraticCheck}) by restricting to $\text{sf}$, i.e. setting $\rho_{\text{sf}}=0$.
The analysis on $\tilde{\mathcal{R}}^+_{z_0,\partial_+ X}=\beta^{-1}_{sc}(\tilde{\mathcal{R}}^+_{z_0,\partial_+ X})\subset \text{sf}^\circ$ remains exactly the same as before. 

With this we can state the semiclassical sc-b radial point estimates. We use the conventions
\begin{align}
    P_{h,z}= h^2\hat{\Box}_{g_{\bfa,\bfb,\mfrm}}(h^{-1}z),\; h=|\sigma|^{-1}, z=h\sigma,
\end{align}
with $\sgn \text{ Re }\sigma = z_0\in \{-1,1\}$ and 
\begin{align}
    0\leq \text{Im }\sigma \leq \nu. 
\end{align}
with $\nu$ fixed. 
Choose a compact collar $K$ of $\partial_+X$ and $\chi\equiv 1$ near $K$. All $G,B,E\in \Psi_{sc,b,\hbar}^{0,0,0}(X)$ have kernels supported and are chosen in sufficiently small neighbourhoods of  $\tilde{\mathcal{R}}^+_{z_0,\partial X}$ or $\tilde{\mathcal{R}}^-_{\partial X}$.

\begin{proposition}[Incoming estimates]
     (Direct Estimate) We work near $\tilde{\mathcal{R}}^+_{z_0,\partial_+X}\subset \text{sf}^\circ\cap \tilde{\Sigma}_{z_0}$. Assume that
\begin{align}
    \ell_{sc}>\ell_0 >-\frac{1}{2}(1-\theta_{in}), \; s,\ell_b\in\mathbb{R}.
\end{align}
For the direct estimate, there exist $B,G\in \Psi_{sc,b,\hbar}^{0,0,0}(X)$ such that we have the direct estimate
\begin{align}
    \|Bu\|_{H^{s,\ell_{sc},\ell_b}_{sc,b,h}}\leq C\left(h^{-1}\|GP_{h,z} u\|_{H^{s-2,\ell_{sc+1},\ell_b+1}_{sc,b,h}}+ h^N\|\chi u\|_{H^{-N,\ell_0,\ell_b}}\right)
\end{align}
    
    and for the adjoint estimate there exist $B,E,G \in \Psi^{0,0,0}_{sc,b,h}(X)$ such that 
    \begin{align}
        \|Bv\|_{H^{2-s,-\ell_{sc}-1,-\ell_b-1}_{sc,b,h}} \leq C&\left( h^{-1}\|G P^*_{h,z}v\|_{H^{-s,-\ell_{sc},-\ell_b}_{sc,b,h}}\right.\\
        &\left. \|Ev\|_{H^{2-s,-\ell_{sc}-1,-\ell_b-1}_{sc,b,h}}+h^N\|\chi v\|_{H^{-N,-N,-N}_{sc,b,h}} \right)
    \end{align}   
\end{proposition}

\begin{proposition}[Outgoing estimates]
    We work near $\tilde{\mathcal{R}}^-_{\partial_+X}\subset \text{ff}\cap \tilde{\Sigma}_{z_0} $. Assume 
    \begin{align}
        \ell_b<-\frac{1}{2}+\underline{S},\;\;\;\; s,r_{sc}\in\mathbb{R}.
    \end{align}
    For the direct estimate there exist $B,E,G\in \Psi^{0,0,0}_{sc,b,\hbar}(X)$ such that 
    \begin{align}
        \|Bu\|_{H_{sc,b,h}^{s.\ell_{sc},\ell_b}}\leq C&\left(h^{-1}\|G P_{h,z}\|_{H_{sc,b,h}^{s-2,\ell_{sc}+1,\ell_b+1}} +\|Eu\|_{H_{sc,b,h}^{s,\ell_{sc},\ell_b}}\right.\\
        &\left. +h^N \|\chi u\|_{H^{-N,-N,\ell_b}_{sc,b,h}} \right).
    \end{align}
    For the adjoint estimate there exist $B,G \in \Psi^{0,0,0}_{sc,b,\hbar}(X) $ such that 
    \begin{align}
        \|Bv\|_{H_{sc,b,h}^{s-2,-\ell_{sc}-1,-\ell_b-1}}\leq C\left(h^{-1}\|GP^*_{h,z} v\|_{H_{sc,b,h}^{-s,-\ell_{sc},-\ell_b}}+h^N \|\chi v\|_{H_{sc,b,h}^{-N,-N,-\ell_b-1}}\right)
    \end{align}
\end{proposition}

We combine this with semiclassical $sc-b$ and elliptic estimates.
\begin{proposition}(Semiclassical $sc,b$ ellptic estimates) Let 
\begin{align}
    P_h\in \Psi^{m,l_{sc},l_b}_{sc,b,h}(X),\;\;B_h,G_h \in \Psi^{0,0,0}_{sc,b,h}(X),
\end{align}
and assume 
\begin{align}
    \WF'_{sc,b,h}(B_h) \subset \Ell_{sc,b,h}(P_h)\cap \Ell_{sc,b,h}(G_h).
\end{align}
Then for sufficiently large $N$, there exist $C, h_1>0$ such that 
\begin{align}
    \|B_hu\|_{H_{sc,b,h}^{s,\ell_{sc},\ell_b}}\leq C&\left( h^{-1}\|GP_{h,z} u\|_{H^{s-m,\ell_{sc}-l_{sc}+1,\ell_b-l_b}_{sc,b,h}}+h^N\|\chi u\|_{H^{-N,-N,l_b}_{sc,b,h}}\right)
\end{align}

\end{proposition}

\begin{proposition}(Principal type propagation on $\tilde\Sigma_{z_0} \backslash (\tilde{\mathcal{R}}^+_{z_0,\partial X}\cup \tilde{\mathcal{R}}^-_{\partial X})$). Assume $P_{z,h}$ admits a real representative $p_h$ satisfying 
\begin{align}
    P_h-\text{Op}_{sc,b,h}(p_h) \in h\Psi^{m-1,l_{sc}-1,l_b}_{sc,b,h}(X),
\end{align}
with
\begin{align}
    \tilde{p}_h=\varrho_{\text{fi}}^m\rho_{\text{sf}}^{l_{sc}}\rho_{\text{ff}}^{l_b} p_h \in C^\infty({}^{sc,b,h}\overline{T^*}X),
\end{align}
and set 
\begin{align}
    {}^{sc,b,h}\hat{H}_{p_0}=\varrho^{m-1}_{\text{fi}}\rho_{\text{sf}}^{l_{sc}-1}\rho_{\text{ff}}^{l_b}H_{p_0}.
\end{align}
Let $B,E,G\in \Psi^{0,0,0}_{sc,b,h}(X)$, with $G$ elliptic on a bicharacteristic segment in $\tilde{\Sigma}_z$, $E$ elliptic at the endpoint and every backward bicharacteristic from $\WF'(B)$ reach $\Ell(E)$ whilst remaining in $\Ell(G)$. It then holds that 
\begin{align}
    \|Bu\|_{H^{s,\ell_{sc},\ell_b}_{sc,b,h}} \leq C\left(\|Eu\|_{H_{sc,b,h}^{s,\ell_{sc},\ell_b}} +h^{-1}\|G P_{h,z}u\|_{H^{s-1,\ell_{sc}+1,\ell_b+1}_{sc,b,h}} + h^N \|u\|_{H^{-N,-N,\ell_b}_{sc,b,h}} \right)
\end{align}
    
\end{proposition}

\begin{theorem}\label{thm:SemiclasscialHighEnergyEstimate}
    Let $(M,g_{\bfa,\bfb,\mfrm})$ be a Myers-Perry spacetime with spatial manifold $X$ as defined by Definition \ref{def:RadiallyCompactifiedR4}, and let $(\bfa,\bfb) \in D_{\mfrm}$. Let $P_{h,z}:=h^2 \hat{\Box}_{g_{\bfa,\bfb,\mfrm}}(h^{-1}z)$ as before and fix strip height $\nu> 0$, $0<\epsilon<-\gamma_-$, with $\gamma_-$ defined in (\ref{eq:DefinitionoflittleGammas}). Choose $s$ such that it satisfies the horizon condition 
    \begin{align}
        s>s_0>\frac{1}{2}-\frac{\nu}{\kappa}
    \end{align}
    at $\mathcal{H}^+$. Let also $\ell_{sc}>-\frac{1}{2}$ and $\ell_b <-\frac{1}{2}$. Then for for $z\in \Omega_h$ with
    \begin{align}\label{eq:OmegaDomain}
         \Omega_h:=\{z\in \mathbb{C}\;|\; |z|=1, 0\leq \text{Im }z\leq h\nu\},
    \end{align}

    then there exists a $C_\epsilon$ s.t. for every $0<h<h_0$ we have that 
    \begin{align}
        \|u\|_{\overline{H}^{s,\ell_{sc},\ell_b}_{sc,b,h}(X)}\leq C_\epsilon h^{-1-\delta_\Gamma(\text{Im}(h^{-1}z)-\epsilon)} \|P_{h,z} u\|_{\overline{H}^{s-1,\ell_{sc}+1,\ell_b+1}_{sc,b,h}(X)},
    \end{align}
    where $\delta_{\Gamma}$ is the semiclassical loss $\delta_\Gamma$ defined in (\ref{eq:Semiclassicallossfunction}). Likewise we have the adjoint estimate
    \begin{align}
        \|v\|_{\dot{H}^{1-s,-\ell_{sc}-1,-\ell_b-1}}\leq C_\epsilon h^{-1-\delta_\Gamma(\text{Im}(h^{-1}z)-\epsilon)}\|P^*_{h,z}v\|_{\dot{H}_{sc,b,h}^{-s,-\ell_{sc},-\ell_b}}
    \end{align}
\end{theorem}
\begin{remark}
    We conclude from the above that the only semiclassical loss in the global estimate \textit{is} the loss occured through $\delta_\Gamma$ at the trapped set.
\end{remark}

We now use (\ref{eq:TranslationtobSpaces}) to immediately deduce the following corollary, which gives us high energy estimates and thus control over the wave operator $\Box_g(\sigma)^{-1}$ for $|\text{Re }\sigma|>>1$.
\begin{corollary}\label{cor:HighEnergyEstimate}
    Let $(M,g_{\bfa,\bfb,\mfrm})$ be a Myers-Perry spacetime with spatial manifold $X$ as defined by Definition \ref{def:RadiallyCompactifiedR4}, and let $(\bfa,\bfb) \in D_{\mfrm}$. Let $P_{h,z}:=h^2 \hat{\Box}_{g_{\bfa,\bfb,\mfrm}}(h^{-1}z)$. Let $s>s_0>\frac{1}{2}-\frac{\nu}{\kappa}$, $\ell<-\frac{1}{2}$ and $s+\ell>-\frac{1}{2}$. Furthermore fix an $0<\epsilon<-\gamma_-$. 
    
    We then have for $z\in \Omega_h$ as defined in (\ref{eq:OmegaDomain}) that there exists a $C_\epsilon$ such that for all $0<h<h_0$,
    \begin{align}
        \| u\|_{\overline{H}^{s,\ell}_{b,h}(X)}\leq C_\epsilon h^{-1-\delta_\Gamma(\text{Im }(h^{-1}z)-\epsilon)}\|P_{h,z}u\|_{\overline{H}^{s,\ell+1}_{b,h}(X)}.
    \end{align}
    and for the adjoint, we have 
    \begin{align}
        \|v\|_{\dot{H}^{-s,-\ell-1}_{b,|\sigma|^{-1}}(X)}\leq C_\epsilon h^{-1-\delta_\Gamma(\text{Im } h^{-1}z -\epsilon)} \|P_{h,z}^* v\|_{\dot{H}^{-s,-\ell}_{b,|\sigma|^{-1}}}
    \end{align}
        
    Rephrased in the spectral parameter $\sigma\in \mathbb{C}\backslash \{0\}$ we have that for all 
    \begin{align}
        \sigma \in \{\sigma\in \mathbb{C}\;|\; 0\leq \text{Im }\sigma \leq \nu, |\sigma|>h_0^{-1}\}
    \end{align}
    that there exists a $C_\epsilon>0$ such that
    \begin{align}\label{eq:DirectHighEnergyEstimate}
        \| u\|_{\overline{H}^{s,\ell}_{b,|\sigma|^{-1}}(X)}\leq C_\epsilon |\sigma|^{-1+\delta_{\Gamma}(\text{Im }(\sigma)-\epsilon)}\|\hat{\Box}_{g_{\bfa,\bfb,\mfrm}}(\sigma) u\|_{\overline{H}^{s,\ell+1}_{b,|\sigma|^{-1}}(X)}.
    \end{align}
    Similarly, we have the adjoint estimate 
    \begin{align}\label{eq:AdjointHighEnergyEstimate}
        \|v\|_{\dot{H}^{-s,-\ell-1}_{b,|\sigma|^{-1}}(X)}\leq C_\epsilon |\sigma|^{-1+\delta_{\Gamma}(\text{Im }(\sigma) -\epsilon)}\|\hat{\Box}_{g_{\bfa,\bfb,\mfrm}}(\sigma)^* u\|_{\dot{H}^{-s,-\ell}_{b,|\sigma|^{-1}}(X)}.
    \end{align}
\end{corollary}

\newpage
\section{Part III: Resolvent Expansion and Price's Law.}\label{sec:PriceLaw}

In this part, we shall use the high energy estimates from Theorem \ref{thm:SemiclasscialHighEnergyEstimate} and Corollary \ref{cor:HighEnergyEstimate}, to derive Price's law for the wave equation on $(M,g_{\bfa,\bfb,\mfrm})$ in the full extremal range $(|\bfa|+|\bfb|)^2< 2\mfrm$. The approach incorporates a resolvent expansion of $\hat{\Box}_{g_{\bfa,\bfb,\mfrm}}(\sigma)$ which is executed over the full indicial gap range $\alpha\in (0,2)$ and closely follows \cite{hintz2022sharp} and \cite{hintz2023linear}.\\
The key concept in applying the resolvent expansion to the spectral family $\hat{\Box}_g(\sigma)$ of a wavetype operator on an asymptotically flat space time $(M,g)$ is identifying a set admissibility criteria of $(M,g)$ need to satisfy. \cite{hintz2022sharp} divides these these in $2$ categories, a geometric admissiblity criterion on $g$ and a spectral admissibility criterion for the corresponding operator $\hat{\Box}_g(\sigma)$. The purpose of these is two-fold:
\begin{enumerate}
    \item They provide sufficient conditions for the existence of the resovlent 
    \begin{align}
        \hat{\Box}_{g_{\bfa,\bfb,\mfrm}}(\sigma)^{-1}:\mathcal{Y}^{s,r}\to \mathcal {X}^{s',r'},
    \end{align}
    for suitable $\sigma$-independent Sobolev spaces $\mathcal{Y}^{s,r}$ and $\mathcal {X}^{s',r'}$ to justify the Payley-Wiener contourshifting argument producing the solution to the forward solved wave problem
    \begin{align}
        \phi(t_*,x) = \int_{\mathbb{R}}e^{-i\sigma t_*}\hat{\Box}_{g_{\bfa,\bfb,\mfrm}}(\sigma)^{-1}\hat{f}(\sigma,x).
    \end{align}
    \item They provide necessary and sufficient conditions for the refinement of the output of $\hat{\Box}_{g_{\bfa,\bfb,\mfrm}}(\sigma)^{-1}$ on suitable conormal spaces $\mathcal{A}^s(X)$. The resolvent expansion can then proceed over the indicial range of the b-normal operator ${}^bN(\hat{\Box}_{g_{\bfa,\bfb,\mfrm}}(\sigma))$ with $\alpha \in (0,2)$ for sufficiently regular forcing data $f\in \mathcal{A}^{4+\alpha-}(X)$. Once reaching threshold decay this then produces the polynomial late-time tail for $\phi(t_*,x)$ and hence also completes the proof of stability.
\end{enumerate}

\subsection{Geometric and spectral admissibility}

The notations of $(M,g_{\bfa,\bfb,\mfrm})$ and the corresponding compactified spatial manifold $X$ and the related $b$- and $sc$-algebras from Part I and Part II go through unaltered.\\
In this section we identify the geometric admissibility criterion that $(M,g_{\bfa,\bfb,\mfrm})$ needs to satisfy for the contourshifting and resolvent expansion arguments. We will also refine the structure of $\hat{\Box}_{g_{\bfa,\bfb,\mfrm}}(\sigma)$ in terms of its $\text{Diff}_b$ and $\text{Diff}_{sc}$ structure. The following is a direct analogue of Definition 2.3 in \cite{hintz2022sharp}.

\begin{definition}\label{def:AdmissibleMetric}[Geometric Admissiblity Criterion]
    We say that an asymptotically flat Lorentzian metric with mass $\mfrm\in\mathbb{R}$ on $M^{\circ}=\mathbb{R}_{t_*}\times X^\circ$ is admissible if $\partial_{t_*}$ is a Killing vector field (i.e. it is stationary), and if moreover
    \begin{enumerate}
        \item $dt_*$ is everywhere future timelike, i.e. $g^{00}=g(dt_*,dt_*)<0$;
        \item The coefficients of the dual metric 

        \begin{align}
            g^{-1}=g^{00}\partial^2_{t^*}+2\partial_{t_*}\otimes_s g^{0X} + g^{XX}
        \end{align}
        satisfy
        \begin{align}
            g^{00}&\in \rho^2C^\infty(X), \label{eq:g00Component} \\
            g^{0X}&\in \left\langle \rho^2\partial_\rho + \rho^2C^\infty(X;{}^{sc}TX)\right\rangle,\label{eq:g0XComponent}\\
            g^{XX}&\in \left\langle(1-2\mfrm \rho^2)(\rho^2\partial_\rho)^2+\rho^2\cancel{g}^{-1}+\rho^2C^\infty(X;S^2 {}^{sc}TX) \right\rangle \label{eq:gXXComponent}
        \end{align}
    \end{enumerate}
\end{definition}

The first requirement is necessary to provide an actual forward solution in $t_*$ to $\Box_{g_{\bfa,\bfb,\mfrm}}\phi(t_*,x)=f(t_*,x)$ via a Payley-Wiener argument. As it turns out standard $t_*$ provided by standard star coordinates (\ref{eq:MPStartcoordinate2})-(\ref{eq:MPStarpsicoordinate2}) do not automatically satisfy the first requirement due to the ergoregion $\mathcal{E}$. 
Therefore let us define the following coordinates, where we use that $r_+\geq \sqrt{2\mfrm}$:
\begin{align}
    t_* &:= t + \int_{\sqrt{3\mfrm}}^r \frac{(\tilde{r}^2 + \bfa^2)(\tilde{r}^2+\bfb^2)}{\tilde{r}^2 \Delta_{\tilde{r}}} \left( \chi_+(\tilde{r}) - \chi_{\partial}(\tilde{r})\right)\:d\tilde{r} + F(r)\label{eq:MPStartcoordinate},\\
    \phi_* &:= \phi + \int^r_{\sqrt{3\mfrm}} \frac{\bfa(\tilde{r}^2 +\bfb^2)}{\tilde{r}^2\Delta_{\tilde{r}}} \chi_+(\tilde{r})\:d\tilde{r},\label{eq:MPStarphicoordinate} \\
     \psi_* &:= \psi + \int^r_{\sqrt{3\mfrm}} \frac{\bfb(\tilde{r}^2 +\bfa^2)}{\tilde{r}^2\Delta_{\tilde{r}}} \chi_+(\tilde{r})\:d\tilde{r}\label{eq:MPStarpsicoordinate},    
\end{align}
where the cut-off functions $\chi_+,\chi_\partial \in C^\infty(\mathbb{R}_{\geq 0})$ are given by 
\begin{align}
    &\chi_+ \equiv 1 \; \text{ on } \; [0,\sqrt{3\mfrm}) \text{ and } \chi_+\equiv 0 \; \text{ on } [\sqrt{4\mfrm},\infty) \\
    &\chi_\partial \equiv 0 \; \text{ on } \; [0,\sqrt{3\mfrm}] \text{ and } \chi_\partial \equiv 1 \; \text{ on } [\sqrt{6\mfrm},\infty)    
\end{align}
and $F\in C^\infty(\mathbb{R}_{\geq 0})$ to be chosen below. Under these transformations we get 
\begin{align}\label{eq:TrueRescaledDualMetric}
    &\bfr^2  g^{-1}_{\mfrm,\bfa,\bfb} = -\frac{1-\chi_+^2}{r^4 \bfDelta}\left[(r^2+\bfa^2)(r^2 +\bfb^2) \partial_{t_*} + \bfa (r^2 + \bfb^2) \partial_{\phi_*} + \bfb(r^2 + \bfa^2) \partial_{\psi_*} \right]^2  \nonumber \\
    &+ \bfDelta \partial_r^2 +\partial_\theta^2 + \bfDelta\left[-\frac{(r^2 + \bfa^2)(r^2 + \bfb^2)}{r^2\bfDelta} \chi_\partial + F' \right]^2 \partial^2_{t_*}\\
    & +2 \left[\left(\frac{(r^2+\bfa^2)(r^2+\bfb^2)}{r^2}(\chi_+ - \chi_\partial) +\bfDelta F' \right) \partial_{ t_*} + \frac{\bfa(r^2+\bfb^2)}{r^2} \chi_+ \partial_{\phi_*} + \frac{\bfb(r^2+\bfa^2)}{r^2}\chi_+ \partial_{\psi_*}\right]\partial_r \nonumber \\
    & +2\chi_+\left[ \frac{(r^2 + \bfa^2)(r^2+\bfb^2)}{r^2} \partial_{t_*} + \frac{\bfa(r^2 + \bfb^2)}{r^2} \partial_{\phi_*} + \frac{\bfb(r^2 + \bfa^2)}{r^2}\partial_{\psi_*}\right]\left[-\frac{(r^2+\bfa^2)(r^2+\bfb^2)}{r^2\bfDelta}\chi_\partial + F'\right]\partial_{t_*} \nonumber\\
    &+ \frac{\sin^2\theta\cos^2\theta}{\bfp^2}\left[(\bfa^2 - \bfb^2) \partial_{t_*} + \frac{\bfa}{\sin^2\theta} \partial_{\phi_*} -\frac{\bfb}{\cos^2\theta}\partial_{\psi_*}  \right]^2 + \frac{\bfr^2}{r^2 \bfp^2}\left[\bfa\bfb\partial_{t_*} + \bfb \partial_{\phi_*} + \bfa \partial_{\psi_*}\right]^2.
\end{align}
Now we require the inner product for $dt_*$ to be future timelike everywhere, i.e. $\bfr^2g^{-1}_{\mfrm,\bfa,\bfb}(dt_*,dt_*)<0$. For $r\geq r_0>> \sqrt{6 \mfrm}$, we compute that 
\begin{align}
    \bfr^2 g_{\mfrm,\bfa,\bfb}^{-1}(dt_*,dt_*) &= \Delta_r F'(r)^2 - 2\frac{(r^2 + \bfa^2)(r^2 + \bfb^2)}{r^2} F'(r) + \frac{\sin^2\theta \cos^2\theta}{\bfp^2}(\bfa^2-\bfb^2) + \frac{\bfa^2\bfb^2}{r^2 \bfp^2}\nonumber\\
    &=\Delta_r F'(r)^2-2(\Delta_r+2\mfrm)F'(r)+\frac{\sin^2\theta\cos^2\theta}{\bfp^2}(\bfa^2-\bfb^2)+\frac{\bfa^2\bfb^2}{r^2\bfp^2}.
\end{align}
Choosing $F'(r)$ such that $F'(r)=\frac{2\mfrm}{\Delta_r+2\mfrm}$ for $r\geq r_0>>\sqrt{6\mfrm}$ we see that 
\begin{align}\label{eq:Timelikedt*}
    \bfr^2 g_{\mfrm,\bfa,\bfb}^{-1}(dt_*,dt_*) &=-4\mfrm +\underbrace{\bfa^2\sin^2\theta+\bfb^2\cos^2\theta}_{\leq(|\bfa|+|\bfb|)^2} +\frac{\bfa^2\bfb^2}{r^2}+\frac{4\mfrm^2\Delta_r}{(\Delta_r+2\mfrm)^2}\leq -2\mfrm +\mathcal{O}(r^{-2}).
\end{align}
when $r$ is large. On the compact set $[r_X,r_0]$, $F(r)$ can then also be chosen to be such that $\bfr^{-2}g^{-1}_{\mfrm,\bfa,\bfb}(dt_*,dt_*)<0$. Moreover, by a tedious direct inspection one can verify equations (\ref{eq:g00Component})-(\ref{eq:gXXComponent}). We have thus shown.

\begin{proposition}
    Let $(M,g_{\bfa,\bfb,\mfrm})$ be a Myers-Perry spacetime. Then $M$ admits coordintates $t_*,r,\phi_*,\psi_*,\theta$ on $M^\circ$, such that $dt_*$ is timelike on all of $M$, $g^{-1}_{\bfa,\bfb,\mfrm}(dt_*,dt_*)<0$ and (\ref{eq:g00Component})-(\ref{eq:gXXComponent}) are satisfied.
\end{proposition}
\begin{remark}
    We remark that the coordinate changes (\ref{eq:MPStartcoordinate})-(\ref{eq:MPStarpsicoordinate}) on $M^\circ$ are harmless with regard to all our previous arguments. When only considerations near $r=r_+$ are relevant, we will often opt for the transformation setting $\chi_+\equiv 1$, $\chi_\partial\equiv 0$ and $F\equiv 0$ everywhere in (\ref{eq:MPStartcoordinate}), (\ref{eq:MPStarphicoordinate}) and (\ref{eq:MPStarpsicoordinate}) which retrieves equations (\ref{eq:MPStartcoordinate2})-(\ref{eq:MPStarpsicoordinate2}). We also remark here that (\ref{eq:TrueRescaledDualMetric}) is directly equivalent to equation (4.3) in \cite{hintz2022sharp}.
\end{remark}

In Proposition \ref{prop:RoughCharacterization} we already had a rough characterization of $\hat{\Box}_{g_{\bfa,\bfb,\mfrm}}(\sigma)$ as both a $b$- and a $sc$-operator on $X$, the following Lemma provides a sharpening of this statement under the assumptions of Definition \ref{def:AdmissibleMetric}, formulated on $T^*M^\circ$. This also justifies the claim made in equation (\ref{eq:CrudeWaveEquationindDimensions}).

\begin{lemma}\label{lem:AdmissibleWaveEquation}
    Let $g$ be a metric that is admissible in the sense of Definition \ref{def:AdmissibleMetric}. The wave operator $\Box_g$ defined as $\Box_g=-|g|^{-1/2}\partial_\mu(|g|^{1/2}g^{\mu\nu}\partial_\nu)\in \text{Diff}^2(M^\circ)$ is given by 
    \begin{align}\label{eq:AdmissibleMetricExpressoin}
        \Box_g=-2\rho\partial_{t_*}Q + \hat{\Box}_g(0)-g^{00}\partial_{t_*}^2,
    \end{align}
    where $Q\in \text{Diff}_b^1(X)$, and $\hat{\Box}_g(0)\in \rho^2\text{Diff}_b^2(X)$. Near the boundary at $\partial X$ these operators are of the form 
    \begin{align}
        Q &= Q_0 + \tilde Q && \tilde Q\in \rho^2\text{Diff}^1_b(X),\label{eq:QForm}\\
        \rho^{-2}\hat{\Box}_g(0)&= L_0+\tilde L, && \tilde L\in \rho^2\text{Diff}^2_b(X),\label{eq:WaveEqDecomposition}
    \end{align}
    where the $\rho$-dilation-invariant operators $Q_0,L_0,L_1$ are given by 
    \begin{align}
        Q_0&=\rho\partial_\rho -\frac{3}{2}, && \text{(indicial operator)} \label{eq:IndicialOperator}\\
        L_0&=-(\rho\partial_\rho)^2 + 2\rho\partial_\rho +\cancel{\Delta}, && \text{ (normal operator)}\label{eq:NormalOperator}
    \end{align}
    where $\cancel{\Delta}=\frac{1}{\cos\theta\sin\theta}\partial_\theta\left(\cos\theta\sin\theta\,\partial_\theta \right)+\frac{1}{\cos^2\theta}\partial_\phi^2 + \frac{1}{\cos^2\theta}\partial_\psi^2$ is the spherical Laplacian on $\mathbb{S}^3$ (in Hopf coordinates).
\end{lemma}
\begin{remark}
    Note that this is the analogue of the admisibility Lemma 2.7 in \cite{hintz2022sharp}. A crucial difference here is that we did not list the mass term in (\ref{eq:WaveEqDecomposition}). In full we could write (\ref{eq:WaveEqDecomposition})
    \begin{align}
    \rho^{-2}\hat{\Box}_g(0)&= L_0+\rho^2L_2 +\tilde L, & \text{ with } L_2=2\mfrm (\rho \partial_\rho)^2 \;\;\;\;\;\;\;\text{ (mass term)}.
    \end{align}
    However due to $\rho^2$-decay, the mass-term here can be fully absorbed in to $\tilde{L}$. It turns out that this feature, combined with the fact that the indicial roots to (\ref{eq:IndicialOperator}) is here fractional, are what will set apart the late-time asymptotics of wave decay in the $1+4$-dimensional Myers-Perry case from the Kerr case. Within the context of resolvent expansions one can indicate this as the absence of a strong Huygens principle.
\end{remark}
\begin{proof}
    We note that 
    \begin{equation}
        \Box_g= - g^{\mu\nu}\partial_{\mu}\partial_\nu -\frac{\partial_\mu(g^{\mu\nu}|g|^{1/2})}{|g|^{1/2}}\partial_\nu,
    \end{equation}
    hence the coefficients of the 2nd order derivatives of $\Box_g$ are exactly equal to those of $g^{\mu\nu}$. This verifies the first order term in the indicial operator and the 2nd order terms in the normal operator $L_0$ and mass term $L_2$. Focusing now on the lower order terms, we find (using Hopf-coordinates $(\phi,\psi,\theta)$ on $\mathbb{S}^3$) that up to leading order     
    \begin{align}\label{eq:DeterminantAsymptotics}
        |g|^{\pm \frac{1}{2}} \in \rho^{\mp 3}(\cos\theta\sin\theta)^{\pm1}(1+\rho^3 C^\infty). 
    \end{align}
     We now track the $t_*,X$-block, let us define this component of $\Box_g$ as 
     \begin{align}
         T&:= -|g|^{-1/2}\partial_{t_*}(|g|^{1/2} g^{0\nu}\partial_\nu)-|g|^{-1/2}\partial_{\nu}(|g|^{1/2}g^{\nu 0}\partial_{t_*})\\
         &=-g^{0\nu}\partial_{t_*}\partial_\nu -|g|^{-1/2}\partial_\nu (|g|^{1/2}g^{\nu 0}\partial_{t_*}),
     \end{align}
     where we used the fact that $\partial_{t_*}$ is Killing in the second line.
     Let us analyze this per variable. The only non-Killing variables of interest are $\rho$ and $\theta$. From definition \ref{def:AdmissibleMetric} we know that $g^{t_*r}\propto -1+\rho^2C^\infty$. Thus the principal 2nd order part of $T$ is given as
     \begin{align}
         -2\rho\partial_{t_*}(\rho\partial_\rho).
    \end{align}
     The second lower order piece comes from 
    \begin{align}
        -|g|^{-1/2}\partial_r (|g|^{1/2}g^{rt_*})\partial_{t_*}&\propto \frac{\rho^2\partial_\rho (\rho^{-3}(\cos\theta\sin\theta)(1+\rho^3C^\infty)(-1+\rho^2C^\infty)}{\rho^{-3}(\cos\theta\sin\theta)(1+\rho^3C^\infty)}\partial_{t_*}\\
        &\propto (+3\rho^5 \rho^{-4} +\mathcal{O}(\rho^3))\partial_{t^*} = (3\rho+\mathcal{O}(\rho^3))\partial_{t_*}.
    \end{align}
    This exactly gives the the $(0,1)$ $b$-contribution to $Q_0$ whilst the $\mathcal{O}(\rho^3)$ remainder can be absorbed into $\tilde{Q}$. Also the $t_*\phi$ $t_*\psi$ and $t_*\theta$ contributions can be absorbed in the remainder term $\tilde{Q}$. \\
    Now to analyze the zero energy operator we need to look at the $\partial_r^2$ operator. In the radial part we have from definition \ref{def:AdmissibleMetric}, that 
    \begin{equation}
        g^{rr}= 1 -2\mfrm \rho^2+\rho^2C^\infty=1+\rho^2C^\infty,
    \end{equation}
    where here, in contrast to Kerr we absorbed the mass into the remainder term.
    Now using again (\ref{eq:DeterminantAsymptotics}) we then see that 
    \begin{align}
        |g|^{1/2}g^{rr}&\propto\rho^{-3}(1+\rho^3 C^\infty)(1+\rho^2C^\infty)=\rho^{-3}+C^\infty.
    \end{align}
     
Now we see, writing out that 
\begin{align}
    -|g|^{-1/2}\partial_r (|g|^{1/2}g^{rr})\partial_r&\propto -\rho^3\partial_r(\rho^{-3}(1+\rho^2C^\infty)\partial_r)\\
    &=-\rho^3(\rho^2\partial_\rho (\rho^{-3}(1+\rho^2C^\infty)\rho^2\partial_\rho)\\
    &=-\rho^4(\rho \partial_\rho (\rho^{-4}(1+\rho^2C^\infty)(\rho\partial_\rho))\\
    &\propto -(\rho\partial_\rho)^2 - [ \rho^5(-3\rho^{-4} \partial_\rho +\rho^{-4}\partial_\rho ] + \rho^2\text{Diff}^2_b(X)\\
    &=-(\rho\partial_\rho)^2+2\rho\partial_\rho +\rho^2\text{Diff}^2_b(X).
\end{align}
Combining this with the final estimates for the Killing vector fields (i.e. we get a contribution to the wave equation of $\rho^2\cancel{\Delta}$ originating from $\mathbb{S}^3$, we see that the expressions $L_0$ and $\tilde{L}$ for $\hat{\Box}_g$ are exactly as stated above.
\end{proof}

From this we move to spectral admissibility. As a corollary of Lemma \ref{lem:AdmissibleWaveEquation}, we immediately have the following statement for $\hat{\Box}_{g_{\bfa,\bfb,\mfrm}(\sigma)}$ on $X$.

\begin{corollary}\label{cor:SpectralFormAdmissibleWaveOpertor}
    Assume $(M,g)$ is a spacetime with a metric that is admissible in the sense of Definition \ref{def:AdmissibleMetric}. We then have
    \begin{align}
    \hat{\Box}_g(\sigma):=e^{i\sigma t_*}\circ\Box_g\circ e^{-i\sigma t}=2i\sigma \rho Q +  \hat{\Box}_g(0) + \sigma^2 g^{00}\in \rho\text{Diff}^2_b(X).\label{eq:GeneralSpectralFamily},
\end{align}
with $Q\in \text{Diff}_b^1(X)$ given by (\ref{eq:QForm}) and (\ref{eq:IndicialOperator}).
\end{corollary}

We then identify the following spectral admissibility criterion, following Definition
\cite{hintz2022sharp}.

\begin{definition} \label{def:SpectrallyAdmissbleMetric}
    (Spectral admissibility criterion) The metric $g$ is spectrally admissibile if the following conditions are stasfied:
    \begin{enumerate}
        \item (Absence of bound states) The nullspace $\hat{\Box}_g(0)$ on $\mathcal{A}^1(X)$ is trivial. 
        \item (Mode stability.) The nullspace of $\hat{\Box}_g(\sigma)$ over $\mathcal{A}^{1}(X)$ is trivial for all $\sigma \in \overline{\mathbb{H}}$. 
        \item (High energy estimates.) There exists $\delta\in \mathbb{R}$ such that for $s\in \mathbb{R}$, $\ell<-\frac{1}{2}, s+\ell>-\frac{1}{2}$ there exists $C>0$ such that the estimate 
        \begin{align}
            \|u\|_{H^{s,\ell}_{b,|\sigma|^{-1}}(X)}\leq C|\sigma|^{-1+\delta}\|\hat{\Box}(\sigma)u\|_{H^{s,\ell+1}_{b,|\sigma|^{-1}}(X)}, \text{ for all }\sigma\in T_{2C}\subset\mathbb{C}, 
        \end{align}
        where $T_{2C}$ truncated strip $T_{2C}\subset \overline{\mathbb{H}}$ with horizontal gapsize $2C$, i.e.
        \begin{align}\label{eq:Truncated2CStrip}
            T_{2C}:=\{\sigma\in \overline{\mathbb{H}}\;|\; |\text{Re }\sigma|\geq C, \text{ Im }\sigma\in [0,1)\}.
        \end{align}
    \end{enumerate}
\end{definition}

For $(M,g_{\bfa,\bfb,\mfrm})$, the last condition in \ref{def:SpectrallyAdmissbleMetric} follows for all $\bfa,\bfb\in \mathbb{R}^2$ satisfying $(|\bfa|+|\bfb|)^2\leq 2\mfrm$ from the semiclassical high energy results in Theorem \ref{thm:SemiclasscialHighEnergyEstimate}, and Corollary \ref{cor:HighEnergyEstimate}. The mode stability assumption for $\overline{\mathbb{H}}\backslash \{0\}$ has been proven for $(M,g_{\bfa,\bfb,\mfrm})$ for all $\sigma\in \mathbb{C}$ with $\text{Im } \sigma>0$ whilst for $\sigma \in \mathbb{R}\backslash\{0\}$ it remains a conditional assumption. Finally we mention condition (1) in Definition \ref{def:SpectrallyAdmissbleMetric} seperately as it serves a special role in the contourshifting argument. \\
In what follows we will subdivide into first proving that the resolvent exists on suitable Sobolev spaces $\hat{\Box}_{g_{\bfa,\bfb,\mfrm}}(\sigma)^{-1}:\mathcal{Y}\to \mathcal{X}$ for all $\sigma \in \overline{\mathbb{H}}$ and use it to produce the forward solution. We will then refine its output at finite $\sigma$ and near $\sigma =0$ to prepare for the resolvent expansion.

\subsubsection{Existence of the resolvent and forward solution.}

The following proposition immediately proves the absence of bound states over $\mathcal{A}^{1}(X)$ for the full subextremal range $(|\bfa|+|\bfb|)^2<2\mfrm$, i.e. condition (1) in Definition \ref{def:SpectrallyAdmissbleMetric}.

\begin{proposition}\label{prop:NoBoundStates}
    Let $(M,g_{\bfa,\bfb,\mfrm})$ be a subextremal Myers-Perry spacetime satisfying $(|\bfa|+|\bfb|)^2<2\mfrm$. We have that $\ker_{\mathcal{A}^1(X)}\hat{\Box}_{g_{\bfa,\bfb,\mfrm}}(0)=\{0\}$.
\end{proposition}
\begin{remark}
    For an analogous proof for the Kerr black hole see Lemma 3.19 in \cite{hintz2024mode}.
\end{remark}
\begin{proof}
    Starting point here is equation (\ref{eq:FullWaveEquation}). Making the Ansatz 
    \begin{align}
        u=e^{-i\sigma t}e^{i(m\phi+k\psi)}R(r) S(\theta),
    \end{align}
    we read of its radial and angular parts as 
    \begin{align}
        \frac{1}{r}\partial_r(r\Delta_r\partial_rR) + \left[\frac{(r^2+\bfa^2)^2(r^2+\bfb^2)^2}{r^4\Delta_r}\left(\frac{m\bfa}{r^2+\bfa^2}+\frac{k\bfb}{r^2+\bfb^2}\right)^2-\frac{1}{r^2}(m\bfb+k\bfa)^2-\ell(\ell+2)\right]R=0
    \end{align}
    and the angular equation reads (in Hopf coordinates)
    \begin{align}
        \frac{1}{\cos\theta\sin\theta}\partial_\theta(\cos\theta\sin\theta \partial_\theta S_{kml}) +\left[\ell(\ell+2) -\frac{m^2}{\sin^2\theta}-\frac{k^2}{\cos^2\theta}\right]S_{kml}=0.
    \end{align}
    Multiplying the angular equation by $\overline{S}\sin\theta\cos\theta$. 
    \begin{align}
        \ell(\ell+2)\|S_{kml}\|^2 = \int^{\pi/2}_0\left( |S'_{kml}(\theta)|^2+\frac{m^2}{\sin^2\theta}|S_{kml}(\theta)|^2+\frac{k^2}{\cos^2\theta}|S_{kml}(\theta)|^2\right)\sin\theta\cos\theta d\theta,
    \end{align}
    $\ell(\ell+2)=0$ only for the constant angular mode $\ell = m = k =0$. In fact imposing regularity at $\theta=0$, $\theta=\frac{\pi}{2}$ easily shows that $\ell=|m|+|k|+2n$ with $n\in \mathbb{N}_0$. \\
    We note that the radial equation can be recast in the form 
    \begin{align}
        \frac{1}{r}\partial_r(r\Delta_r\partial_rR) + \left[\frac{(m\bfa(r^2+\bfb^2)+k\bfb(r^2+\bfa^2))^2}{r^4\Delta_r}-\frac{(m\bfb+k\bfa)^2}{r^2}-\ell(\ell+2)\right]R=0 \label{eq:RewrittenRadialEquation}
    \end{align}
    Because we are dealing with real coefficients, we know that the Wronskian current 
    \begin{align}\label{eq:WronksianCurrent}
        J(r):= \frac{r^2\Delta_r}{r}\text{Im}(\overline{R}(r) R'(r)),
    \end{align}
    should be constant, i.e. $J'(r) = 0$. Furthermore, we know that the admissible indicial roots to (\ref{eq:MinkowskiModel}) are $\rho^{-\ell}=r^l$ and $\rho^{\ell+2}=r^{-\ell-2}$ so for a decaying equation we know that
    \begin{align}
        R(r)\sim \mathcal{O}(r^{-\ell-2})
    \end{align}
    and hence it follows that $J(\infty)=0$.
    Switching to starlike coordinates (\ref{eq:MPStartcoordinate2})-(\ref{eq:MPStarpsicoordinate2}) we see that 
    \begin{align}
        e^{im\phi + k\psi}=e^{i (m \phi_* +k\psi_*)}\exp\left(-i\int^\infty_{r_+} \frac{m\bfa (r^2+\bfb^2)+k\bfb (r^2+\bfa^2)}{r^2\Delta_r}dr\right)
    \end{align}
    so smoothness across $\mathcal{H}^+$ imposes 
    \begin{align}
        R(r) = \exp\left(i\int^\infty_{r_+} \frac{m\bfa (r^2+\bfb^2)+k\bfb (r^2+\bfa^2)}{r^2\Delta_r}dr\right)\tilde{R}(r)\label{eq:RadialSolutionHorizonCondition}
    \end{align}
    with $\tilde R(r)$ smooth at $\{r=r_+\}$. 
    Plugging (\ref{eq:RadialSolutionHorizonCondition}) back into (\ref{eq:WronksianCurrent}) we get 
    \begin{align}
        \text{Im}(\overline R R') =\frac{m\bfa(r^2+\bfb^2)+k\bfb(r^2+\bfa^2)}{r^2\Delta_r}|\tilde R|^2+\text{Im}(\overline{\tilde R}\tilde R')
    \end{align}
    and thus 
    \begin{align}
        J(r) = \frac{m\bfa(r^2+\bfb^2)+k\bfb(r^2+\bfa^2)}{r}|\tilde R|^2+\frac{r^2\Delta_r}{r}\text{Im}(\overline{\tilde R}\tilde R')
    \end{align}
    from which we then obtain
    \begin{align}\label{eq:WronksianCurrentatHorizon}
        J(r_+) = \frac{m\bfa (r_+^2+\bfb^2)+k\bfb(r_+^2+\bfa^2)}{r_+}|\tilde R(r_+)|^2
    \end{align}
    as $r^2\Delta_r = (r^2-r_+^2)(r^2-r_-^2)$.
    Finally let us introduce the horizon angular velocities (already seen as coefficients in the expression of the horizon Killing field (\ref{eq:HorizonKillingField})) as 
    \begin{align}
        \Omega_\bfa = \frac{\bfa}{r_+^2+\bfa^2}, \text{ and }, \Omega_\bfb=\frac{\bfb}{r_+^2+\bfb^2}.
    \end{align}
    We then see that (\ref{eq:WronksianCurrentatHorizon}) can be rewritten as 
    \begin{align}
        J(r_+)=\frac{(r^2+\bfa^2)(r^2+\bfb^2)(m\Omega_\bfa+k\Omega_\bfb)}{r_+}|\tilde R(r_+)|^2.
    \end{align}
    Consequently of $m\Omega_\bfa+k\Omega_\bfb\neq 0$, conservation of $J$ together with $J(\infty)=0$ forces
    \begin{align}
        \tilde R(r_+)=0.
    \end{align}
    and thus that $R \equiv 0$ by regular-singular uniqueness for the horizon-smooth branch.\\
    It remains to also exclude the case where $m\Omega_\bfa+k\Omega_\bfb = 0$. Assuming that this holds we still wish to show that $R\equiv 0$. From the current conservation of $J(r)$ we know that 
    \begin{align}
        \int^\infty_{r_+}r\Delta_r|R'(r)|dr +\int^\infty_{r_+}\left[\frac{1}{r}\left((m\bfb+k\bfa)^2-\frac{m\bfb(r^2+\bfa)^2+k\bfa(r^2+\bfb^2)}{r^2\Delta_r}\right)+\ell(
        \ell+2)r\right]|R|^2 dr =0,
    \end{align}
    thus the claim would follow if we can show that 
    \begin{align}\label{eq:FinalInequalityNullEnergy}
    (m\bfb+k\bfa)^2-\frac{m\bfb(r^2+\bfa)^2+k\bfa(r^2+\bfb^2)}{r^2\Delta_r}\geq 0.
    \end{align}
    This follows immediately when we use that $r_+^2r_-^2=\bfa^2\bfb^2$ (\ref{eq:VietaRelation}) and recast (\ref{eq:FinalInequalityNullEnergy}) in the following form 
    \begin{align}
        (m\bfb+k\bfa)^2-\frac{m\bfb(r^2+\bfa)^2+k\bfa(r^2+\bfb^2)}{r^2\Delta_r}= \frac{(m\bfb+k\bfa)^2r^2(r^2_+-r_-^2)}{r_+^2(r^2-r_-^2)}\geq 0.
    \end{align}  
\end{proof}
We have the following theorem due to Vasy (\cite{vasy2021resolvent}, Theorem 1.1). The framework proving it is closely related to the $sc-b$-framework we employed in Section \ref{subsec:scbcalculus} coming from its companion paper \cite{vasy2021limiting}, and it is known as the \textit{resolved} $sc-b$-calculus. It currently lies outside of the scope of this work, but it roughly involves turning calculus on ${}^{sc,b}\overline{T^*}$ with a new calculus by blowing up ${}^{sc,b}\overline{T^*}X\times \overline{\mathbb{H}}_{\sigma}$ near $\sigma =0$. Later in this work we will work with the resolved space $X_{res}^+$ (see Definition \ref{def:PositiveResolvedSpace} which can be loosely interpreted of the base projection of the sc-b resolved cotangent bundle ${}^{sc,b,res}\overline{T^*}X$. This result enables one to provide uniform estimates on $\sigma$-independent $b$-spaces near $\sigma=0$.

\begin{theorem}\label{thm:SpectralEstimate}
    For $s,\ell,\nu\in \mathbb{R}$ with $s>\frac{1}{2}$ $\ell<-\frac{1}{2}$, $s+\ell>-\frac{1}{2}$, $\ell-\nu\in (-2,0)$ the bound
    \begin{equation}
        \|(\rho+|\sigma|)^\nu u\|_{H_b^{s,\ell}(X)} \leq C\|(\rho+|\sigma|)^{\nu-1}\hat{\Box}(\sigma)u\|_{H_b^{s,\ell+1}(X)},\label{eq:SemiFredholmEstimate}
    \end{equation}
    holds for $\text{Im}\,\sigma\geq 0$, with $|\sigma|\leq \sigma_0 <<1$.
\end{theorem}
\begin{remark}
    We did did have to weaken Theorem \ref{thm:SpectralEstimate} slightly with the horizon threshold condition $s>\frac{1}{2}$ coming in from interior propagation, as Theorem 1.1 in \cite{vasy2021resolvent} assumes a no horizon condition. The condition $\ell-\nu \in (-2,0)$ comes from the indicial gap of the Mellin transformed $b$-normal operator of $\hat{\Box}(0)$, i.e. $\mathscr{M} {}^bN(\hat{\Box}(0))$ which is shifted from $(-1,0)$ on Kerr to $(-2,0)$ to $1+4$-dimensional Myers-Perry.  
\end{remark}

We now prove that $\hat{\Box}_g(\sigma)$ Fredholm on suitable Sobolev spaces. The key ingredient for its proof is modeled after \cite{hintz2026non} Corollary 9.43 using augmented wave operators.

\begin{proposition}\label{prop:IndexZeroResolvent}(Existence of the Resolvent)
    Let $(\bfa,\bfb)\in D_{\mfrm}$ and fix 
    \begin{align}
        s>\frac{1}{2}, \;\ell<-\frac{1}{2},\;s+\ell>-\frac{1}{2},
    \end{align}
    and define 
    \begin{align}
        \mathcal{X}^{s,\ell}_\sigma:=\{u\in \overline{H}^{s,\ell}_b(X)\;|\; \hat{\Box}(\sigma)u\in \overline{H}^{s,\ell+1}_b(X)\},
    \end{align}
     equipped with the graph norm
    \begin{align}
        \|u\|_{\mathcal{X}^{s,\ell}_\sigma}:=\|u\|_{\overline{H}^{s,\ell}_b}+\|\hat{\Box}(\sigma)u\|_{\overline{H}^{s,\ell+1}_b(X)}.
    \end{align}

    Then
    \begin{align}
        \hat{\Box}(\sigma):=\hat{\Box}_{g_{\bfa,\bfb,\mfrm}}(\sigma):\mathcal{X}^{s,\ell}_\sigma\to\overline{H}^{s,\ell+1}_b(X)
    \end{align}
    is Fredholm of index zero on all for all of $\sigma \in \overline{\mathbb{H}}\backslash$ and under the mode stability assumptions of Definition \ref{def:SpectrallyAdmissbleMetric} invertible. 
    At $\sigma =0$, for $s>\frac{1}{2}$ and $-2<\ell_0<0$, define 
    \begin{align}
        \tilde{\mathcal{X}}_0^{s,\ell_0}:=\{u\in H^{s,\ell_0}_b(X)\;|\;\hat{\Box}(0) u \in H^{s,\ell_0+2}_b(X)\},
    \end{align}
    equipped with its graph norm. Then
    \begin{align}
        \hat{\Box}(0):X_0^{s,\ell_0}\to H_b^{s,\ell_0+2}(X)
    \end{align}
    is invertible 
    \begin{align}
        \|\hat{\Box}(0)^{-1}f\|_{H^{s,\ell_0}_b}\leq C\|f\|_{H^{s,\ell_0+2}_b}.
    \end{align}

\end{proposition}

\begin{proof}
     By Corollary (\ref{cor:FredholmFiniteSigmabStatement}), equations (\ref{eq:FredholmbStatementFiniteSigma}) and (\ref{eq:FredholmbStatementFiniteSigmaAdjoint}) $\hat{\Box}(\sigma)$ and its adjoint satisfy Fredholm estimates with compact remainders in compact subsets of $\overline{\mathbb{H}}\backslash \{0\}$. We thus know that for each $\sigma \in \overline{\mathbb{H}}\backslash \{0\}$, $\dim \ker (\hat{\Box}(\sigma))<\infty$ and $\dim \text{coker} (\hat{\Box}(\sigma))<\infty$. We claim that the function 
     \begin{align}
         \sigma \mapsto \text{Ind} \hat{\Box}(\sigma) = \dim \ker \hat{\Box}(\sigma) -\dim \text{coker} \hat{\Box}(\sigma),
     \end{align}
    is locally constant on $\mathbb{H}\backslash \{0\}$. To see this fix a $\sigma_0\in \overline{\mathbb{H}}\backslash \{0\} $ first. Because $\hat{\Box}(\sigma)$ is Fredholm, let us set
    \begin{align}
        n =\dim \ker (P(\sigma_0)), \text{ and } m=\dim \text{coker} P(\sigma_0).
    \end{align}
    Now choose finite-rank maps with smooth data 
    \begin{align}
        R_-:\mathbb{C}^m\to H^{s,\ell+1}_b(X) \text{ and } R_+:\overline{H}_b^{s,\ell}(X)\to \mathbb{C}^n,
    \end{align}
    such that $\text{Ran } R_-$ satisfies $\text{Ran } R_- \oplus \text{Ran }(\hat{\Box}(\sigma)) = \overline{H}^{s,\ell+1}_b(X)$ and such that $R_+|_{\ker (P(\sigma_0))}:\ker P(\sigma_0)\to \mathbb{C}^n$ is an isomorphism. Then consider the augmented operator
    \begin{align}
        \tilde{P}(\sigma) :=\begin{pmatrix}
            \hat{\Box}(\sigma) & R_-\\
            R_+ & 0
        \end{pmatrix}\to \mathcal{X}_\sigma^{s,\ell}\oplus \mathbb{C}^m \to \overline{H}^{s,\ell+1}_b(X)\oplus \mathbb{C}^n,
    \end{align}
    which is invertible at $\sigma_0$ by construction. We claim that $\tilde{P}$ remains invertible nearby $\sigma_0$. Assume by contradiction that first that injectivity fails arbitrarily close to $\sigma_0$, then we can find a sequence $\sigma_j \in $ $0\neq(u_j,c_j) \in \ker \tilde{P}(\sigma_j)\subset  \mathcal{X}_\sigma^{s,\ell} \oplus \mathbb{C}^m$.Normalize these elements by 
    \begin{align}
        \|u_j\|_{\overline{H}^{-N,\ell-\delta}_b}+|c_j|=1.
    \end{align}
    Now (\ref{eq:FredholmbStatementFiniteSigma}) together with $\tilde{P}(\sigma_j)u_j=-R_-c_j$ bounds implies that $\|u_j\|_{\overline{H}^{s,\ell}_b}\leq C$ for some $C>0$ uniformily. After passing to a subsequence we may assume $c_j\to c \in \mathbb{C}^m$ strongly and $u_j \rightharpoonup u \in \overline{H}^{s,\ell}_b$ weakly and $u_j\to u$ strongly in $\overline{H}^{-N,\ell-\delta}_b(X)$, the latter through the compact inclusion $ \overline{H}^{s,\ell}_b(X) \hookrightarrow \overline{H}^{-N,\ell-\delta}_b(X)$. Taking distributional limits then gives 
    \begin{align}
        \hat{\Box}(\sigma_0) u +R_-c=0, \;R_+u=0.
    \end{align}
    As $u\in \mathcal{X}^{s,\ell}_{\sigma_0}$ this contradicts invertibility of $\hat\Box(\sigma_0)$. For the assumption that surjectivity fails arbitrarily close to $\sigma_0$ the corresponding analogous argument for the adjoint $\tilde{P}(\sigma)^*$ for the annihilator of $\text{Ran }\tilde P(\sigma_0)$ in $\dot{H}^{-N,-\ell-1-\delta}_b(X)$ gives a contradiction as well.
    We conclude that $\tilde{P}(\sigma)$ is invertible near $\sigma_0$, and that 
    \begin{align}
        0=\ind (\tilde{P}(\sigma))= \ind(\hat{\Box}(\sigma))+m-n.
    \end{align}
    We conclude 
    \begin{align}
        \ind \hat{\Box}(\sigma) =n-m=\ind \hat{\Box}(\sigma_0).
    \end{align}
    Local constancy and connectedness of the punctured closed upper half-plane now imply 
    \begin{align}
        \ind \hat{\Box}(\sigma ) = 0 \text{ for all } \sigma \in \overline{\mathbb{H}}\backslash\{0\}
    \end{align}
    To determine the index choose a frequency $\sigma_*\in \overline{\mathbb{H}}\backslash\{0\}$ large enough in a fixed horizontal strip. The direct estimate (\ref{eq:DirectHighEnergyEstimate}) of Corollary \ref{cor:HighEnergyEstimate} implies that 
    \begin{align}
        \ker P(\sigma_*)=\{0\}
    \end{align}
    Similarly, we know that for $v\in \text{Ann}(\text{Ran }  \hat{\Box}(\sigma))$ that $v\in \dot{H}^{-s,-\ell-1}_b(X)$, and satisfies the supported adjoint equation 
    \begin{align}
        P(\sigma_*)^*v=0.
    \end{align}
    Thus  $\text{Ann}(\text{Ran }  \hat{\Box}(\sigma))$ is trivial by (\ref{eq:AdjointHighEnergyEstimate}). As $\text{Ran }P(\sigma_*)$ is closed, $P(\sigma_*)$ is an isomorphism and hence $\text{ind }P(\sigma_*) =0 $. By local constancy this we conclude that 
    \begin{align}
        \ind \hat{\Box}(\sigma) = 0 \text{ for all } \sigma \in \overline{\mathbb{H}}\backslash \{0\}.
    \end{align}
    It remains to show the statement at $\sigma = 0$ as well. Fix $-2<\ell_0<0$ and choose $\ell=-1$, $\nu = -1 -\ell_0$. These choices satisfy its hypotheses $s-1>-\frac{1}{2}$ and $\ell-\nu =\ell_0\in (-2,0)$. The preceding argument gives inveritibility of $P(i\epsilon)$ on $\mathcal{X}^{s,-1}_{i\epsilon}$ for arbitrarily small $i\epsilon$. Hence given a $f\in H^{s,\ell_0+2}_b(X)\subset H^{s,0}_b(X)$, let 
    \begin{align}
    u_\epsilon= \hat{\Box}(i\epsilon)^{-1}f,\; v_\epsilon = (\rho+\epsilon)^\nu u_\epsilon.
    \end{align}
    Then Theorem \ref{thm:SpectralEstimate} yields 
    \begin{align}
        \|v_\epsilon\|_{H^{s,-1}_b}\leq \|(\rho+\epsilon)^{-\ell_0-2}f\|_{\overline{H}^{s,0}_b}\leq C'\|f\|_{\overline{H}^{s,\ell_0+2}_b},\label{eq:WeakLowerSemiContinuityEstimate}
    \end{align}
    where the last bound is uniform due to 
    \begin{align}
        (\rho+\epsilon)^{-\ell_0-2}f=\left(\frac{\rho}{\rho+\epsilon}\right)^{\ell_0+2}\rho^{-\ell_0-2}f,
    \end{align}
    and $\left(\frac{\rho}{\rho+\epsilon}\right)$ factor has uniformly bounded $b$-derivatives. Passing to a weakly convergent subsequence, 
    \begin{align}
        v_\epsilon \rightharpoonup v_0 \text{ in } H^{s,-1}_b(X), \; u_0=\rho^{-\nu}v_0\in H^{s,\ell_0}_b(X).
    \end{align}
    On compact subsets subsets away from $\partial_+X$ we have that the multipliers converge smoothly. Hence we can take $u_\epsilon \to u_0$ distributionally and thus 
    \begin{align}
       f = \hat\Box(i\epsilon)u_\epsilon \to \hat{\Box}(0)u_0,
    \end{align}
    showing that $\hat{\Box}(0):\mathcal{X}_0^{s,\ell}\to H^{s,\ell}_b(X)$ is surjective.    
    Weak lower semicontinuity in (\ref{eq:WeakLowerSemiContinuityEstimate}) then gives 
    \begin{align}
        \|u_0\|_{H^{s,\ell_0}_b}\leq C'\|f\|_{\overline{H}_b^{s,\ell_0+2}}. 
    \end{align}
    Finally horizon and zero-frequency indicial regularity improve every solution to $\hat{\Box}(0)u=0$ in $H^{s,\ell_0}_b(X)$ with $-2<\ell_0<0$ to $\mathcal{A}^1(X)$ as in the proof of Proposition \ref{prop:NoBoundStates} and thus $u=0$. We thus find that 
    \begin{align}
        \ker \hat\Box(0) =\{0\}.
    \end{align}
    and hence that $\hat\Box(0):\mathcal{\tilde{X}}^{s,\ell}_0\to H^{s,\ell+1}_b(X)$ is invertible.
\end{proof}

The main problem with our index proposition in Proposition \ref{prop:IndexZeroResolvent}, is that in the context of the Payley-Wiener theorem the resolvent $\hat{\Box}(\sigma)^{-1}:\overline{H}^{s,\ell+1}_b(X)\to \mathcal{X}_\sigma^{s,\ell}$ maps into spaces depending on $\sigma\in \mathbb{C}\backslash \{0\}$. Moreover at $\sigma=0$ the domain space also shifts to $H^{s,\ell_0+2}_b(X)$ with $\ell_0$ chosen in the indicial gap, $\ell_0\in (-2,0)$. Fortunately, the $\sigma=0$ part can largely be ignored due to the high energy estimates giving bounds on arbitrarily small neigbhourhoods near $0$, upon which one can deform the contour slightly around this point. The aim is now to find coarser spaces that are independent of $\sigma\in \mathbb{C}$ for which the resolvent exists and is holomorphic for $\sigma\in \overline{\mathbb{H}}\backslash \{0\}$.

\begin{lemma}\label{lem:HolomorphyLemma}
    Let $s>\frac{1}{2}$. Let $\hat{\Box}(\sigma)^{-1}=\hat{\Box}^{-1}_{g_{\bfa,\bfb,\mfrm}}(\sigma)$ be defined for $\sigma\in \mathbb{C}\backslash\{0\}$ as in Proposition \ref{prop:IndexZeroResolvent}. Then $\hat{\Box}(\sigma)^{-1}$ is holomorphic for $\sigma \in \overline{\mathbb{H}}\backslash \{0\}$ with values in 
    \begin{align}
        \mathcal{L}(H^{s,0}_b(X),H^{s,-1}_b(X)).
    \end{align}
    Consequently, it is also holomorphic for  $\sigma \in \overline{\mathbb{H}}\backslash \{0\}$ with values in 
     \begin{align}\label{eq:ExtrasourceEstimate1}
        \mathcal{L}(H^{s,1}_b(X),H^{s,-1}_b(X)).
    \end{align}
    As a consequence $\overline{H}\backslash \{0\}\ni\sigma \mapsto \hat{\Box}^{-1}(\sigma)\in \mathcal{L}(H^{s,0}_b(X),H^{s,-1}_b(X))$ is continuous up to $\mathbb{R}\backslash \{0\}$ and for $0<|\sigma|<\sigma_0$, $\hat{\Box}(\sigma)^{-1}$ it satisfies uniform bounds
    \begin{align}\label{eq:ExtraSourceEstimate2}
        \sup_{\substack{\sigma<|\sigma|\leq \sigma_0\\\text{Im }\sigma\geq 0}}\|\hat{\Box}(\sigma)^{-1}\|_{\mathcal{L}(H^{s,1}_b(X),H^{s,-1}_b(X))}\leq C'.
    \end{align}
    
\end{lemma}
\begin{proof}
    For bookkeeping let us define 
    \begin{align}
        \mathcal{X}^p := H^{p,-1}_b(X),\;\text{ and } \;\mathcal{Y}^p:=H^{p,0}_b(X).
    \end{align}
    We will work with orders $p=s,s+1,s+2$. 
    For every compact $K\subset\{\text{Im }\sigma >0\}$, we have 
    \begin{align}\label{eq:ComplexSupremum}
        \sup_{\sigma\in K}\|\hat{\Box}(\sigma)^{-1}\|_{\mathcal{L}(\mathcal{Y}^p,\mathcal{X}^p)}\leq C_{K,p}.
    \end{align}
    This follows from (\ref{eq:FredholmbStatementFiniteSigma}) and injectivity of $\hat{\Box}(\sigma)$. If it failed we would find again 
    \begin{align}\label{eq:UpgradetoLowerOrder}
        \sigma_j \to \sigma _0 \in K, \; \|u_j\|_{\mathcal{X}^p}=1, \text{ and } P(\sigma_j)u_j\to 0 \text{ in }\mathcal{Y}^p.
    \end{align}
    After passing to a subsequence, we can get $u_j\rightharpoonup u$ weakly in $\mathcal{X}^p$ and strongly in the compact remainder space $H^{-N,-1,-\delta}_b(X)$. Taking limits we get $P(\sigma_0) u =0$ hence $u=0$, contradicting $\|u_j\|=1$ for all $j$.
    We will now establish a resolvent identity for compositions. From (\ref{eq:GeneralSpectralFamily}) we have that 
    \begin{align}
        \hat{\Box}(\sigma) -\hat{\Box}(\sigma')=(\sigma-\sigma')A_{\sigma,\sigma'},\; \text{ with } A_{\sigma,\sigma'}=2i\rho Q+(\sigma+\sigma')g^{00}.
    \end{align}
    Since $Q\in \text{Diff}_b^1(X)$ and $g^{00}\in \rho^2C^\infty(X)$ we have that 
    \begin{align}
        A_{\sigma,\sigma'}:\mathcal{X}^{p+1}\to \mathcal{Y}^p,
    \end{align}
   is locally uniformily bounded in $\sigma$ and $\sigma'$.
    Take $f\in \mathcal{Y}^{p+1}$ and set $u_{\sigma'} =\hat{\Box}(\sigma')^{-1}f$. Then $u_{\sigma'} \in \mathcal{X}^{p+1}$ and 
    \begin{align}
        \hat{\Box}(\sigma)u_{\sigma'}= f+(\sigma-\sigma')A_{\sigma,\sigma'}u_{\sigma'}\in \mathcal{Y}_{\sigma'}.
    \end{align}
    We thus see that $u_{\sigma'} \in \mathcal{X}^p$. Then applying $\hat{\Box}^{-1}(\sigma)$ we obtian 
    \begin{align}\label{eq:ComplexDifferentiabilityEquation1}
        \hat{\Box}(\sigma)^{-1}f-\hat{\Box}(\sigma')^{-1}f=-(\sigma-\sigma')\hat{\Box}(\sigma)^{-1}A_{\sigma,\sigma'}\hat{\Box}(\sigma')^{-1}f.
    \end{align}
   From this we conclude that 
   \begin{align}\label{eq:ComplexDifferentiabilityEquation2}
       \|\hat{\Box}(\sigma)^{-1}-\hat{\Box}(\sigma')^{-1}f\|_{\mathcal{X}^p}\leq C_K |\sigma -\sigma'|\|f\|_{\mathcal{Y}^{p+1}},\;\;\sigma,\sigma'\in K.
   \end{align}
    Refining to a dense set of more regular data, let us take $f\in \mathcal{Y}^{s+2}$, $\sigma\in \mathbb{H}\backslash \{0\}$ and then apply (\ref{eq:ComplexDifferentiabilityEquation2}) with $p=s+1$ we get 
    \begin{align}
        \hat{\Box}(\sigma_0+z)^{-1}f\to \hat{\Box}(\sigma_0)^{-1}f \text{ in } \mathcal{X}^{s+1}
    \end{align}
   when taking the limit $z\to 0$.Then applying (\ref{eq:ComplexDifferentiabilityEquation1}) at $\sigma_0$ on the left
   \begin{align}\label{eq:ComplexDifferentiabilityEquation3}
       \frac{\hat{\Box}(\sigma_0+z)^{-1}f-\hat{\Box}(\sigma)^{-1}f}{z}=-\hat{\Box}(\sigma_0)^{-1}A_{\sigma_0+z,\sigma_0}\hat{\Box}(\sigma_0+z)^{-1}f.
   \end{align}
   Oberserve that 
   \begin{align}
       A_{\sigma_0+z,\sigma_0}=\partial_\sigma\hat{\Box}(\sigma)|_{\sigma_0}+zg^{00}\to \partial_\sigma\hat{\Box}(\sigma)|_{\sigma_0} \text{ in } \mathcal{L}(\mathcal{X}^{s+1},\mathcal{Y}^s),
   \end{align}
from which we take the limit in (\ref{eq:ComplexDifferentiabilityEquation3}) to get 
\begin{align}
    \frac{\hat{\Box}(\sigma_0+z)^{-1}f-\hat\Box(\sigma_0)^{-1}f}{z}\to - \hat{\Box}(\sigma_0)^{-1} \partial_\sigma\hat{\Box}(\sigma)|_{\sigma=\sigma_0} \hat{\Box}(\sigma_0)^{-1}f \text{ in } \mathcal{X}^{s}.
\end{align}
and we have that
\begin{align}
    \|\hat{\Box}^{-1}(\sigma_0) \partial_\sigma\hat{\Box}(\sigma)|_{\sigma=\sigma_0} \hat{\Box}^{-1}(\sigma_0)f\|_{\mathcal{X}^s}\leq C|z|\|f\|_{\mathcal{Y}^{s+2}}.
\end{align}
Now we extend by density, of $\mathcal{Y}^{s+2}$ lies dense in $\mathcal{Y}^s$. Thus for $f\in \mathcal{Y}^s$, choose $f_j\in \mathcal{Y}^{s+2}$ with $f_j\to f$ in $\mathcal{Y}^s$, we then get by (\ref{eq:ComplexSupremum}) that 
\begin{align}
    \sup_{\sigma\in K}\|\hat{\Box}(\sigma)^{-1}f_j-\hat{\Box}(\sigma)^{-1}f\|_{\mathcal{X^s}}\leq C_{K,s}\|f_j-f\|_{\mathcal{Y}_s}
\end{align}
Finally we proceed to upgrade these uniform bounds to operator-norm holomorphy. For $\sigma \in \mathbb{H}$ choose a closed $\overline{D}(\sigma_0,r)\subset\{\text{Im }\sigma>0\}$, and let $C_{\sigma_0,r}$ be such that $\|\hat{\Box}^{-1}(\zeta)\|_{\mathcal{L}(\mathcal{Y}^s,\mathcal{X}^s)}\leq C_{\sigma_0,r}$ for all $z\in \partial \overline{B(\sigma_0,r)}=S(\sigma_0,r)$
Defining $B_n:\mathcal{Y^s\to \mathcal{X}^s}$
\begin{align}
    B_n f:=\frac{1}{2\pi i}\int_{S(\sigma_0,r)} \frac{R(\zeta) f}{(\zeta-\sigma_0)^{n+1}}d\zeta
\end{align}
we then have $\|B_n\|_{\mathcal{L}(\mathcal{Y}^s,\mathcal{X}^s)}\leq C_{\sigma_0,r}r^{-n}$.
Then by the vector valued Cauchy formula we have that 
\begin{align}
    \hat{\Box}(\sigma)^{-1}f =\sum_{n=0}^\infty B_n f (\sigma-\sigma_0)^n \text{ in } \mathcal{L}(H^{s,0}_b(X),H^{s,-1}_b(X)).
\end{align}
Finally, to prove continuity down to the real line notice that injectivity and the compact Fredholm estimates still give (\ref{eq:ComplexSupremum}) for $\sigma\in \mathbb{R}\backslash \{0\}$, from which we get
\begin{align}
    \|(\hat{\Box}(\sigma)^{-1}-\hat{\Box}(\sigma')^{-1})f\|_{\mathcal{X^s}}\leq C_K |\sigma -\sigma'|\|f\|_{\mathcal{Y}^{s+1}}.
\end{align}
then for $f\in \mathcal{Y^{s}}$ choose $f_j \in \mathcal{Y}^{s+1}$ with $f_j\to f$ in $\mathcal{Y}^s$ to get 
\begin{align}
    \|(\hat{\Box}(\sigma)^{-1}-\hat{\Box}(\sigma')^{-1})f\|_{\mathcal{X^s}}\leq 2C_K \|f-f_j\|_{\mathcal{Y}^s}+C_K|\sigma-\sigma_0|\|f_j\|_{\mathcal{Y}^{s+1}},
\end{align}
which goes to $0$ when taking $\sigma\to \sigma_0$, $f_j\to f$.
For the uniform bounds at $\sigma =0$, we invoke Theorem \ref{thm:SpectralEstimate}, with $\ell=-1$, $\nu=0$, which are admissible because $s+\ell=s-1>-\frac{1}{2}$.
Applying the estimate (\ref{eq:SemiFredholmEstimate}) $u:=\hat{\Box}(\sigma)^{-1}f$ then gives
\begin{align}
    \|\hat{\Box}(\sigma)^{-1}f\|_{H^{s,-1}_b}\leq C\|(\rho+|\sigma|)^{-1}f\|_{H^{s,0}_b},
\end{align}
uniformily for $0<|\sigma|\leq \sigma_0$. Now writing again
\begin{align}
    (\rho+|\sigma|)^{-1}f =\frac{\rho}{\rho+|\sigma|}\rho^{-1}f.
\end{align}
We then immediately conclude
\begin{align}
    \sup_{\substack{0<|\sigma|\leq \sigma_0\\\text{Im } \sigma\geq 0}}\|\hat{\Box}(\sigma)^{-1}\|_{\mathcal{L}(H^{s,1}_b(X),H_b^{s,-1}(X))}<C'.
\end{align}

\end{proof}
\begin{remark}
    Note that we included the extra source weight in (\ref{eq:ExtrasourceEstimate1}) to conclude the uniform bounds in (\ref{eq:ExtraSourceEstimate2}).
\end{remark}

\begin{proposition}(Polynomial bounds on the resolvent.) \label{prop:PolynomialEnergyBounds}Choose $s>\frac{1}{2}$, and let $\hat{\Box}(\sigma)^{-1}=\hat{\Box}_{g_{\bfa,\bfb,\mfrm}}(\sigma)^{-1}:H^{s,1}_b(X)\to H_b^{s,-1}(X)$. Then for every $s>\frac{1}{2}$ and $\nu>0$ there exists $C_{s,\nu}>0$ such that 
\begin{align}
    \|\hat{\Box}(\sigma)^{-1}\|_{\mathcal{L}(H^{s,1}_b(X),H_b^{s,-1}(X)}\leq C_{\nu,s}\langle \sigma\rangle^{\lceil s \rceil}, \;0\leq \text{Im }\sigma\leq \nu.
\end{align}
Moreover, for $s=1$, there exist $\gamma_E\geq 1$ and $C>0$ such that
\begin{align}
    \|\hat{\Box}(\sigma)^{-1}f\|_{H^{1,-1}_b(X)}\leq \frac{C}{\text{Im }\sigma}\|f\|_{H_b^{0,1}}, \;\text{Im }\sigma\geq \gamma_E,
\end{align}
for every $f\in H^{1,1}_b(X)$.
Consequently there exists a $C>0$ such that we have the polynomial bound 
\begin{align}\label{eq:GlobalUpperHalfPlaneEstimate}
    \|R(\sigma)\|_{\mathcal{L}(H^{1,1}_b(X),H^{1,-1}_b(X)}\leq C\langle \sigma \rangle,\; \text{Im }\sigma\geq 0.
\end{align}

\end{proposition}
\begin{proof}
    For the first statement, we use the high energy estimates in Corollary \ref{cor:HighEnergyEstimate} with $\ell=-1$ and $h=|\sigma|^{-1}$. We have for $0<\epsilon<-\gamma_-$ that 
    \begin{align}
        \|\hat{\Box}(\sigma)^{-1}f\|_{H^{s,-1}_{b,h}}\leq C_{s,\nu,\epsilon}|\sigma|^{-1+\delta_\Gamma(-\epsilon)}\|f\|_{H^{s,0}_{b,h}},
    \end{align}
    for sufficiently large $|\sigma|$, uniformily on $0\leq \text{Im }\sigma \leq \nu$.
    For $s\geq 0$ and $0<h\leq 1$ the ordinary and semiclassical norms satisfy
    \begin{align}
        \|v\|_{H_b^{s,\ell}}\leq Ch^{-s}\|v\|_{H^{s,\ell}_{b,h}},\; \|v\|_{H^{s,\ell}_{b,h}}\leq \|v\|_{H^{s,\ell}_b},
    \end{align}
    and we get 
    \begin{align}
        \|\hat{\Box}(\sigma)^{-1}f\|_{H^{s,-1}_b(X)}&\leq C_{s,\nu,\epsilon}|\sigma|^{s-1+\delta_{\Gamma}(-\epsilon)}\|f\|_{H^{s,0}_b(X)}\\
        &\leq C_{s,\nu,\epsilon}|\sigma|^{s-1+\delta_{\Gamma}(-\epsilon)}\|f\|_{H^{s,1}_b(X)}\\
        &\leq C_{s,\nu,\epsilon}\langle\sigma\rangle^{\lceil s\rceil}\|f\|_{H^{s,1}_b(X)}.
    \end{align}
    For the second claim let us choose $K=\partial_{t_*}$. We work with stationary spatial density $ dg_{\bfa,\bfb,\mfrm}|_X$ on $X$ (recall that $dg_{\bfa,\bfb,\mfrm}|_{X}=r\bfr^2 dr\cos\theta\sin\theta d\theta d\phi d\psi$ on the interior $X^\circ$ in Hopf coordinates). Choose a stationary future timelike vector field $N$ which equals $T$ for sufficiently large $r$\footnote{Such an $N$ exists: on a compact region use $\nabla t_*$ and interpolate to $T$ when $T$ is future timelike using convexity}. Since $T$ is a Killing vector field, the deformation tensor 
    \begin{align}
        \pi_N =\frac{1}{2}\mathcal{L}_N g
    \end{align}
    has compact spatial support. Take $f$ smooth up to $\partial_-X$ and vanishing near spatial infinity. Set 
    \begin{align}
        u=\hat{\Box}(\sigma)^{-1}f,\; U(t_*,x)=e^{-i\sigma t_*}u(x), \; \gamma =\text{Im }\sigma.
    \end{align}
    As in the proof of Lemma \ref{lem:HolomorphyLemma}, especially the argument after (\ref{eq:UpgradetoLowerOrder}) the inverses at different orders $s$ agree; hence the following computations are justified on spatial slabs 
    \begin{align}
        [0,1]_{t_*}\times \Omega_R, \text{ with } \Omega_R= [r_X,R]_r\times \mathbb{S}^3, \text{ and } R>r_+.
    \end{align}
    For the complex scalar field, define its stress-energy tensor through
    \begin{align}
        T[U]_{\alpha\beta}=\text{Re}(\partial_\alpha U \overline{\partial_\beta U})-\frac{1}{2}g_{\alpha\beta}g^{\mu\nu}\partial_\mu U \overline{\partial_\nu U}.
    \end{align}
    Define the energy $E_\sigma(u)$ through the energy density $dg_{\bfa,\bfb,\mfrm}|_X=r\bfr^2 dr dg_{\mathbb{S}^3}$ (see equation (\ref{eq:StationarySplitting}))
    \begin{align}
        e_\sigma(u)=T[U](N,-\nabla t_*)|_{t_*=0},\; E_\sigma(u)=\int_X e_\sigma(u) \,dg_{\bfa,\bfb,\mfrm}|_X
    \end{align}
    Near $\partial_+X$ we have $N=K$, and the mixed terms cancel:
    \begin{align}\label{eq:EnergyIntegrand}
        e_\sigma(u)=\frac{1}{2}\left(-g^{00}|\sigma u|^2+g^{ij}\partial_iu\overline{\partial_ju}\right).
    \end{align}
    where $i,j$ are the spatial indices.
    From (\ref{eq:TrueRescaledDualMetric}) and (\ref{eq:Timelikedt*}) we have that
    \begin{align}\label{eq:EnergyEquivalence}
        -g^{00}\sim \rho^2, \; \text{ and } g^{ij}\partial_iu\overline{\partial_ju}\asymp |\partial_ru|^2+r^{-2}\|d_\omega u\|^2_{g^{-1}_{\mathbb{S}^3}}
    \end{align}
   where we defined
    
    \begin{align}
        d_\omega u = (\partial_\theta u)d\theta +(\partial_{\phi_*}u)d\phi_* +(\partial_{\psi_*}u)d\psi_*.
    \end{align}
    On the remaining compact region $N$ and $-\nabla t_*$ are both timelike by construction, so we have normal coercivity of the stress tensor. Consequentely we have for uniform $|\sigma|\geq 1$ that 
    \begin{align}\label{eq:EnergyEquivalence}
        E_{\sigma}(u)^{1/2}\asymp \|u\|_{H^{1,-1}_b}+|\sigma|\|u\|_{H^{0,-1}_b}.
    \end{align}
 Now we apply the divergence identity with exponential weights using $\Box_g = -\text{div}_g\nabla$, using currents $J^\alpha$ and bulk term $\mathcal{K}_N$ and we define
 \begin{align}
     J^\alpha=-T[U]^\alpha{}_\beta N^\beta, \; \mathcal{K}_N[U]=T[U]^{\alpha\beta}(\pi_N)_{\alpha\beta}.
 \end{align}
And we find 
\begin{align}
    \text{div}_g (e^{-2\gamma t_*}J) =\text{Re}((\Box_g U) \overline{NU})-\mathcal{K}[U].
\end{align}
Then using $J(t) = T[U](N,-\nabla t_*)$ it follows that 
\begin{align}
    \text{div}_g(e^{-\gamma t}J)= e^{-2\gamma t}\left[\text{Re}((\Box_g U)\overline{NU}-\mathcal{K}_N[U]-2\gamma T[U](N,-\nabla t_*)\right].
\end{align}

Now integrating over the slabs $[0,1]_{t_*} \times \Omega_R$, because $U=e^{-i\sigma t_*}u$ the weighted fluxes through $\{t_*=0\}$ and $\{t^*=1\}$ agree and cancel out. Now let
\begin{align}
    N_\sigma u=-i\sigma N^t u +N^j\partial_j u,
\end{align}
and define the truncated energy 
\begin{align}
    E_{\sigma,R}(u)=\int_{\Omega_R} e_\sigma(u) \,dg_{\bfa,\bfb,\mfrm}|_X
\end{align}

we obtain for fixed $R>r_+$,
\begin{align}\label{eq:FluxIdentity}
    2\gamma E_{\sigma,R}+\Phi_{in} +\Phi_R = \text{Re}\int_{\Omega_R}f\overline{N_\sigma u} \,dg_{\bfa,\bfb,\mfrm}|_X-\int_{\Omega_R} \mathcal{K}_\sigma(u)\,dg_{\bfa,\bfb,\mfrm}|_X,
\end{align}
where we set $\mathcal{K}_\sigma(u)=\mathcal{K}_N[U]|_{t_*=0}$ and let $\Phi_{in}$ be the flux through the inner boundary $\{r=r_X\}$ and $\Phi_R$ the flux through the outer boundary $\{r=R\}$.
Since $\pi_N$ has compact support, positivity on $\Omega_R$ gives 
\begin{align}
    |\mathcal{K}_\sigma(u)|\leq C_0 e_\sigma(u).
\end{align}
Moreover, (\ref{eq:EnergyIntegrand}) and (\ref{eq:EnergyEquivalence}) imply 
\begin{align}
    \rho |N_\sigma u|\leq Ce_\sigma(u)^{1/2},
\end{align}
from which we get
\begin{align}\label{eq:FluxTerm1}
    \left|\int_{\Omega_R} f\overline{N_\sigma u}\, dg_{\bfa,\bfb,\mfrm}|_X\;\right|\leq C\|f\|_{H_b^{0,1}}E_{\sigma,R}(u)^{1/2}.
\end{align}
Now note that the fluxes have the correct sign, i.e. $\Phi_{in}\geq 0$ and $\Phi_R=0$. To see the former, note that at $r=r_X$, $\nabla r$ is future timelike, since this boundary lies in the black-hole interior and is a future outflow boundary. Its energy flux is given by 
\begin{align}\label{eq:FluxTerm2}
    -T[U]^r{}_\beta(-N^\beta) = T[U](N,\nabla r)\geq 0, 
\end{align}
i.e. $-J^r\geq 0$ and therefore $\Phi_{in}\geq 0$.
For the outer boundary, take $R$ large enough such that $N\sim \partial_{t_*}$, thus the weighted flux $\Phi_R$ is given here by 
\begin{align}\label{eq:FluxTerm3}
    \Phi_R=\int_{\{r=R\}}[-g^{rt}|\sigma u|^2-\text{Re }(g^{rj}\partial_j u - \overline{i\sigma u})\,r\bfr^2 dg_{\mathbb{S}^3}.
\end{align}
The outgoing asymptotics give $-g^{rt}\geq c >0$. Completing the square shows
\begin{align}
    \Phi_R \geq -C \underbrace{\int_{\{r=R\}}\left(|\partial_ru|^2+r^{-2}|d_\omega u|^2\right) r\bfr^2dg_{\mathbb{S}^3}}_{:=A(R)},
\end{align}
where the r.h.s is finite as $u\in H^{1,-1}_b(X)$, hence taking $R_j\to \infty$ we have that $A(R_j)\to 0$ and hence $\Phi_R\geq 0$.
Taking this limit $R_j\to \infty$ in (\ref{eq:FluxIdentity}) and using (\ref{eq:FluxTerm1}),(\ref{eq:FluxTerm2}) and (\ref{eq:FluxTerm3}) we find
\begin{align}\label{eq:FluxLimit}
    (2\gamma - C_0)E_\sigma(u)\leq C\|f\|_{H^{0,1}_b}E_\sigma(u)^{1/2}.
\end{align}
Thus choose $\gamma_E\geq 1$ sufficiently large such that $2\gamma -C_0 \geq\gamma$ for $\gamma \geq \gamma_E$ in (\ref{eq:FluxLimit}) and using the equivalence (\ref{eq:EnergyEquivalence}) we see that
\begin{align}\label{eq:FinalNormMajoration}
    \|u\|_{H^{1,-1}_b(X)}+|\sigma|\|u\|_{H^{0,-1}_b(X)}\leq \frac{C}{\gamma}\|f\|_{H^{0,1}_b(X)},
\end{align}
Finally approximate any $f\in H^{1,1}_b(X)$ with smooth functions $f_j \in \overline{H}_b^{1,1}(X)$ and set $u_j := \hat{\Box}(\sigma)^{-1}f_j$ then passing to the limit again (\ref{eq:FinalNormMajoration}) gives
\begin{align}
    \|\hat{\Box}(\sigma)^{-1}f\|_{H^{1,-1}_b(X)}\leq \frac{C}{\text{Im }
    \sigma}\|f\|_{H^{0,1}_b(X)}.
\end{align}

\end{proof}
\begin{remark}
   The first assertion follows from the semiclassical high-energy estimates of Corollary \ref{cor:HighEnergyEstimate} together with Lemma \ref{lem:HolomorphyLemma}. For analogous high energy estimates in fixed strips see \cite{hafner2021linear} Theorem 4.3. The second part of the proof uses an exponentially weighted stress-energy argument as in Proposition 9.35 (which is in fact our Proposition \ref{prop:HighEnergyEstimate}) of \cite{hintz2026non}, and related energy estimates near null infinity in \cite{hintz2026non} Lemma 8.1. Another closely related reference is \cite{hintz2025introduction} Chapter 9. 
\end{remark}

With the existence and holomorphy properties of $\hat{\Box}(\sigma)^{-1}:H^{s,1}_b(X)\to H^{s,-1}_b(X)$ Lemma \ref{lem:HolomorphyLemma} and Proposition \ref{prop:PolynomialEnergyBounds} we can then contour-shift on $\mathbb{H}\backslash\{0\}$ to the real axis to produce a forward solution.

\begin{theorem}
    (Forward wave solution) Assume the mode stability condition of Definition \ref{def:SpectrallyAdmissbleMetric}, and let 
    \begin{align}
        f\in C_c^\infty(\mathbb{R}_{t_*};H^{\infty,1}(X)).
    \end{align}
    Choose $T_-\in \mathbb{R}$ such that $f(t_*)=0$ for $t_*<T_-$. Then there exist 
    \begin{align}
        \phi \in C^\infty(\mathbb{R}_{t_*};H^{\infty,-1}_b(X))
    \end{align}
    satisfying 
    \begin{align}
        \Box_g\phi, \; \phi(t_*)=0 \text{ for }  t_*<T_-.
    \end{align}
    Moreover, for every $s>\frac{1}{2}$ and $j\in \mathbb{N}_0$, such that
    \begin{align}
        \sup_{t_*\in \mathbb{R}}\|\partial^j_{t_*}\phi(t_*)\|_{H^{s,-1}_b(X)}<\infty.
    \end{align}
    Writing 
    \begin{align}
        \hat{f}(\sigma)= \int_{\mathbb{R}}e^{i\sigma t_*}f(t_*) dt_*, 
    \end{align}
    the solution is given, for every $\eta>0$ by 
    \begin{align}
        \phi(t_*)&=\frac{1}{2\pi}\int_{\mathbb{R}+i\eta}e^{-i\sigma t_*}\hat{\Box}_{g_{\bfa,\bfb,\mfrm}}(\sigma)^{-1}\hat{f}(\sigma) d\sigma
        \\&=\frac{1}{2\pi}\int_{\mathbb{R}}e^{-i\sigma t_*}\hat{\Box}_{g_{\bfa,\bfb,\mfrm}}(\sigma)^{-1}\hat{f}(\sigma) d\sigma.
    \end{align}
    Both integrals converge absolutely in $H^{s,-1}_b(X)$ for every $s>\frac{1}{2}$. The solution is unique in the class of distributional solutions $C^0(\mathbb{R}_{t_*};H^{1,-1}_b(X))$, which vanish for sufficiently large $t_*$ and have at most exponential growth in $H^{1,-1}_b(X)$ as $t_*\to \infty.$
\end{theorem}
\begin{remark}
    Obviously, we will ultimately show that exponential growth cannot occur and that solutions $\phi$ exhibit $t_*^{-3}$ decay.
\end{remark}
\begin{proof}
    Repeated integration by parts gives for every $s>\frac{1}{2}$, $N\in \mathbb{N}_0$, 
    \begin{align}
        \|\hat{f}(\sigma)\|_{H^{s,1}_b(X)\leq C_{s,N}}e^{-T_- \text{Im }\sigma}\langle\sigma \rangle^{-N},\;\text{Im }\sigma\geq 0.
    \end{align}
    For $\eta>0$, we define 
    \begin{align}
        \phi_\eta(t_*)=\frac{1}{2\pi}\int_{\mathbb{R}+i\eta}e^{-i\sigma t_*}\hat{\Box}(\sigma)^{-1}\hat{f}(\sigma) d\sigma.
    \end{align}
    Now the finite strip bounds of Proposition \ref{prop:PolynomialEnergyBounds} (originally stemming from the high energy estimates in Corollary \ref{cor:HighEnergyEstimate}) we have absoltute convergence in $H^{s,-1}_b(X)$, of $\phi_\eta(t_*)$ and all its time derivatives $t_*$. Holomorphy from Lemma \ref{lem:HolomorphyLemma} and Cauchy's theorem on rectangles shows they are independent of $\eta>0$, i.e. completing the contour of the integrals between two lines at $\text{Im }\sigma = i\eta_1$ and $\text{Im }\sigma=i\eta_2$ have vanishing side integrals. We thus write $\phi=\phi_\eta$.\\
    Using the global upper half-plane estimate (\ref{eq:GlobalUpperHalfPlaneEstimate}) for $N>2$, we get 
    \begin{align}
    \|\phi(t_*) \|_{H^{1,-1}_b(X)}&\leq C_N e^{\eta(t_*-T_-)}\int_{\mathbb{R}}\langle\xi+i\eta\rangle^{1-N}d\xi\\
    &\leq C'_N e^{\eta(t_*-T_-)}\langle \eta\rangle^{2-N}.
    \end{align}
    Letting $\eta\to\infty$, we conclude $\phi(t_*)=0$ for $t_*<T_-$.
    Now, lowering the contour to zero, uniform boundedness and finite strip estimates allow lowering to the real axis by $\eta\to 0$ by dominated convergence. We thus get
    \begin{align}
        \phi(t_*)=\frac{1}{2\pi}\int_{\mathbb{R}}e^{-i\sigma t_*}\Box_{g_{\bfa,\bfb,\mfrm}}(\sigma)^{-1}\hat{f}(\sigma) d\sigma,
    \end{align}
    from Fourier inversion we conclude $\Box_{g_{\bfa,\bfb,\mfrm}}(\sigma)\phi=f$. Finally uniqueness follows from mode stability. Let $w\in C^0(\mathbb{R}_{t_*};H^{1,-1}_b(X))$ be a homogeneous solution, then for sufficiently large $\text{Im }\sigma$ its Laplace transform exists in $H^{1,-1}_b(X)$ and satisfies
    \begin{align}
        \hat{\Box}_g(\sigma)\hat{w}(\sigma)=0.
    \end{align}
    We see that $\hat{w}(\sigma)\in \mathcal{X}^{1,-1}_\sigma$ and injectivity from Proposition \ref{prop:IndexZeroResolvent} implies $\hat{w}(\sigma)=0$. Fourier transforming back then gives $w=0$.
\end{proof}

\subsubsection{Resolvent regularity at nonzero frequencies.}

We will now refine the input of the resolvent operator $\hat{\Box}(\sigma)^{-1}$ following \cite{hintz2022sharp}. One important notational distinction should be made  In \cite{hintz2022sharp}, the index $\alpha$ in \cite{hintz2022sharp} lies within the indicial gap $(0,1)$ and is therefore also used in the eventual leading asymptotics formulations as the next order decay term will lie at least one integer away. In our case the indicial gap is $(0,2)$ but not every argument goes through unchanged. We will take $\alpha \in (0,1)$ whereever needed but try to prove resolvent outputs for an as larger $\alpha$ when possible.

\begin{lemma}\label{lem:LimitingAbsorptionFredholm}
Let us fix $R_0,R_1>0$ and let us consider $\sigma$ in the ''horseshoe'' $A_{1/2}=\{\sigma\in\mathbb{H}\;|\;R_0<|\sigma|<R_1\}\subset \mathbb{C}$,
and set 
\begin{align}
    s>\frac{1}{2}, \;\ell<-\frac{1}{2}, \;s+\ell>-\frac{1}{2}
\end{align}

\begin{align}
    \|u\|_{H_b^{s,\ell}(X)}\leq C \|\hat{\Box}(\sigma)u\|_{H_b^{s,\ell+1}(X)}, \text{ for all } \sigma \in \overline{A_{1/2}},\label{eq:FredholmAfterCompactEmbedding}
\end{align}
where $C>0$.
\end{lemma}
\begin{proof}
    Write $P(\sigma)=\hat{\Box}_g(\sigma)$. From Corollary \ref{cor:FredholmFiniteSigmabStatement} we have 
    \begin{align}
        \|u\|_{\overline{H}^{s,\ell}_b}\leq C_0\left(\|P_\sigma u\|_{\overline{H}^{s,\ell+1}_b}+\|u\|_{\overline{H}_b^{-N,\ell-\delta}}\right),
    \end{align}
    where we note that $\inf_{\sigma\in \mathcal{K}} |\text{Im } \sigma|=0$.
    Assume that the estimate (\ref{eq:FredholmAfterCompactEmbedding}) fails. Then we can select sequences $\{\sigma_j\}_{j\in \mathbb{N}}\in \overline{A_{1/2}}$, and $\{u_j\}_{j\in \mathbb{N}}$, 
    \begin{align}
        \|u_j\|_{\overline{H}^{s,\ell}_b}=1,\;\;\|P_{\sigma_j}u_j\|_{\overline{H}_b^{s,\ell+1}}\to 0.
    \end{align}
    After passing to a subsequence compactness of $\overline{A_{1/2}}$ and the compact inclusion
    \begin{align}
        \overline{H}^{s,\ell}_b(X)\hookrightarrow\overline{H}^{-N,-\ell-\delta}_b(X), 
    \end{align}

    gives $\sigma_j \to \sigma _\infty$, $u_j \rightharpoonup u $ weakly in $\overline{H}^{s,\ell}_b(X)$ and $u_j\to u$ in $\overline{H}^{-N,\ell-\delta}_b(X)$.
    As $P(\sigma)$ is analytic in $\sigma$, $P_{\sigma_\infty}u=0$ distributionaly. As the r.h.s. of $P(\sigma)u=0$ is smooth, we have by homogeneous regularity that $u\in \overline{H}^{\infty,\ell}_b(X)$ for any $\ell<-\frac{1}{2}$, hence by Sobolev embedding (\ref{eq:SobolevEmbedding}) $u\in \mathcal{A}^{\frac{3}{2}-}(X)\subset \mathcal{A}^1(X)$ which is a contradiction with 
    $\ker_{\mathcal{A}^1(X)} \hat{\Box}_{g_{\bfa,\bfb,\mfrm}}(\sigma)=0$.
\end{proof}

\begin{lemma}\label{lem:SigmaDerivativeControl}
    Let $m\in \mathbb{N}_{\geq 0}$ $s,\ell,R_0,R_1$ and $A_{1/2}$ as in lemma \ref{lem:LimitingAbsorptionFredholm}, then we have for the operator $\partial^m_\sigma\hat{\Box}(\sigma)^{-1}:H_b^{s,\ell+1}(X)\to H_b^{s-m,\ell}(X)$ that
    \begin{align}
        \|u\|_{H_{b}^{s-m,\ell}(X)}\leq C \|\partial^m_\sigma \hat{\Box}(\sigma)^{-1} u\|_{H_b^{s,\ell+1}(X)} \text{ for all } \sigma \in \overline{A_{1/2}}.  
    \end{align}
    Semiclassically, a similar estimate holds for the operator \begin{equation*}
    \partial^m_\sigma \hat{\Box}(\sigma)^{-1}:H^{s,\ell+1}_{b,|\sigma|^{-1}}(X)\to |\sigma|^{-1+(m+1)\delta}H^{s-m,\ell}_{b,|\sigma|^{-1}}(X),
    \end{equation*} that is there exists a $C>0$
    \begin{align}
        \|u\|_{H^{s,\ell+1}_{b,|\sigma|^{-1}}(X)}\leq C |\sigma|^{-1+(m+1)\delta} \|u\|_{H^{s-m,\ell}_{b,|\sigma|^{-1}}}(X), \text{ for all } \sigma\in T_{2C}.
    \end{align}
\end{lemma}
\begin{proof}
    This is Proposition 12.10 in \cite{hafner2021linear} or Lemma $2.1$ in \cite{hintz2022sharp}. 
    Notice that by (\ref{eq:GeneralSpectralFamily}) we have that
    \begin{align}
        \partial_\sigma \hat{\Box}(\sigma)&= 2i\rho Q+ 2\sigma g^{00},
    \end{align}
    and hence by equations (\ref{eq:g00Component}) and (\ref{eq:QForm}) we have that $\partial_\sigma \hat{\Box}(\sigma) \in \rho\text{Diff}^1_b(X) + \sigma\rho^2C^\infty (X)$. Now using that $\partial_\sigma\hat{\Box}(\sigma)$ is bounded in norm we have 
    \begin{align}
        0&=\partial_\sigma (\hat{\Box}(\sigma)\hat{\Box}(\sigma)^{-1})\\
        &=\partial_\sigma\hat{\Box}(\sigma)\hat{\Box}(\sigma)^{-1}+\hat{\Box}(\sigma)\partial_\sigma(\hat{\Box}(\sigma)^{-1}),
    \end{align}
    which yields the factorization
    \begin{align}\label{eq:Factorization}
    \partial_\sigma\hat{\Box}(\sigma)=-\hat{\Box}(\sigma)^{-1}\circ \partial_\sigma\hat{\Box}(\sigma)\circ \hat{\Box}(\sigma)^{-1}.
    \end{align}
    Keeping track of this in a chain of mapping spaces we see that
    \begin{align}
        H^{s,\ell+1}_{b,|\sigma|^{-1}}(X) \xrightarrow{\hat{\Box}(\sigma)^{-1}}|\sigma|^{-1+\delta}H^{s,\ell}_{b,|\sigma|^{-1}}(X)\xrightarrow{\partial_\sigma\hat\Box(\sigma)}|\sigma|^\delta H^{s-1,\ell+1}_{b,|\sigma|^{-1}}(X)\xrightarrow{\hat{\Box}(\sigma)^{-1}}|\sigma|^{-1+2\delta}H^{s-1,\ell}_{b,|\sigma|^{-1}}(X).
    \end{align}
    Applying this sequence inductively for $m\in\mathbb{N}$ yields the claim.
\end{proof}

\subsubsection{Resolvent at zero energy.}

The Sobolev embedding (\ref{eq:SobolevEmbedding}) combined with the fact that the size of the indicial gap of $\hat{\Box}(0)$ is of size an increased gapsize $2$\footnote{In contrast to $\mathbb{R}^3$} by Lemma 

\begin{lemma}\label{lem:SobolevEmbeddingClosedGraph}
We have that for $\alpha\in (0,2)$, and $\sigma = 0$ that the output of the resolvent operator $\hat{\Box}(0)^{-1}$ is given by 
\begin{align}
    \hat{\Box}(0)^{-1}:\mathcal{A}^{2+\alpha}(X) \to \mathcal{A}^{\alpha-}(X).
\end{align}
\end{lemma}
\begin{proof}
    Foix $0<\alpha<2$, $f\in \mathcal{A}^{2+\alpha}(X)$. Fix $0<\epsilon<\alpha$ and set $\ell_0 = \alpha-2-\epsilon$. Then $-2<\ell_0<0$ with, and so by Proposition \ref{prop:IndexZeroResolvent} we have that $f\in \mathcal{A}^{2+\alpha}(X)\subset H^{\infty,\alpha-\epsilon}_b(X)=H_b^{\infty,\ell_0+2}(X)$. Then by Theorem \ref{thm:SpectralEstimate} we have $u=\hat{\Box}(0)^{-1} f \in H^{\infty,\ell_0}_b(X)\subset \mathcal{A}^{2+\ell_0}(X)=\mathcal{A}^{\alpha-}(X)$.
\end{proof}

To describe the output of the resolent $\hat{\Box}(\sigma)^{-1}$, $\rho=\sigma=0$, we introduce a resolvent space. This allows us to deal with $\hat\Box_g(\sigma)$ transitioning form a principally $\text{Diff}^2_{sc}(X)$ operator into a $\text{Diff}_b^2(X)$ operator in the limit $\sigma\to 0$ in a uniform manner.

\begin{definition}\label{def:PositiveResolvedSpace}
    The resolved space (for positive frequencies) $X_{res}^+$ is the blow-up of $X\times [0,1)_\sigma$ at $\rho=0,\sigma=0$, i.e. 
    \begin{equation}
        X^+_{res}:= [X\times [0,1)_\sigma; \partial X\times\{0\}_\sigma].
    \end{equation}
    with blowdown map $\beta:X_{res}^+\to [0,1)_\sigma \times X$. Within this context we identify
    \begin{enumerate}
        \item the transition face (tf), lift of $\{0\}_\sigma\times\partial X$, i.e. $\overline{\beta^{-1}(\{0\}_\sigma\times\partial X)}$, with boundary defining function $\rho_{\text{tf}}:=\rho+\sigma$;\\
        \item the b-face (bf), the lift of $[0,1)_\sigma\times \partial X$, i.e. $\overline{\beta^{-1}((0,1)_\sigma\times \partial X)}$, with boundary defining function $\rho_{bf}:=\frac{\rho}{\rho+\sigma}$;\\
        \item the zero face (zf), the lift of $\{0\}\times X$, i.e. $\overline{\beta^{-1}(\{0\}_\sigma\times X^\circ}$, with boundary defining function $\rho_{zf}:=\frac{\sigma}{\sigma+\rho}$.
    \end{enumerate}
\end{definition}
Note that if we are away from the zero face i.e. $\sigma>0$, we can use projective coordinates $\hat{\rho}=\frac{\rho}{\sigma}$ as a boundary defining function for the $b$-face and vice versa for the zero-face we can use boundary defining coordinates $\hat{r}=\frac{\sigma}{\rho}$ if we are away from the $b$-face.\\
As $X_{res}$  is a manifold with corners we will now have to work with conormal spaces that are weighted with respect to all boundaries, i.e. we introduce 
\begin{align}
    \mathcal{A}^{\alpha,\beta,\gamma}(X^+_{res})= \rho^\alpha_{\text{bf}}\rho^\beta_{\text{tf}}\rho^\gamma_{\text{zf}}\mathcal{A}^{0,0,0}(X^+_{res}),
\end{align}
where a $u\in \mathcal{A}^{0,0,0}(X_{res}^+)$ only if it is conormal with respect to all $3$-faces of $X_{res}^+$ that is
\begin{align}
    \mathcal{A}^{0,0,0}(X^+_{res}):=\left\{u\in L^{\infty}(X^+_{res})\;|\; (\rho_\bullet\partial_{\rho_\bullet}) u \in L^\infty(X^+_{res}) \text{ where } \bullet\in\{\text{tf,bf,zf} \} \right\}.
\end{align}

Notice that with the above one can immediately derive a useful general  transformation rule regarding $\sigma$ dependent conormal functions on $X$. 
\begin{lemma}\label{lem:TransformationRuleXres}
  Let $\beta,\gamma\in \mathbb{R}$, we have that 
\begin{align}
    \mathcal{A}^{\beta}([0,1)_\sigma,\mathcal{A}^\gamma(X))=\mathcal{A}^{\gamma,\beta+\gamma,\beta}(X_{res}^+).
\end{align}
\end{lemma}

\begin{proof}
    This follows immediately from the following equality:
    \begin{align}
        \rho_{bf}^\gamma\cdot \rho_{tf}^{\beta+\gamma}\cdot \rho_{zf}^\beta = \left(\frac{\rho}{\rho+\sigma}\right)^{\gamma}(\rho+\sigma)^{\beta+\gamma}\left(\frac{\sigma}{\rho+\sigma}\right)^\beta=\sigma^\beta\rho^\gamma.
    \end{align}
\end{proof}

As we will run into logarithmic outputs of the the resolvent $\hat{\Box}_g(0)$ later on, we will have to allow in our definition of conormal spaces functions that exhibit partial polyhomogeneous expansions with respects to the boundaries $\text{zf,bf,tf}$.\\
To this ends we introduce first index sets.

\begin{definition}
   A set $\mathcal{E}\subset\mathbb{C}\times\mathbb{N}_{\geq 0}$ is an index set if 
   \begin{enumerate}
       \item $(z,k) \in \mathcal{E}\implies (z+1,k)\in \mathcal{E}$, and additionally, for all $k\geq 1$ it holds that $(z,k)\in \mathcal{E}\implies (z,k-1)\in \mathcal{E}$.
       \item For each $C\in \mathbb{R}$ we have that 
       \begin{align}
           |\mathcal{E}\cap \{(z,k)\in \mathbb{C}\;|\; \text{Re }z<C\}|<\infty \label{eq:InfimumIndexSet}.
       \end{align}
   \end{enumerate}
   Finally we define $\inf \text{Re }\mathcal{E} := \inf \{\text{ Re } z\;|\; (z,k)\in \mathcal{E} \}$, which is finite because of (\ref{eq:InfimumIndexSet}).
\end{definition}
This allows us to generalize our notion of conormal spaces to include poly-homogeneous expansions.

\begin{definition}
    \begin{enumerate}
        \item Let $\mathcal{E}$ be an index set, let us set $\alpha\in \mathbb{R}$ $\alpha_0:=\inf\text{Re }\mathcal{E}$ and set $\beta=\min \{\alpha,\alpha_0\}$. Then $\mathcal{A}^{(\mathcal{E},\alpha)}(X)\subset\mathcal{A}^{\beta-}(X)$ consists of functions $u$ that are smooth in $X^{\circ}$ and near $\partial X$ has a partial expansion
        \begin{align}
            u=\sum_{(z,k)\in\mathcal{E}, \text{Re }z\leq \alpha } u_{z,k}(\omega)\rho^z(\log \rho)^k +\mathcal{R}, \text{ where } \mathcal{R}\in\mathcal{A}^{\alpha}(X)
        \end{align}
        for some $u_{z,k}\in C^\infty (\partial X)$, i.e. $u_{z,k}$ is a smooth function on $\mathbb{S}^3$.\\
        \item The notion under 1. can be generalized, that is, for index sets $\mathcal{E}_{\bullet}$ with $\bullet = \text{bf,tf,zf}$ and $\alpha_{0,\bullet}=\inf\text{Re }\mathcal{E}_\bullet$ and $\alpha_\bullet\in \mathbb{R}$, set $\beta_\bullet = \min\{\alpha_{0,\bullet},\alpha_\bullet\}$, contains all $u\in \mathcal{A}^{(\mathcal{E}_{\text{bf}},\alpha_{\text{bf}}),(\mathcal{E}_{\text{tf}},\alpha_{\text{tf}}),(\mathcal{E}_{\text{zf}},\alpha_{\text{zf}})}(X^+_{res})$ that have partial polyhomogenous expansions with conormal remainders to each respective boundary. Explicitly we have
        \begin{itemize}
            \item In a collar neighbourhood of \text{zf}, i.e. $[0,\epsilon)_{\rho_{\text{zf}}}\times \text{zf}$, there exists $u_{\text{zf},(z,k)}\in\mathcal{A}^{(\mathcal{E}_{\text{tf}},\alpha_{\text{tf}})}(\text{zf})$ s.t.
        \begin{align}
            u = \sum_{\substack{(z,k)\in \mathcal{E}_{\text{zf}}\\\text{Re }z\leq\alpha_{\text{zf}}}} u_{\text{tf},(z,k)}(x)\rho_{\text{tf}}^z(\log\rho_{\text{tf}})^k + \mathcal{R}, \text{ where }\mathcal{R}\in \mathcal{A}^{\beta_{\text{tf}}-,\alpha_{\text{zf}}}([0,\epsilon)_{\rho_{\text{zf}}}\times \text{zf}).
        \end{align}\\
        \item In a collar neighbourhood of $\text{tf}$ near $\text{tf}\cap\text{zf}$, i.e. $[0,\epsilon)_{\rho_{\text{tf}}}\times \text{tf}$, there exist $u_{\text{tf},(z,k)}\in\mathcal{A}^{(\mathcal{E}_{\text{zf}},\alpha_{\text{zf}})}(\text{tf})$, s.t. 
              \begin{align}
            u = \sum_{\substack{(z,k)\in \mathcal{E}_{tf}\\\text{Re }z\leq\alpha_{tf}}} u_{\text{tf},(z,k)}(x)\rho_{\text{tf}}^z(\log\rho_{\text{tf}})^k + \mathcal{R}, \text{ where }\mathcal{R}\in \mathcal{A}^{\alpha_{\text{tf}},\beta_{\text{zf}}-}([0,\epsilon)_{\rho_{\text{tf}}}\times \text{tf}),
        \end{align}
    whilst on $\text{tf}$ near $\text{tf}\cap\text{bf}$, there exist $\tilde{u}_{\text{tf},(z,k)}\in \mathcal{A}^{(\mathcal{E}_{\text{bf}},\alpha_{\text{bf}})}(tf)$, s.t. 
    \begin{align}
            u = \sum_{\substack{(z,k)\in \mathcal{E}_{tf}\\\text{Re }z\leq\alpha_{tf}}} \tilde{u}_{\text{tf},(z,k)}(x)\rho_{\text{tf}}^z(\log\rho_{\text{tf}})^k + \mathcal{R}, \text{ where }\mathcal{R}\in \mathcal{A}^{\beta_{\text{bf}}-,\alpha_{\text{tf}}}([0,\epsilon)_{\rho_{\text{tf}}}\times \text{tf}).
        \end{align}\\
        \item In a collar neighbourhood of $\text{bf}$, i.e. $[0,\epsilon)_{\rho_{\text{bf}}}\times\text{bf} $ there exist $u_{\text{bf},(z,k)}\in \mathcal{A}^{(\mathcal{E}_{\text{tf}},\alpha_{\text{tf}})}(\text{bf})$, s.t.
        \begin{align}
            u = \sum_{\substack{(z,k)\in \mathcal{E}_{bf}\\\text{Re }z\leq\alpha_{bf}}} u_{\text{bf},(z,k)}(x)\rho_{\text{bf}}^z(\log\rho_{\text{bf}})^k + \mathcal{R}, \text{ where }\mathcal{R}\in \mathcal{A}^{\alpha_{\text{bf}},\beta_{\text{tf}}-}([0,\epsilon)_{\rho_{\text{tf}}}\times \text{tf}).
        \end{align}    
    \end{itemize}
    \item Partial polyhomogeneous spaces, such as $\mathcal{A}^{\alpha_{\text{bf}},\alpha_{\text{tf}},(\mathcal{E}_{\text{zf}},\alpha_{\text{zf}})}(X^+_{res})$ only have partial expansions at the boundary for which the index set is given (in this case only at the zero face).
    \end{enumerate}
\end{definition}
 
\begin{proposition}\label{prop:FirstMappingProposition}
    Let $\alpha \in (0,1)$ and suppose we are given $\sigma$-independent radial data $f\in \mathcal{A}^{2+\alpha}(X)$, then we have that 
    \begin{align}\label{eq:FirstMappingEquation1}
        \hat{\Box}(\sigma)^{-1}f\in \mathcal{A}^{\alpha-,\alpha-,((0,0),\alpha-)}(X^+_{res}).
    \end{align}
    For $\sigma$-dependent inputs $f\in \mathcal{A}^\beta([0,1)_\sigma;\mathcal{A}^{2+\alpha}(X))$ with $\beta\in\mathbb{R}$ we have
    \begin{align}\label{eq:FirstMappingEquation2}
        \hat{\Box}(\sigma)^{-1}f\in \mathcal{A}^{\alpha-,\alpha+\beta-,\beta}(X_{res}).    \end{align}
\end{proposition}

\begin{proof}
    By Theorem \ref{thm:SpectralEstimate} and Sobolev embedding in Lemma \ref{lem:SobolevEmbedding} we have that 
    \begin{equation}
        u(\sigma):=\hat{\Box}(\sigma)^{-1}f\in \mathcal{A}^{\alpha -}(X),
    \end{equation}
    which at $\sigma = 0$ follows from a factorization similar to (\ref{eq:Factorization}), in this case by setting 
    \begin{align}
        \partial_\sigma (\sigma \hat{\Box}(\sigma) \circ\hat{\Box}(\sigma)^{-1})=\mathds{1},
    \end{align}
    from which one obtains
    \begin{align}
        \sigma \partial_\sigma\hat{\Box}(\sigma)^{-1}=-\hat{\Box}(\sigma)^{-1}\circ \sigma \partial_\sigma \hat{\Box}(\sigma) \circ \hat{\Box}(\sigma)^{-1}.
    \end{align}
      Tracking this chain of maps we get
    \begin{align}\label{eq:SecondFactorization}
        \mathcal{A}^{2+\alpha}(X) \xrightarrow{\hat{\Box}(\sigma)^{-1}} \mathcal{A}^{\alpha-}(X) \xrightarrow{\sigma \partial_\sigma\hat{\Box}(\sigma)}|\sigma|\mathcal{A}^{1+\alpha -}(X) \xrightarrow{\hat{\Box}(\sigma)^{-1}}\mathcal{A}^{\alpha-}(X); 
    \end{align}
    Oberserve that the boundedness of $\sigma \partial_\sigma \hat{\Box}(\sigma): \mathcal{A}^{\alpha-}(X)\to |\sigma|\mathcal{A}^{1+\alpha-}(X)$ follows directly from (\ref{eq:GeneralSpectralFamily}), i.e. we have that 
    \begin{align}
       \sigma \partial_\sigma\hat{\Box}(\sigma)=2i\sigma\rho Q + 2\sigma^2 g^{00},
    \end{align}
    with $Q\in \text{Diff}^1_b(X)$.
   To prove the boundedness of $\hat{\Box}(\sigma)^{-1}:|\sigma|\mathcal{A}^{1+\alpha-}(X) \rightarrow \mathcal{A}^{\alpha-}(X)$ let us start out with $v\in \mathcal{A}^{1+\alpha-}(X)$, we then aim to show that 
   \begin{align}\label{eq:RenormalizedSigmaMapping}
       \hat{\Box}(\sigma)^{-1}(|\sigma|v) \in \mathcal{A}^{\alpha-}(X)
   \end{align}
   uniformly for small $|\sigma|$.
  Let us set $v:=\hat{\Box}(\sigma)^{-1}(|\sigma| w)$
 and apply Theorem \ref{thm:SpectralEstimate} with $\nu=0$ and $\ell = -2 + \alpha -\delta$ with $0<\delta < \alpha$ arbitrary, to give the estimate 
    \begin{align}
        \|v\|_{H_b^{s, \alpha-2 -\delta}(X)}\leq C\|w\|_{H^{s,\alpha-1 -\delta}_b(X)},
    \end{align}
    where we also used $\frac{|\sigma|}{\rho+|\sigma|}\leq 1$. Now by the Sobolev embedding Lemma \ref{lem:SobolevEmbedding}, we get that for large $s$, we have that 
    \begin{align}
        v\in H^{s,-2+\alpha-\delta}_b(X) \implies v\in \mathcal{A}^{\alpha-\delta-}(X),
    \end{align}
    as $\delta>0$ can be arbitrary, it thus follows $v \in \mathcal{A}^{\alpha-}(X)$.
    
    Similarly, we have that 
    \begin{align}
        w\in \mathcal{A}^{1+\alpha-}(X) \implies w \in H_b^{s,\alpha-1-\delta}(X).
    \end{align}
    Finally observe that we can establish exactly the same mapping properties for repeated applications of the $b$-operator (in $X^+_{res}$) of $\sigma\partial_\sigma$. More precisely, we have that for $k\in \mathbb{N}_0$, that 
    \begin{align}\label{eq:RepeatedSigmaDerivations}
        (\sigma\partial_\sigma)^k u(\sigma) = (\sigma\partial_\sigma)^k \hat{\Box}(\sigma)^{-1}f \in \mathcal{A}^{\alpha-}(X),
    \end{align}
    this follows from applying the mapping property in equation (\ref{eq:SecondFactorization}) inductively, as we have that
    \begin{align}\label{eq:BinomialExpansionSigmaDerivations}
        \hat{\Box}(\sigma) (\sigma\partial_\sigma)^k u=-\sum^k_{j=1}\begin{pmatrix}
            k\\j 
        \end{pmatrix}
        (\sigma\partial_\sigma)^j \hat{\Box}(\sigma) (\sigma\partial_\sigma)^{k-j}u,
    \end{align}
    with
    \begin{align}\label{eq:RepeatedSigmaDerivations}
        (\sigma\partial_\sigma)^j \hat{\Box}(\sigma)= 2i\sigma\rho Q + 2^j\sigma^2g^{00}
    \end{align}
    so $ (\sigma\partial_\sigma)^j \hat{\Box}(\sigma)$ has the mapping property of (\ref{eq:SecondFactorization}). If we assume that $(\sigma \partial_\sigma)^m u \in \mathcal{A}^{\alpha_-}(X)$ for every $m<k$ then the r.h.s belongs to $|\sigma|\mathcal{A}^{1+\alpha-}(X)$. Applying (\ref{eq:RenormalizedSigmaMapping}) then proves the claim.    
    From this we conclude that $u(\sigma)$ is conormal of order $0$ w.r.t. $\sigma\partial_\sigma$ at zf, i.e. we have that 
    \begin{align}\label{eq:InitialuCharacterization}
        u \in \mathcal{A}^0([0,1)_\sigma; \mathcal{A}^{\alpha-}(X)) \subset \mathcal{A}^{\alpha-,\alpha-,0}(X^+_{res}).
    \end{align}
   
    To finish up the lemma, we still need to improve regularity at the zf-face slightly more. This time we use a sharper two parameter inequality in $\rho$ and $\sigma$, more specifically $\frac{|\sigma|}{\rho+|\sigma|} \leq |\sigma|^{\alpha-\delta}\rho^{-\alpha+\delta}$, which is valid for $0<\alpha-\delta<1$. With exactly analaogous steps (again using Theorem \ref{thm:SpectralEstimate} and Lemma \ref{lem:SobolevEmbedding}) one then obtains a chain 
    \begin{align}\label{eq:ThirdFactorization}
         \mathcal{A}^{2+\alpha}(X) \xrightarrow{\hat{\Box}(\sigma)^{-1}} \mathcal{A}^{\alpha-}(X) \xrightarrow{\sigma \partial_\sigma\hat{\Box}(\sigma)}|\sigma|\mathcal{A}^{1+\alpha -}(X) \xrightarrow{\hat{\Box}(\sigma)^{-1}}|\sigma|^{\alpha-\delta}\mathcal{A}^{\delta-}(X). 
    \end{align}
    Applying Lemma \ref{lem:TransformationRuleXres} then gives 
    \begin{align}
        \sigma\partial_\sigma u\in \mathcal{A}^{\alpha-\delta}([0,1)_\sigma;\mathcal{A}^{\delta-}(X))\subset \mathcal{A}^{\delta-,\alpha-,\alpha-\delta}(X_{res}^+).
    \end{align}
     For each $k\geq 1$ the r.h.s. of (\ref{eq:BinomialExpansionSigmaDerivations}) belongs to $|\sigma|\mathcal{A}^{1+\alpha-})(X)$ by (\ref{eq:RepeatedSigmaDerivations}). Applying the refined final mapping in (\ref{eq:ThirdFactorization}) thus gives
    \begin{align}\label{eq:SigmaDerviativeArgument1}
        (\sigma\partial_\sigma)^k u\in |\sigma|^{\alpha-\delta}\mathcal{A}^{\delta-}(X).
    \end{align}
    To finish the proof, switch to projective coordinates $\hat{r}=\frac{\sigma}{\rho}$ at $\text{zf}$. We have exactly $\sigma \partial_\sigma =\hat{r}\partial_{\hat{r}}$ at fixed $\rho$. For $0<\epsilon<\delta<\alpha$ estimate (\ref{eq:SigmaDerviativeArgument1}) therefore gives 
    \begin{align}
        (\sigma\partial_\sigma) u \in \hat{r}^{\alpha -\delta}\rho^{\alpha-\epsilon}\mathcal{A}^{0,0}(U), 
    \end{align}
    where $U$ is a small enough neighbourhood of $\text{zf}\cap\text{tf}$.
    Indeed notice that 
    \begin{align}
        \sigma^{\alpha-\delta}\rho^{\delta-\epsilon}=\hat{r}^{\alpha-\delta}\rho^{\alpha-\epsilon}.
    \end{align}
    Since $\alpha-\delta>0$, integration in $\hat{r}$ therefore guarantees that $u$ has a limit $u_0(\rho,\omega)$ as $\hat{r}\to 0$ and 
    \begin{align}
        u(\rho,\hat{r},\omega)-u_0(\rho,\omega)=\int^{\hat{r}}_0 (\hat{r}\partial_{\hat{r}}) u(\rho,s,\omega) \frac{ds}{s}.
    \end{align}
    This argument can be repeated upon repeated application of $b$-derivatives, and logrithmic derivatives in $\hat{r}$ are controlled by estimates for $(\hat{r}\partial_{\hat{r}})^ku$. We conclude 
    \begin{align}
        u-u_0 \in \hat{r}^{\alpha-\delta}\rho^{\alpha-\epsilon}\mathcal{A}^{0,0}(U),\;\; u_0\in \mathcal{A}^{\alpha-}(X),
    \end{align}
    where away from $\text{tf}$ the same integration argument applies with $\sigma$ as boundary defining function, as $\rho$ is bounded from below.
    Since $\delta$ and $\epsilon$ can be taken to be arbitrarily small this give sthe claimed partial expansion at $\text{zf}$. Together with (\ref{eq:InitialuCharacterization}) away from $\text{zf}$ this proves \ref{eq:FirstMappingEquation1}.
    For the $\sigma$-dependent statement, simply set $h=\sigma^{-\beta}f$ and $v=\sigma^{-\beta} u$, then $\hat{\Box}(\sigma)v = h$ with 
    \begin{align}
        h\in \mathcal{A}^{0}([0,1)_\sigma;\mathcal{A}^{2+\alpha}(X)).
    \end{align}
    The same differentiated-equation argument, now including $(\sigma\partial_\sigma)^kh$ on the right-hand side then gives
    \begin{align}
        v\in \mathcal{A}^{0}([0,1)_\sigma;\mathcal{A}^{\alpha-}(X)),
    \end{align}
    and thus multiplication by $\sigma^\beta$ and Lemma \ref{lem:MappingIntoXres} yield (\ref{eq:FirstMappingEquation2}).
\end{proof}

\begin{lemma}\label{lem:ResolvedonResolvent}
    Let $\alpha \in (0,1)$ and $f\in \mathcal{A}^{2+\alpha,2+\alpha,\alpha}(X^+_{res})$, then
    \begin{equation*}
        \hat{\Box}(\sigma)^{-1}f\in \mathcal{A}^{\alpha-,\alpha-,\alpha-}(X_{res}^+).
    \end{equation*}
\end{lemma} 

\begin{proof}
    First, let us observe that $f\in \mathcal{A}^{2+\alpha,2+\alpha,\alpha}(X^+_{res})$ implies that 
    \begin{align}
        f&= \rho_{bf}^{\alpha+2}\rho_{tf}^{\alpha+2}\rho _{zf}^\alpha\cdot g
    \end{align}
    with $g \in \mathcal{A}^0(X^+_{res})$ smooth and conormal. Untangling the decay factors this means that 
    \begin{align}
        f&= \frac{\rho^{\alpha+2}}{(\rho+\sigma)^{\alpha+2}}(\rho+\sigma)^{\alpha+2} \frac{\sigma^\alpha}{(\rho+\sigma)^{\alpha}}g = \frac{\rho^{\alpha+2}\sigma^\alpha}{(\rho+\sigma)^{\alpha}}g
    \end{align}
    To control the spatial decay uniformly in $\sigma$, in particular near \text{bf}, we observe that 
    \begin{align}
        \frac{f}{\rho^{2+\alpha}}=\left(\frac{\sigma}{\rho+\sigma}\right)^\alpha g\qquad \text{ with } g \in \mathcal{A}^{0,0,0}(X_{res}^+)
    \end{align}
    To additionally capture the vanishing towards \text{zf}, we factor out $\sigma^{\alpha-\delta}$, where $0<\delta<\alpha$, obtaining
    \begin{align}
        \frac{f}{\rho^{2+\delta}\alpha^{\alpha-\delta}}=\left(\frac{\rho}{\rho+\sigma}\right)^{\alpha-\delta}\left(\frac{\sigma}{\rho+\sigma}\right)^\delta g\qquad \text{ with } g \in \mathcal{A}^{0,0,0}(X_{res}^+)
    \end{align}
    
    Combining the above two statements, we see that if $f\in \mathcal{A}^{2+\alpha,2+\alpha,\alpha}(X^+_{res})$ we must have 
    \begin{align}
        f\in \mathcal{A}^{0}([0,1)_\sigma ; \mathcal{A}^{2+\alpha}(X))\; \cap \;\mathcal{A}^{\alpha-\delta}([0,1)_\sigma;\mathcal{A}^{2+\delta}(X))].
    \end{align}
    Now with the help of equation (\ref{eq:FirstMappingEquation2}) from Proposition \ref{prop:FirstMappingProposition} we can simply compute the action of $\hat{\Box}(\sigma)^{-1}$ on both conormal spaces. To summarize, from $f\in \mathcal{A}^{0}([0,1)_\sigma ; \mathcal{A}^{2+\alpha}(X))$ it follows (using $\beta=0$ in (\ref{eq:FirstMappingEquation2})) that 
    \begin{align}\label{eq:DomainSPace1}
        \hat{\Box}(\sigma)^{-1}f \in \mathcal{A}^{\alpha-,\alpha-,0}(X_{res}^+).
    \end{align}
    Similarly for from $f\in \mathcal{A}^{\alpha-\delta}([0,1)_\sigma;\mathcal{A}^{2+\delta}(X))$ we get from (\ref{eq:FirstMappingEquation2}) using $\beta=\alpha-\delta$ and substituting $\delta$ for $\alpha$ that 
    \begin{align}\label{eq:DomainSpace2}
        \hat{\Box}(\sigma)^{-1}f \in \mathcal{A}^{\delta-,\alpha-,\alpha-\delta}(X_{res}^+).
    \end{align}
    The statement now follows from the fact that $\hat{\Box}(\sigma)^{-1}f$ should therefore lie in the intersection  $\mathcal{A}^{\alpha-,\alpha-,0}(X_{res}^+)\;\cap \mathcal{A}^{\delta-,\alpha-,\alpha-\delta}(X_{res}^+)$ and $\text{bf}\cap \text{zf}=\emptyset$ therefore picking up the largest decay order in each factor of the b.d.f's and the fact that we may take $\delta>0$ arbitrarily small.
\end{proof}

\subsection{Action of the resolvent on large inputs.}

We can still prove this resolvent output with $\alpha$ in the ranging over the full indicial gap.

\begin{lemma}\label{lem:MappingIntoXres}
    Let $\alpha \in (0,2)$. Suppose $f\in \mathcal{A}^{2-\alpha}(X)$. Then 
    \begin{equation}\label{eq:ResolventOutput}
        \hat{\Box}(\sigma)^{-1}f\in \mathcal{A}^{1-\alpha-,-\alpha-,-\alpha-}(X).
    \end{equation}
    The same conclusion holds if, more generally, $f\in \mathcal{A}^0([0,1)_\sigma,\mathcal{A}^{2-\alpha}(X))$.
    
\end{lemma}
\begin{remark}
    We can still prove this resolvent output with $\alpha$ in the ranging over the full indicial gap $\alpha\in (0,2)$ as we have more freedom to choose the parameters $\ell$ and $\nu$ to satisfy Theorem \ref{thm:SpectralEstimate}.
\end{remark}
\begin{proof}
    Put $u=\hat{\Box}(\sigma)^{-1}f$, choose a small $\delta>0$ and set 
    \begin{align}
        \ell = -1-\alpha-2\delta,\;\qquad \nu=1-\alpha-3\delta.
    \end{align}
Then $\ell<-\frac{1}{2}$ and $\ell-\nu=-2+\delta\in (-2,0)$ so Theorem \ref{thm:SpectralEstimate} applies and for sufficiently large $s$,
\begin{align}\label{eq:FirstClaim}
    \|(\rho+|\sigma|)^{1-\alpha-3\delta}u\|_{H^{s,-1-\alpha-2\delta}_b}\lesssim \|(\rho+|\sigma|)^{-\alpha-3\delta}f\|_{H_b^{s,-\alpha-2\delta}}.
\end{align}

We now claim that 
\begin{align}\label{eq:SecondClaim}
    \|(\rho+|\sigma|)^{-\alpha-3\delta}f\|_{H_b^{s,-\alpha-2\delta}}\lesssim |\sigma|^{-\alpha-2\delta}\|f\|_{H_b^{s,-\alpha-\delta}}.
\end{align}

The norm on the right is finite because we have 
\begin{align}
    f\in \mathcal{A}^{2-\alpha}\subset H_b^{s-\alpha-\delta}(X).
\end{align}

Thus if we can show that the pointwise weighting function relating these two norms defined by
\begin{align}
    M(\rho,\sigma):=\frac{\rho^{\alpha+2\delta}(\rho+|\sigma|)^{-\alpha-3\delta}}{|\sigma|^{-\alpha-2\delta}\rho^{\alpha+\delta}}= \left(\frac{|\sigma|}{\rho+|\sigma}\right)^{\alpha+2\delta}\left(\frac{\rho}{\rho+|\sigma|}\right)^\delta
\end{align}
is uniformily bounded in $0<\rho\leq 1$ near $\text{bf}\cap \text{tf}$ and $|\sigma|\leq |\sigma_0|$ near $\text{tf}\cap \text{zf}$ the claim (\ref{eq:SecondClaim}) would follow. For the region $\rho\geq |\sigma|$ we have the $(\rho +|\sigma|)\sim \rho$ and thus 
\begin{align} 
    M(\rho,\sigma) \simeq \left( \frac{|\sigma|}{\rho}\right)^{\alpha+2\delta}\leq  1.
\end{align}
Similarly, if $|\sigma|\geq \rho$ we get that $(\rho+|\sigma|)\sim \sigma$ and thus
\begin{align}\label{eq:CursedEquation}
    M(\rho,\sigma)\simeq\left(\frac{\rho}{|\sigma|}\right)^{\delta}\leq 1.
\end{align}
The multiplier $M$ is uniformily bounded bounded with respect to all spatial deriviatives, so multiplication by $M$ is uniformily bounded on $H^{s}_b(X_{res}^+)$ proving (\ref{eq:SecondClaim}).
Combining (\ref{eq:FirstClaim}) and (\ref{eq:SecondClaim}) we then find
\begin{align}
    u\in |\sigma|^{-\alpha-2\delta}\rho^{1-\alpha-2\delta}(\rho+|\sigma|)^{\alpha-1+3\delta}\mathcal{A}^{0,0,0}(X_{res}^+).
\end{align}
for positive frequencies, we have that 
\begin{align}
    |\sigma|^{-\alpha-2\delta}\rho^{1-\alpha-2\delta}(\rho+|\sigma|)^{\alpha-1+3\delta}=\rho_{\text{bf}}^{1-\alpha-2\delta}\rho_{\text{tf}}^{-\alpha-\delta}\rho_{\text{zf}}^{-\alpha-2\delta},
\end{align}
and thus
\begin{align}
    u\in \mathcal{A}^{1-\alpha-2\delta,-\alpha-\delta,-\alpha-\delta}(X_{res}^+)
\end{align}
and upon taking intersections for smaller and smaller $\delta>0$ the claim follows. Note that conormality w.r.t. $\sigma\partial_\sigma$ follows from the same argument as in the proof of Proposition \ref{prop:FirstMappingProposition}, equation (\ref{eq:RepeatedSigmaDerivations})
\end{proof}

\begin{lemma}\label{lem:KernelElements}
    We have 
    \begin{align}
        \mathcal{A}^0(X)\cap\ker\hat{\Box}(0)=\mathbb{C}u_{(0)},\;\;\;\;\;\;\;\;\ \dot{H}^{s,\ell}_b(X)\cap\ker\hat{\Box}(0)^*=\mathbb{C}u_{(0)}^*, \;s<\frac{1}{2}, -3<\ell<-2.
    \end{align}\footnote{The general relation on $\ell$ depends on the dimension $d$ through $-1<\ell+\frac{d}{2}<0$.}
for states $\chi_\infty u_{(0)},\chi_{\infty}u_{(0)}^*\in \mathcal{A}^{((0,0),1-)}(X)$, which are uniquely determined by their leading order behaviour $\chi_\infty u_{(0)},\chi_\infty u_{(0)}^*\in 1+\mathcal{A}^{1-}(X)$ for a cut-off $\chi_\infty\in C_c^\infty(X)$ with $\chi_\infty\equiv 1$ near $\partial_+X$.
\end{lemma}

\begin{proof}
    This proof is essentially the same as the proof of Lemma 2.19 in \cite{hintz2022sharp}. We first prove existence of $u_{(0)}$. Let $\chi_{\partial}\in C^\infty(X)$ be a cutoff function with $\chi_{\partial}\equiv 1$ near $\partial X$. Now using Lemma \ref{lem:AdmissibleWaveEquation} we see by freezing coefficients at $\rho=0$ that the $b$-normal operator at $\partial X$ of $\rho^{-2} \hat{\Box}(0)$ is given by $L_0=-(\rho\partial_{\rho})^2+2\rho\partial_{\rho} + \cancel{\Delta}$. As $L_0$ annihilates constants, we have $L_0 \chi_{\partial} \in \rho C^\infty(X)$ and therefore that $e:=-\hat{\Box}(0)(\chi_{\partial})\in \rho^3 C^\infty(X)$, which improves the decay of one order relative to the usual mapping $\hat{\Box}(0):C^\infty(X) \to \rho^2 C^\infty(X)$. By Lemma \ref{lem:SobolevEmbeddingClosedGraph} we have that $\hat{\Box}(0):\mathcal{A}^{\alpha-}(X)\to \mathcal{A}^{2+\alpha}$ is surjective for the entire range $\alpha \in (0,2)$. Thus for $\alpha = 1$ there must exist a $\tilde u \in \mathcal{A}^{1-}(X)$ with $\hat{\Box}(0)\tilde u = e$. Let us then set $u_{(0)}:=\chi_\partial +\tilde u$. Because of the spectral admissibility assumption 1. in Definition \ref{def:SpectrallyAdmissbleMetric}, we know that $\hat{\Box}(0)$ is injective on $\mathcal{A}^{1-}(X)$ hence for any other admissible energy state $\tilde{u}_{(0)}$ we have $u_{(0)}-\tilde{u}_{(0)}\in \ker \hat{\Box}(0) \cap \mathcal{A}^{1-}(X)=\{0\}$. This establishes uniqueness. For $u^*_{(0)}$ the claim follows by duality.
 \end{proof}
Throughout what follows we will keep track of $u_{(0)}$ and $u_{(0)}^*$ in the resolvent expansions. However, for Myers-Perry, these functions can be specified explicitly. 

\begin{corollary}\label{cor:KernelElementsMP}
    Let $(\bfa,\bfb)\in D_\mfrm$ and set $\hat{\Box}(0):=\hat{\Box}_{g_{\bfa,\bfb,\mfrm}}(0)$. Then we have that the $u_{(0)},u_{(0)}^*$ of Lemma \ref{lem:KernelElements} are given by
\begin{align}
    u_{(0)}=1,\;\text{ and } u_{(0)}^* = H(r-r_+),
\end{align}
where $H:\mathbb{R}\to\mathbb{R}$ is the heavyside step function.
\end{corollary}
\begin{proof}
    First observe that because $\hat{\Box}(0)1=0$, uniqueness of the normalized direct state gives $u_{(0)}=1$.
    With respect to $dg_{\bfa,\bfb,\mfrm}|_X$ we have that $\hat{\Box}(0)=\hat{\Box}(0)^*$ as formal adjoints. From the dual metric (\ref{eq:MPStarDualMetric2}) we obtain 
    \begin{align}\label{eq:VanishingHeavyside}
        [\hat{\Box}(0),H(r-r_+)]=-\frac{2\delta(r-r_+)}{r^2\bfr^2}(\bfa(r^2+\bfb^2)\partial_{\phi_*} +\bfb(r^2+\bfa^2)\partial_{\psi_*}).
    \end{align}
    Indeed the purely angular terms commute with $H(r-r_+)$ while the radial divergence term contrubtes nothing we have that 
    \begin{align}
        [\partial_r(r\Delta_r\partial_r,H(r-r_+)]v &=\partial_r(r\Delta_r\delta(r-r_+) v)| + r\Delta_r\delta(r-r_+)\partial_rv|=0
    \end{align}
    as $\Delta_r|_{r=r_+}=0$. Thus only mixed radial-angular terms contribute to the commutator.
    As $u_{(0)}=1$ is independent of both angular variables we have 
    \begin{align}
        \hat{\Box}(0)^*(H(r-r_+)u_{(0)}=H(r-r_+)\hat{\Box}(0)u_{(0)}+[\hat{\Box}(0),H(r-r_+)]u_{(0)}.
    \end{align}
    Moreover $H(r-r_+)u_{(0)}=H(r-r_+)$ vanishes near the artifical inner boundary and belongs to $\dot{H}^{s,\ell}_b(X)$ for $s<\frac{1}{2}$, for $-3<\ell <-2$ and equals $1$near infinity. Uniqueness of the normalized adjoint then gives
    \begin{align}
        u_{(0)}^*=H(r-r_+)u_{(0)}=H(r-r_+).
    \end{align}
\end{proof}
\begin{remark}
    The preceding proof was on purpose kept in line with Lemma 4.4 in \cite{hintz2022sharp}. The fastest argument to verify the above claim is to simply verify that $[\hat{\Box}(0),H(r-r_+)]1=\hat{\Box}(0)H(r-r_+)$ and to observe that $\hat{\Box}(0)H(r-r_+)=0$ directly.
\end{remark}

\begin{definition}
    Let $\hat{\rho}=\frac{\rho}{\sigma}$. The model operator $\rho=\sigma=0$ is
    \begin{equation}
        \tilde{\Box}(1)=2i\hat{\rho}\left(\hat{\rho}\partial_{\hat\rho}-\frac{3}{2}\right)+\hat{\rho}^2(-(\hat{\rho}\partial_{\hat\rho})^2+2\hat{\rho}\partial_{\hat{\rho}}+\cancel{\Delta})\in \text{Diff}^2((0,\infty)_{\hat{\rho}}\times \mathbb{S}^3).
    \end{equation}
\end{definition}

The following is a direct consequence of Lemma \ref{lem:AdmissibleWaveEquation} and writing out the model operator in different regions of $X_{res}^+$.

\begin{lemma}\label{lem:NormalOperatorOnEachFace}
    The operator $\hat{\Box}(\sigma)$, as a second order differential operator on $X_{res}^+$ is a $b$-differential operator of class $\hat{\Box}(\sigma)\in \rho_{b}\rho_{tf}^2\text{Diff}_b^2(X_{res}^+)$. It's $b$-normal operators on each face can be given as:
    \begin{itemize}
        \item $2i\sigma^2\hat{\rho}\left(\hat{\rho}\partial_{\hat{\rho}}-\frac{3}{2}\right)$ at $\text{bf}$, i.e. we have 
        \begin{align}
            2i\sigma^2\hat{\rho}\left(\hat{\rho}\partial_{\hat{\rho}}-\frac{3}{2}\right)-\hat{\Box}(\sigma) \in \rho_{\text{b}}^2\rho_{\text{tf}}^2\text{Diff}_b^2(X);
        \end{align}\\
        \item $\sigma^2\tilde{\Box}(1)$ at $\text{tf}$, i.e. 
        \begin{align}
            \hat{\Box}(\sigma)-\sigma^2\tilde{\Box}(1) \in \rho_{\text{b}}\rho_{\text{tf}}^3\text{Diff}_b(X);
        \end{align}
        \item $\hat{\Box}(0)$ at \text{zf}, i.e. we have
        \begin{align}
            \hat{\Box}(\sigma)-\hat{\Box}(0) \in \rho_{\text{b}}\rho_{\text{tf}}^2\rho_{\text{zf}}\text{Diff}^2_b(X).
        \end{align}
        In fact, $\hat{\Box}(\sigma)-\sigma^2\tilde{\Box}(1)\in \rho_{\text{b}}^2\rho^3_{\text{tf}}\text{Diff}^2_b(X)$.
    \end{itemize}
\end{lemma}

On $tf$ we have the following solvability theorem for the normal operator of $\hat{\Box}(\sigma)$, directly obtained from \cite{vasy2021resolvent} Proposition 5.4 after adjusting for the central weight range for $\nu$.

\begin{theorem}\label{thm:VasyTransitionFaceTheorem}
    For $s\in \mathbb{R}$, $\ell<-\frac{1}{2}$, $s+\ell>-\frac{1}{2}$, and $\nu\in(0,2)$ the operator 
    \begin{align}
        \tilde{\Box}(1):\{u\in H^{s,\ell,\nu}_b(\text{tf})\;|\; \tilde{\Box}(1)u\in H_b^{s,\ell+1,\nu-2}(\text{tf})\}\to H^{s,\ell+1,\nu-2}_b(\text{tf}),
    \end{align}
is invertible. In particular, $\tilde{\Box}(1):\{u\in\mathcal{A}^{\beta-,\gamma-}(\text{tf})\;|\;\tilde\Box(1) u \in \mathcal{A}^{\beta+1,\gamma-2}(tf)\}\to \mathcal{A}^{\beta+1,\gamma-2}(\text{tf})$ is an isomorphism $\beta<1$ and $\gamma\in (-2,0)$.    
\end{theorem}

\begin{lemma}\label{lem:tfResolventOutput}
    Let $\tilde{f}\in C^\infty(\partial X)$, the unique solution $\tilde{u}\in \mathcal{A}^{1-,0-}(\text{tf})$ of 
    \begin{equation}
        \tilde{\Box}(1)\tilde{u}=\hat{\rho}^2\tilde{f}\in\mathcal{A}^{2,-2}(\text{tf}),
    \end{equation}
    lies in the space $\mathcal{A}^{1-,((0,1),1-)}(\text{tf})$; the leading term at $\text{zf}$ is $-\frac{1}{4\pi^2}(\int_{\partial X}\tilde{f} dg_{\mathbb{S}^3})\log(\hat{r}).$
\end{lemma}

\begin{proof}
    We only need to analyze $\tilde u$ near the intersection $tf\cap zf$, i.e. near $\hat{r}=\frac{\sigma}{\rho}=0$. Let $\psi = \psi(\hat{r})\in C_c^{\infty}([0,\frac{1}{2}))$ be a cutoff with $\psi\equiv 1$ near $\hat{r}=0$ and set $v:=\psi \tilde{u}$. Then on $[0,1)_{\hat{r}}\times \partial X$, we have 
    \begin{align}
        \tilde{\Box}(1) v = \hat{r}^{-2}\psi\tilde f + [\tilde{\Box}(1),\psi]\tilde u \in \hat{r}^{-2}C^\infty(X).
    \end{align}
    The $b$-normal operator of $\tilde{\Box}(1)$ is $\hat{r}^{-2} L_0$ in the notation of (\ref{eq:NormalOperator}) (under change of variables of $\rho=\frac{\sigma}{\hat{r}}$), so $v\in \mathcal{A}^{0-}(X)$ solves 
    \begin{align}
        L_0v = h := \psi \tilde f(\omega) +\hat{r}^{1-}g(\hat{r},\omega)\in \psi \tilde f+\mathcal{A}^{1-},\;\;\;\; L_0= -(\hat{r}\partial_{\hat{r}})^2-2\hat{r}\partial_{\hat{r}}+\cancel{\Delta},
    \end{align}
    where we note that $[\tilde{\Box}(1),\psi]$ once again picks up one extra order of decay in $\hat{r}$ due to the cutoff function.
   We then treat the above equations by Mellin transforming in $\hat{r}$ defined 
   \begin{equation}
       \hat{v}(\lambda):=\int^\infty_0 \hat{r}^{-i\lambda}v(\hat{r})\frac{d\hat{r}}{\hat{r}},
   \end{equation}

   where we note that $\widehat{\hat{r}\partial_{\hat{r}}}v(\lambda)=i\lambda \hat{v}(\lambda)$. Also recall that as $\hat{r}$ is the b.d.f. for the $tf$ face at $zf\cap tf$ a Mellin transform of a function $u\in \hat{r}^\nu L^2(\partial X)$ implies that $\psi\hat{u}(\lambda)$ is exactly holomorphic for $\text{ Im } (\lambda) > -\nu$. In this setting holomorphy of $\hat{u}(\lambda)$ is determined by the $\nu$ implicitly chosen in Theorem \ref{thm:VasyTransitionFaceTheorem}. \\
   We thus see that 
   \begin{align}\label{eq:NormaltfzfOperator}
       \widehat{L_0}(\lambda) \hat{v}(\lambda) = \hat{h}(\lambda), \;\;\text{ and }\;\;\; \widehat{L_0}:=-(i\lambda)^2-2i\lambda + \cancel{\Delta}.
   \end{align}
    We recall that $\hat{v}(\lambda)$ is holomorphic in $\text{Im}(\lambda)>0$ and therefore satisfies estimates 
    \begin{align}
        \|\hat{v}(\lambda)\|_{H^s(\partial X)} \leq C_{s,N,\epsilon} \langle\lambda\rangle^{-N}, \;\;\;\;\; \epsilon < \text{Im}( \lambda)<1,
    \end{align}
    for all $s,N\in \mathbb{R}$ and $\epsilon >0$. Similarly the above estimate holds for $\hat{h}(\lambda)$ in the region $-1<\text{Im } (\lambda)< 1$, $|\lambda|>\epsilon$ but $\hat{h}(\lambda)$ has a simple pole by equation (\ref{eq:NormaltfzfOperator}) at $\lambda=0$ with residue $\text{Res}_{\lambda=0}(\hat{h})=i\tilde{f}$.\\
    Let us expand $v$ and $h$ in hyperspherical harmonics $S_{kml}$, $k=0,1,2,...$, $l=0,1...,k$, $|m|\leq l$ with $\cancel{\Delta}S_{klm}= l(l+2)$ and let $\mathcal{Y}_{kl}=\text{Span}\{S_{klm}\;|\; |m|\leq l\}$. Then restricting to $\mathcal{Y}_{00}$ the inverse $\left.\widehat{L_0(\lambda)}^{-1}\right|_{\mathcal{Y}_{kl}}=(-(i\lambda)^2-2i\lambda+(l(l+2))^{-1}$ is meromorphic with simple poles at $i\lambda = -l-2,l$.  Let $h= h_0 + h_1$ with $h_0$ the orthogonal projection onto the s-wave $Y_{000}$ by spherically averaging, i.e.
    \begin{align}
        h_0(\hat{r}) =\frac{1}{\text{vol}(\mathbb{S}^3)}\int_{\mathbb{S}^3} h(\hat{r},\omega) dg_{\mathbb{S}^3}, 
    \end{align}
    and also define 
    \begin{align}
        h_0(0) =\frac{1}{\text{Vol}(\mathbb{S}^3)}\int_{\mathbb{S}^3}\tilde f\, dg_{\mathbb{S}^3},\qquad h_1(0,\omega)=\tilde{f}(\omega)-h_0(0)
    \end{align}
    
    Note then that for $\hat{v}_0(\lambda):=\hat{L}_0(\xi)^{-1}\widehat{h_0}(\lambda)$ ($l=0$) 
    and for $\hat{v}_1(\lambda):=\hat{L}_0(\lambda)^{-1}\hat{h}_1(\lambda)$ $(l>0)$ we have at most a double, resp. single pole, as $\widehat{L_0}(\lambda)$ has no roots at $i\lambda = 0$ on the complement of $\mathcal{Y}_{00}$ and are meromorphic for $\text{ Im } \lambda> -1$. In the inverse Mellin transform 
    \begin{align}
        v(\hat{r},\omega) =\frac{1}{2\pi}\int_{\text{Im}(\lambda)=\epsilon} \hat{r}^{i\lambda}\widehat{L_0}(\lambda)^{-1}\hat{h}(\lambda) d\lambda
    \end{align}
   we then shift the contour down to $\text{Im}(\lambda)=-1+\epsilon$ and we get 
   \begin{align}
       v(\hat{r},\omega)= -2\pi i\cdot \text{Res}_{\lambda=0}\left(\frac{1}{2\pi} \hat{r}^{i\lambda}\widehat{L_0}((\lambda)^{-1}\hat{h}(\lambda)) \right)+ \frac{1}{2\pi}\int_{\text{Im}(\lambda)=-1+\epsilon} \hat{r}^{i\xi}\widehat{L_0}(\lambda)^{-1}\hat{h}(\lambda) d\lambda
   \end{align}
    As the new contour on the r.h.s. will have higher $\hat{r}$-decay due to the shifting, we now only need to analyze the resulting residue at $\lambda = 0$. Let us decompose $h_1$ into spherical modes of degree $l$, that is
    \begin{align}
        h_1 = \sum_{l=1}^\infty h_{1,l}, \quad \Delta_{\mathbb{S}^3}h_{1,l}=l(l+2)h_{1,l}
    \end{align}
    We then have that 
    \begin{align}
        \hat{r}^{i\lambda}\widehat{L_0}(\lambda)^{-1}\hat{h}(\lambda)&=\hat{r}^{i\lambda}\overbrace{ \frac{\hat{h}_0(\lambda)}{\lambda^2-2i\lambda}}^{\hat{v}_0(\lambda)} + \sum_{l=1}^\infty\hat{r}^{i\lambda}\overbrace{\frac{\hat{h}_{1,l}(\lambda)}{-(i\lambda-l)(i\lambda + (l+2))}}^{\hat{v}_1(\lambda)}.
    \end{align}
  Analyzing the first part in the sum and expanding $\hat{r}^{i\lambda}\sim 1+\log(\hat{r})i\lambda +\mathcal{O}(\lambda^2)$, we see that 
   \begin{align}
       \hat{v}_0(\lambda)\sim \left(\frac{-1}{2i\lambda}+\mathcal{O}(1) \right)\left(\frac{-h_0(0)}{i\lambda}+\mathcal{O}(1)\right) = \frac{h_0(0)}{2(i\lambda)^2}+\mathcal{O}((i\lambda)^{-1}),
   \end{align}
   so 
   \begin{align}
      \hat{r}^{i\lambda} \hat{v}_0(\lambda)\sim  \frac{h_0(0)}{2}\left(\frac{1}{(i\lambda)^2}+\frac{\log(\hat{r})+const}{i\lambda}+\mathcal{O}(1)\right),
   \end{align}
   hence the residue of this part is $\frac{h_0(0)}{2i}\log(\hat{r})+\text{const}$. Likewise we can compute for the single pole part $\hat{v}_1(\lambda)$ that
   \begin{align}
       \hat{r}^{i\lambda}\hat{v}_1(\lambda)=\left(\frac{C_1(\omega)}{i\lambda}+C_1(\omega)\log(\hat{r}) +\mathcal{O}(1)\right), \text{ with } C_1(\omega)=-\sum^\infty_{l=1} \frac{h_{1,l}(0,\omega)}{l(l+2)},
   \end{align}
   hence here the residue only contributes a power of $\hat{r}^0$.

    We conclude therefore that 
    \begin{align}\label{eq:GeneralLogofV}
        v(\hat{r})= c\log (\hat r) +\tilde v +\mathcal{A}^{1-}, 
    \end{align}
    where $c\in \text{Span}\{Y_{000}\}\cong\mathbb{C}\cdot \mathbf{1}$ and $\tilde {v}\in C^\infty(\partial X)$. 
    Finally to compute $c$ explicitly note that 
    \begin{align}\label{eq:PreSphericallyAveraged}
        L_0 v(\hat{r}) = c L_0(\log(\hat{r})) +\Delta_{\mathbb{S}^3}\tilde v +\mathcal{A}^{1-},
    \end{align}
    now given that we know that $L_0 v = h_0 + h_1$ where $h_0$ is the \textit{unique} projection onto $Y_{000}$ we conclude that we must have $-2c = h_0$ after spherical averaging of \ref{eq:PreSphericallyAveraged}  and thus that
      $c= - \frac{1}{2}h_0(0) = -\frac{1}{4\pi^2}\int_{\mathbb{S}^3}\tilde f dg_{\mathbb{S}^3}.$ 
   
\end{proof}

\begin{lemma}\label{lem:ThresholdOutput}
    Let $u_{(0)} \in \mathcal{A}^0(X)$ with $\hat{\Box}(0)u_{(0)} = 0$ be as in Lemma \ref{lem:KernelElements}. Let $u^{(2)}:= \tilde{\Box}(1)^{-1}(\hat{\rho}^2)$ as computed by Lemma \ref{lem:tfResolventOutput}, i.e. by setting $\tilde{f}\equiv 1$. Then 
    \begin{align}
        \hat{\Box}(\sigma)^{-1}\rho^2 \in \mathcal{A}^{1-,((0,0),1-),((0,1),1-)}(X_{res}^+).   \end{align}
    We then have that up to leading order
    \begin{align}
        \hat{\Box}(\sigma)^{-1}\rho^2\sim \tilde{u}^{(2)} \text{ at tf} \;\;\;\;\;\text{ and } \;\;\;\;\;\hat{\Box}(\sigma)^{-1}\rho^2\sim -\frac{1}{2}\log(\frac{\sigma}{\rho}) u_{(0)} \text{ at zf}.
    \end{align}
\end{lemma}

\begin{proof}
    Let $\chi_\partial \in C^\infty(X)$ be a cutoff supported near $\partial X$. Note that as $zf$ is the lift of $\{0\}_{\sigma} \times X^\circ$ that the lift of $\chi_\partial$ to $X_{res}^+$ is then supported \textit{away} from zf and we may moreover assume that $\chi_\partial\equiv 1$ in some open neighbourhood of $tf$. Rewriting $\hat{\Box}(\sigma)^{-1}\rho^2$ we get that 
    \begin{align}\label{eq:OriginalErrorTerm}
        \hat{\Box}(\sigma)^{-1}\rho^2:=\chi_\partial \tilde u ^{(2)} + \hat{\Box}(\sigma)^{-1}\tilde{e}, \text{ with } \tilde e := \rho^2-\hat{\Box}(\sigma)(\chi_\partial \tilde u^{(2)}).
    \end{align}
    Now note that by constuction $\chi_\partial \tilde{u}^{(2)}$ will be smooth conormal w.r.t. $\rho_{tf}$ on $X^+_{res}$.  Combined with the $bf$ and $zf$ output given by Lemma \ref{lem:tfResolventOutput} we then see that $\chi_\partial \tilde u ^{(2)} \in \mathcal{A}^{1-,(0,0),(0,1)+1-}(X_{res}^+)$.\\
    For the next step, we will show that the remainder $\tilde {e}$ exhibits leading order cancellations at both $bf$ and $tf$, thereby improving conormal regularity at those faces. The Abschnitt will be that $\tilde{e}\in\mathcal{A}^{3-,3,0}(X)$ improving over $\rho^2\in \mathcal{A}^{2,2,0}(X_{res}^+)$  At $tf$ note that we can expand 
    \begin{align}
        \hat{\Box}(\sigma)(\chi_\partial\tilde u ^{(2)})&= \chi_\partial \sigma^2\tilde \Box(1) \tilde u^{(2)} +[\hat{\Box}(\sigma),\chi_\partial]\tilde u^{(2)}+R(\chi_\partial \tilde u ^{(2)}),
    \end{align}
    where $R\in \rho_b\rho_{tf}^3\text{Diff}_b^2(X)$ by Lemma \ref{lem:NormalOperatorOnEachFace}. Now note that we can expand $\tilde{e}$ at tf as
    \begin{align}
        \tilde e = \underbrace{\rho^2-\chi_\partial\sigma^2\tilde \Box(1) \tilde u^{(2)}}_{\tilde e_1} -\underbrace{[\sigma^2\tilde{\Box}(1),\chi_\partial]\tilde u^{(2)}}_{\tilde{e}_2}-\underbrace{R(\chi_\partial \tilde u^{(2)})}_{\tilde e_3}.
    \end{align}
    Now at $tf$ where $\chi_\partial \equiv 1$, we see that $\tilde e_1$ vanishes up to one extra order in $\rho_{tf}$  as $\sigma^2\tilde\Box(1)\tilde u^{(2)}=\rho^2$ a priori. Moreover the commutator $\tilde{e}_2$ picks up an extra order of $\rho_{tf}$ on top of the fact that $\sigma^2\tilde\Box(1)$ already gives $\rho_{tf}^2$-decay (using also again that $\tilde u^{(2)}$ is smooth conormal on $X_{res}^+\backslash zf$). Furthermore $\tilde{e}_3$ is also of order $\rho_{tf}^3$ due to $R\in \rho_b\rho_{tf}^3\text{Diff}^2_b(X)$. For $bf$ we can gain a similar improvement, note that here also $\chi_\partial\equiv 1$ hence that 
    \begin{align}
        \tilde e =\rho^2-\hat{\Box}(\sigma)\tilde u^{(2)},
    \end{align}
where we emphasize here that $\tilde u^{(2)}$ is smoothly defined in the coordinate $\hat{\rho} =\frac{\rho}{\sigma}$ on $X^+_{res}\backslash zf$. Expanding using the normal operator at $bf$ we have 
\begin{align}
    \tilde e&= \rho^2-N_{bf}(\hat{\Box}(\sigma))\tilde u^{(2)}-R_{bf}\tilde u^{(2)} \text{ with } R_{bf}\in \rho_{bf}^2\rho_{tf}^2\text{ Diff }_b^2(X).
\end{align}
Now observe that $N_{bf}(\hat\Box(\sigma))=2i\sigma^2\hat \rho (\hat\rho \partial_{\hat\rho}-\frac{3}{2})$ so $N_{bf}(\hat\Box(\sigma))-\sigma^2\tilde\Box(1)=\sigma^2\hat\rho^2(-(\hat\rho\partial_{\hat\rho})^2+2\hat\rho\partial_{\hat\rho}+\Delta_{\mathbb{S}^3})$ hence we see that 
\begin{align}
    \tilde {e} = \sigma^2\hat\rho^2-\underbrace{(N_{bf}(\hat\Box(\sigma))- \sigma^2\tilde \Box(1))\tilde u^{(2)}}_{\mathcal{O}(\sigma^2\hat{\rho}^3)}-\underbrace{\sigma^2\tilde\Box(1)\tilde u^{(2)}}_{:=\sigma^2\hat\rho^2}-R(\tilde u^{(2)}).
\end{align}
Hence there is also a leading order cancellation in $\tilde{e}$ at $bf$. Combined with the statement from Lemma ... that $\tilde{u}^{(2)}$ is of $1-$ order at bf we conclude that $\tilde{e}\in \mathcal{A}^{3-,3,0-}(X_{res}^+)$ indeed improves over the standard $\rho^2\in \mathcal{A}^{2,2,0}(X_{res}^+)$.\\
In the above we have been careful not to declare any positive vanishing of $\tilde e$ at zf as we know that $N_{zf}(\hat{\Box}(\sigma)=\hat{\Box}(0)$ has a zero resonance, hence we could have boundedness or a logarithmic singularity at zf. By Lemma~\ref{lem:tfResolventOutput}, with $\widetilde f\equiv 1$,
we have near $\mathrm{zf}$, in the coordinate $\widehat r=\sigma/\rho$ using (\ref{eq:GeneralLogofV}) that
\begin{equation}\label{eq:Tildeu2Expansion}
    \tilde u^{(2)}\sim c\log(\hat r) + v +\mathcal{A}^{1-}(X)
\end{equation}
where $c=-\frac{1}{2}$ due to $L_0(\log\hat r) = -\frac{1}{2}$. Thus, using Lemma \ref{lem:NormalOperatorOnEachFace} and equation (\ref{eq:GeneralSpectralFamily})  we see that by expanding at \text{zf} 
\begin{align}
    \hat{\Box}(\sigma)(\chi_\partial \tilde u^{(2)})=\underbrace{N_{zf}(\hat{\Box}(\sigma))}_{\hat\Box(0)}(\chi_{\partial}\tilde u^{(2)}) + R(\chi_\partial \tilde u^{(2)}) \; \text{ with }\; R\in \rho_{bf}\rho_{tf}^2\rho_{zf}\text{Diff}_b^2(X_{res}^+)
\end{align}
Looking at the leading order approximate solution we then see by plugging in (\ref{eq:Tildeu2Expansion})
\begin{align}
    \hat{\Box}(0)(\chi_\partial \tilde u^{(2)}) &= \hat{\Box}(0)(-\frac{1}{2}\log(\hat{r})\chi_\partial)+\mathcal{O}(1)\nonumber\\
    &=-\frac{1}{2}\log(\hat{r})\hat{\Box}(0)\chi_\partial +\underbrace{[\hat{\Box}(0),\frac{1}{2}\log(\hat r)](-\chi_\partial)}_{\mathcal{O}(1)}+\mathcal{O}(1).
\end{align}

From this we conclude that $\tilde{e}\sim \rho^2-\hat{\Box}(0)(\chi_\partial\tilde u^{(2)})$ at $zf$ must satisfy
\begin{align}
    \tilde e \in \mathcal{A}^{3-,3,(0,1)+1-}(X_{res}^+),
\end{align}

with leading order term at zf given by $\frac{1}{2}\log(\hat{r})\hat{\Box}(0)(\chi_\partial)$. \\
Now let us rewrite $\hat{\Box}(\sigma)^{-1}\tilde e$ algebraically once more:
    \begin{align}\label{eq:TildeEAlgebraically}
        \hat{\Box}(\sigma)^{-1}\tilde e = \tilde v -\sigma \hat{\Box}(\sigma)^{-1}[\sigma^{-1}(\hat{\Box}(\sigma)-\hat{\Box}(0))]\tilde v, \text{ with } \;\tilde v:= \hat{\Box}(0)^{-1}\tilde e.
    \end{align}

    Now we recall that near zf, working with the coordinates $\sigma = \rho \hat{r}$ and $\hat{r}=\frac{\sigma}{\rho}$ that we have 
    \begin{align}
        \sigma \partial_\sigma\sim \hat{r}\partial_{\hat{r}}.
    \end{align}
    Thus we have that $(\sigma\partial_{\sigma})^2\log(\hat{r})=(\hat{r}\partial_{\hat{r}})^2\log(\hat{r})=0$, i.e. the logarithmic decay vanishes and we have 
    $(\sigma\partial_\sigma)^2\tilde{e}\in \mathcal{A}^{3-,3-,1-}(X_{res}^+)$, and thus by Lemma \ref{lem:ResolvedonResolvent} we have for $\alpha =1$ that 
    \begin{align}{\label{eq:DoubleDerofTildevisBounded}}
        (\sigma\partial_\sigma)^2 \tilde v\in \mathcal{A}^{1-,1-,1-}(X_{res}^+).
    \end{align}
    On the other hand, using only the rougher estimate $\tilde{e}\in \mathcal{A}^{3-,3-,0-}(X_{res}^+)\subset\mathcal{A}^{0-}([0,1)_\sigma;\mathcal{A}^{3-}(X))$ we defined earlier we get again with Proposition \ref{prop:FirstMappingProposition}, that 
    \begin{align}\label{eq:TildevisBounded}
    \tilde{v} \in \mathcal{A}^{0-}([0,1)_\sigma;\mathcal{A}^{1-}(X))\subset \mathcal{A}^{1-,1-,0-}(X_{res}^+). 
    \end{align}
    These two facts tell us that if we were to integrate with regard to $s=\log(\sigma)$ twice (equivalently the $b$-density $\frac{d\sigma}{\sigma}$) from some fixed basepoint $\sigma_0>0$ that by (\ref{eq:DoubleDerofTildevisBounded}) that $\tilde v$ is at most affine in $s=\log(\sigma)$, i.e.
    \begin{align}
        \tilde v(\sigma,\rho,\omega) = A(\rho,\omega) + B(\rho,\omega)\log(\sigma) + R(\sigma,\rho,\omega), 
    \end{align}
    where $A,B\in \mathcal{A}^{1-,1-,0}(X_{res}^+)$ because of the bound (\ref{eq:TildevisBounded}) and 
    \begin{align}
    R(\sigma,\rho,\omega)=\int^{\log(\sigma)}_{\log(\sigma_0)}(\log(\sigma)-\log(t))  (t\partial_t)^2\tilde v (t,\rho,\omega) \frac{dt}{t}\in \mathcal{A}^{1-,1-,1-}(X_{res}^+).
    \end{align}
    Certainly near zf, we have the polyhomogeneous expansion (in $\hat{r}$) given by
    \begin{equation}
        \tilde v \sim c_0 + c_1 \log(\hat{r}) + \mathcal{O}(\hat{r}^{1-\epsilon}),
    \end{equation}
    whence $\tilde v\in \mathcal{A}^{1-,1-,(0,1)+1-}(X_{res}^+)$.
    Now again using Lemma \ref{lem:tfResolventOutput}, we see see that the logarithmic term for $\tilde{v}$ can be directly computed from the one we found for $\tilde{e}$, i.e. 
    \begin{align}
    \frac{1}{2}\log(\hat{r})\hat{\Box}(0)^{-1}(\hat{\Box}(0)\chi_\partial).
    \end{align}
    Now by Lemma \ref{lem:KernelElements}, we know that the unique solution $\tilde{v}_{0,1}$ (modulo kernel) of $\hat{\Box}(0)\tilde{v}_{0,1}=\hat{\Box}(0)\chi_\partial$ is given a priori by $\tilde{v}_{0,1}=\chi_\partial - c u_{(0)}$, where $u_{(0)}\in 1+\mathcal{A}^{1-}(X)$ is given by Lemma \ref{lem:KernelElements}. However for $\tilde{v}_{0,1}$ to lie in the proper domain for $\hat{\Box}(0)^{-1}$, i.e. $\mathcal{A}^{1-}(X)$ we must have that the constant term vanishes, which forces $c=1$, hence we conclude that logarithmic term of $\tilde{v}$ at zf is therefore $\frac{1}{2}\log(\hat{r})(\chi_\partial - u_{(0)})$. \\
    For the second term (\ref{eq:TildeEAlgebraically}) we want to show that $\sigma \hat{\Box}(\sigma)^{-1}[\sigma^{-1}(\hat{\Box}(\sigma)-\hat{\Box}(0))]\tilde v \in \mathcal{A}^{1-,1-,1-}(X_{res}^+)$. However this follows readily from (\ref{eq:GeneralSpectralFamily}) and Lemma \ref{lem:AdmissibleWaveEquation}. More specifically, recall that $\hat{\Box}(\sigma)-\hat{\Box}(0)=2i\sigma\rho Q + \sigma ^2g^{00}$ with $g^{00}\in \rho^2C^\infty(X)$ and $Q\in \text{Diff}_b(X)$. On $X_{res}^+$ this then implies that if $\tilde{v}\in \mathcal{A}^{1-,1-,0}(X_{res}^+)$ then $2i\sigma\rho Q \tilde {v} \in \mathcal{A}^{2-,2-,0}(X_{res}^+)$ as we only multiply by $\rho$, picking up an order at $bf$ and $tf$ and multiplying by $\sigma$ only improves the weights at tf and zf. Trivially we also have that $\sigma^2g^{00}\tilde v \in \mathcal{A}^{3-,3-,0}(X_{res}^+)$ so overall we conclude that 
    \begin{align}
        \sigma^{-1}(\hat{\Box}(\sigma)-\hat{\Box}(0)):\mathcal{A}^{1-,1-,0}(X_{res}^+)\to \mathcal{A}^{2-,2-,0}(X_{res}^+).
    \end{align}
    Using again that $\mathcal{A}^{2-,2-,0-}(X_{res}^+)\subset |\sigma|^{-\delta}\mathcal{A}^0([0,1)_\sigma;\mathcal{A}^{2-\delta-}(X))$ for any $\delta>0$. Using Lemma \ref{lem:MappingIntoXres} we then have that $\hat{\Box}(\sigma)^{-1}$ maps $\mathcal{A}^{2-\delta}(X)$ into $\mathcal{A}^{2-\delta,-\delta-,-\delta-}(X_{res}^+)$. Multiplication by $\sigma$ then improves $zf$ and $tf$ decay and we get that $\sigma\hat{\Box}(\sigma)^{-1}$ maps $\mathcal{A}^{2-,2-,0}(X_{res}^+)$ into $|\sigma|^{1-\delta}\mathcal{A}^{1-\delta,-\delta-,-\delta-}(X_{res}^+)=\mathcal{A}^{1-,1-,1-}(X_{res}^+)$.\\
    From this we conclude that the only logarithmic contributions to $\hat{\Box}(\sigma)^{-1}\rho^2$ are the ones coming from $\chi_\partial \tilde{u}^{(2)}$ in equation (\ref{eq:OriginalErrorTerm}) and from the first term $\tilde{v}$ in equation (\ref{eq:TildeEAlgebraically}). Explcitly they are given by 
    \begin{align}
        \frac{1}{2}\log(\hat{r})(-\chi_\partial)+\frac{1}{2}\log(\hat{r})(\chi_\partial-u_{(0)}) = -\frac{1}{2}\log(\hat{r})u_{(0)}.
    \end{align}
\end{proof}

\begin{remark}
    It should be noted here that $\tilde u^{(2)}$ can again be evaluated explicitly as was done in the Kerr case in \cite{hintz2022sharp}. Changing variables from $\hat{\rho}$ to $\hat{r}=\hat{\rho}^{-1}$ we see that $\tilde\Box(1)$ can be written as 
    \begin{align}
        \tilde\Box(1)= -2 i \hat{r}^{-1}(\hat{r}\partial_{\hat{r}}+\frac{3}{2}) + \hat{r}^{-2}(-(\hat{r}\partial_{\hat{r}})^2 - 2\hat{r}\partial_{\hat r}+\Delta_{\mathbb{S}^3}),
    \end{align}
    and we wish to solve $\tilde\Box(1)u(\hat{r})=\hat{r}^{-2}$. Note that this amounts to solving this regular-signular ODE is equivalent to solving $\hat{r}^2\tilde\Box(1)=1$ i.e.
    \begin{align}
        (-(\hat{r}\partial_{\hat{r}})^2 -2\hat{r}\partial_{\hat{r}}+\Delta_{\mathbb{S}^3})-\underbrace{2i\hat{r}(\hat{r}\partial_{\hat{r}}+\frac{3}{2})}_{\mathcal{O}(\hat{r})}u(\hat{r})=1.
    \end{align}
    Note that the indicial equation has roots $\alpha=0,-2$ for solutions $\hat{r}^\alpha$. The zero-branch corresponds to a $u_{hom}(\hat{r})=c_0+c_1\log(\hat{r})+\mathcal{O}(r^\epsilon)$. Solving the indicial equations yields $c_0 =0$, $c_1=-\frac{1}{2}$.
    Unfortunately we cannot factorize the full ODE neatly into two Euler equations but what we can do is remove the linear terms.
    First we set $u(\hat{r}):=\hat{r}^{-3/2} y(\hat{r})$. We then get that 
    \begin{align}
        y''(\hat{r})+2iy'(\hat{r}) +\frac{3i}{4\hat{r}^2}=-\hat{r}^{-1/2}.
    \end{align}
    We then subistitute $y(\hat{r})=e^{-i\hat{r}}w(\hat{r})$ to get 
    \begin{align}
        w''(\hat{r})+\left(1-\frac{3}{4\hat{r}^2} \right)w(\hat{r})=-\hat{r}^{-1/2}e^{i\hat{r}}.\label{eq:BesselSchroedingerForm}
    \end{align}
    Now we recognize this as a Bessel-Schr\"oedinger ODE with forcing term. If we impose the Sommerfeld condition $(\partial_{\hat{r}}-i)\hat{r}^{3/2}u(\hat{r})=0$ as $r\to \infty$ or equivalently $u(r)\sim r^{-3/2}e^{ir}$, then its particular solution is a linear $\hat{r}$-dependent combination of 
    \begin{align}
        w_1(\hat{r})=\sqrt{\hat{r}}J_1(\hat{r}),\;\;\; w_2(\hat{r}) =\sqrt{r}H^{(1)}_1(\hat{r}),
    \end{align}
    where $J_1$ and $H^{(1)}_1$ are the Bessel function of the first kind and an outgoing Hankel function, i.e. $H^{(1)}_1=J_1+iY_1$ with $Y_1$ a Bessel function of the second kind. The equation has constant Wronksian $W(w_1,w_2)(s)\equiv\frac{2i}{\pi}$. 
    Finally we vary parameters and set $f(\hat{r}) = -\hat{r}^{1/2}e^{i\hat{r}}$.The unique solution regular at $\hat{r}=0$ and outgoing at $\infty$ is then given by
    \begin{align}
        w(\hat{r})=w_1(\hat{r}) \int^\infty_{\hat{r}} \frac{w_2(s)f(s)}{W}ds + w_2(\hat{r}) \int^{\hat r}_0 \frac{w_1(s)f(s)}{W}ds.
    \end{align}
    or writing out in full:
    \begin{align}
        w(\hat{r}) =\frac{\pi}{2i} \left[ \sqrt{\hat{r}}J_1(\hat{r}) \int^\infty_{\hat{r}}H^{(1)}_1(s)e^{is} ds - \sqrt{\hat r}H^{(1)}_1(\hat{r}) \int^{\hat{r}}_0 J_1(s) e^{is}ds \right].
    \end{align}
    Hence substituting back we get 
    \begin{align}\label{eq:FullFormOfU}
        u(\hat{r})=e^{-i\hat{r}}\frac{\pi}{2i\hat{r}} \left[ J_1(\hat{r}) \int^\infty_{\hat{r}}H^{(1)}_1(s)e^{is} ds - H^{(1)}_1(\hat{r}) \int^{\hat{r}}_0 J_1(s) e^{is}ds \right],
    \end{align}
    where the first integral term is understood as the limiting absorption outgoing condition (to guarantee convergence)
    \begin{align}
        \int^\infty_{\hat{r}}e^{is}H^{(1)}_1 ds = \lim_{\epsilon\searrow 0}\int^{\infty}_{\hat{r}}e^{i(s+i\epsilon)}H^{(1)}_1(s) ds.
    \end{align}
    
    Note that we have 
    \begin{align}
        &\frac{d}{ds}\left(e^{is}(is H^{(1)}_0 + s H^{(1)}_1(s)-H^{(1)}_0(s)) \right)=e^{is}H^{(1)}_1(s)\\
        &\int^{\hat{r}}_0 J_1(s) e^{is}ds = 1 - e^{i\hat{r}}+\hat{r}e^{i\hat{r}}J_1(\hat{r}) + i \hat{r}J_0(\hat{r}),
    \end{align}
    and with this we can write $u_{part}(\hat{r})$ in the fully closed form 
    \begin{align}
        u_{part}(\hat{r}) =\frac{\pi e^{-i\hat{r}}}{2i\hat{r}}\left[-J_1(\hat{r})(i\hat{r}H^{(1)}_{0}(\hat{r})+\hat{r}H^{(1)}_1(\hat{r}) - H^{(1)}_0(\hat{r}))- H_1^{(1)}(\hat{r}) (1 - e^{i\hat{r}}+\hat{r}e^{i\hat{r}}J_1(\hat{r}) + i \hat{r}J_0(\hat{r}))
        \right]
    \end{align}
    Now that we have the particular solution $u_{part}$ we note that the total outgoing solution can be writen as 
    \begin{align}
        \tilde{u}^{(2)}(\hat{r})=u_{part}(\hat{r}) + C u_{hom}(\hat{r}).
    \end{align}
    Note that regular and singular branch to the homogeneous radial equation can be found by solving\footnote{We neglect the $2i u'(\hat{r})$ and $\frac{3i}{\hat{r}}u(\hat{r})$ here which are more regular in $\hat{r}$. }
    \begin{align}
        u''(\hat{r})+\frac{3}{\hat{r}}u'(\hat{r})=0,
    \end{align}
    which for $u(\hat{r})=\hat{r}^\alpha$ has indicial roots $\alpha =0$ (regular) and $\alpha =-2$ (singular). This motivates us to require that $u_{hom}(\hat{r})\in \mathcal{A}^{-2+}(X)$.
    Furthermore, the homogeneous solution corresponding to the Sommerfeld radiation condition is given by 
    \begin{align}
     u_{hom}(\hat{r}) = \frac{e^{-i\hat{r}}}{\hat{r}}H^{(1)}_1(\hat{r}),
    \end{align}
    and given that $H^{(1)}_1(\hat{r}) \sim -\frac{2i}{\pi\hat{r}}+\mathcal{O}(\hat{r}\log \hat{r})$ we see that this forces $C=0$.
    Finally, using that $Y_0(\hat{r})\sim \frac{2}{\pi}(\log\frac{\hat{r}}{2}+\gamma) +\mathcal{O}(\hat{r}^2\log\hat{r})$, $Y_1(\hat{r})\sim -\frac{2}{\pi\hat{r}}+\frac{\hat{r}}{\pi}(\log\frac{\hat{r}}{2}+\gamma-\frac{1}{2})+\mathcal{O}(\hat{r}^3\log\hat{r})$, $J_0(\hat{r)}\sim 1-\frac{\hat{r}^2}{4}+...$, $J_1(\hat{r})\sim \frac{\hat{r}}{2}+...$ we can expand $u(\hat{r})$ at $\hat{r}=0$ to find
    \begin{align}
    \tilde u^{(2)}(\hat{r}) \sim-\frac{1}{2}\log \hat{r}-\frac{1}{2}(\gamma-\log 2) - \frac{3}{4}+\frac{i\pi}{4}+\mathcal{O}(\hat{r}\log \hat{r)},
    \end{align}
    where $\gamma$ is the Euler-Mascheroni constant. For a full solution to the model equation at $\text{tf}$ see also Appendix A and Theorem \ref{thm:GlobalAsymptotics}.
\end{remark}

\subsection{Resolvent expansion.}

With the refinement of the output of $\hat{\Box}_{g_{\bfa,\bfb,\mfrm}}(\sigma)^{-1}$ on conormal inputs established in the previous section, we can now state the first theorem on the resolvent expansion.

\begin{theorem}\label{thm:SingularRegularDecomposition}
    Let $f=f(x)\in \mathcal{A}^{4+\alpha}(X)$ where again $\alpha\in (0,1)$. For positive frequencies, the resolvent acting on $f$ is then the form 
    \begin{align}
        \hat{\Box}(\sigma)&= u_{\text{reg}}(\sigma) + u_{\text{sing}}(\sigma),
    \end{align}
    where the singular term is given by 
    \begin{align}\label{eq:SingularTerm}
        u_{\text{sing}}(\sigma) \in \sigma^2 \mathcal{A}^{1-,((0,0),1-),((0,1),1-)}(X_{res}^+)
    \end{align}
    and the regular term is given by  \begin{align}
        u_{reg}(\sigma)\in C^\infty([0,1)_{\sigma};\mathcal{A}^{1-}(X)))+\mathcal{A}^{\alpha-,2+\alpha-,((2,0),2+\alpha-)}(X_{res}^+).
    \end{align}
    Here the leading terms of $\sigma^{-2}u_{sing}(\sigma)$ at zf $-\frac{1}{2}\log(\frac{\sigma}{\rho})c_X(f)u_{(0)}$ and $c_{X}(f)\tilde u^{(2)}$ where 
    \begin{align}
       c_X(f)=\frac{1}{4\pi^2}\langle f,u_{(0)}^*\rangle,
    \end{align}
    and at tf $-c_{(0)}\tilde u^{(2)}$ with $\tilde u^{(2)}=\tilde\Box(1)^{-1}(\hat\rho^2)$ and $u_{(0)}^*$ given by Lemma \ref{lem:KernelElements}.
\end{theorem}
\begin{remark}
    We emphasize here that the singular term $c_X(f)$ does not depend on the mass $\mfrm>0$, in contrast to the Kerr case (compare \cite{hintz2022sharp}) Theorem 3.1. Once we have established the global asymptotics on $(M,g_{\bfa,\bfb,\mfrm})$ we will link this back to the absence of the strong Huygen's principle.
\end{remark}

\textbf{Outline of the algorithm.} We briefly sketch the algorithm here using the propositions from the preciding algorithm. We start with conormal data $f\in \mathcal{A}^{4+\alpha}(X)$, where we take $0<\alpha<1$. We repeatedly apply the resolvent
\begin{align}
    \hat{\Box}(\sigma)^{-1}f= \hat{\Box}(0)^{-1}f -\sigma \hat{\Box}(\sigma)^{-1}(\sigma^{-1}(\hat{\Box}(\sigma)-\hat{\Box}(0))\hat{\Box}(0).
\end{align}
Then by Corollary \ref{cor:SpectralFormAdmissibleWaveOpertor}, we have that 
\begin{align}
    \sigma^{-1}(\hat{\Box}(\sigma)-\hat\Box(0))= 2i \rho (Q_0 +\tilde Q) +\sigma g^{00},\qquad Q_0=\rho\partial_\rho-\frac{3}{2},
\end{align}
where $\tilde Q \in \rho^2\text{Diff}^{1}_b(X)$ and $g^{00}\in \rho^2C^\infty(X)$. In particular,
\begin{align}
    Q_0\rho^2=\frac{1}{2}\rho^2,\qquad Q_0\rho=-\frac{1}{2}\rho. \label{eq:Q0Identities}
\end{align}

\begin{enumerate}
    \item \textbf{First zero-energy inversion.} Lemma \ref{lem:ZeroEnergyResolventOutputAsymptotics} gives 
    \begin{align}
        u_0:=\hat{\Box}(0)^{-1}f=c_{(0)}\rho^2 +\tilde{u}_0, \qquad \tilde{u}_0 \in \mathcal{A}^{2+\alpha-}(X),
    \end{align}
    where we find back the leading coefficient 
    \begin{align}
        c_{(0)}=\frac{1}{4\pi^2}\langle f,u_{(0)}^*\rangle.
    \end{align}
    \item \textbf{Input of the next iteration.} Define the remainder first iteration data through 
    \begin{align}
          f_{1}(\sigma) = -\sigma^{-1}(\hat{\Box}(\sigma)-\hat{\Box}(0)).
    \end{align}
    Then, the identity $Q_0\rho^2=\frac{1}{2}\rho^2$ in (\ref{eq:Q0Identities}) gives
    \begin{align}
        f_1(\sigma)=-ic_{(0)}\rho^3 +\mathcal{A}^{3+\alpha-}(X) +\sigma \mathcal{A}^4(X),
    \end{align}
    and the resolvent identity becomes 
    \begin{align}
        \hat{\Box}(0)^{-1}f=u_0 + \sigma \hat{\Box}(\sigma)^{-1}f_1(\sigma).
    \end{align}
    \item \textbf{Second zero-energy inversion.} Lemma \ref{lem:FirstIteration} then shows that 
    \begin{align}
        u_1=\hat{\Box}(0)^{-1}f_1(0) = -ic_{(0)}\rho +\mathcal{A}^{1+\alpha-}(X),
    \end{align}
    where the leading coefficient follows again from (\ref{eq:Q0Identities}) with
    \begin{align}
        L_0\rho=\rho, \qquad L_0=-(\rho\partial_\rho)^2+2\rho\partial_\rho +\Delta_{\mathbb{S}^3}.
    \end{align}
    \item \textbf{The borderline input.} Now we set 
    \begin{align}
        f_2(\sigma):=\sigma^{-1}(\hat{\Box}(\sigma)-\hat{\Box}(0))u_1,
    \end{align}
    and since we have 
    \begin{align}
        2i\rho Q_0(-ic_{(0)}\rho)=-c_{(0)}\rho^2,
    \end{align}
    again by (\ref{eq:Q0Identities}), and we obtain
    \begin{align}
        f_2(\sigma)=\underbrace{-c_{(0)}\rho^2}_{f_{2,0}},\qquad \tilde f_{2}(\sigma)\in \mathcal{A}^{2+\alpha-}(X)+\sigma\mathcal{A}^3(X).
    \end{align}
    \item \textbf{The logarithmic singularity and remaining terms.} To account for the frequency dependence of $f_1(\sigma)$ put 
    \begin{align}
        f_2:=\sigma^{-1}(f_1(\sigma)-f_1(0))=-g^{00} u_0 \qquad u_2:=\hat{\Box}(0)^{-1}f_2 \in\mathcal{A}^{1-}(X)
    \end{align}
    yielding the expansion 
    \begin{align}
        \hat{\Box}(\sigma)^{-1}= u_0+\sigma u_1 +\sigma^2 u_2 +c_{(0)}\sigma^2\hat{\Box}(\sigma)^{-1}\rho^2 +\mathcal{A}^{\alpha-,2+\alpha-,((2,0),2+\alpha-)}(X_{res}^+)
    \end{align}
    The singular contirubtiion yielding the late time asymptotics is then given by 
    \begin{align}
        u_{sing}(\sigma):=c_{(0)}\sigma^2\hat{\Box}(\sigma)^{-1}\rho^2,
    \end{align}
    as by Lemma \ref{lem:ThresholdOutput} the leading terms of $\sigma^{-2}u_{sing}$ are 
    \begin{align}
        -\frac{1}{2}c_{(0)}\log\left(\frac{\sigma}{\rho}\right)u_{(0)} \text{ at zf},\qquad c_{(0)}\tilde u_{(2)} \text{ at tf,} 
    \end{align}
    which produces the $\sigma^2\log \sigma$ singularity on compact spatial sets. 
\end{enumerate}
\begin{remark}
    We emphasize here that this $\sigma^2\log(\sigma)$ is again the same as found for the Kerr spacetime \cite{hintz2022sharp}. It is however caused by different reasons, as the indicial gap size is increased from $(0,1)$ to $(0,2)$ but the roots of the transport $Q_0$ are also fractional $(3/2)$, hence $Q_0\rho=-\frac{1}{2}\rho$ whereas in Kerr $Q_0\rho=0$. Thus the leading borderline source is already generated by the flat leading operators, without an additional mass prefactor. The shifting of the indicial roots finally then causes the output on radial data in $\rho^3$ to have no resonance output of decay $\rho\log\rho$ in Lemma \ref{lem:FirstIteration}. Thus the Myers-Perry is the $d=4$ case of the mechanism described in \cite{hintz2023linear} for generally flat asymptotic spacetimes in odd dimensions (particularly Proposition 6.4).
\end{remark}
\begin{table}[htbp]
    \centering
    
    \renewcommand{\arraystretch}{1.5}
    \begin{tabular}{lcc}
        \hline
        Stage & Kerr & Myers--Perry \\
        \hline
        Relevant terms in $u_0$
        & $c_K\rho + \mfrm c_{(0)}\rho^2$
        & $c_{{{(0){}}}}\rho^2$ \\

        First transport action
        & $Q_0\rho = 0$
        & $Q_0\rho^2 = \frac12\rho^2$ \\

        Leading spherical term in $f_1(0)$
        & $-2imc_{(0)}\rho^3$
        & $-ic_{{(0)}}\rho^3$ \\

        Corresponding term in $u_1$
        & $2i\mfrm c_{(0)}\rho\log\rho$
        & $-ic_{(0)}\rho$ \\
        \hline
    \end{tabular}
    \caption{Comparison of the resolvent iterations for Kerr and
    $1+4$-dimensional Myers--Perry. Only the spherically symmetric
    contributions relevant to the leading singularity are displayed.
    Here $c_{(0)}$ denotes the respective leading zero-energy coeffcient respective to each spacetime.}
    \label{tab:kerr-myers-perry-resolvent}
\end{table}

\begin{lemma}\label{lem:ZeroEnergyResolventOutputAsymptotics}
    Let $\alpha \in (0,1)$. For $f\in \mathcal{A}^{4+\alpha}(X)$ we have 
    \begin{align}
        \hat{\Box}(0)^{-1} f = c_{(0)} \rho^2 + \tilde{u},
    \end{align}
    where $c_{(0)}=(4\pi^2)^{-1}\langle f,u^*_{(0)}\rangle$,  $\tilde{u}\in \mathcal{A}^{2+\alpha-}(X).$
\end{lemma}

\begin{proof}
    Clearly, we have $u:=\hat{\Box}(0)^{-1}f\in \mathcal{A}^{1-}(X)$. From this we have that 
    \begin{align}\label{eq:MinkowskiLaplacianinTermsofBox}
        L_0 u =\rho^{-2}f - \tilde L u \in \mathcal{A}^{2+\alpha-}(X)+\mathcal{A}^{3-}(X).
    \end{align}
    Let $\chi_\partial\in C_c^\infty(X)$ be a cutoff supported near $\partial X$, then $u_\partial:=\chi_\partial u \in \mathcal{A}^{1-}(X)$ satisfies an equation 
    \begin{align}
        L_0(u_\partial)=f_\partial \in \mathcal{A}^{2+\alpha-}(X) \label{eq:MinkowskiLaplacianCutOff}
    \end{align}
    We now improve the decay of $u_\partial$ using a $b$-normal operator argument similar to Lemma \ref{lem:tfResolventOutput}, this time by Mellin transforming in the variable $\rho \sim \hat{r}^{-1}$, i.e.
    \begin{align}
        \widehat{u_\partial}(\lambda):= \int^\infty_0 \rho ^{-i\lambda}u_\partial(\rho) \frac{d\rho}{\rho},
    \end{align}
    and (\ref{eq:MinkowskiLaplacianCutOff}) becomes 
    \begin{align}
        \widehat{L_0}(\lambda)\widehat{u_\partial}(\lambda) = \widehat{f_\partial}(\lambda), \; \;\;\;\; \widehat{L_0}(\lambda):=-(i\lambda)^2 + 2 i\lambda +\cancel{\Delta}.
    \end{align}
    Again by the same arguments as before, we know that $\widehat{u_\partial}$ is holomorphic in $\text{ Im }\lambda >-1$ and satisfies 
    \begin{align}
        \|\widehat{u_\partial}(\lambda)\|_{H^s(\mathbb{S}^3)} \leq C_{sN\epsilon}\langle \lambda \rangle^{-N}, \;\;\;\; -1 + \epsilon< \text{Im }\lambda<0, 
    \end{align}
    for all $s,N\in \mathbb{R}$ and $\epsilon >0$ by analytic continuation. Since $f_\partial\in \mathcal
    A^{2+\alpha-}(X)$ the same estimate holds for $\widehat{f_\partial}(\lambda)$ but in the region $-2-\alpha < \text{ Im } \lambda$ due to its decay. Once again expanding $\widehat{u_\partial}$ into spherical harmonics, we have $\widehat{u_\partial}=\sum_{mlk}u_{mlk}(\lambda)Y_{mlk}$ and  $\widehat{f_\partial}=\sum_{mlk}f_{mlk}(\lambda)Y_{mlk}$. Once again we have that 
    \begin{align}
        \left.\widehat{L_0}(\lambda)\right|_{\mathcal{Y}_{kl}}=-(i\lambda)^2+2i\lambda + l(l+2),
    \end{align}
    which is invertibile for $i\lambda \neq -l, l+2$. Thus, for $u_{mlk}(\lambda)= \left(\left.\widehat{L_0}(\lambda)\right|_{\mathcal{Y}_{kl}}\right)^{-1}f_{mlk}(\lambda)$ the only relevant poles would be at $l=0$ $i=-2$. Hence in the inverse Mellin transform 
    \begin{align}
        u_\partial(\rho) = \frac{1}{2\pi}\int_{\text{Im} \lambda=-\epsilon}\rho^{i\lambda}\widehat{L_0}(\lambda)^{-1}\widehat{f_\partial}(\lambda) d\lambda,
    \end{align}
    we can shift the contour of integration through the pole at $\lambda = -2i$ to $\text{ Im }\lambda =-2-\alpha+\epsilon$. The residue theorem then gives $u_\partial(\rho)\in c_{(0)}\rho^2 + \mathcal{A}^{2+\alpha-}(X)$,
    \begin{align}\label{eq:PartialExpansionofu}
        u(\rho)\sim c_{(0)}\rho^2 +\tilde u_0, \text{ with } \tilde u_0\in\mathcal{A}^{2+\alpha-}(X)
    \end{align}
    The constant $c_{(0)}$ can be evaluated. Let $\chi_{\epsilon}(\rho)=\chi(\rho/\epsilon)$ where $\chi \in C^\infty([0,\infty))$ is a function supported away from $0$ and identically $1$ for $\rho/\epsilon>1$, we have 
    \begin{align}
         0&=\langle u, \hat{\Box}(0)^*u^*_{(0)}\rangle = \lim_{\epsilon\to 0} \langle \chi_\epsilon u,\hat{\Box}(0)^*u^*_{(0)}\rangle =\lim_{\epsilon\to 0}\langle \hat{\Box}(0)\chi_\epsilon u,u_{(0)}^*\rangle\\
         &=\langle f, u^*_{(0)}\rangle +\lim_{\epsilon\to 0}\langle[\hat{\Box}(0),\chi_\epsilon]u,u_{(0)}^*\rangle,\label{eq:IntegrationCommutator}
    \end{align}
   where we used again the identity $\hat{\Box}(0) \chi_\epsilon u=\chi_\epsilon\hat{\Box}(0)u+ [\hat{\Box}(0),\chi_\epsilon]u$. $[\hat{\Box}(0),\chi_\epsilon]\in \rho^2\text{Diff}_b^1(X)$ converges to $0$ strongly as an operator  $\mathcal{A}^{\beta}(X)\to \mathcal{A}^{\beta+2}(X)$ for any $\beta\in \mathbb{R}$ and this implies that the formal paring $\langle[\hat{\Box}(0),\chi_\epsilon]\tilde{u}_0,u^*_{(0)}\rangle$ vanishes in the limit $\epsilon \to 0$ and we conclude that the second term in the sum vanishes in (\ref{eq:IntegrationCommutator}). Note that $dg_{\bfa,\bfb,\mfrm}|X$ is a smooth positive multiple of $r^3 dr dg_{\mathbb{S}^3} = \rho^{-4} \frac{d\rho}{\rho}dg_{\mathbb{S}^3}$ hence convergence to $0$ of the commutator in $\mathcal{A}^{4+\delta}(X)$ suffice for convergence of the inner product to $0$. For the same reason subleading terms in the limit do not contribute. Thus using that $u_{(0)}^*\sim 1 + \mathcal{A}^{1-}(X)$ we see that using that $\hat{\Box}(0)\sim\rho^2L_0$ that 
   \begin{align}
       c_{(0)}&\lim_{\epsilon\to 0}\langle[\hat{\Box}(0),\chi_\epsilon]\rho^2,1\rangle\\
       &=c_{(0)} \text{vol}(\mathbb{S}^3)\lim_{\epsilon\to0}\int[\rho^2(-\rho\partial_\rho[\rho\partial_\rho,\chi_\epsilon]-[\rho\partial_\rho,\chi_\epsilon](\rho\partial\rho-2))\rho^2\cdot \rho^{-4}\frac{d\rho}{\rho}],
   \end{align}
   which after writing this term out completely and substituting $\rho=\epsilon x$ we find that 
   \begin{align}
       c_{(0)}&\lim_{\epsilon\to 0}\langle[\hat{\Box}(0),\chi_\epsilon]\rho^2,1\rangle\\
       &=-\text{vol}(\mathbb{S}^3)\int^\infty_0 3\chi'(x) + x\chi''(x) dx\\
       &=-\text{vol}(\mathbb{S}^3)\left([x\chi'(x)]_{0}^\infty+\int^\infty_02\chi'(x) dx\right)\\
       &=-2\text{vol}(\mathbb{S}^3)c_{(0)}=-4\pi^2c_{(0)},
   \end{align}
   from which we conclude that 
   \begin{align}
       c_{(0)}=\frac{1}{4\pi^2}\langle f,u_{(0)}^*\rangle.
   \end{align}
\end{proof}

Now, let us denote the output of Lemma \ref{lem:ZeroEnergyResolventOutputAsymptotics} by 
\begin{align}
    u_0:=\hat{\Box}(0)^{-1}f=c_{(0)}\rho^2 + \tilde u, \text{ with } \tilde u \in \mathcal{A}^{2+\alpha-}(X),
\end{align}

and by Lemma \ref{lem:AdmissibleWaveEquation}, we now wish to analyze the output of
\begin{align}\label{eq:DifferentiatedOperator}
    \sigma^{-1}(\hat{\Box}(\sigma)-\hat{\Box}(0))=-2i\rho(Q_0+\tilde{Q})+\sigma g^{00}.
\end{align}
We note that we have $Q_0 \rho^2 = \frac{1}{2}\rho^2$ so unlike in Kerr, there is no cancellation when the operator (\ref{eq:DifferentiatedOperator}) acts on the first order term of the expansion of $u$. We then see that 
\begin{align}
    f_1(\sigma)&:=-\sigma^{-1}(\hat{\Box}(\sigma)-\hat{\Box}(0))u_0\\
    &=-i c_{(0)}\rho^3+\mathcal{A}^{4+\alpha-}(X) + \sigma\mathcal{A}^{4+\alpha-}(X)
\end{align}
The aim is now to analyze the expansion 
\begin{align}
    \hat{\Box}(\sigma)^{-1}f&=u_0+\sigma\hat{\Box}(\sigma)^{-1}f_1(\sigma)\\
    &=u_0+\sigma \hat{\Box}(0)f_1(\sigma)-\sigma\hat
    \Box(\sigma)^{-1}(\hat{\Box}(\sigma)-\hat\Box(0))\hat{\Box}(0)^{-1}f_1(\sigma)
\end{align}

To this end we now set 
\begin{align}
    u_1:=\hat{\Box}(0)^{-1}f_1(0) \in \mathcal{A}^{1-}(X)
\end{align}

\begin{lemma}\label{lem:FirstIteration}
    Let $\alpha\in (0,1)$. We have that
    \begin{align}
        u_1=-ic_{(0)}\rho+\mathcal{A}^{1+\alpha-}(X).
    \end{align}
\end{lemma}

\begin{proof}
   We again consider 
   \begin{align}
       L_0 u_1= \rho^{-2}f_1(0) - \tilde{L}u_1\in \mathcal{A}^1(X)+\mathcal{A}^{3-}(X)=\mathcal{A}^1(X).
   \end{align}
We Mellin transform in $\rho$ once again with a cutoff $\chi_\partial$ that is constant to $1$ near $\partial X$. Let us set $u_{1,\partial}:=\chi_\partial u_1$, and Let us set $f_{1,\partial}:= \chi_\partial L_0 u_1 +[L_0,\chi_\partial]u_1$. We then have that 
\begin{align}
    L_0 u_{1,\partial}=f_{1,\partial} =- i c_{(0)}(\rho\chi_\partial)+\mathcal{A}^{1+\alpha-}(X).
\end{align}
Applying the Mellin transform in $\rho$ we see that 
\begin{align}
    \widehat{u_{1,\partial}}(\lambda)&=-ic_{(0)}\left.\frac{\widehat{\rho\chi_\partial}(\lambda)}{-(i\lambda)^2+2i\lambda+l(l+2)}\right|_{l=0}+ R(\lambda)\\
    &=-ic_{(0)}\widehat{L}_0(\lambda)^{-1}\widehat{\rho\chi_\partial}(\lambda)+R(\lambda).\label{eq:MellinInverseIndicialFamily}
\end{align}
where $R(\lambda)$ has no poles in the strip $-1-\alpha < \text{Im}\,\lambda<0$ and we evaluate on $l=0$ as the displayed leading term $u_1$ is the constant spherical mode. Thus, the Mellin transform of $\widehat{\rho\chi_\partial}(\lambda)$ has a simple pole at $\lambda=-i$, whilst the inverse indicial family $\widehat{L}_0(\lambda)$ has no poles in $-2<\text{Im }\lambda<0$. We therefore shift the inverse Mellin contour to $\text{Im }\lambda=-1-\alpha+\epsilon$, for $0<\epsilon<\alpha$, crossing only the pole at $-i$ whilst remaining above the indicial pole at $-2i$. As $\widehat{L}_0(-i) = 1$, we find the leading coefficient to be $-ic_{(0)}$
from which we conclude
\begin{align}
    u_1(\rho)=-ic_{(0)}\rho+\mathcal{A}^{1+\alpha-}(X).
\end{align}
\end{proof}
\begin{remark}
    Note that in the Kerr case (see \cite{hintz2022sharp}, Lemma 3.3) the indicial family is given by $\widehat{L}_0(\lambda)=-(i\lambda)^2+i\lambda +l(l+1)$. Hence it is at this step that the contourshifting of the equivalent of (\ref{eq:MellinInverseIndicialFamily}) encounters a double pole for the $l=0$ mode, leading to an additional $\rho\log\rho$ term in this step of the expansion.
\end{remark}

Now we plug this back we are set to develop the expension. To this end we do not only develop around $\hat{\Box}(0)^{-1}$ but also around $f_1(0)$. To this end let us set $u'_1=\hat{\Box}(0)^{-1}(f'_1)\in \mathcal{A}^{1-}(X)$ where $f'_1=\sigma^{-1}(f_1(\sigma)-f_1(0))\in \mathcal{A}^{4+\alpha-}(X)$.  We then get in the first step of the expansion that 
\begin{align*}
    \hat{\Box}(\sigma)^{-1}f&=\underbrace{\hat{\Box}(0)^{-1}f
}_{u_0}+(\hat{\Box}(\sigma)^{-1}-\hat{\Box}(0)^{-1})f\\
    &=u_0+\sigma\hat{\Box}(\sigma)^{-1}\underbrace{[\sigma^{-1}(\hat{\Box}(0)-\hat{\Box}(\sigma))\hat{\Box}(0)^{-1}f]}_{f_1(\sigma)}\\
    &=u_0+\sigma\hat\Box(\sigma)^{-1}f_1(\sigma)\\
    &=u_0+\sigma[(\hat\Box(\sigma)^{-1}-\hat{\Box}(0)^{-1})f_1(\sigma)+\underbrace{\hat{\Box}(0)^{-1}(f_1(\sigma)-f_1(0))}_{\sigma u_1'}+\underbrace{\hat{\Box}(0)^{-1}f_1(0)}_{u_1}]\\
    &=u_0+\sigma u_1 +\sigma^2u_1' + \sigma(\hat{\Box}(\sigma)^{-1}-\hat{\Box}(0)^{-1})(f_1(\sigma)-f_1(0))+\sigma(\hat{\Box}(\sigma)^{-1}-\hat{\Box}(0)^{-1})f_1(0)]\\
    &=u_0+\sigma u_1+\sigma^2u_1' + \sigma^2\hat{\Box}(\sigma)^{-1}(\hat{\Box}(0)-\hat{\Box}(\sigma))\underbrace{\hat{\Box}(0)^{-1}f_1'}_{u_1'}+\sigma\hat{\Box}(\sigma)^{-1}(\hat{\Box}(0)-\hat{\Box}(\sigma))\underbrace{\hat{\Box}^{-1}(0)f_1(0)}_{u_1}\\
    &=u_0+\sigma u_1 + \sigma^2 u_1'-\sigma\hat{\Box}(\sigma)^{-1}(\hat\Box(\sigma)-\hat\Box(0))u_1 -\sigma^2\hat\Box(\sigma)^{-1}(\hat\Box(\sigma)-\hat\Box(0))u'_1.
\end{align*}

We now analyze the $4$-th term in the expansion. Once again we have 
\begin{align}
    \sigma^{-1}(\hat{\Box}(\sigma)-\hat{\Box}(0))u_1&=(-2i\rho(Q_0+\tilde Q)+\sigma g^{00})u_1:=f_{2,0}+\tilde{f}_2(\sigma)\\
    &\text{ with } f_{2,0}=-c_{(0)}\rho^2, \tilde{f}_2(\sigma)\in \mathcal{A}^{2+\alpha-}(X)+\sigma\mathcal{A}^{3+\alpha-}(X).
\end{align}
When it comes to this term, we only neeed to analyze the output of $\hat{\Box}(\sigma)^{-1}f_{2,0}$ and $\hat{\Box}(\sigma)^{-1}\tilde f_2(\sigma)$. By Lemma \ref{lem:ThresholdOutput} we have that 
\begin{align}
    \hat{\Box}(\sigma)^{-1}f_{2,0}\in \mathcal{A}^{1-,((0,0),1-),((0,1),1-)}(X_{res}^+),
\end{align}
with leading order at zf equal to $\log(\frac{\sigma}{\rho})c_{(0)}u_{(0)}$ and at tf equal to $-c_{(0)}\tilde u^{(2)}$, which is exactly the singular contribution. To analyze $\sigma^2\hat{\Box}(\sigma)^{-1}\tilde f_2(\sigma)$ we notice first that by Proposition \ref{prop:FirstMappingProposition}, that 
\begin{align}
    \sigma^2\hat{\Box}(\sigma)^{-1}\tilde f_2(0) \in \sigma^2\mathcal{A}^{\alpha-,\alpha-,((0,0),\alpha-)}(X_{res^+})=\mathcal{A}^{\alpha-,2+\alpha-,((2,0),2+\alpha-)}(X_{res}^+),
\end{align}
and 
\begin{align}
    \sigma^2\hat{\Box}(\sigma)^{-1}(\tilde f_2(\sigma)-\tilde f_2(0))\in \sigma^3\mathcal{A}^{1-,0-,0-}(X_{res}^+)=\mathcal{A}^{1-,3-,3-}(X_{res}^+). \label{eq:ConormalRemainder}
\end{align}
To see (\ref{eq:ConormalRemainder}) more precisely, note that $\tilde f_2(\sigma)-\tilde f_2(0)\in \sigma\mathcal{A}^{3-}(X)$, hence it lies certainly also in $\sigma \mathcal{A}^{2-\alpha}(X)$ for any small $\alpha>0$, hence the statement follows from Lemma \ref{lem:MappingIntoXres} and the fact that multiplying by $\sigma^3$ increases tf and zf regularity by $3$. Finally, we note that $\sigma^2\hat{\Box}(\sigma)^{-1}(\hat{\Box}(\sigma)-\hat{\Box}(0))u_1'\in 
\sigma \mathcal{A}^{2-}(X)+\sigma^3\mathcal{A}^{3-}(X)$ as we have that term $(\hat{\Box}(\sigma)-\hat{\Box}(0))u_1'\in \mathcal{A}^{2-}(X)+\sigma \mathcal{A}^{3-}(X)$. As we also have $\sigma \mathcal{A}^{2-}(X)+\sigma^3\mathcal{A}^{3-}(X)\subset \mathcal{A}^{1-,3-,3-}(X_{res}^+)$ this final term requires no futher consideration. Plugging all statements back into the expansion we get that 
\begin{align}
    \hat{\Box}(\sigma)^{-1}f& \in \mathcal{A}^{(2,0),2+\alpha-}(X) +\sigma \mathcal{A}^{(1,0),1+\alpha-}(X)+\sigma^2\mathcal{A}^{1-}(X)\nonumber\\&+\sigma^2\mathcal{A}^{1-,((0,0),1-),((0,1),1-)}(X_{res}^+) +\mathcal{A}^{\alpha-,2+\alpha-,((2,0),2+\alpha-)}(X_{res}^+).\label{eq:FullResolventExpansion}
\end{align}

\subsection{Asymptotic behaviour of waves on $1+4$ dimensional Myers-Perry}

\begin{theorem}\label{thm:DecayonCompactSubsets}
    Let $\alpha \in (0,2)$ and let $f=f(t_*,x)\in C^\infty_c(\mathbb{R}_{t_*};\mathcal{A}^{4+\alpha}(X))$. Then the unique forward solution of $\Box \phi = f $ satisfies 
    \begin{align}
        |\phi(t_*,x)-c_{M}(f)u_{(0)}t_*^{-3}|\leq Ct_*^{-3-\alpha+},
    \end{align}
    for $x\in K\Subset X^\circ$ with $K$ compact. In the above, $c_M(f)$ is given by 
    \begin{align}
        c_M(f)=-\frac{1}{4\pi^2}\int_{\mathbb{R}_{t_*}\times X}f(t_*,x)u^*_{(0)}(x) dg_{\bfa,\bfb,\mfrm}.
    \end{align}

    Consequently, derivatives of $\phi$ decay as 
    \begin{align}
        |\partial_{t_*}^j\partial_x^\beta(\phi-c_M(f) u_{(0)}t_*^{-3})|\leq C_{j\beta}t_*^{-3-\alpha -j+}\;\;\;\;\forall j\in\mathbb{N}_0, \beta \in \mathbb{N}^3_0.
    \end{align}
\end{theorem}

\begin{proof}
 We start out with the Fourier transform $\hat{\phi}(\sigma,x)=\int_{\mathbb{R}}e^{i\sigma t_*}\phi(t_*,x) dt_* $ we then can express $\phi$ as
 \begin{align} 
 \phi(t_*,x)=\mathcal{F}^{-1}(\hat{\Box}(\sigma)^{-1}\hat{f}=\frac{1}{2\pi}\int_{\mathbb{R}}e^{-i\sigma t_*} \hat{\Box}(\sigma)^{-1}\hat{f}(\sigma,x)d\sigma,
 \end{align}
    which gives the unique forward solution from the Payley-Wiener theorem upon deforming the integration contour $iC+\mathbb{R}$ and then taking the limit $C\to \infty$. It is here that the mode stability assumption in Definition \ref{def:SpectrallyAdmissbleMetric} comes into play.
    Let now $\chi \in C_c^\infty((-1,1))$ be a bump function supported near $\sigma = 0$. We split $\phi = \phi_0+\phi_1+\phi_2$ with 
    \begin{align}
        \widehat{\phi_0}(\sigma)&:= \chi(\sigma)\hat{\Box}(\sigma)^{-1}\hat{f}(0),\\
        \widehat{\phi_1}(\sigma)&:=\sigma\chi(\sigma)\hat{\Box}(\sigma)^{-1}(\sigma^{-1}(\hat{f}(\sigma)-\hat{f}(0))),\\
        \widehat{\phi_2}(\sigma)&:=(1-\chi(\sigma))\hat{\Box}(\sigma)^{-1}\hat{f}(\sigma).
    \end{align}
    Now we know that $1-\chi(\sigma)$ vanishes to arbitary order at $\sigma =0$. Therefore $\widehat{\phi_2}$ allows us to consider high energy estimates for $\hat{\Box}(\sigma)^{-1}$ away from the singular structure at $\sigma = 0$. In particular we have that $\hat{\Box}(\sigma)^{-1}$ and its $\sigma$-derivatives are smoothly controlled by Lemma \ref{lem:SigmaDerivativeControl}. Now as $f$ is also compactly supported in $t_*$ we know that $\tilde{f}$ is Schwarz in $\sigma$. Thus, upon taking the inverse Fourier transform we have that 
    \begin{align}\label{eq:phi2Schwarz}
        \phi_2(t_*,x)\in \mathscr{S}(\mathbb{R}_{t_*},\mathcal{A}^{1-}(X)).
    \end{align}
    Now we consider the low energy piece $\phi_0$. We are restricting to $\rho>0$ (recall that we consider $x\in K\Subset X^\circ$), therefore by Theorem \ref{thm:SingularRegularDecomposition} we have that $\phi_0$ is of class 
    \begin{align}
        C^\infty([0,1)_\sigma;C^\infty(X^\circ))+\mathcal{A}^{2+\alpha-}([0,1)_\sigma;C^\infty(X^\circ)) -\sigma^2(\frac{1}{2}\log(\sigma)-\frac{i\pi}{4})c_{(0)}u_{(0)}+C^\infty(X^\circ).
    \end{align}
    The first two terms are coming from the conormal regularity of $u_{reg}$ in Theorem \ref{thm:SingularRegularDecomposition} upgraded to smoothness on the interior of $X$, anFd the last part is the contribution from $u_{sing}$. As in the case of Kerr, the smooth pieces of the expansion (\ref{eq:FullResolventExpansion}) fit together for $\sigma>0$ and $\sigma<0$ at $\sigma=0$. For real $\sigma$ we have $\hat{\Box}(\sigma)^{-1} f = \overline{\hat{\Box}(-\sigma)^{-1}\overline{f}}$, but for $\sigma<0$ with a branch cut of log around $-i[0,\infty)$, with 
    \begin{align}
        \overline{\frac{1}{2}\log(-\sigma)-\frac{i
        \pi}{4}}=\frac{1}{2}\log(\sigma+i0)-\frac{i
        \pi}{4},
    \end{align}
    hence the logarithmic terms glue at $\pm \sigma>0$ to form a $\sigma^2\log(\sigma+i0^+)$ term. Absorbing $-\frac{i\pi}{4}$ into the smooth parts we find that 
    \begin{align}
        \widehat{\phi_0}(\sigma)\in C^\infty([0,1)_\sigma;C^\infty(X^\circ))+\mathcal{A}^{2+\alpha-}(\pm[0,1)_\sigma;C^\infty(X^\circ)) -\frac{1}{2}\sigma^2\log(\sigma+i0^+)c_{(0)}u_{(0)}+C^\infty(X^\circ),
    \end{align}
    where we note that the entire piece $\widehat{\phi_1}(\sigma)$ can be absorbed into the remainder space of $\widehat{\phi_0}$ due to it having an extra factor of $\sigma$.
    Now, as in Kerr, we argue that the Fourier transform of the smooth term contributes terms that are Schwarz in $t_*$, the Fourier transform of the conormal function class have a decay of $t_*^{-3-\alpha+}$ in all $\partial_x, t_*\partial_{t_*}$-derivatives. So the main decay contribution to $\phi_0$ comes from the logarithmic singularity.\\
    As this term is again of the kind $\sigma^2\log(\sigma+i0^+)$ the rest of this theorem proceeds entirely analogous to Kerr. In short, the inverse Fourier transform of $(\sigma+i0)^z$ for $z\in \mathbb{C}$ is $\mathcal{F}^{-1}((\sigma+i0)^z)=e^{i\pi z/2}\chi_+^{-z-1}(t_*)$ where $\chi^z_+(t_*)$ is the holomorphic continuation from $\text{Re }z>-1$ to $\mathbb{C}$ of $(t_*)/\Gamma(z+1)$ with $(t_*)_+ =0$, and $\Gamma$ the Euler gamma function. Therefore for $t_*>0$,
    \begin{align}
        \mathcal{F}^{-1}(\sigma^k\log(\sigma+i0))=\partial_z\mathcal{F}^{-1}((\sigma+i0)^z)|_{z=k}=-e^{i\pi k/2}\partial_z(\chi_+^z(t_*))|_{z=-k-1}.
    \end{align}
    To evaluate this, we note that in $t_*>0$, 
    \begin{align}
        \partial_z\chi_+^z(t_*)|_{z=-k-1}&=\left.\partial_z\left(\frac{(z+1)\dots(z+k+1)}{\Gamma(z+k+2)}t_*^z\right)\right|_{z=-k-1}\\
        &=\frac{(z+1)\dots(z+k+1)}{\Gamma(z+k+2)}t_*^{-k-1
        }\\
        &=(-1)^kk!t_*^{-k-1}.
    \end{align}
    In summary we find that 
    \begin{align}
        \mathcal{F}^{-1}(\sigma^2\log(\sigma+i0))=\begin{cases}
            &2t_* \; \text{ for }\; t_*>0\\
            &0 \; \text{ for } t_*\leq 0
        \end{cases}
    \end{align}
    giving the $t_*^{-3}$ leading order term with the stated constant.
\end{proof}

In the above we used the result for inverse Fourier transforms on inverse conormal functions, found directly in \cite{hintz2022sharp} Lemma 3.6.

\begin{lemma}[Inverse Fourier transform of conormal functions]
\label{lem:InverseFourierConormal}
Let $\beta>-1$, and let
$\widehat{\phi}\in\mathcal A^\beta([0,1)_\sigma)$
have compact support in $[0,1)$. Extend $\widehat{\phi}$ by
zero to $\mathbb R_\sigma$, and define
\[
    \phi(t_*)
    :=\frac{1}{2\pi}\int_0^1
      e^{-i\sigma t_*}\widehat{\phi}(\sigma)\,d\sigma.
\]
Then $\phi\in C^\infty(\mathbb R)$ is a symbol of order
$-1-\beta$. More precisely, for every $j\in\mathbb N_0$,
\[
    |\partial_{t_*}^j\phi(t_*)|
    \le C_j\langle t_*\rangle^{-1-\beta-j}.
\]
In particular,
\[
    |(t_*\partial_{t_*})^j\phi(t_*)|
    \le C_j'\langle t_*\rangle^{-1-\beta}.
\]
The same conclusion holds for negative-frequency support,
upon replacing $\widehat{\phi}(\sigma)$ by
$\widehat{\phi}(-\sigma)$ in the hypotheses.
\end{lemma}

\begin{proof}
Conormality implies that, for every $k\in\mathbb N_0$,
\[
    |\partial_\sigma^k\widehat{\phi}(\sigma)|
    \le C_k\sigma^{\beta-k},
    \qquad 0<\sigma<1.
\]
Since $\beta>-1$, the inverse Fourier integral and all its
$t_*$-derivatives converge absolutely. Thus
$\phi\in C^\infty(\mathbb R)$.

We first prove the estimate for $j=0$. It suffices to
consider $T:=|t_*|\ge2$. Choose
$\chi\in C_c^\infty([0,\infty))$ with $\chi=1$ on $[0,1]$
and $\operatorname{supp}\chi\subset[0,2)$, and write
\[
    2\pi\phi(t_*)=I(t_*)+II(t_*),
\]
where
\begin{align*}
    I(t_*)
    &=\int_0^1 e^{-i\sigma t_*}
      \chi(T\sigma)\widehat{\phi}(\sigma)\,d\sigma,\\
    II(t_*)
    &=\int_0^1 e^{-i\sigma t_*}
      (1-\chi(T\sigma))\widehat{\phi}(\sigma)\,d\sigma.
\end{align*}
The first term satisfies
\[
    |I(t_*)|
    \le C\int_0^{2/T}\sigma^\beta\,d\sigma
    \le C'T^{-1-\beta}.
\]

For the second term, choose an integer $N>\beta+1$.
Integrating by parts $N$ times gives
\[
    II(t_*)
    =(it_*)^{-N}\int_0^1 e^{-i\sigma t_*}
      \partial_\sigma^N
      \bigl((1-\chi(T\sigma))\widehat{\phi}(\sigma)\bigr)
      \,d\sigma.
\]
There are no boundary terms, since the integrand vanishes
near both endpoints. Derivatives of $\chi(T\sigma)$ are
supported where $\sigma\asymp T^{-1}$, so the product rule
and conormality give
\[
    \left|
      \partial_\sigma^N
      \bigl((1-\chi(T\sigma))\widehat{\phi}(\sigma)\bigr)
    \right|
    \le C_N\sigma^{\beta-N}
\]
on its support, which is contained in $\{\sigma\ge T^{-1}\}$.
Consequently,
\[
    |II(t_*)|
    \le C_NT^{-N}\int_{T^{-1}}^1
      \sigma^{\beta-N}\,d\sigma
    \le C_N'T^{-1-\beta}.
\]
This proves the estimate for $\phi$.

For $j\ge0$, apply the same argument to
$(-i\sigma)^j\widehat{\phi}(\sigma)
 \in\mathcal A^{\beta+j}([0,1))$, whose inverse Fourier
transform is $\partial_{t_*}^j\phi$. This yields
\[
    |\partial_{t_*}^j\phi(t_*)|
    \le C_j\langle t_*\rangle^{-1-\beta-j}.
\]
The estimates for $(t_*\partial_{t_*})^j\phi$ follow by
expanding this operator as a linear combination of
$t_*^k\partial_{t_*}^k$, $0\le k\le j$.
Negative-frequency support is handled by reflection.
\end{proof}

\begin{remark}\label{rem:AmplitudeRemark}
The proof also applies to amplitudes with values in
$\mathcal A^0(X)$, with all estimates understood in each
conormal seminorm.
\end{remark}

\subsection{Global asymptotics}

\begin{definition}
    Let $X=\overline{\mathbb{R}^4}$, and let us compactify the nonnegative half line the same way we $\mathbb{R}^4$ i.e. by defining $[0,\infty]_{t_*}:=([0,\infty)_{t_*}\sqcup [0,\infty)_{\tau})/\sim$ where $\sim$ glues $(0,\infty)_{t_*}$ and $(0,\infty)_{\tau}$ via $\tau = t_*^{-1}$. This setup allows us to define the compactification of the causal future of $t_*=0$ inside of $M^\circ =(\mathbb{R}_{t_*}\times X^\circ,g)$ is defined as the blow-up 
    \begin{align}
        M_+:= [[0,\infty]_{t_*}\times X;\{\infty\}\times\partial X],
    \end{align}
    which contains $[0,\infty)\times X^\circ \subset M^{\circ}$ as an open dense subset. Four boundary hypersurfaces $M_+$ will be important in the analysis of the global asymptotics. They are given by
    \begin{enumerate}
        \item the Cauchy surface $\Sigma=t_*^{-1}(0)\cong X$;
        \item $\mathscr{J}^+$: the lift of $[0,\infty]_{t_*}\times\partial X$;
        \item Minkowski future time like infinty $I^+$: the front face;
        \item spotatially compact infinity $\mathcal{K}^+$: the lift of $\tau^{-1}(0)$.
    \end{enumerate}
    \end{definition}
    With respect to these boundaries, we can again introduce weights.
    \begin{definition}
        Let $\rho_{\mathscr{J}^+},\rho_{I^+},\rho_{\mathcal{K}^+}\in C^\infty(M_+)$ be boundary defining functions of $\mathscr{J}^+,I^+,\mathcal{K}^+$ respectively. We define the space of conormal functions on $M_+$ as 
        \begin{align}
            \mathcal{A}^{\alpha,\beta,\gamma}(M_+):=\rho^\alpha_{\mathscr{J}^+}\rho^\beta_{I^+}\rho^\gamma_{\mathcal{K}^+}\mathcal{A}^{0,0,0}(M_+),
        \end{align}
        which are smooth at $\Sigma$ and remain bounded upon application of any finite number of smooth vector fields on $M_+$ that are tangent to $\mathscr{J}^+$
    \end{definition}

    Notice that we can simplify our local coordinates to some extent depending what part of $M_+$ we considering. Away from $\mathcal{K}^+$, we have that 
    \begin{align}
        \rho_{\mathscr{J}^+,1}&:=\rho/\tau =\rho t_*=t_*/r, \; \;\rho_{I^+,1}:=\tau=t_*^{-1},
    \end{align}
    are defining function for $\mathscr{J}^+$ and $I^+$ respecitively. On the other hand away from $\mathscr{J}^+$ we have that 
    \begin{align}
        \rho_{\mathcal{K}^+,0}:=\tau/\rho=(\rho t_*)^{-1}=r/t_*, \; \rho_{I^+,0}:=\rho=r^{-1},
    \end{align}
      
    in particular we have,
    \begin{align}
        I^+=[0,\infty]_{\rho_{\mathscr{J}^{+,1}}} \times \partial X,
    \end{align}
    with $I^+\cap \rho^{-1}_{\mathscr{J}^{+},1}(0)$ the future boundary of $\mathscr{J^+}$. For a diagram respresenting the relevant geometry we refer back to Figure \ref{fig:myers-perry-forward} in Section \ref{subsec:ExtensionCompactification}.

    \begin{theorem}\label{thm:GlobalAsymptotics}
        (Global asymptotics) Let $\alpha\in (0,2)$, $(\bfa,\bfb)\in D_{\mfrm}$ and $f\in C^\infty_c((0,\infty)_{t_*};\mathcal{A}^{4+\alpha}(X))$. Denote with $u$ the unique forward solution to $\Box_{g_{\bfa,\bfb,\mfrm}} u =f$. Then 
        \begin{align}
            \phi\in \mathcal{A}^{((3/2,0),3+\alpha),((3,0),3+\alpha-),((3,0),3+\alpha-)}(M_+).
        \end{align}
        Setting $c_M(f)=-\frac{1}{4\pi^2}\langle f, u_{(0)}^*\rangle$, the leading order terms of $u$ are 
        \begin{enumerate}
            \item at $\mathcal{K}^+$: $c_M(f) u_{(0)} t_*^{-3}$.\\
            \item at $I^+$, $c_M(f) u^+(\rho t_*)t_*^{-3}$ \text{ where } 
            \begin{align}
                u^+(v) :=\left(\frac{v}{v+2}\right)^{3/2}.
            \end{align}
            \item At $\mathscr{J}^+$: $\rho_{\mathscr{J}^+,1}^{3/2}\phi_{rad}(t_*,\omega)$ where the radiation field $\phi_{rad}\in \mathcal{A}^{((3,0),3+\alpha-)}(\mathscr{J}^+)$ satisfies 
            \begin{align}
                \phi_{rad}(t_*,\omega)-\frac{1}{2\sqrt{2}}c_M(f)t_*^{-3} \in \mathcal{A}^{3+\alpha-}(\mathscr{J}^+)
            \end{align}
        \end{enumerate}
    \end{theorem}
\begin{proof}
As in Theorem \ref{thm:DecayonCompactSubsets} we will use the splitting 
\begin{align}
    \phi= \phi_0+\phi_1+\phi_2,
\end{align}
where we know that $\phi_2\in \mathcal{A}^{1-,\infty,\infty}(M_+)$ due to $\phi_2 \in \mathscr{S}(\mathbb{R}_{t_*};\mathcal{A}^{1-}(X))$ by (\ref{eq:phi2Schwarz}).
\begin{itemize}
    \item (Low energy contribution: regular part) We consider $\phi_0$ starting with the contribution from the regular part $u_{reg}(\sigma)$ coming from Theorem \ref{thm:SingularRegularDecomposition}. Writing 
    \begin{align}
        u_{reg}(\sigma) = u_{reg,1}(\sigma)+ u_{reg,2} (\sigma) \; \text{ with } u_{reg,1}\in C^\infty ([0,1)_\sigma;\mathcal{A}^{1-}(X))], u_{reg,2}\in \mathcal{A}^{\alpha-,2+\alpha-,((2,0),2+\alpha-)}(X_{res}^+).
    \end{align}
    we have as before 
    \begin{align}    
        \phi_{1,reg}\in \mathcal{A}^{1-,\infty,\infty}(M).
    \end{align}
    The remaining piece, $u_{reg,2}$ we glue for $\sigma>0$ and $\sigma<0$ to give
    \begin{align}
        u_{reg,2}&= u_{reg,2,0}+\sum_{\pm} u_{reg,2,\pm} \\
        &\text{ with } u_{reg,2,0}\in \mathcal{A}^{\alpha-,2+\alpha-}(X_{res}), \text{ and }\;u_{reg,2,\pm} \in \mathcal{A}^{\alpha-,2+\alpha-,2+\alpha-}(X_{res}^\pm).
    \end{align}

Now since $u_{reg,\pm} \in \mathcal{A}^{2+\alpha-}(\pm [0,1)_\sigma;\mathcal{A}^0(X))$, we have by Lemma \ref{lem:InverseFourierConormal} that  
\begin{align}
    \phi_{0,reg,\pm}:=\mathcal{F}^{-1}u_{reg,2,\pm}\in \mathcal{A}^{3+\alpha-}([0,1)_{\tau};\mathcal{A}^0(X))\subset \mathcal{A}^{0,3+\alpha-,3+\alpha-}(M_+).
\end{align}
Because $\phi_{0,reg,2,0}=\mathcal{F}^{-1}u_{reg,2,0}$ we used that $\phi_{0,,reg,2,0}$ vanishes rapidly in compact spatial sets as $u_{reg,2,0}$ is smooth across $\sigma =0$. We refine this argument at $\rho=\sigma =0$ and set 
\begin{align}
    \psi = \psi(\sigma/\rho) \in C_c^\infty(\mathbb{R}), \; \psi\equiv 1 \text{ near } 0 
\end{align}
is a cutoff then $(1-\psi) u_{reg,2} \in \sum_{\pm}\mathcal{A}^{\alpha-,2+\alpha-,\infty}(X_{res}^\pm)$ and thus we can conclude that a fortriori $\mathcal{F}^{-1}((1-\psi)u_{reg,2})\in\mathcal{A}^{0,3+\alpha-,3+\alpha-}(M_+)$. It remains to see that $\psi u_{reg,2}\in \mathcal{A}^{\infty,2+\alpha-}(X_{res})$ as a function
\begin{align}\label{eq:BoundOnpsiu}
    \psi u_{reg,2}(\rho,\hat{r},\omega)\in \mathcal{A}^{2+\alpha-}([0,1)_\rho; C^\infty_c(\mathbb{R}_{\hat{r}}\times\partial X_\omega), \;\;\;\;\;\;\; \hat{r}=\frac{\sigma}{\rho}.
\end{align}
Its inverse Fourier transform is thus
\begin{align}
    \frac{1}{2\pi}\int^\infty_{-\infty}e^{-i\sigma t_*}\psi u_{reg,2}(\rho,\frac{\sigma}{\rho},\omega)\;d\sigma&=\frac{\rho}{2\pi}\int_{-\infty}^\infty \psi u_{reg,2}(\rho,\hat{r},\omega) d\hat{r}\\
    &=\rho\mathcal{F}^{-1}_{v\to \hat{r}}(\psi u_{reg,2})(\rho,\rho t_*,\omega) \lesssim \rho \rho^{2+\alpha-}\langle\rho t_*\rangle^{-N}.
\end{align}
Where the last bound comes from (\ref{eq:BoundOnpsiu}) whence we connclude that $\rho \mathcal{F}^{-1}_{v\to\hat{r}}(\psi u_{reg,2}) \in \mathcal{A}^{3+\alpha-,3+\alpha-,\infty}(M_+)$. Putting all estimates together we conclude that for the regular contributions coming from $u_{reg,2}$ lead to 
\begin{align}
    \phi_{0,reg,2}=\mathcal{F}^{-1} u_{reg,2}\in \mathcal{A}^{0,3+\alpha-,3+\alpha-}(M_+).
\end{align}
\item (Low energy contribution: singular part): We now forucs on the main contribution from $\phi_0$ at $I^+$ and $\mathcal{K}^+$ comes from the singular part $u_{sing}(\sigma)$ in Theorem \ref{thm:SingularRegularDecomposition}. In particular $u_{sing}$ at $tf$ is $1$-homogeneous in the dilation $(\rho,\sigma)\mapsto(\lambda\rho,\lambda\sigma)$ as in fact $u_{sing}(\sigma,\rho,\omega)\sim \sigma^2c_X(f)\tilde u^{(2)}(\rho/\sigma,\omega)$. For this reason the leading term gets annihilated up to leading order by $\sigma\partial_\sigma+\rho\partial_{\rho}-2$. Thus the $tf$ vanishing order gets improved by $1$ for $\tilde u_{sing}:=(\sigma\partial_\sigma+\rho\partial_\rho-2 )u_{sing}$ and we conclude
\begin{align}
   \tilde u_{sing}\in \sum_{
   \pm } \mathcal{A}^{1-,3-,((2,1),3-)}(X_{res}^\pm).
\end{align}
Let $\psi$ be a cutoff function as before, we see that $(1-\psi)\tilde u_{sing}\in \mathcal{A}^{1-,3-,\infty}(X^\pm_{res})\subset \mathcal{A}^{3-}([0,1)_\sigma;\mathcal{A}^0(X))$ have inverse Fourier transform in $\mathcal{A}^{0,4,4-}(X)$. The piece $\psi\tilde u_{sing}\in \sum_{\pm} \mathcal{A}^{\infty.3-,((2,1),3-)}(X_{res}^\pm)$ is the sum of two terms, $\psi \tilde u_{sing} = \psi\tilde u_{sing,1} +\psi\tilde u_{sing,2}$, where 
\begin{align}\label{eq:SingularContributionSplitup}
    \tilde {u}_{sing,1}&=a(\rho,\omega)\hat{r}^2\log(\hat{r}+i0), \; a\in\mathcal{A}^{3-}(X)\\ \tilde{u}_{sing,2}&\in \mathcal{A}^{\infty,3-,3-}(X_{res}^\pm) \subset\mathcal{A}^{3-}(\pm[0,1)_\sigma;\mathcal{A}^0(X)) \text{ for } \pm\sigma>0.
\end{align}
We note that 
\begin{align}
    \tilde\phi_{sing,1}=\mathcal{F}^{-1}(\psi\tilde u_{sing,1})=\int_{-\infty}^\infty e^{-i\sigma t_*} \psi a(\rho,\omega)\hat{r}^2\log(\hat{r}+i0) d\sigma\lesssim \rho \rho^{3-}(\rho t_*)^{-3}
\end{align}
Hence $\tilde \phi_{sing,1} \in \mathcal{A}^{4-,4-,((3,0),4-}(M_+)$. We also see directly from (\ref{eq:SingularContributionSplitup}) that $\tilde{\phi}_{sing,2}:=\mathcal{F}^{-1}(\psi\tilde u_{sing,2})\in\mathcal{A}^{0,4-,4-}(M_+)$. Taking intersections we see that in $(t_*,\rho)$ coordinates after Fourier transforming back (using $\sigma\partial\sigma \mapsto -t_*\partial_{t_*}-1$) that
\begin{align}
    \tilde{\phi}_{sing}:= \mathcal{F}^{-1}(-t_*\partial_{t_*}+\rho\partial_{\rho}-3)\phi_{sing}\in\mathcal{A}^{0,4-,((3,0),4-)}(M_+), 
\end{align}
with $\tilde\phi_{sing}=\mathcal{F}^{-1}u_{sing}=(-t_*\partial_{t_*}+\rho\partial_{\rho}-3)\phi_{sing} \in \mathcal{A}^{0,4-,((3,0),4-)}(M_+).$
Note that the Euler vector field annihilates $\rho^3_{I^+}$, and indeed we have  
\begin{align}
    -t_*\partial_{t_*}+\rho\partial_{\rho}-3=\rho_{I^+,1}\partial_{\rho_{I^+,1}}-3=\rho_{I^+,0}\partial_{\rho_{I^+,0}}-3
\end{align}

Where the second expression holds at $I^+$ close to $\mathscr{J}^+$ with boundary defining functions $\rho_{I^+,1}=t_*^{-1}$, $\rho_{\mathscr{J}^+,1}=\rho t_*$, and the last expression holds near $\mathcal{K}^+$ with $\rho_{\mathcal{K},0}:=(\rho t_*)^{-1}$, $\rho_{I^+,0}:=\rho$.
Using degree $-2$ homoegeneity under the transformation $(\rho,\sigma)\mapsto(\lambda\rho,\lambda\sigma)$, for $\lambda>0$ we can match the $\rho_{I^+}^3$ leading order term of $\phi$ at $I^+$ to the inverse Fourier transform of $u_{sing}$ at tf, i.e. the inverse Fourier transform of 
\begin{align}
    u_{sing,2}= \sigma^2 c_X(\hat{f}(0))(\tilde{u}^{(2)}(\hat{r}) H(\hat{r}) +\overline{\tilde u^{(2)}}(-\hat{r})H(-\hat{r}),
\end{align}
where as usual $\hat{r}=\sigma\rho^{-1}$. This glues together the leading orders from $X_{res}^\pm$ at $tf$ leading to a logarithmic singularity at the seam $\sigma=0$ of the form $\hat{r}^2\log(\hat{r}+i0)$ as it does in the Kerr case. Where Myers-Perry diverges is the actual inverse Fourier transform, we have
\begin{align}
    c_X(\hat{f}(0))&\mathcal{F}^{-1}(\sigma^2 (\tilde{u}^{(2)}(\hat{r})H(\hat{r})+\overline{\tilde{u}^{(2)}}(-\hat{r})H(-\hat{r})\\
    &=c_{X}\hat{f}(0)\rho^3\frac{1}{2\pi}\int^\infty_0 e^{-i\hat{r}\rho t_*}\hat{r}^2\tilde{u}^{(2)}(\hat{r})+e^{i\hat{r}\rho t_*} \hat{r}^2\tilde u^{(2)}(\hat{r}) d\hat{r}\\
    &=\frac{c_X(\hat{f}(0))\rho^3}{\pi}\text{ Re}\left.\left(\int^\infty_0 e^{-i\hat{r}v} \hat{r}^2 \tilde u^{(2)}(\hat{r})d\hat{r}\right)\right|_{v=\rho t_*}\\
    &=-\frac{c_X(\hat{f}(0))\rho^3}{\pi v^2}\text{Re}\left.\left(\int^\infty_0 e^{-i\hat{r}v} \partial^2_{\hat{r}}\hat{r}^2\tilde u^{(2)}(\hat{r}) \,d\hat{r}\right)\right|_{v=\rho t_*}.\label{eq:IvIntegral}
\end{align}
This latter integral can be evaluated using the fundamental solution for $\tilde u^{(2)}$ and the necessary Laplace transforms for Bessel and Hankel functions. We find that 
\begin{align}\label{eq:IvIntegralSolution}
    \text{ Re}\left.\left(\int^\infty_0 e^{-i\hat{r}v} \partial^2_{\hat{r}}\hat{r}^2\tilde u^{(2)}(\hat{r}) \,d\hat{r}\right)\right|_{v=\rho t_*}=\frac{\pi v^{1/2}}{(v+2)^{3/2}}
\end{align}
Plugging this back in and using $\rho^3 = v^3t_*^{-3}$ we find that 
\begin{align}
    c_X(\hat{f}(0))&\mathcal{F}^{-1}(\sigma^2 (\tilde{u}^{(2)}(\hat{r})H(\hat{r})+\overline{\tilde{u}^{(2)}}(-\hat{r})H(-\hat{r})=c_M(f)t_*^{-3}\left(\frac{v}{v+2}\right)^{3/2},
\end{align}
where we used that $c_M(f)=-c_X(\hat{f}(0))$.
We now conclude that from the above it follows that 
\begin{align}\label{eq:RegularityatM_+}
    \phi \in \mathcal{A}^{0,((3,0),3+\alpha-),((3,0),3+\alpha-)}(M_+).
\end{align}
\item (Decay at null infinity) We integrate along asymptotic characteristics as in \cite{hintz2022sharp} and \cite{hintz2020stability} and use the same notation as in \cite{hintz2022sharp}, i.e. let 
\begin{align}
    v:=\rho_{\mathscr{J}^+,1}=\rho t_*, \;\;\; \tau =\rho_{I^+,1}=t_*^{-1},
\end{align}
near $\mathscr{J}^+\subset M_+$. We then get the following conormal transformations: $\partial_{t_*}=\tau(-\tau\partial_{\tau} + v\partial_v)$ and $\rho\partial_{\rho}=v\partial_v$. On the neighbourhood $M'=[0,1)_\tau \times [0,1)_v\times \partial X$ of $\mathscr{J}^+\cap I^+\subset M_+$ we then see $\Box$ factorizes as follow
\begin{align}
    \Box = 2v\tau^2(A+R),\;\;\;\;A:=(\tau\partial_{\tau}-v\partial_v)(v\partial_v-\frac{3}{2}) \in \text{Diff}^2_b(M'),\; R\in v\text{Diff}_b^2(M').
\end{align} 

By (\ref{eq:RegularityatM_+}) we then know that $\phi\in\mathcal{A}^{0,3}(M')=v^0\tau^3\mathcal{A}^{0,0}(M')$ is a solution of 
\begin{align}
    (A+R)\phi =f':=(2v\tau^2)^{-1}f\in \mathcal{A}^{3+\alpha,\infty}(M'),
\end{align}
which can again be rewritten as $A\phi = f'-R\phi\in \mathcal{A}^{1,3}(M')$. We integrate out $\tau\partial_\tau-v\partial_v$ with initial data set at $v=\frac{1}{2}$ of class $\mathcal{A}^3([0,1)_\tau)$ (which is justified by the analysis at $I^+$, i.e. (\ref{eq:RegularityatM_+})) and we obtain $(v\partial_v-3/2)\phi\in \mathcal{A}^{1,3}(M')$ and thus $\phi \in \mathcal{A}^{1,3}(M')$ as there is no logarithmic loss due to $v^1$ not being an indicial solution to $v\partial_v-3/2$. Now iterating once more we thus find that $A\phi \in \mathcal{A}^{2,3}(M')$ and the inversion of $v\partial_v-3/2$ this time gives a $v^{3/2}$ leading term of $\phi$ coming from the homogeneous solution which dominates near $\mathscr{J}^+$, so we conclude that 
\begin{align}
    \phi \in \mathcal{A}^{((3/2,0),2),3}(M')
\end{align}
To now narrow down the $\tau$ asymptotics even further next to this argument note that $\phi|_{v=1/2} \in \mathcal{A}^{((3,0),3+\alpha-)}([0,1)_\tau)$ has a $\tau^3$ leading order term at $\tau=0$ (again, justified by (\ref{eq:RegularityatM_+}))and $R\phi \in A^{1,((3,0),3+\alpha-)}(M_+)$. Integration gives $\phi\in \mathcal{A}^{1-,((3,0),3+\alpha-)}(M')$, hence we have
\begin{align}\label{eq:NormalOperatorDecayArgument}
    A\phi=f'-R\phi\in \mathcal{A}^{2-,((3,0),3+\alpha-)}.
\end{align}
Integrating this further then shows that the $v^{3/2}$ leading order term (i.e. the radiation field in $1+4$ dimensions), $\phi_{rad}$ at $\mathscr{J}^+$ has a leading order term
\begin{align}
    \phi_{rad}-c_{rad}(f)\tau^3 \in \mathcal{A}^{3+\alpha-}([0,1)_\tau),\;\;\;\; c_{rad}(f) \in C^\infty(\partial X) =C^\infty(\mathbb{S}^3).
\end{align}
By matching at the intersection $I^+\cap \mathscr{J}^+$ we also know that the prefactor of $\phi_{rad}$ needs to match that of $\phi$ found in (\ref{eq:RegularityatM_+}). More specifically, near the edge $v=\tau=0$ we need to have the equality 
\begin{align}
    c_{rad}(f)v^{3/2} \tau^3 =c_M(f) \tau^3
 u^+(v).
 \end{align}
Solving for $c_{rad}(f)$ and taking the limit $v\searrow0$ we find
\begin{align}
    c_{rad}(f) = \lim_{v\searrow0} c_M(f) \frac{1}{(v+2)^{3/2}}= \frac{1}{2\sqrt{2}}c_M(f)
\end{align}
as required. Iterating this process over the indicial gap $\alpha \in(0,2)$ for (\ref{eq:NormalOperatorDecayArgument}) until we hit the $v^{3+\alpha}$ term we can conclude $\phi \in \mathcal{A}^{(3/2,0),3+\alpha),((3,0),3+\alpha-)}(M')$, proving the theorem.
\end{itemize}

\end{proof}

\begin{remark}
    Note that the computation of $I^+$ agrees exactly with the $d=4$ case in \cite{hintz2023linear}. Indeed writing $v=\rho t_*= r^{-1}t_*$ we have that 
    \begin{align}
        t_*^{-3} u^+(v) = t_*^{-3}\left(\frac{v}{v+2}\right)^{3/2}=\bigl[t_*(t_*+2r)\bigr]^{-3/2},
    \end{align}
    which is the corresponding Minkowski profile. \cite{hintz2023linear} computes this profile using the Minkowski fundamental solution, with an alternative derivation through the distributional equation (see 4.66a-4.66b in \cite{hintz2023linear}). In our proof we recover the profile directly from the inverse Fourier transform of the transition-face model. Specifically the representation of $\tilde{u}^{(2)}$ in terms of Bessel and Hankel functions allows \ref{eq:IvIntegral} directly to be evaluated using Laplace transforms and Abel regularization, yielding (\ref{eq:IvIntegralSolution}).   
\end{remark}

We will show that that initial value problems can be reduced to forced wave-problems in the present setting as in \cite{hintz2022sharp} Corollary 3.11.
\begin{corollary}\label{cor:InitialDataCorollary}
    Let $(M,g_{\bfa,\bfb,\mfrm})$ be a Myers-Perry spacetime with $(|\bfa|+|\bfb|)^2< 2\mfrm$ and spatial manifold $X$. Then the solution of 
    \begin{align}
        \Box_{g_{\bfa,\bfb,\mfrm}}\phi=0,\;\phi|_{t_*=0}=\phi_0,\; \partial_{t_*}\phi|_{t_*=0}=\phi_1
    \end{align}
    with $\phi_0,\phi_1\in C_c^\infty(X^\circ)$ has the asymptotic behaviour as stated in Theorem \ref{thm:GlobalAsymptotics} for every $0<\alpha<1$, with $c_M(f)$ replaced by
    \begin{align}
        c_M(\phi_0,\phi_1) :=\frac{1}{4\pi^2}\left\langle2\rho Q\phi_0 +g^{00}\phi_1, u_{(0)}^*\right\rangle, 
    \end{align}
    where the pairing is taken with respect to the stationary spatial density $dg_{\bfa,\bfb,\mfrm}$ and $Q$ and $g^{00}=g^{-1}(dt_*,dt_*)$ are as in Lemma \ref{lem:AdmissibleWaveEquation}.\\
    Equivalently, after extending $\phi_0$ independently of $t_*$, we have that 
    \begin{align}\label{eq:CMphi0phi1Factor}
        c_M(\phi_0,\phi_1) =\frac{1}{4\pi^2}\left\langle-[\Box_{g_{\bfa,\bfb,\mfrm}},t_*]\phi_0|_{t_*=0}+g^{-1}(dt_*,dt_*)\phi_1,u_{(0)}^* \right\rangle.
    \end{align}
    In particular the leading terms at are given by
    \begin{align}
        c_M(\phi_0,\phi_1)u_{(0)}t_*^{-3} \text{ at } \mathcal{K}^+ \text{ and } c_M(\phi_0,\phi_1)\left(\frac{\rho t_*}{\rho t_*+2}\right)^{3/2}t_*^{-3} \text{ at } I^+,
    \end{align}
    while the radiation coefficient satisfies 
    \begin{align}
        \phi_{rad}(t_*,\omega)-\frac{c_M(\phi_0,\phi_1)}{2\sqrt{2}}t_*^{-3}\in \mathcal{A}^{3+\alpha-}(\mathscr{J}^+).
    \end{align}
\end{corollary}

\begin{proof}
    By finite propagation speed, there exist $\epsilon >0$ and a compact $K\Subset X^\circ$ such that 
    \begin{align}
        \supp \,\phi(t_*,\cdot)\subset K, \;\;\; 0\leq t_* \leq 3\epsilon.
    \end{align}
    Choose $\chi\in C^\infty(\mathbb{R})$ which is forward supported in $\mathbb{R}^+$, that is 
    \begin{align}
        \chi \equiv 0 \text{ for } t_*\leq \epsilon, \text{ and } \chi\equiv 1 \text{ for } t_*\geq 2\epsilon.
    \end{align}
Extending $\phi_+ := \chi\phi$ by $0$ for $t_*<0$ we obtain
\begin{align}
    \Box_{g_{\bfa,\bfb,\mfrm}}\phi_+ = f, \text{ and } f:=[\Box_{g_{\bfa,\bfb,\mfrm}},\chi]\phi \in C^\infty_c((0,\infty)_{t_*}C^\infty_c(X^\circ),
\end{align}
So Theorem \ref{thm:GlobalAsymptotics} applies to $\phi_+$. Since $\phi-\phi_+$ is smooth and supported away from $\mathcal{K}^+,I^+$ and $\mathscr{J}^+$, $\phi$ must have the same asymptotics as $\phi_+$.\\
It remains now to compute $c_M(f)$. We have by Lemma \ref{lem:AdmissibleWaveEquation} that 
\begin{align}
    \Box_{g_{\bfa,\bfb,\mfrm}}= \hat{\Box}_{g_{\bfa,\bfb,\mfrm}}(0)-2\rho Q\partial_{t_*} -g^{00}\partial_{t_*}^2,
\end{align}
and thus we have that 
\begin{align}
    f= -2\chi'\rho Q\phi -g^{00}\chi"\phi - 2g^{00}\chi'\partial_{t_*}\phi.
\end{align}
Integrating the middle term by parts once in $t_*$ gives
\begin{align}
    \int_{\mathbb{R}}\langle f(t_*),u^*_{(0)} \rangle dt_* =-\int_\mathbb{R} \chi'(t_*) \left\langle 2\rho Q\phi(t_*) + g^{00}\partial_{t_*}\phi(t_*), u^*_{(0)} \right\rangle dt_*
\end{align}
We claim that the pairing in the last integrand is constant on $\supp \chi'\subset \{0\leq t_*\leq 3\epsilon\}$. To see this, observe that 
\begin{align}
    \frac{d}{dt_*}\left\langle 2\rho Q\phi + g^{00}\partial_{t_*}\phi,u_{(0)}^* \right\rangle&=\langle \hat{\Box}_{g_{\bfa,\bfb,\mfrm}}(0)\phi, u_{(0)}^*\rangle\\
    &=\langle \phi,\hat{\Box}_{g_{\bfa,\bfb,\mfrm}}(0)^*u_{(0)}^*\rangle =0,
\end{align}
where we used stationarity, Lemma \ref{lem:AdmissibleWaveEquation}, and the fact that in the proof of Corollary \ref{cor:KernelElementsMP} we derived that $\hat{\Box}_{g_{\bfa,\bfb,\mfrm}}(0)^*u_{(0)}^*=0$. Since we can arrange that $\int^{3\epsilon}_0 \chi' dt_*=1$ it follows that 
\begin{align}
    \int_{\mathbb{R}}\langle f(t_*),u_{(0)}^*\rangle dt_*=-\left\langle 2\rho Q\phi + g^{00}\partial_{t_*}\phi,u_{(0)}^* \right\rangle
\end{align}
Finally, (\ref{eq:CMphi0phi1Factor}) follows from the fact that
\begin{align}
    [\hat{\Box}_{g_{\bfa,\bfb,\mfrm}}(0),t_*]=-2\rho Q -2g^{00}\partial_{t_*}.
\end{align}

\end{proof}
\begin{remark}
    As in \cite{hintz2022sharp} the first thing to remark is that the identity $u^*_{(0)}=H(r-r_+)$ shows that only the exterior part of the forcing contributes to $c_M(f)=-\int_{\mathbb{R}_{t_*}\times\{r>r_+\}} f\,dg_{\bfa,\bfb,\mfrm}$ consistent with the causal indepence of the forward exterior solution from such data. In contrast to Kerr however, the absence of an overall mass factor exhibits the absence of the strong Huygens principle and can already be seen on the $1+4$-dimensional Minkowsky space: its retarded fundamental solution has a $t_*^{-3}$ tail see e.g. H\"ormander Chapter 6. Correspondingly, the flat spatial model produces the logarithmic low-frequency term responsible for the $t_*^{-3}$ asymptotics see e.g. \cite{hormander1997lectures} Chapter 6.1. Note that also the radiation field in Theorem \ref{thm:GlobalAsymptotics} exactly carries the $3/2$ decay factor as expected from our early calculations for the Minkowski model in Proposition \ref{prop:IndicialDecayRootsMinkowski} and should also be compared to the differentiated radiation field in Theorem 1.1 and Corollary 1.7 of \cite{baskin2015asymptotics}. We stress that the the mass $\mfrm>0$ still enters the geometry and the Cauchy initial data through (\ref{eq:CMphi0phi1Factor}).
\end{remark}

\begin{theorem}[Price's law on Myers--Perry]
\label{thm:MPPriceLaw}
Let $(M,g_{\bfa,\bfb,\mfrm})$ be a subextremal Kerr spacetime with $\mfrm>0$, let $(\bfa,\bfb)\in D_m$, and and $X$ be its associated spatial manifold.
Fix $0<\alpha<1$, and let
\begin{align}
    f\in C_c^\infty
    \left((0,\infty)_{t_*};\mathcal A^{4+\alpha}(X)\right), 
\end{align}

and denote by $\phi$ the unique forward solution of
\begin{align}
    \Box_{g_{a,b,m}}\phi=f.
\end{align}
Let
\begin{align}
    c_M(f)=-\frac{1}{4\pi^2}\int_0^\infty\int_{r_+}^\infty\int_{\mathbb S^3}f(t_*,r,\phi_*,\psi_*,\theta)\,
    r\,r_{\bfa,\bfb}^{\,2}\cos\theta\sin\theta\,
 \,dr\,dt_* d\phi\,d\psi\,d\theta. 
\end{align}

For every compact $K\Subset X^\circ$, every $j\in\mathbb N_0$,
every spatial multi-index $\beta\in\mathbb N_0^4$, we have in local spatial coordinates $x=(r,\theta,\phi_*,\psi_*)$
\begin{align}
    \left|\partial_{t_*}^{\,j}\partial_x^\beta \bigl(\phi(t_*,x)-c_M(f)t_*^{-3}\bigr)
    \right|\leq C_{K,j,\beta,\epsilon}\, t_*^{-3-\alpha-j+},\qquad t_*>1,\quad x\in K.
\end{align}
More generally, let
$\Omega=(\Omega_1,...,\Omega_3)$ denote the standard rotation
vector fields that span $\mathcal{V}(\mathbb S^3)$ over $\mathbb{S}^3$\footnote{Again we use parallzelizabilit of $\mathbb{S}^3$ here.}. Then, for every $j,k\in\mathbb N_0$, $\gamma\in\mathbb N_0^3$, and $0<\epsilon<\alpha$, we have the global estimate
\begin{align}
    \left|
    (t_*\partial_{t_*})^j
    (r\partial_r)^k
    \Omega^\gamma
    \left(
        \phi(t_*,r,\omega)
        -
        \frac{c_M(f)}
        {[t_*(t_*+2r)]^{3/2}}
    \right)
    \right|
    \leq
    C_{j,k,\gamma,\epsilon}\,
    t_*^{-3-\alpha+\epsilon}
    \left(\frac{t_*}{t_*+r}\right)^{3/2},
\end{align}
which holds for $t_*>1$ and all $(r,\phi_*,\psi_*,\theta)\in X^\circ$.

In particular, the radiation field of $\phi$, defined by

\begin{align}
    F(t_*,\omega)
    :=
    \lim_{r\to\infty}
    r^{3/2}\phi(t_*,r,\omega)
\end{align}
is given by $F(t_*,\omega) = \frac{1}{2\sqrt{2}}c_M(f)t_*^{-3/2} +\mathcal{O}(t_*^{-3/2-\alpha+})$,
and satisfies the pointwise decay estimates
\begin{align}
    \left|
    \partial_{t_*}^{\,j}\Omega^\gamma
    \left(
        F(t_*,\omega)
        -
        \frac{c_M(f)}{2\sqrt{2}}\,t_*^{-3/2}
    \right)
    \right|
    \leq
    C_{j,\gamma,\epsilon}\,
    t_*^{-3/2-\alpha-j+}
\end{align}

for all $j\in\mathbb N_0$, $\gamma\in\mathbb N_0^3$,
$0<\epsilon<\alpha$, and $t_*>1$.
\end{theorem}

\begin{proof}
By Corollary \ref{cor:KernelElementsMP} we have $u_{(0)}=1$ and
$u_{(0)}^*=H(r-r_+)$. Thus the coefficient in
Theorem \ref{thm:GlobalAsymptotics} is precisely $c_M(f)$, and the estimates
on compact spatial sets follow from that theorem.\\
Since
\begin{align}
    t_*^{-3}u^+(\rho t_*)
    =
    \frac{1}{[t_*(t_*+2r)]^{3/2}},
\end{align}
the global expansions of Theorem \ref{thm:GlobalAsymptotics} give
\begin{align}\label{eq:SubstractedProfile}
    \phi-\frac{c_M(f)}{[t_*(t_*+2r)]^{3/2}}
    \in
    \mathcal A^{3/2,\,3+\alpha-,\,3+\alpha-}(M_+).
\end{align}
Here the discrepancy at $\mathcal{K}^+$ between the subtracted
profile and $c_M(f)t_*^{-3}$ is $\mathcal{O}(t_*^{-4})$, which
is absorbed in the remainder for $\alpha<1$. The exact weight $3/2$
at $J^+$ follows from the radiation-coefficient estimate
and interpolation with the higher-order $\mathscr{J}^+$ remainder,
allowing arbitrarily small losses in the temporal decay order.

On $t_*>1$, choose equivalent boundary defining functions for $\mathscr{J}^+$, 
\begin{align}
    \rho_{J^+}=\frac{t_*}{t_*+r},
    \qquad
    \rho_{I^+}=\frac{t_*+r}{t_*r},
    \qquad
    \rho_{K^+}=\frac{r}{t_*+r}.
\end{align}
Then
\begin{align}
    \rho_{J^+}^{3/2}
    \rho_{I^+}^{3+\alpha-\epsilon}
    \rho_{K^+}^{3+\alpha-\epsilon}
    =
    \left(\frac{t_*}{t_*+r}\right)^{3/2}
    t_*^{-3-\alpha+\epsilon}.
\end{align}
Since $t_*\partial_{t_*}$, $r\partial_r$, and the rotation
vector fields lift to b-vector fields on $M_+$, the
conormal estimate for the substracted profile (\ref{eq:SubstractedProfile}) proves the global pointwise bounds.

Finally, the expansion at $\mathscr{J}^+$ implies
\begin{align}
    F=t_*^{3/2}\phi_{\mathrm{rad}},
    \qquad
    F-\frac{c_M(f)}{2\sqrt{2}}\,t_*^{-3/2}
    \in\mathcal A^{3/2+\alpha-}(\mathscr{J}^+).
\end{align}
Conormal differentiation gives the asserted
radiation-field estimates.
\end{proof}

\begin{corollary}[Price's law for Myers-Perry in Boyer-Lindquist coordinates]. Let $(M,g_{\bfa,\bfb,\mfrm})$ be a subextremal Myers-Perry spacetime with $\mfrm>0$, $(|\bfa|+|\bfb|)^2<2\mfrm$. Suppose $\phi_0,\phi_1\in C_c^\infty((r_+,\infty)_r\times\mathbb{S}^3)$. Then the solution to the initial value problem
\begin{align}
    \Box_{g_{\bfa,\bfb,\mfrm}}\phi=0,\;\phi|_{t=0}=\phi_0,\;\partial_t\phi|_{t=0}=\phi_1,
\end{align}
has the asymptotic behaviour as stated in Theorem \ref{thm:MPPriceLaw} for every $\alpha<1$ with the constant $c_M(f)$ replaced by 
\begin{align}\label{eq:BoyerLindquistPrincipalComponent}
    c_M(\phi_0,\phi_1)=-\frac{1}{4\pi^2}\int_{r_+}^\infty\int_{\mathbb{S}^3}\left(r\bfr^2+\frac{2\mfrm(r^2+\bfa^2)(r^2+\bfb^2)}{r\Delta_r}\right)\phi_1(r,\phi,\psi,\theta)\cos\theta\sin\theta\, dr\,d\phi\,d\psi\,d\theta.
\end{align}

\end{corollary}
\begin{proof}
    Choose the admissible time function $t_*$ to agree with Boyer-Lindquist time $t$ near $\supp\, \phi_0 \cup \supp\, \phi_1$. Corollary \ref{cor:InitialDataCorollary} applies for these data extended by zero on the remainder of the spatial slice $\{t_*=0\}$. We then get for its coefficient formula that 
    \begin{align}
        c_M(\phi_0,\phi_1) =\frac{1}{4\pi^2}\int_{r>r_+} \int_{\mathbb{S}^3} \left(-[\Box_{g_{\bfa,\bfb,\mfrm}},t]\phi_0 +g^{tt}\phi_1\right)r\bfr^2\cos\theta\sin\theta drd\phi d\psi d\theta,
    \end{align}
    with $\phi_0$ extended independently of $t$ and $g^{tt}=g^{-1}(dt,dt)$.
    In Boyer-Lindquist coordinates, we have $g^{tr}=g^{t\theta}=0$, and the metric and its density are still independent of the Killing vector fields $\partial_\phi$ and $\partial_\psi$. Consequently,
    \begin{align}
        -[\Box_{g_{\bfa,\bfb,\mfrm}},t]\phi_0= 2g^{t\phi}\partial_\phi\phi_0 +2g^{t\psi}\partial_\psi \phi_0,
    \end{align}
    and the integral of this component vanishes by integration by parts over $\mathbb{S}^3$.
    To evaluate this, we use the Boyer-Lindquist expression in (\ref{eq:FrolovExpressiontt}),
    observe that
    \begin{align}
        g^{tt}&=\frac{1}{\bfr^2}\left[(\bfa^2-\bfb^2)\sin^2\theta-\frac{(r^2+\bfa^2)(r^2\Delta_r +2\mfrm(r^2+\bfb^2)}{r^2\Delta_r} \right]\\
        &=\frac{1}{\bfr^2}\left[(\bfa^2-\bfb^2)\sin^2\theta-(r^2+\bfa^2)-\frac{2\mfrm(r^2+\bfb^2)}{r^2\Delta_r} \right]\\
        &=-1-\frac{2\mfrm(r^2+\bfa^2)(r^2+\bfb^2)}{r^2\Delta_r\bfr^2},
    \end{align}
    where we used $(\bfa^2-\bfb^2)\sin^2\theta -(r^2+\bfa^2) =-\bfr^2$ in the final step, from which the claim follows.
\end{proof}

\begin{remark}
    Note that the mass-dependent time in (\ref{eq:BoyerLindquistPrincipalComponent}) arises from the Boyer-Lindquist term in $g^{tt}$. Moreover as outlined in the proof, the contribution of $\phi_0$ vanishes by angular integration by parts, without any symmetry assumptions on the data. In particular, if $\phi_1=0$, or if the weighted integral in (\ref{eq:BoyerLindquistPrincipalComponent}) vanishes, the preceding estimates give $\phi= \mathcal{O}(t_*^{-4+})$ at $K\Subset X^\circ$ and a radiation field $F = \mathcal{O}(t_*^{-5/2+})$, where $\phi$ denote sthe horizon-regular extension constructed in the proof.
\end{remark}

\section{Appendix A. Solution to the transition face model.}

In this section we will explore the solution to the $tf$ wave equation in terms of it's Bessel and Hankel functions, and show that its Fourier transform gives the desired qualitative decay profile at $I^+$. 

Starting out we are interested in solving 
\begin{align}
    \tilde\Box(1)u=\hat{r}^{-2}
\end{align}
in $\hat{\rho}=\frac{\sigma}{\rho}$ coordinates. Here $\tilde{\Box}(1)$ is given by 
\begin{align}
    \tilde{\Box}(1) = -\hat{r}^{-2}[(\hat{r}\partial_{\hat{r}})^2-2\hat{r}\partial_{\hat{r}} +\Delta_{\mathbb{S}^3}] -2i\hat{r}^{-1}\left(\hat{r}\partial_{\hat{r}}+\frac{3}{2}\right).
\end{align}
Hence solving this equation amounts to solving 
\begin{align}
-[(\hat{r}\partial_{\hat{r}})^2-2\hat{r}\partial_{\hat{r}} +\Delta_{\mathbb{S}^3}]u -2i\hat{r}\left(\hat{r}\partial_{\hat{r}}+\frac{3}{2}\right)u=1.
\end{align}
Unlike the Kerr case, this equation cannot be factored into two Euler equations (see \cite{hintz2022sharp}, Remark 2.25). However, the left hand side of this equation can be reduced to be of Bessel-type following \cite{derezinski2017schrodinger} Chapter 3. Moreover we are fundamentally only interested in the $\ell = 0$ mode. Let us write $D=\hat{r}\partial_{\hat{r}}$ and note that 
\begin{align}
    D^2+2D &u = \hat{r}^2 u''+3\hat{r}u'\\
    \hat{r}\left(D+\frac{3}{2}\right)u&=\hat{r}^2u'+\frac{3}{2}\hat{r}u.
\end{align}
With this we see that we are solving 
\begin{align}
    -(\hat{r}^2 u''+3\hat{r}u')-2i\left(\hat{r}^2u'+\frac{3}{2}\hat{r}u\right)=1,
\end{align}
or
\begin{align}
    \left( u''+\frac{3}{\hat{r}}u'\right)+2i\left(u'+\frac{3}{2\hat{r}}u\right)=-\hat{r}^{-2}.
\end{align}
Now to get this into Bessel form we perform a weighting by the indicial root:
\begin{align}
    w(\hat{r}):=\hat{r}^{3/2}u(\hat{r}),
\end{align}
where we note that
\begin{align}
    u'(\hat{r}) &=-\frac{3}{2}\hat{r}^{-5/2}w(\hat{r})+\hat{r}^{3/2}w'(\hat{r}).\\
    u''(\hat{r})&=\hat{r}^{-3/2}w''(\hat{r})-3\hat{r}^{-5/2}w'(\hat{r})+\frac{15}{4}\hat{r}^{-7/2}w(\hat{r}).
\end{align}
We then see that the equation transforms and we get rid of the transport term through cancellations
\begin{align}
    u''(\hat{r})+3\hat{r}^{-1}u'(\hat{r})&=\hat{r}^{-3/2}\left[w''(\hat{r})-\frac{3}{4}\hat{r}^{-2}w(\hat{r})\right]\\
    2iu'(\hat{r})+\frac{3i}{\hat{r}}u(\hat{r})&=2i\hat{r}^{-3/2}w'(\hat{r}),
\end{align}
whence the ODE reduces after multiplying through by $\hat{r}^{3/2}$
\begin{align}
    w''(\hat{r})+2iw'(\hat{r})-\frac{3}{4\hat{r}^2}w(\hat{r})=-\hat{r}^{-1/2}.
\end{align}
To get rid of the full linear term we now substitute again by a phase factor, 
\begin{align}
    q(r):=e^{i\hat{r}}w(\hat{r})=e^{i\hat{r}}\hat{r}^{3/2}u(\hat{r}).
\end{align}
We then get that 
\begin{align}
    w'(\hat{r})&=e^{-i\hat{r}}q'(\hat{r})-ie^{-i\hat{r}}q(\hat{r})\\
    w''(\hat{r})&=e^{-i\hat{r}}q''(\hat{r})-2ie^{-i\hat{r}}q'(\hat{r})-e^{-i\hat{r}}q(\hat{r}).
\end{align}
We then see that this yields us a reduced inhomogenous Bessel equation: \cite[Eq.~(10.13.1)]{NIST:DLMF}
\begin{align}\label{eq:ReducedBesselEquation}
    q''(\hat{r})+\left(1-\frac{3}{4\hat{r}^2}\right)q(\hat{r})=-\hat{r}^{-1/2}e^{i\hat{r}}.
\end{align}
It is well known that that a general type equation of the type 

\begin{align}
    q''(z)+\left(a-\frac{\nu^2-\frac{1}{4}}{z^2}\right)q(z)=0
\end{align}
transforms into the standard Bessel equation 
\begin{align}
    x^2y''(x)+xy'(x) +(x^2-\nu^2)y=0
\end{align}
under the standard substitutions $x=az$ and $q(z)=\sqrt{z}y(az)$. In our case case $a=1$ and $\nu=\pm 1$. We thus get after also transforming the r.h.s.
\begin{align}
    \hat{r}^2y''(\hat{r})+\hat{r}y'(\hat{r})+(\hat{r}^2-1)y(\hat{r}) = -\hat{r}^{-1}e^{i\hat{r}}.
\end{align}
The basis of solutions to the homogeneous equation to (\ref{eq:ReducedBesselEquation}) is given by 
\begin{align}
    q_0(\hat{r})=\hat{r}J_1(\hat{r}), \text{ and } q_+(\hat{r})=\hat{r}H^{(1)}_1(\hat{r}).
\end{align}
where $J_1$ and $H^{(1)}_1$ are first order Bessel and Hankel functions respectively. Recall that in particular 
\begin{align}
    H^{(1)}_1(x) = J_1(x) + i Y_1(x),
\end{align}
where $Y_1(x)$ is the first order Neumann function.
Their Wronksian is convienently well-known and constant:
\begin{align}
    \mathscr{W}(q_0,q_+)(\hat{r})&\equiv \frac{2i}{\pi}.
\end{align}
The outgoing solution (i.e. $q(\hat{r})\sim \hat{r}^{-1/2}$ as $\hat{r}\to\infty$) solution that is regular at $\hat{r}=0$ can be constructed using variation of parameters:
\begin{align}
    q(\hat{r}) &= q_+(\hat{r}) \int^{\hat{r}}_0 \frac{q_0(s)(-s^{-1/2}e^{is})}{\mathscr{W}(q_0,q_+)(s)}+q_0(\hat{r}) \int^\infty_0 \frac{q_+(s)(-s^{-1/2}e^{is})}{\mathscr{W}(q_0,q_+)(s)}\\
    &=\frac{\pi i}{2}\left[\sqrt{\hat{r}}H^{(1)}_1(\hat{r}) \int^{\hat{r}}_0 e^{is}J_1(s)ds + \sqrt{\hat{r}}J_1(s)\int^\infty_{\hat{r}}e^{is} H^{(1)}_1(s)ds\right],
\end{align}

and thus substituting back $u(\hat{r})=\hat{r}^{-3/2}e^{-i\hat{r}}q(\hat{r})$ we find
\begin{align}
    u(\hat{r})=\frac{\pi i}{2\hat{r}}e^{-i\hat{r}}\left[H^{(1)}_1(\hat{r}) \int^{\hat{r}}_0 e^{is}J_1(s)ds + J_1(s)\int^\infty_{\hat{r}}e^{is} H^{(1)}_1(s)ds\right],
\end{align}

\label{eq:FullVoPSolution}
Notice that we understand the second integral of (\ref{eq:FullVoPSolution}) to be convergent under standard Abel regularization, that is by including a small factor $e^{-\epsilon s}$ in the integrand and letting $\epsilon\to 0$ in the evaluation.
To massage this equation further we make use of the following derivative identities:
\begin{align}
    e^{ix}J_1(x)&= \frac{d}{dx}\left(1+e^{ix}((ix-1)J_0(x) + xJ_1(x) \right)\\
    e^{is}H^{(1)}_1(x)&=-\frac{d}{dx}(e^{ix}((ix-1)H_0^{(1)}(x)+xH^{(1)}_1(x)).
\end{align}
Combined with the Wronksian identity $J_1(x)H^{(1)}_0(x)-J_0(x)H^{(1)}_1(x)=\frac{2i}{\pi x}$ we see that we can collapse the above to 
\begin{align}\label{eq:FullyRedcedSolution}
    u(\hat{r})=\frac{i\hat{r}-1}{\hat{r}^2}+\frac{\pi i}{2}\frac{e^{-i\hat{r}}H^{(1)}_1(\hat{r})}{\hat{r}}.
\end{align}
Note that $\hat{r}^{-1}e^{-i\hat{r}}H_1^{(1)}(\hat{r})\sim \hat{r}^{-3/2}$ as $\hat{r}\to \infty$ and that $u(\hat{r})$ is regular in $\hat{r}=0$. To see the later notice that the small $\hat{r}$ we will expand for small $\hat{r}$. Fist recall that 
\begin{align}
    H^{(1)}_1(x) =J_1(x) + i Y_1(x) \sim -\frac{2i}{\pi x}+\frac{x}{2}+\frac{ix}{\pi}\left(\log\left(\frac{x}{2}\right)+\gamma - \frac{1}{2}\right)+\mathcal{O}(x^3\log(x)),
\end{align}
where $\gamma$ is the Euler-Maschorini constant.
Now also using $e^{-ix}\sim 1-ix+\mathcal{O}(x^2)$, we see that 
\begin{align}
    \frac{\pi i }{2}\frac{e^{-ix}H^{(1)}_1(x)}{x}=\frac{1}{x^2}-\frac{i}{x}-\frac{1}{2}+\frac{\pi i}{4}-\frac{1}{2}\left(\log\left(\frac{x}{2}\right)+\gamma - \frac{1}{2}\right)+x\left(\frac{\pi}{4}+\frac{i}{2}\left(\log\left(\frac{x}{2}\right)+\gamma-\frac{1}{2}\right)\right)+\mathcal{O}(x^2\log(x)).
\end{align}
Now the crucial insight is that the first two terms of this expansion are exactly cancelling out against 
\begin{align}
\frac{ix-1}{x^2}=\frac{i}{x}-\frac{1}{x^2}.
\end{align}
We conclude that indeed 
\begin{align}
    u(\hat{r})\sim -\frac{1}{2}\log\left(\frac{\hat{r}}{2}\right)-\frac{\gamma}{2}-\frac{1}{4}+\frac{i\pi}{4}+\mathcal{O}(\hat{r}\log\hat{r)}.
\end{align}

\subsection{Fourier transforming the fundamental solution.}

We are now interested in solving the inverse Fourier transform of $\partial_{\hat{r}}\hat{r}^2\tilde{u}^{(2)}(\hat{r})$ in $v$ at $I^+$, where $\tilde{u}^{(2)}$ is our tf fundamental solution to $\tilde \Box (1) \tilde u^{(2)}=\hat{r}^{-2}$ at tf, i.e. (\ref{eq:IvIntegral}) in Theorem \ref{thm:GlobalAsymptotics} and Let us set
\begin{align}
   I(v):=\left. \int^\infty_0 e^{-i\hat{r}v}\partial_{\hat{r}}^2\hat{r}^2\tilde u^{(2)}(\hat{r}) d\hat{r}\,\right|_{v=\rho t_*},
\end{align}
which with Abel regularization we would interpret as 
\begin{align}
    I_{\epsilon}(v)&= \left.\int^{\infty}_0 e^{-iv\hat{r}-\epsilon\hat
    r} F''(\hat{r}) d\hat{r} \right|_{v=\rho t_*},
\end{align}
where we also defined
\begin{align}
    F(\hat{r}):=\hat{r}^2u(\hat{r}).
\end{align}
We know from (\ref{eq:FullyRedcedSolution}) that we have 
\begin{align}
    F(\hat{r})= i\hat{r}-1 +\frac{\pi i}{2}\hat{r}e^{-i\hat{r}}H^{(1)}_1(\hat{r}).
\end{align}
We now note three more recurrence identities for Hankel functions (see \cite[10.5]{NIST:DLMF}, $10.5$), in this case
\begin{align}
    \frac{\partial}{\partial x}(xH^{(1)}_1(x))&= xH^{(1)}_0(x)\\ \frac{\partial}{\partial x}H^{(1)}_0(x)&=-H^{(1)}_1(x)\\
     \frac{\partial}{\partial x}(xH^{(1)}_0(x))&= H^{(1)}_0(x) - xH^{(1)}_1(x).
\end{align}

From this we can directly compute:
\begin{align}
    F'(\hat{r})=i+\frac{\pi i}{2}\hat{r}e^{-i\hat{r}}(H^{(1)}_0(\hat{r})-i H^{(1)}_1(\hat{r})).
\end{align}

and thus we see that 
\begin{align}
    F''(\hat{r})&=\frac{\pi i}{2}e^{-i\hat{r}} \left(-i\hat{r}H^{(1)}_0(\hat{r})-\hat{r}H^{(1)}_1(\hat{r}) + (H^{(1)}_0(\hat{r})-\hat{r}H^{(1)}_1(\hat{r})-i\hat{r}H^{(1)}_1(\hat{r})\right)\\
    &=\frac{\pi i}{2}e^{-i\hat{r}}\left((1-2i\hat{r})H^{(1)}_0(\hat{r})-2\hat{r}H^{(1)}_1(\hat{r})\right).
\end{align}

The above is exaclty the function that we want to Fourier transform in $\hat{r}$. Let us define a new conjugate variable as
\begin{align}
    p:=\epsilon+i(v+1),
\end{align}
we then see that that $I_\epsilon$ can be interpreted as a Laplace transform:
\begin{align}
    I_{\epsilon}(p)=\frac{\pi i}{2}\int^{\infty}_0 e^{-p\hat{r}}((1-2i\hat{r}) H^{(1)}_0(\hat{r}) -2\hat{r}H^{(1)}_1(\hat{r})) d\hat{r}.
\end{align}
The Laplace transforms for Bessel $J_\nu$ and Neuman functions $Y_\nu$ for $\text{ Re } p>0$ are well known, see (\cite{Polyanin2005Laplace} and \cite{GradshteynRyzhik2015} 6.611, using the logarithmic formula for $\text{arcsinh}(p)$).
\begin{align}
    \int^\infty_0 e^{-px}J_0(x) dx &= \frac{1}{\sqrt{1+p^2}}\\
    \int^\infty_0e^{-px}Y_0(x)dx&=-\frac{2}{\pi}\frac{\text{arcsinh}(p)}{\sqrt{1+p^2}},
\end{align}
from this one can again construct the Laplace transform for $H^{(1)}_0(x)$, using that $H^{(1)}_0(x)=J_0(x)+iY_0(x)$, and we have 
\begin{align}
    I^{(1)}_0(p):=\int^\infty_0 e^{-px}H^{(1)}_0(x) dx =\frac{1-\frac{2i}{\pi}\text{arcsinh}(p)}{\sqrt{1+p^2}}.
\end{align}

For the $xH^{(1)}_1(x)$ term we can proceed by integration by parts and using again that $(xH^{(1)}_1)'=xH^{(1)}_0(x)$. We see that 
\begin{align}
    \int^\infty_0 e^{-ipx}\underbrace{(xH^{(1)}_1(x))'}_{= xH^{(1)}_0(x)} dx &= \left[e^{-px}xH^{(1)}_1(x)\right]^\infty_0 + p \underbrace{ \int^\infty_0e^{-px} xH^{(1)}_1(x) dx}_{:=I^{(1)}_1(p)}\\
    -\frac{d}{dp}I^{(1)}_0(p)&=\left(0-\left(-\frac{2i}{\pi}\right)\right) + p\cdot I^{(1)}_1(p).
\end{align}
Rearranging we see that 
\begin{align}
    I_1^{(1)}(p)=-\frac{\frac{d}{dp}I^{(1)}_0(p)+\frac{2i}{\pi}}{p}.
\end{align}
Finally, the expression for $\frac{d}{dp}I^{(1)}_0(p)$ can be evaluated directly and we get
\begin{align}
    \frac{d}{dp} I^{(1)}_0(p)&=-\frac{2i}{\pi(1+p^2)}-\frac{p(1-\frac{2i}{\pi}\text{arcsinh}(p)}{(1+p^2)^{3/2}}\\
    &=-\frac{p(1-\frac{2i}{\pi}\text{arcsinh}(p))\sqrt{1+p^2}+\frac{2i}{\pi}\sqrt{1+p^2}}{(1+p^2)^{3/2}}.
\end{align}

Putting it all back together we see that
\begin{align}
    I_\epsilon(p)&=\frac{\pi i}{2}\left(I^{(1)}_0(p)-2i\frac{d}{dp}I^{(1)}_0(p)+\frac{2\frac{d}{dp}I^{(1)}_0(p)+\frac{4i}{\pi}}{p} \right)\\
    &=\frac{i}{2}\frac{(\pi-2i\text{ arcsinh}(p))\sqrt{1+p^2}+4(ip -1)}{(p+i)^2}
\end{align}
Now letting $\epsilon \searrow 0$ we see that as $v>0$, 
\begin{align}
    p\to i(v+1), \text{ and } v+1>1,
\end{align}
from which we find that 
\begin{align}
    \sqrt{1+p^2}&\to i\sqrt{v(v+2)}\\
    \text{arcsinh(p)}&\to \text{arcosh}(v+1)+i\frac{\pi}{2}.
\end{align}
Putting all of this in we conclude that the integral $I(v)$ is given by
\begin{align}
    I(v)=\frac{\pi\sqrt{v(v+2)}}{(v+2)^2}-i \frac{\sqrt{v(v+2)}\text{arcosh}(v+1)-2(v+2)}{(v+2)^2}.
\end{align}
And thus that 
\begin{align}
    \text{ Re } I(v)= \frac{\pi\sqrt{v}}{(v+2)^{3/2}} \text{ for } v>0.
\end{align}
We recall that this comes in with a prefactor $(\rho t_*)^{-2}\rho^3$. Using that $v=\rho t_*$ we see that this gives 
\begin{align}
    (\rho t_*)^{-2}\rho^2\frac{\sqrt{v}}{(v+2)^{3/2}}=t_*^{-3}\left(\frac{v}{v+2}\right)^{3/2}. 
\end{align}
On the other hand looking at the asymptotic profile in Proposition 6.4 in \cite{hintz2023linear} it is given as
\begin{align}
    (t_*(t_*+2\rho^{-1}))^{-3/2}&=t_*^{-3} (1+2 t_*/\rho)^{-3/2}\\
    &=t_*^{-3}(1+2v^{-1})^{-3/2}\\
    &=t_*^{-3}\left(\frac{v}{v+2}\right)^{3/2}
\end{align}
as desired.

\section{Appendix B. Visualization of the reduced quartic}

The mathematica code below vinishes the comparison of the rescaled polynomial $\overline{p}(x)$ versus $\overline{q}(x)$ in the proof of Lemma \ref{lem:NoLocalMinimum}. Note that values for $(|\bfa|+|\bfb|)^2 > 2\mfrm $ and $r<r_+$ are not prohibited by the Manipulate function. One can check if one is in the right regime by checking if the value of $r_+$ (in the inset) is real.

\label{app:reduced-quartic}

\label{app:mathematica}

\begin{lstlisting}[style=MPMathematica]
Δ[r_, a_, b_, m_] =
  r^2 + a^2 + b^2 + a^2 b^2 r^(-2) - 2 m

δΔ[r_, a_, b_, m_] =
  2 (r - a^2 b^2 r^(-3))

δδΔ[r_, a_, b_, m_] =
  2 + 6 a^2 b^2 r^(-4)

EventHorizon[r_, a_, b_, m_] =
  1/Sqrt[2] Sqrt[
    2 m - a^2 - b^2 +
    Sqrt[(2 m - a^2 - b^2)^2 - 4 a^2 b^2]
  ]

EventHorizon2[a_, b_, m_] =
  1/Sqrt[2] Sqrt[
    2 m - a^2 - b^2 +
    Sqrt[(2 m - a^2 - b^2)^2 - 4 a^2 b^2]
  ]

ExtremalityValue[r_, a_, b_, m_] =
  a^2 + b^2 + 2 Abs[a*b]

ExtremalityHorizon[r_, a_, b_, m_] =
  Sqrt[Abs[a*b]]

MonicRescaledSACmap[r_, a_, b_, m_] =
  8 a*b/(r*Sqrt[Δ[r, a, b, m]])
  (* Defines d(r)/a(r) *)

MonicRescaledSC2map[r_, a_, b_, m_] =
  r*δΔ[r, a, b, m]/Δ[r, a, b, m] - 4
  (* Defines c(r)/a(r) *)

MonicRescaledSA2map[r_, a_, b_, m_] =
  r^2 δΔ[r, a, b, m]^2/(4 Δ[r, a, b, m]^2) +
  r*δΔ[r, a, b, m]/(2*Δ[r, a, b, m]) -
  r^2 δδΔ[r, a, b, m]/(2 Δ[r, a, b, m])
  (* Defines e(r)/a(r) *)

DiscrDiff[r_, a_, b_, m_] =
  r^2 Δ[r, a, b, m]^3 (
    MonicRescaledSACmap[r, a, b, m]^2 -
    4*(MonicRescaledSC2map[r, a, b, m] + 2)*
      (MonicRescaledSA2map[r, a, b, m] - 1)
  )
  (* Computes 3.57 *)

FullSimplify[
  r^2 Δ[r, a, b, m]^3 (
    MonicRescaledSACmap[r, a, b, m]^2 -
    4*(MonicRescaledSC2map[r, a, b, m] + 2)*
      (MonicRescaledSA2map[r, a, b, m] - 1)
  )
]
(* Fully simplifies to 3.57 *)

Manipulate[
  Plot[
    {
      x^4 +
      MonicRescaledSC2map[r, a, b, m] x^2 +
      MonicRescaledSACmap[r, a, b, m] x +
      MonicRescaledSA2map[r, a, b, m],
      x^4 - 2 x^2 + 1
    },
    {x, -xbound, xbound},
    Epilog -> {
      Text[
        ExtremalityValue[r, a, b, m] "=Extr.",
        Scaled[{0.15, 0.9}]
      ],
      Text[
        2 m "=2m",
        Scaled[{0.15, 0.8}]
      ],
      Text[
        EventHorizon[r, a, b, m] "=r+",
        Scaled[{0.15, 0.7}]
      ],
      Text[
        ExtremalityHorizon[r, a, b, m] "=|ab|",
        Scaled[{0.15, 0.6}]
      ],
      Text[
        a "=a",
        Scaled[{0.15, 0.3}]
      ],
      Text[
        b "=b",
        Scaled[{0.15, 0.2}]
      ]
    },
    PlotRange -> {-ybound, ybound}
  ],
  {a, -4, 4},
  {b, -4, 4},
  {m, 0.1, 2},
  {r, 0.1, 50},
  {xbound, 1, 100},
  {ybound, 1, 100}
]
\end{lstlisting}

\bibliographystyle{alpha}
\bibliography{Bibliography}
\end{document}